\input amstex
\documentstyle{amsppt}
\nopagenumbers
\nologo
\magnification\magstephalf
%----------------------------------------------------------------
% Title:     Course of analytical geometry
% Author:    Ruslan Sharipov
% Comments:  AmSTeX, 226 pages, amsppt style
% MSC-class: 51N20, 51-01
% \\
% This book is a regular textbook of analytical geometry covering 
% vector algebra and its applications to describing straight lines,
% planes, and quadrics in two and three dimensions. The stress is 
% made on vector algebra by using skew-angular coordinates and 
% introducing some notations and prerequisites for understanding 
% tensors. The book is addressed to students specializing in mathematics, 
% physics, engineering, and technologies and to students of other 
% specialities where educational standards require learning this 
% subject. 
% \\
%---------------------------------------------------------------
\catcode`@=11
\redefine\output@{%
  \def\break{\penalty-\@M}\let\par\endgraf
  \global\voffset=-20pt
  \ifodd\pageno\global\hoffset=1pt\else\global\hoffset=-38pt\fi  
  \shipout\vbox{%
    \ifplain@
      \let\makeheadline\relax \let\makefootline\relax
    \else
      \iffirstpage@ \global\firstpage@false
        \let\rightheadline\frheadline
        \let\leftheadline\flheadline
      \else
        \ifrunheads@ %\let\makefootline\relax
        \else \let\makeheadline\relax
        \fi
      \fi
    \fi
    \makeheadline \pagebody \makefootline
  }%
  \advancepageno \ifnum\outputpenalty>-\@MM\else\dosupereject\fi
}
%---------------------------------------------------------------
\font\cpr=cmr7
\newcount\xnumber
\footline={\xnumber=\pageno
\divide\xnumber by 7
\multiply\xnumber by -7
\advance\xnumber by\pageno
\ifnum\xnumber>0\hfil\else\vtop{\vskip 0.5cm
\noindent\cpr CopyRight \copyright\ Sharipov R.A.,
2010.}\hfil\fi}
%---------------------------------------------------------------
\def\setfirstpage{\global\firstpage@true}
\def\myglue{\hskip 0pt plus 0.5pt minus 0.5pt}
\newbox\tmponebox
\newbox\tmptwobox
\catcode`@=\active
%---------------------------------------------------------------
%\def\startpage#1{\pageno=#1}
%---------------------------------------------------------------
\fontdimen3\tenrm=3pt
\fontdimen4\tenrm=0.7pt

%---------------------------------------------------------------
\def\leaderfill{\leaders\hbox to 0.3em{\hss.\hss}\hfill}
\font\tvbf=cmbx12
\font\tvrm=cmr12
\font\etbf=cmbx8

% Cyrillic fonts definition
\font\tencyr=wncyr10
\font\tencyrit=wncyi10
\font\tencyrbf=wncyb10
\font\eightcyr=wncyr8
\font\ninecyr=wncyr9
\font\ninecyrbf=wncyb9
%---------------------------------------------------------------
\def\negskp{\hskip -2pt}

\def\sign{\operatorname{sign}}

\def\Sopfr{\operatorname{Sopfr}}

\chardef\No=125
\def\compos{\,\raise 1pt\hbox{$\sssize\circ$} \,}
\redefine\cot{\operatorname{ctg}}
\redefine\arctan{\operatorname{arctg}}
\accentedsymbol\updownarrows{\,\uparrow\hskip-1.6pt\downarrow\,}
\def\msum#1{\operatornamewithlimits{\sum^#1\!{\ssize\ldots}\!\sum^#1}}
%\def\LaTeX{L\kern-.36em\raise.35ex\hbox{\Smallcaps a}\kern-.15em\TeX} 
%\def\LaTeX{L\kern-.36em\raise.3ex\hbox{\sc \uppercasesc a}\kern-.15em\TeX}
%---------------------------------------------------------------
\Monograph
\loadbold
\TagsOnRight
\newcount\chapternum
\newcount\partno
\accentedsymbol\tbb{\kern-3pt\tilde{\kern 3pt\bold b}}
\def\blue#1{#1}
\def\red#1{#1}
\catcode`#=11\def\diez{#}\catcode`#=6
\catcode`_=11\catcode`_=8
\catcode`~=11\catcode`~=\active
\chardef\wnIshort='022
\let\wnsoft=\volna
\def\mycite#1{\cite{\blue{#1}}\immediate\special{ps:
     ShrHPSdict begin /ShrBORDERthickness 0 def}}

\def\mytag#1{%
    \tag#1}
\def\mythetag#1{\thetag{\blue{#1}}\immediate\special{ps:
     ShrHPSdict begin /ShrBORDERthickness 0 def}}
\def\mythetagchapter#1#2{\thetag{\blue{#1}}\immediate\special{ps:
     ShrHPSdict begin /ShrBORDERthickness 0 def}}
\def\myrefno#1{\no#1}
\def\mylitno#1{\no#1}
\def\myhref#1#2{\blue{#2}\immediate\special{ps:
     ShrHPSdict begin /ShrBORDERthickness 0 def}}
\def\myEarXivlink{\myhref{http://arXiv.org}{http:/\negskp/arXiv.org}}
\def\mytheorem#1{\csname proclaim\endcsname{Theorem #1}}
\def\mytheoremwithtitle#1#2{\csname proclaim\endcsname{Theorem #1#2}}
\def\mythetheorem#1{\blue{#1}\immediate\special{ps:
     ShrHPSdict begin /ShrBORDERthickness 0 def}}
\def\mythetheoremchapter#1#2{\blue{#1}\immediate\special{ps:
     ShrHPSdict begin /ShrBORDERthickness 0 def}}
\def\mylemma#1{\csname proclaim\endcsname{Lemma #1}}
\def\mylemmawithtitle#1#2{\csname proclaim\endcsname{Ëåììà #1#2}}
\def\mythelemma#1{\blue{#1}\immediate\special{ps:
     ShrHPSdict begin /ShrBORDERthickness 0 def}}

\def\myproposition#1{\csname proclaim\endcsname{Proposition #1}}
\def\mypropositionwithtitle#1#2{\csname proclaim\endcsname{Proposition #1#2}}

\def\mycorollary#1{\csname proclaim\endcsname{Corollary #1}}

\def\myexercise#1{\csname proclaim\endcsname{Exercise #1}}
\def\mytheexercise#1{\blue{#1}\immediate\special{ps:
     ShrHPSdict begin /ShrBORDERthickness 0 def}}

\def\mydefinition#1{\definition{Definition #1}}
\def\mythedefinition#1{\blue{#1}\immediate\special{ps:
     ShrHPSdict begin /ShrBORDERthickness 0 def}}
\def\mythedefinitionchapter#1#2{\blue{#1}\immediate\special{ps:
     ShrHPSdict begin /ShrBORDERthickness 0 def}}
\def\myaxiom#1{\csname proclaim\endcsname{\red{Axiom #1}}}
\def\myaxiomwithtitle#1#2{\csname proclaim\endcsname{\red{Axiom #1}#2}}

\def\SectionNum#1#2{\S
     \,#1.}
%---------------------------------------------------------------
\pagewidth{10cm}
\pageheight{15.4cm}
%---------------------------------------------------------------
\fontdimen3\tenrm=3pt
\fontdimen4\tenrm=0.7pt

%---------------------------------------------------------------
\document

\chapternum=1
\vbox to\vsize{\centerline{\etbf MINISTRY OF EDUCATION AND SCIENCE}
\centerline{\etbf OF THE RUSSIAN FEDERATION}
\medskip 
\centerline{\etbf BASHKIR STATE UNIVERSITY}
\vskip 3cm
\centerline{RUSLAN A. SHARIPOV}
\vskip 1.5cm
\centerline{\tvbf COURSE \ OF \ ANALYTICAL \ GEOMETRY}
\vskip 1.3cm
\centerline{\tvrm The textbook: second English edition}
\vskip 5.0cm
\centerline{UFA 2013}

\vss}
\newpage
\vbox to 13.5cm{
{\tencyr UDK 514.123}\par
{\tencyr BBK 22.151}\par
{\tencyr \qquad SH25}\par
\vskip 0.5cm
Referee:\ \ \ \
\vtop{\hsize 7.5cm\noindent The division of Mathematical Analysis 
of Bash\-kir State Pedagogical University ({\tencyr BGPU}) named after 
Miftakhetdin Akmulla, Ufa.}
\vskip 0.5cm
     {\bf Course of analytical geometry}, second English edition: corrected 
with the use of the errata list by Mr\.~\'Eric~Larouche, 
Universit\'e du Qu\'ebec \`a Chicoutimi, Qu\'ebec, Canada. 
In typesetting this second edition the \AmSTeX\ package was used.\par
\vfill
\noindent
\vbox{\hsize=9.5cm\nointerlineskip\baselineskip=10pt
\noindent\kern 1.0cm {\ninecyrbf Sharipov R\.\,A\.}\par
\noindent\kern 0.1cm{\ninecyr Sh25\kern 0.4cm Kurs analitichesko{\ae}
geometrii:}\par
\noindent
\kern 1.0cm\vbox{\hsize=8.1cm\noindent\ninecyr Uchebnoe posobie
/ R\. A\. Sharipov. --- Ufa: RIC Bash\-GU, 2010. --- 228 s.}\par
\vskip 0.1cm
\noindent\kern 1.25cm ISBN 978-5-7477-2574-4\par
\vskip 0.2cm
\noindent\kern 0.8cm\hfill\vbox{\hsize=8.5cm\ninecyr Uchebnoe posobie
po kursu analitichesko{\ae} geometrii adresovano studentam matematikam,
fizikam, a takzhe studentam injenerno-tehnicheskih, tehnologicheskih 
i\linebreak inyh
spetsial\wnsoft noste{\ae}, dlya kotoryh gosudarstvennye 
ob\-razovatel\wnsoft nye standarty predusmatrivayut izuchenie dannogo 
predmeta.\par}\par
\hfill\vbox{\hsize 2.5cm\noindent\ninecyr  UDK 514.123
\newline BBK 22.151\par}
}
\vskip 1cm
\line{ISBN 978-5-7477-2574-4\hss\copyright\ Sharipov R.A., 2010}
\line{English Translation\hss\copyright\ Sharipov R.A., 2011}
\line{Second Edition\hss\copyright\ Sharipov R.A., 2013}\par
\vfill}
\newpage
\ \bigskip\medskip
\leftheadtext{CONTENTS.}
\rightheadtext{CONTENTS.}
\setfirstpage
\centerline{\bf CONTENTS.}
\bigskip
\vskip 0pt plus 1pt minus 1pt
\line{CONTENTS.\ \leaderfill\ 3.}
\medskip
\line{PREFACE.\ \leaderfill\ \myhref{\diez pg7}{7}.}
\medskip
\line{CHAPTER~\uppercase\expandafter{\romannumeral 1}.
VECTOR ALGEBRA.\ \leaderfill\ 
\myhref{\diez s1pg9}{9}.}
\medskip\baselineskip=12pt plus 0.05pt minus 0.05pt
\line{\ \S\,1. Three-dimensional Euclidean space.
Acsiomatics\hss}
\line{\ \qquad and visual evidence.\ \leaderfill\ 
\myhref{\diez s1pg9}{9}.}
\line{\ \S\,2. Geometric vectors. Vectors bound to 
points.\ \leaderfill\ \myhref{\diez s2pg11}{11}.}
\line{\ \S\,3. Equality of vectors.\ \leaderfill\ 
\myhref{\diez s3pg13}{13}.}
\line{\ \S\,4. The concept of a free vector.\ \leaderfill\ 
\myhref{\diez s4pg14}{14}.}
\line{\ \S\,5. Vector addition.\ \leaderfill\ 
\myhref{\diez s5pg16}{16}.}
\line{\ \S\,6. Multiplication of a vector by a number.\ \leaderfill\ 
\myhref{\diez s6pg18}{18}.}
\line{\ \S\,7. Properties of the algebraic operations with vectors.\ 
\leaderfill\ \myhref{\diez s7pg21}{21}.}
\line{\ \S\,8. Vectorial expressions and their transformations.\ 
\leaderfill\ \myhref{\diez s8pg28}{28}.}
\line{\ \S\,9. Linear combinations. Triviality, non-triviality,\hss}
\line{\ \qquad and vanishing.\ \leaderfill\ 
\myhref{\diez s9pg32}{32}.}
\line{\S\,10. Linear dependence and linear independence.\ \leaderfill\ 
\myhref{\diez s10pg34}{34}.} 
\line{\S\,11. Properties of the linear dependence.\ \leaderfill\ 
\myhref{\diez s11pg36}{36}.} 
\line{\S\,12. Linear dependence for $n=1$.\ \leaderfill\ 
\myhref{\diez s12pg37}{37}.} 
\line{\S\,13. Linear dependence for $n=2$. Collinearity\hss}
\line{\ \qquad of vectors.\ \leaderfill\ 
\myhref{\diez s13pg38}{38}.} 
\line{\S\,14. Linear dependence for $n=3$. Coplanartity\hss}
\line{\ \qquad of vectors.\ \leaderfill\ 
\myhref{\diez s14pg40}{40}.} 
\line{\S\,15. Linear dependence for $n\geqslant 4$.\ 
\leaderfill\ \myhref{\diez s15pg42}{42}.} 
\line{\S\,16. Bases on a line.\ \leaderfill\ 
\myhref{\diez s16pg45}{45}.} 
\line{\S\,17. Bases on a plane.\ \leaderfill\ 
\myhref{\diez s17pg46}{46}.} 
\line{\S\,18. Bases in the space.\ \leaderfill\ 
\myhref{\diez s18pg48}{48}.} 
\line{\S\,19. Uniqueness of the expansion \pagebreak of a vector 
in a basis.\ \leaderfill\ \myhref{\diez s19pg50}{50}.} 
\line{\S\,20. Index setting convention.\ \leaderfill\ 
\myhref{\diez s20pg51}{51}.} 
\line{\S\,21. Usage of the coordinates of vectors.\ \leaderfill\ 
\myhref{\diez s21pg52}{52}.} 
\line{\S\,22. Changing a basis. Transition formulas\hss}
\line{\ \qquad and transition matrices.\ \leaderfill\ 
\myhref{\diez s22pg53}{53}.} 
\line{\S\,23. Some information on transition matrices.\ \leaderfill\ 
\myhref{\diez s23pg57}{57}.} 
\line{\S\,24. Index setting in sums.\ \leaderfill\ 
\myhref{\diez s24pg59}{59}.}
\line{\S\,25. Transformation of the coordinates of vectors\hss}
\line{\ \qquad under a change of a basis.\ \leaderfill\ 
\myhref{\diez s25pg63}{63}.}
\line{\S\,26. Scalar product.\ \leaderfill\ 
\myhref{\diez s26pg65}{65}.}
\line{\S\,27. Orthogonal projection onto a line.\ \leaderfill\ 
\myhref{\diez s27pg67}{67}.}
\line{\S\,28. Properties of the scalar product.\ \leaderfill\ 
\myhref{\diez s28pg73}{73}.}
\line{\S\,29. Calculation of the scalar product through \hss}
\line{\ \qquad  the coordinates of vectors in a skew-angular basis.\ 
\leaderfill\ \myhref{\diez s29pg75}{75}.}
\line{\S\,30. Symmetry of the Gram matrix.\ \leaderfill\ 
\myhref{\diez s30pg79}{79}.}
\line{\S\,31. Orthonormal basis.\ \leaderfill\ 
\myhref{\diez s31pg80}{80}.}
\line{\S\,32. The Gram matrix of an orthonormal basis.\ 
\leaderfill\ \myhref{\diez s32pg81}{81}.}
\line{\S\,33. Calculation of the scalar product through\hss} 
\line{\ \qquad the coordinates of vectors in an orthonormal basis.\ 
\leaderfill\ \myhref{\diez s33pg82}{82}.}
\line{\S\,34. Right and left triples of vectors.\hss} 
\line{\ \qquad The concept of orientation. \ \leaderfill\ 
\myhref{\diez s34pg83}{83}.}
\line{\S\,35. Vector product.\ \leaderfill\ 
\myhref{\diez s35pg84}{84}.}
\line{\S\,36. Orthogonal projection onto a plane.\ 
\leaderfill\ \myhref{\diez s36pg86}{86}.}
\line{\S\,37. Rotation about an axis.\ \leaderfill\ 
\myhref{\diez s37pg88}{88}.}
\line{\S\,38. The relation of the vector product\hss} 
\line{\ \qquad with projections and rotations.\ \leaderfill\ 
\myhref{\diez s38pg91}{91}.}
\line{\S\,39. Properties of the vector product.\ \leaderfill\ 
\myhref{\diez s39pg92}{92}.}
\line{\S\,40. Structural constants of the vector product.\ \leaderfill\ 
\myhref{\diez s40pg95}{95}.}
\line{\S\,41. Calculation of the vector product through\hss} 
\line{\ \qquad the coordinates of vectors in a skew-angular basis.\  
\leaderfill\ \myhref{\diez s41pg96}{96}.}
\line{\S\,42. Structural constants of the vector product\hss}
\line{\ \qquad  in an orthonormal basis.\ 
\leaderfill\ \myhref{\diez s42pg97}{97}.}
\line{\S\,43. Levi-Civita symbol.\ \leaderfill\ 
\myhref{\diez s43pg99}{99}.}
\line{\S\,44. Calculation of the vector product through the co-\hss} 
\line{\ \qquad ordinates of vectors in an orthonormal basis.\ 
\leaderfill\ \myhref{\diez s44pg102}{102}.}
\line{\S\,45. Mixed product.\ \leaderfill\ 
\myhref{\diez s45pg104}{104}.}
\line{\S\,46. Calculation of the mixed product through the co-\hss}
\line{\ \qquad ordinates of vectors in an orthonormal basis.\ 
\leaderfill\ \myhref{\diez s46pg105}{105}.}
\line{\S\,47. Properties of the mixed product.\ \leaderfill\ 
\myhref{\diez s47pg108}{108}.}
\line{\S\,48. The concept of the oriented volume.\ \leaderfill\ 
\myhref{\diez s48pg111}{111}.}
\line{\S\,49. Structural constants of the mixed product.\ \leaderfill\ 
\myhref{\diez s49pg113}{113}.}
\line{\S\,50. Calculation of the mixed product through the co-\hss}
\line{\ \qquad ordinates of vectors in a skew-angular basis.\ 
\leaderfill\ \myhref{\diez s50pg115}{115}.}
\line{\S\,51. The relation of structural constants of the vectorial\hss} 
\line{\ \qquad and mixed products.\ \leaderfill\ 
\myhref{\diez s51pg116}{116}.}
\line{\S\,52. Effectivization of the formulas for calculating\hss} 
\line{\ \qquad vectorial and mixed products.\ \leaderfill\ 
\myhref{\diez s52pg121}{121}.}
\line{\S\,53. Orientation of the space.\ \leaderfill\ 
\myhref{\diez s53pg124}{124}.}
\line{\S\,54. Contraction formulas.\ \leaderfill\ 
\myhref{\diez s54pg125}{125}.}
\line{\S\,55. The triple product expansion formula and\hss} 
\line{\ \qquad the Jacobi identity.\ \leaderfill\ 
\myhref{\diez s55pg131}{131}.}
\line{\S\,56. The product of two mixed products.\ \leaderfill\ 
\myhref{\diez s56pg134}{134}.}
\bigskip
\line{CHAPTER~\uppercase\expandafter{\romannumeral 2}.
GEOMETRY OF LINES\hss}
\line{\qquad AND SURFACES.\ \ \leaderfill\ 
\myhref{\diez s1pg139}{139}.}
\medskip
\line{\ \S\,1. Cartesian coordinate systems.\ \leaderfill\ 
\myhref{\diez s1pg139}{139}.}
\line{\ \S\,2. Equations of lines and surfaces.\ \leaderfill\ 
\myhref{\diez s2pg141}{141}.}
\line{\ \S\,3. A straight line on a plane.\ \leaderfill\ 
\myhref{\diez s3pg142}{142}.}
\line{\ \S\,4. A plane in the space.\ \leaderfill\ 
\myhref{\diez s4pg148}{148}.}
\line{\ \S\,5. A straight line in the space.\ \leaderfill\ 
\myhref{\diez s5pg154}{154}.}
\line{\ \S\,6. Ellipse. Canonical equation of an ellipse.\ 
\leaderfill\ \myhref{\diez s6pg160}{160}.}
\line{\ \S\,7. The eccentricity and directrices of an ellipse.\hss} 
\line{\ \qquad The property of directrices.\ \leaderfill\ 
\myhref{\diez s7pg165}{165}.}
\line{\ \S\,8. The equation of a tangent line to an 
ellipse.\ \leaderfill\ 
\myhref{\diez s8pg167}{167}.}
\line{\ \S\,9. Focal property of an ellipse.\ \leaderfill\ 
\myhref{\diez s9pg170}{170}.}
\line{\S\,10. Hyperbola. Canonical equation of 
a hyperbola.\ \leaderfill\ 
\myhref{\diez s10pg172}{172}.}
\line{\S\,11. The eccentricity and directrices of a hyperbola.\hss} 
\line{\ \qquad The property of directrices.\ \leaderfill\ 
\myhref{\diez s11pg179}{179}.}
\line{\S\,12. The equation of a tangent line to 
a hyperbola.\ \leaderfill\ 
\myhref{\diez s12pg181}{181}.}
\line{\S\,13. Focal property of a hyperbola.\ \leaderfill\ 
\myhref{\diez s13pg184}{184}.}
\line{\S\,14. Asymptotes of a hyperbola.\ \leaderfill\ 
\myhref{\diez s14pg186}{186}.}
\line{\S\,15. Parabola. Canonical equation of 
a parabola.\ \leaderfill\ 
\myhref{\diez s15pg187}{187}.}
\line{\S\,16. The eccentricity of a parabola.\ \leaderfill\ 
\myhref{\diez s16pg190}{190}.}
\line{\S\,17. The equation of a tangent line to 
a parabola.\ \leaderfill\ 
\myhref{\diez s17pg190}{190}.}
\line{\S\,18. Focal property of a parabola.\ \leaderfill\ 
\myhref{\diez s18pg193}{193}.}
\line{\S\,19. The scale of eccentricities.\ \leaderfill\ 
\myhref{\diez s19pg194}{194}.}
\line{\S\,20. Changing a coordinate system.\ \leaderfill\ 
\myhref{\diez s20pg195}{195}.}
\line{\S\,21. Transformation of the coordinates of a point\hss}
\line{\ \qquad under a change of a coordinate system.\ \leaderfill\ 
\myhref{\diez s21pg196}{196}.}
\line{\S\,22. Rotation of a rectangular coordinate system\hss}
\line{\ \qquad on a plane. The rotation matrix.\ \leaderfill\ 
\myhref{\diez s22pg197}{197}.}
\line{\S\,23. Curves of the second order.\ \leaderfill\ 
\myhref{\diez s23pg199}{199}.}
\line{\S\,24. Classification of curves of the second order.\ \leaderfill\ 
\myhref{\diez s24pg200}{200}.}
\line{\S\,25. Surfaces of the second order.\ \leaderfill\ 
\myhref{\diez s25pg206}{206}.}
\line{\S\,26. Classification of surfaces of the second order.\ \leaderfill\ 
\myhref{\diez s26pg207}{207}.}
\bigskip 
\line{REFERENCES.\ \ \leaderfill\ \myhref{\diez pg216}{216}.}
\bigskip 
\line{CONTACT INFORMATION.\ \ \leaderfill\ \myhref{\diez pg217}{217}.}
\bigskip
\line{APPENDIX.\ \ \leaderfill\ \myhref{\diez pg218}{218}.}
\newpage
\newtoks\truehead
\truehead=\headline
\leftheadtext{PREFACE.}
\rightheadtext{PREFACE.}
\setfirstpage
%\headline{\hfill}
\centerline{\bf PREFACE.
}
\baselineskip=12pt
\medskip
\bigskip
The elementary geometry, which is learned in school, deals with
basic concepts such as {\it a point, a straight line, a segment}.
They are used to compose more complicated concepts: {\it a polygonal 
line, a polygon, a polyhedron}. Some curvilinear forms are also
considered: {\it a circle, a cylinder, a cone, a sphere, a ball}.\par
     The analytical geometry basically deals with the same geometric
objects as the elementary geometry does. The difference is in a method
for studying these objects. The elementary geometry relies on visual
impressions and formulate the properties of geometric objects in its
axioms. From these axioms various theorems are derived, whose 
statements in most cases are also revealed in visual impressions. 
The analytical geometry is more inclined to a numeric description 
of geometric objects and their properties.\par
     The transition from a geometric description to a numeric
description becomes possible due to coordinate systems. Each coordinate
system associate some groups of numbers with geometric points, these
groups of numbers are called coordinates of points. The idea of using
coordinates in geometry belongs French mathematician Rene Descartes.
Simplest coordinate systems suggested by him at present time are 
called {\it Cartesian coordinate systems}.\par
     The construction of Cartesian coordinates particularly and the 
analytical geometry in general are based on the concept of {\it a 
vector}. The branch of analytical geometry studying vectors is called 
{\it the vector algebra}. The vector algebra constitutes the first 
chapter of this book. The second chapter explains {\it the theory
of straight lines and planes} and {\it the theory of curves of the
second order}. In the third chapter {\it the theory of surfaces 
of the second order} is explained in brief.\par
     The book is based on lectures given by the author during several
years in Bashkir State University. It was planned as the first book in 
a series of three books. However, it happens that the second and the 
third  books in this series were written and published before the first 
book. These are 
\roster
\rosteritemwd=-5pt
\item"--" {\tencyr\char '074}Course of linear algebra and 
multidimensional
geometry{\tencyr\char '076} \mycite{1};
\item"--" {\tencyr\char '074}Course of differential geometry{\tencyr\char 
'076} \mycite{2}. 
\endroster
Along with the above books, the following books were written:
\roster
\rosteritemwd=-5pt
\item"--" {\tencyr\char '074}Representations of finite 
group{\tencyr\char '076} \mycite{3};
\item"--"\parshape 2 0cm 9.5cm 0.2cm 8cm {\tencyr\char '074}Classical 
electrodynamics and theory of relativity{\tencyr\char '076} \mycite{4};
\item"--" {\tencyr\char '074}Quick introduction to tensor 
analysis{\tencyr\char '076} \mycite{5}.
\item"--" {\tencyr\char '074}Foundations of geometry for university
students and high school students{\tencyr\char '076} \mycite{6}.
\endroster
The book \mycite{3} can be considered as a continuation of the book 
\mycite{1} which illustrates the application of linear algebra to another
branch of mathematics, namely to the theory of groups. The book \mycite{4} 
can be considered as a continuation of the book \mycite{2}. It illustrates
the application of differential geometry to physics. The book \mycite{5} 
is a brief version of the book \mycite{2}. As for the book \mycite{6},
by its subject it should precede this book. It could br recommended to the 
reader for deeper logical understanding of the elementary geometry.\par
     I am grateful to Prof\.~R\.~R\.~Gadylshin and Prof\.~D\.~I\.~Borisov
for reading and refereeing the manuscript of this book and for valuable 
advices. I am grateful to Mr\.~\'Eric~Larouche who has read the book and 
provided me the errata list for the first English edition of this book.
\bigskip\bigskip
\line{\qquad June, 2013.\hss R.~A.~Sharipov.\qquad}
\newpage
%---------------------------------------------------------------
\setfirstpage
\topmatter
\title\chapter{1}
Vector Algebra.
\endtitle
\endtopmatter
\leftheadtext{Chapter~\uppercase\expandafter{\romannumeral 1}.
VECTOR ALGEBRA.}
\document
\headline=\truehead
\head
\SectionNum{1}{9} Three-dimensional Euclidean space.
Acsiomatics and visual evidence.
\endhead
\rightheadtext{\S\,1. Three-dimensional Euclidean space.}
Like the elementary geometry explained in the book \mycite{6}, 
the analytical geometry in this book is a geometry of three-dimensional 
space $\Bbb E$. We use the symbol $\Bbb E$ for to denote the space
that we observe in our everyday life. Despite being seemingly simple,
even the empty space $\Bbb E$ possesses a rich variety of properties.
These properties reveal through the properties of various geometric
forms which are comprised in this space or potentially can be 
comprised in it.\par
      Systematic study of the geometric forms in the space $\Bbb E$ was
initiated by ancient Greek mathematicians. It was Euclid who succeeded 
the most in this. He has formulated the basic properties of the space
$\Bbb E$ in five postulates, from which he derived all other properties
of $\Bbb E$. At the present time his postulates are called axioms. 
On the basis of modern requirements to the rigor of mathematical reasoning
the list of Euclid's axioms was enlarged from five to twenty. These
twenty axioms can be found in the book \mycite{6}. In favor of Euclid the 
space that we observe in our everyday life is denoted by the symbol
$\Bbb E$ and is called the three-dimensional Euclidean space.\par
     The three-dimensional Euclidean point space $\Bbb E$ consists of 
points. All geometric forms in it also consist of points. 
subsets of the space $\Bbb E$. Among subsets of the space $\Bbb E$ 
straight lines and planes (see Fig\.~1.2) play an especial role. They are 
used in the statements of the first eleven Euclid's axioms. On the base
of these axioms the concept of a segment (see Fig\.~1.2) is introduced.
The concept of a segment is used in the statement of the twelfth axiom. 
\par
     The first twelve of Euclid's axioms appear to be sufficient to
define the concept of a ray and the concept of an angle between two
rays outgoing from the same point. \vadjust{\vskip 5pt\hbox to 0pt{\kern 
-10pt\includegraphics{angemeng01.eps}\hss}\vskip 125pt}The
concepts of a segment and an angle along with the concepts of a straight
line and a plane appear to be sufficient in order to formulate the 
remaining eight Euclid's axioms and to build the elementary geometry in
whole.\par
     Even the above survey of the book \mycite{6}, which is very short, 
shows that building the elementary geometry in an axiomatic way on the
basis of Euclid's axioms is a time-consuming and laborious work. However,
the reader who is familiar with the elementary geometry from his school 
curriculum easily notes that proof of theorems in school textbooks are
more simple than those in \mycite{6}. The matter is that proofs in 
school textbooks are not proofs in the strict mathematical sense. Analyzing
them carefully, one can find in these proofs the usage of some non-proved 
propositions which are visually obvious from a supplied drawing since we 
have a rich experience of living within the space $\Bbb E$. Such proofs
can be transformed to strict mathematical proofs by filling omissions, 
i\.\,e\. by proving visually obvious propositions used in them.\par
     Unlike \mycite{6}, in this book I do not load the reader by absolutely
strict proofs of geometric results. For geometric definitions, constructions,
and theorems the strictness level is used which is close to that of school 
textbooks and which admits drawings and visually obvious facts as arguments.
Whenever it is possible I refer the reader to strict statements and proofs 
in \mycite{6}. As far as the analytical content is concerned, i\.\,e\. in 
equations, in formulas, and in calculations the strictness level is applied 
which is habitual in mathematics without any deviations.\par
\head
\SectionNum{2}{11} Geometric vectors. Vectors bound to points.
\endhead
\rightheadtext{\S\,2. Geometric vectors.}
\mydefinition{2.1} A geometric vectors $\overrightarrow{AB\,\,}\!$ is
a straight line segment in which the direction from the point $A$ to
the point $B$ is specified. The point $A$ is called the {\it initial point\/} 
of the vector $\overrightarrow{AB\,\,}\!\!$, while the point $B$ is called 
its {\it terminal point}.
\enddefinition
\parshape 3 0cm 10cm 0cm 10cm 5cm 5cm 
     The direction of the vector $\overrightarrow{AB\,\,}\!$ in drawing is marked
by an arrow (see Fig\.~2.1). \vadjust{\vskip 5pt\hbox to 0pt{\kern 5pt
\includegraphics{angemeng02.eps}\hss}\vskip -5pt}For this reason vectors
sometimes are called {\it directed segments}.\par
\parshape 1 5cm 5cm 
Each segment $[AB]$ is associated with two different vectors: 
$\overrightarrow{AB\,\,}\!$ and $\overrightarrow{BA\,\,}\!$. The vector 
$\overrightarrow{BA\,\,}\!$ is usually called the {\it opposite vector } for the
vector $\overrightarrow{AB\,\,}\!$.\par
\parshape 2 5cm 5cm 0cm 10cm
       Note that the arrow sign on the vector $\overrightarrow{AB\,\,}\!$ and
bold dots at the ends of the segment $[AB]$ are merely symbolic signs used to 
make the drawing more clear. When considered as sets of points the vector 
$\overrightarrow{AB\,\,}\!$ and the segment $[AB]$ do coincide.\par
     A direction on a segment, which makes it a vector, can mean different 
things in different situations. For instance, drawing a vector 
$\overrightarrow{AB\,\,}\!$ on a geographic map, we usually mark the 
displacement of some object from the point $A$ to the point $B$. However, 
if it is a weather map, the same vector $\overrightarrow{AB\,\,}\!$ can mean 
the wind direction and its speed at the point $A$. In the first case the 
length of the vector $\overrightarrow{AB\,\,}\!$ is proportional to the 
distance between the points $A$ and $B$. In the second case the length of 
$\overrightarrow{AB\,\,}\!$ is proportional to the wind speed at the point 
$A$.\par
     There is one more difference in the above two examples. In the first
case the vector $\overrightarrow{AB\,\,}\!$ is bound to the points $A$ and
$B$ by its meaning. In the second case the vector $\overrightarrow{AB\,\,}\!$
is bound to the point $A$ only. The fact that its arrowhead end is at the
point $B$ is a pure coincidence depending on the scale we used for 
translating the wind speed into the length units on the map. According to 
what was said, geometric vectors are subdivided into two types:
\roster
\item"1)" purely geometric;
\item"2)" conditionally geometric.
\endroster\par
     Only displacement vectors belong to the first type; they actually bind 
some two points of the space $\Bbb E$. The lengths of these vectors are always
measured in length units: centimeters, meters, inches, feets etc.\par
     Vectors of the second type are more various. These are velocity vectors,
acceleration vectors, and force vectors in mechanics; intensity vectors for
electric and magnetic fields, magnetization vectors in magnetic materials 
and media, temperature gradients in non-homogeneously heated objects et al.
Vectors of the second type have a geometric direction and they are bound to
some point of the space $\Bbb E$, but they have not a geometric length. 
Their lengths can be translated to geometric lengths only upon choosing
some scaling factor.\par 
     Zero vectors or null vectors hold a special position among geometric 
vectors. They are defined as follows. 
\mydefinition{2.2} A geometric vector of the space $\Bbb E$ whose initial
and terminal points do coincide with each other is called a {\it zero 
vector\/} or a {\it null vector}.
\enddefinition
     A geometric null vector can be either a purely geometric vector or a
conditionally geometric vector depending on its nature.
\head
\SectionNum{3}{13} Equality of vectors.
\endhead
\rightheadtext{\S\,3. Equality of vectors.}
\mydefinition{3.1} Two geometric vectors $\overrightarrow{AB\,\,}\!$ and
$\overrightarrow{CD\,\,}\!$ are called {\it equal\/} if they are equal in length
and if they are {\it codirected}, i\.\,e\. $|AB|=|CD|$ and $\overrightarrow{AB\,\,}
\!\!\upuparrows\overrightarrow{CD\,\,}\!$.
\enddefinition
     The vectors $\overrightarrow{AB\,\,}\!$ and $\overrightarrow{CD\,\,}\!$ are
said to be codirected if they lie on a same line as shown in Fig\.~3.1 of if they
lie on parallel lines as shown in Fig\.~3.2. 
\vadjust{\vskip 5pt\hbox to 0pt{\kern 5pt\includegraphics{angemeng03.eps}\hss}\vskip 115pt}In both cases they should be pointing in the same 
direction. Codirectedness of geometric vectors and their equality are that very 
visually obvious properties which require substantial efforts in order to derive 
them from Euclid's axioms (see \mycite{6}). Here I urge the reader not to focus
on the lack of rigor in statements, but believe his own geometric intuition.\par
     Zero geometric vectors constitute a special case since they do not fix any 
direction in the space.
\mydefinition{3.2} All null vectors are assumed to be codirected to each other
and each null vector is assumed to be codirected to any nonzero vector.
\enddefinition
     The length of all null vectors is zero. However, depending on the physical 
nature of a vector, its zero length is complemented with a measure unit. In the 
case of zero force it is zero newtons, in the case of zero velocity it is zero 
meters per second. For this reason, testing the equality of any two zero vectors, 
one should take into account their physical nature.
\mydefinition{3.3} All null vectors of the same physical nature are assumed to
be equal to each other and any nonzero vector is assumed to be not equal to any
null vector. 
\enddefinition
     Testing the equality of nonzero vectors by means of the 
definition~\mythedefinition{3.1}, one should take into account its physical
nature. The equality $|AB|=|CD|$ in this definition assumes not only the equality
of numeric values of the lengths of $\overrightarrow{AB\,\,}\!$ and
$\overrightarrow{CD\,\,}\!$, but assumes the coincidence of their measure units
as well.\par
     {\bf A remark}. Vectors are geometric forms, i\.\,e\. they are sets of points
in the space $\Bbb E$. However, the equality of two vectors introduced in the
definition~\mythedefinition{3.1} differs from the equality of sets.\par
\head
\SectionNum{4}{14} The concept of a free vector.
\endhead
\rightheadtext{\S\,4. The concept of a free vector.}
     Defining the equality of vectors, it is convenient to use parallel translations.
Each parallel translation is a special transformation of the space $p\!:\,\Bbb E\to
\Bbb E$ under which any straight line is mapped onto itself or onto a parallel line
and any plane is mapped onto itself or onto a parallel plane. When applied to 
vectors, parallel translation preserve their length and their direction, i\.\,e\.
they map each vector onto a vector equal to it, but usually being in a different place 
in the space. The number of parallel translations is infinite. As appears, the
parallel translations are so numerous that they can be used for testing the equality
of vectors. 
\mydefinition{4.1} A geometric vector $\overrightarrow{CD\,\,}\!$ is called equal
to a geometric vector $\overrightarrow{AB\,\,}\!$ if there is a parallel translation
$p\!:\,\Bbb E\to\Bbb E$ that maps the vector $\overrightarrow{AB\,\,}\!$ onto the
vector $\overrightarrow{CD\,\,}\!$, i\.\,e\. such that $p(A)=C$ and $p(B)=D$.
\enddefinition
     The definition~\mythedefinition{4.1} is equivalent to the 
definition~\mythedefinition{3.1}. I do not prove this equivalence, relying on its
visual evidence and assuming the reader to be familiar with parallel translations
from the school curriculum. A more meticulous reader can see the theorems 8.4 and
9.1 in Chapter~\uppercase\expandafter{\romannumeral 6} of the book \mycite{6}.\par
\mytheorem{4.1} For any two points $A$ and $C$ in the space $\Bbb E$ there is
exactly one parallel translation $p\!:\,\Bbb E\to\Bbb E$ mapping the point $A$ onto
the point $C$, i\.\,e\. such that $p(A)=C$. 
\endproclaim
     The theorem~\mythetheorem{4.1} is a visually obvious fact. On the other hand
it coincides with the theorem 9.3 from Chapter~\uppercase\expandafter{\romannumeral 6} 
in the book \mycite{6}, where it is proved. For these two reasons we exploit the
theorem~\mythetheorem{4.1} without proving it in this book.\par
     Lei's apply the theorem~\mythetheorem{4.1} to some geometric vector
$\overrightarrow{AB\,\,}\!$. Let $C$ be an arbitrary point of the space $\Bbb E$ 
and let $p$ be a parallel translation taking the point $A$ to the point $C$. The
existence and uniqueness of such a parallel translation are asserted by the 
theorem~\mythetheorem{4.1}. Let's define the point $D$ by means of the formula
$D=p(B)$. Then, according to the definition~\mythedefinition{4.1}, we have
$$
\overrightarrow{AB\,\,}\!=\overrightarrow{CD\,\,}\!.
$$
These considerations show that each geometric vector $\overrightarrow{AB\,\,}\!$ 
has a copy equal to it and attached to an arbitrary point $C\in\Bbb E$. 
In the other words, by means of parallel translations each geometric vector
$\overrightarrow{AB\,\,}\!$ can be replicated up to an infinite set of vectors
equal to each other and attached to all points of the space $\Bbb E$. 
\mydefinition{4.2} A {\it free vector}\/ is an infinite collection of geometric 
vectors which are equal to each other and whose initial points are at all points 
of the space $\Bbb E$. Each geometric vector in this infinite collection is called
a {\it geometric realization} of a given free vector. 
\enddefinition
     Free vectors can be composed of purely geometric vectors or of conditionally 
geometric vectors as well. For this reason one can consider free vectors of 
various physical nature.\par
     In drawings free vectors are usually presented by a single geometric
realization or by several geometric realizations if needed.\par
     Geometric vectors are usually denoted by two capital letters: 
$\overrightarrow{AB\,\,}\!$, $\overrightarrow{CD\,\,}\!$,
$\overrightarrow{EF\,\,}\!$ etc. Free vectors are denoted by single lowercase 
letters: $\vec{\bold a}$,  $\vec{\bold b}$, $\vec{\bold c}$ etc. Arrows over 
these letters are often omitted since it is usually clear from the context
that vectors are considered. Below in this book I will not use arrows in denoting
free vectors. However, I will use boldface letters for them. In many other books,
but not in my book \mycite{1}, this restriction is also removed. 
\head
\SectionNum{5}{16} Vector addition.
\endhead
\rightheadtext{\S\,5. Vector addition.}
     Assume that two free vectors $\bold a$ and $\bold b$ are given. Let's choose
some arbitrary point $A$ and consider the geometric realization of the vector
$\bold a$ with the initial point $A$. Let's denote through $B$ the terminal point
of this geometric realization. As a result we get $\bold a=\overrightarrow{AB\,\,}\!$. 
Then we consider the geometric realization of the vector $\bold b$ with initial
point $B$ and denote through $C$ its terminal point. This yields 
$\bold b=\overrightarrow{BC\,\,}\!$. 
\mydefinition{5.1} The geometric vector $\overrightarrow{AC\,\,}\!$ connecting
the initial point of the vector $\overrightarrow{AB\,\,}\!$ with the terminal
point of the vector $\overrightarrow{BC\,\,}\!$ is called the {\it sum\/} of the 
vectors $\overrightarrow{AB\,\,}\!$ and $\overrightarrow{BC\,\,}\!$:
$$
\hskip -2em
\overrightarrow{AC\,\,}\!=\overrightarrow{AB\,\,}\!
+\overrightarrow{BC\,\,}\!.
\mytag{5.1}
$$ 
\enddefinition
The vector $\overrightarrow{AC\,\,}$ constructed by means of the vectors 
$\bold a=\overrightarrow{AB\,\,}\!$ and $\bold b=\overrightarrow{BC\,\,}\!$ 
can be replicated up to a free vector $\bold c$ by parallel translations to
all points of the space $\Bbb E$. Such a vector $\bold c$ is naturally
called the sum of the free vectors $\bold a$ and $\bold b$. For this vector
we write $\bold c=\bold a+\bold b$. The correctness of such a definition is 
guaranteed by the following lemma.
\mylemma{5.1} The sum $\bold c=\bold a+\bold b=\overrightarrow{AC\,\,}\!$ of 
two free vectors\linebreak $\bold a=\overrightarrow{AB\,\,}\!$ and\/ 
$\bold b=\overrightarrow{BC\,\,}\!$ expressed by the formula \mythetag{5.1} does 
not depend on the choice of a point $A$ at which the geometric realization 
$\overrightarrow{AB\,\,}\!$ of the vector $\bold a$ begins. 
\endproclaim
\demo{Proof} In addition to $A$, let's choose another initial point $E$. Then
\vadjust{\vskip 5pt\hbox to 0pt{\kern 0pt
\includegraphics{angemeng04.eps}\hss}\vskip 115pt}in the above 
construction of the sum $\bold a+\bold b$ the vector $\bold a$ has two geometric
realizations $\overrightarrow{AB\,\,}\!$ and $\overrightarrow{EF\,\,}\!$. 
The vector $\bold b$ also has two geometric realizations 
$\overrightarrow{BC\,\,}\!$ and $\overrightarrow{FG\,\,}\!$ (see Fig\.~5.1). 
Then 
$$
\xalignat 2
&\hskip -2em
\overrightarrow{AB\,\,}\!=\overrightarrow{EF\,\,}\!,
&&\overrightarrow{BC\,\,}\!=\overrightarrow{FG\,\,}\!.
\mytag{5.2}
\endxalignat 
$$
Instead of \mythetag{5.1} now we have two equalities
$$
\xalignat 2
&\hskip -2em
\overrightarrow{AC\,\,}\!=\overrightarrow{AB\,\,}\!
+\overrightarrow{BC\,\,}\!,
&&\overrightarrow{EG\,\,}\!=\overrightarrow{EF\,\,}\!
+\overrightarrow{FG\,\,}\!.\quad
\mytag{5.3}
\endxalignat 
$$\par
Let's denote through $p$ a parallel translation that maps the point $A$ to the
point $E$, i\.\,e\. such that $p(A)=E$. Due to the theorem~\mythetheorem{4.1} such
a parallel translation does exist and it is unique. From $p(A)=E$ and from the 
first equality \mythetag{5.2}, applying the definition~\mythedefinition{4.1}, we
derive $p(B)=F$. Then from $p(B)=F$ and from the second equality \mythetag{5.2},
again applying the definition~\mythedefinition{4.1}, we get $p(C)=G$. As a result
we have
$$
\xalignat 2
&\hskip -2em
p(A)=E,
&&p(C)=G.
\mytag{5.4}
\endxalignat 
$$
The relationships \mythetag{5.4} mean that the parallel translation $p$ maps the
vector $\overrightarrow{AC\,\,}\!$ to the vector $\overrightarrow{EG\,\,}\!$. Due
to the definition~\mythedefinition{4.1} this fact yields
$\overrightarrow{AC\,\,}\!=\overrightarrow{EG\,\,}\!$. Now from the equalities 
\mythetag{5.3} we derive
$$
\hskip -2em
\overrightarrow{AB\,\,}\!+\overrightarrow{BC\,\,}\!
=\overrightarrow{EF\,\,}\!+\overrightarrow{FG\,\,}\!.
\mytag{5.5}
$$
The equalities \mythetag{5.5} complete the proof of the lemma~\mythelemma{5.1}. 
\qed\enddemo
     The addition rule given by the formula \mythetag{5.1} is called the 
{\it triangle rule}. It is associated with the triangle $ABC$ in Fig\.~5.1.\par
\head
\SectionNum{6}{18} Multiplication of a vector by a number.
\endhead
\rightheadtext{\S\,6. Multiplication of a vector by a number.}
      Let $\bold a$ be some free vector. Let's choose some arbitrary point 
$A$ and consider the geometric realization of the vector $\bold a$ with initial
point $A$. Then we denote through $B$ the terminal point of this geometric 
realization of $\bold a$.
\vadjust{\vskip 5pt\hbox to 0pt{\kern -5pt\includegraphics{angemeng05.eps}\hss}\vskip 75pt}Let $\alpha$ be some number. It can be either 
positive, negative, or zero.\par
     Let $\alpha>0$. In this case we lay a point $C$ onto the line $AB$ so that
the following conditions are fulfilled:
$$
\xalignat 2
&\hskip -2em
\overrightarrow{AC\,\,}\!\!\upuparrows\overrightarrow{AB\,\,}\!,
&&|AC|=|\alpha|\cdot|AB|.
\mytag{6.1}
\endxalignat 
$$
As a result we obtain the drawing which is shown in Fig\.~6.1.\par
     If $\alpha=0$, we lay the point $C$ so that it coincides with the
point $A$. In this case the vector $\overrightarrow{AC\,\,}\!$ appears to be
zero as shown in Fig\.~6.2 and we have the relationship
$$
\hskip -2em
|AC|=|\alpha|\cdot|AB|.
\mytag{6.2}
$$\par
     In the case $\alpha<0$ we lay the point $C$ onto the line $AB$ so that
the following two conditions are fulfilled:
$$
\xalignat 2
&\hskip -2em
\overrightarrow{AC\,\,}\!\!\updownarrows\overrightarrow{AB\,\,}\!,
&&|AC|=|\alpha|\cdot|AB|.
\mytag{6.3}
\endxalignat 
$$
This arrangement of points is shown in Fig\.~6.3.\par
\mydefinition{6.1} In each of the three cases $\alpha>0$, $\alpha=0$, and
$\alpha<0$ the geometric vector $\overrightarrow{AC\,\,}\!$ defined through
the vector $\overrightarrow{AB\,\,}\!$ according to the drawings in Fig\.~6.1, 
in Fig\.~6.2, and in Fig\.~6.3 and according to the formulas \mythetag{6.1}, 
\mythetag{6.2}, and \mythetag{6.3} is called the product of the vector
$\overrightarrow{AB\,\,}\!$ by the number $\alpha$. This fact is expressed 
by the following formula:
$$
\hskip -2em
\overrightarrow{AC\,\,}\!=\alpha\cdot\overrightarrow{AB\,\,}\!.
\mytag{6.4}
$$
\enddefinition
     The case $\bold a=\bold 0$ is not covered by the above drawings in 
Fig\.~6.1, in Fig\.~6.2, and in Fig\.~6.3. In this case the point $B$ coincides 
with the points $A$ and we have $|AB|=0$. In order to provide the equality
$|AC|=|\alpha|\cdot|AB|$ the point $C$ is chosen coinciding with the point $A$. 
Therefore the product of a null vector by an arbitrary number is again a 
null vector.\par 
     The geometric vector $\overrightarrow{AC\,\,}$ constructed with the use
of the vector $\bold a=\overrightarrow{AB\,\,}\!$ and the number $\alpha$ can be
replicated up to a free vector $\bold c$ by means of the parallel translations
to all points of the space $\Bbb E$. Such a free vector $\bold c$ is called
the product of the free vector $\bold a$ by the number $\alpha$. For this vector
we write $\bold c=\alpha\cdot\bold a$. The correctness of this definition of
$\bold c=\alpha\cdot\bold a$ is guaranteed by the following lemma.
\mylemma{6.1} The product $\bold c=\alpha\cdot\bold a=\overrightarrow{AC\,\,}\!$ 
of a free vector $\bold a=\overrightarrow{AB\,\,}\!$ by a number $\alpha$
expressed by the formula \mythetag{6.4} does not depend on the choice of a
point $A$ at which the geometric realization of the vector $\bold a$ is 
built. 
\endproclaim
\demo{Proof} Let's prove the lemma for the case $\bold a\neq\bold 0$ and 
$\alpha>0$. In addition to $A$ we choose another initial point $E$. 
\vadjust{\vskip 5pt\hbox to 0pt{\kern 0pt
\includegraphics{angemeng06.eps}\hss}\vskip 80pt}Then
in the construction
of the product  $\alpha\cdot\bold a$ the vector $\bold a$ gets two geometric
realizations $\overrightarrow{AB\,\,}\!$ and $\overrightarrow{EF\,\,}\!$
(see Fig\.~6.4). Hence we have
$$
\hskip -2em
\overrightarrow{AB\,\,}\!=\overrightarrow{EF\,\,}\!.
\mytag{6.5}
$$\par
     Let's denote through $p$ the parallel translation that maps the point $A$ 
to the point $E$, i\.\,e\. such that $p(A)=E$. Then from the equality \mythetag{6.5},
applying the definition~\mythedefinition{4.1}, we derive $p(B)=F$. The point $C$ 
is placed on the line $AB$ at the distance $|AC|=|\alpha|\cdot |AB|$ from the
point $A$ in the direction of the vector $\overrightarrow{AB\,\,}\!$. Similarly,
the point $G$ is placed on the line $EF$ at thew distance $|EG|=|\alpha|\cdot |EF|$ 
from the point $E$ in the direction of the vector $\overrightarrow{EF\,\,}\!$. 
From the equality \mythetag{6.5} we derive $|AB|=|EF|$. Therefore $|AC|=|\alpha|
\cdot |AB|$ and $|EG|=|\alpha|\cdot |EF|$ mean that $|AC|=|EF|$. Due to $p(A)=E$ 
and $p(B)=F$ the parallel translation $p$ maps the line $AB$ onto the line $EF$. 
It preserves lengths of segments and maps codirected vectors to codirected ones. 
Hence $p(C)=G$. Along with $p(A)=E$ due to the definition~\mythedefinition{4.1} 
the equality $p(C)=G$ yields $\overrightarrow{AC\,\,}\!=\overrightarrow{EG\,\,}\!$, 
i\.\,e\.
$$
\alpha\cdot\overrightarrow{AB\,\,}\!=\alpha\cdot\overrightarrow{EF\,\,}\!.
$$
The lemma~\mythelemma{6.1} is proved for the case $\bold a\neq\bold 0$ and
$\alpha>0$. Its proof for the other cases is left to the reader as an 
exercise.
\qed\enddemo
\myexercise{6.1} Consider the cases $\alpha=0$ and $\alpha<0$ for 
\pagebreak $\bold a\neq
\bold 0$ and consider the case $\bold a=\bold 0$. Prove the 
lemma~\mythelemma{6.1} 
for these cases and provide your proof with drawings analogous to that of 
Fig\,~6.4.
\endproclaim
\head
\SectionNum{7}{21} Properties of the algebraic operations with vectors. 
\endhead
\rightheadtext{\S\,7. Properties of the algebraic operations \dots}
     The addition of vectors and their multiplication by numbers are two
basic algebraic operations with vectors in the three-dimensional Euclidean 
point space $\Bbb E$. Eight basic properties of these two algebraic operations 
with vectors are usually considered. The first four of these eight properties
characterize the operation of addition. The other four characterize the
operation of multiplication of vectors by numbers and its correlation with the 
operation of vector addition.
\mytheorem{7.1} The operation of addition of free vectors and the operation of
their multiplication by numbers possess the following properties:
\roster
\rosteritemwd=0pt
\item"1)" commutativity of addition: $\bold a+\bold b=\bold b+\bold a$;
\item"2)" associativity of addition: $(\bold a+\bold b)+\bold c=
      \bold a+(\bold b+\bold c)$;
\item"3)" the feature of the null vector: $\bold a+\bold 0=\bold a$;
\item"4)" the existence of an opposite vector: for any vector $\bold a$ there 
      is an opposite vector $\bold a'$ such that $\bold a+\bold a'=\bold 0$;
\item"5)" distributivity of multiplication over the addition of vectors:
      $k\cdot (\bold a+\bold b)=k\cdot\bold a+k\cdot\bold b$;
\item"6)" distributivity of multiplication over the addition of numbers:
      $(k+q)\cdot\bold a=k\cdot\bold a+q\cdot\bold a$;
\item"7)" associativity of multiplication by numbers: $(k\,q)
      \cdot\bold a=k\cdot(q\cdot\bold a)$;
\item"8)" the feature of the numeric unity: $1\cdot\bold a=\bold a$.
\endroster
\endproclaim
     Let's consider the properties listed in the theorem~\mythetheorem{7.1} one
by one. Let's begin with the commutativity of addition. The sum
$\bold a+\bold b$ in the left hand side of the equality $\bold a+\bold b
=\bold b+\bold a$ is calculated by means of the triangle rule upon choosing 
some geometric realizations $\bold a=\overrightarrow{AB\,\,}\!$ and
$\bold b=\overrightarrow{BC\,\,}\!$ as shown in Fig\.~7.1.\par
\parshape 3 0cm 10cm 0cm 10cm 5cm 5cm 
     Let's draw the line parallel to the line $BC$ and passing through the point
$A$. Then we draw the line parallel to the line $AB$ and passing through the point
$C$. Both of these lines are in the plane of the triangle 
$ABC$. For this reason they intersect at some point $D$. The segments $[AB]$, 
$[BC]$, $[CD]$, and $[DA]$ form a parallelogram.\par
\parshape 3 5cm 5cm 5cm 5cm 0cm 10cm 
     Let's mark the vectors $\overrightarrow{DC\,\,}\!$ and 
$\overrightarrow{AD\,\,}\!$ on the segments $[CD]$ and $[DA]$. 
\vadjust{\vskip 5pt\hbox to 0pt{\kern 5pt\includegraphics{angemeng07.eps}\hss}\vskip -5pt}It is easy to see
that the vector $\overrightarrow{DC\,\,}\!$ is produced from the vector 
$\overrightarrow{AB\,\,}\!$ by applying the parallel translation from the point
$A$ to the point $D$. Similarly the vector $\overrightarrow{AD\,\,}\!$ is produced
from the vector $\overrightarrow{BC\,\,}\!$ by applying the parallel translation
from the point $B$ to the point $A$. Therefore $\overrightarrow{DC\,\,}\!
=\overrightarrow{AB\,\,}\!=\bold a$ and $\overrightarrow{BC\,\,}\!
=\overrightarrow{AD\,\,}\!=\bold b$. Now the triangles $ABC$ and $ADC$ yield
$$
\hskip -2em
\aligned
&\overrightarrow{AC\,\,}\!=\overrightarrow{AB\,\,}\!
+\overrightarrow{BC\,\,}\!=\bold a+\bold b,\\
&\overrightarrow{AC\,\,}\!=\overrightarrow{AD\,\,}\!
+\overrightarrow{DC\,\,}\!=\bold b+\bold a.
\endaligned
\mytag{7.1}
$$
From \mythetag{7.1} we derive the required equality $\bold a+\bold b
=\bold b+\bold a$.\par
     The relationship $\bold a+\bold b=\bold b+\bold a$ and Fig\.~7.1 yield
another method for adding vectors. It is called the {\it parallelogram rule}. 
In order to add two vectors $\bold a$ and $\bold b$ their geometric
realizations $\overrightarrow{AB\,\,}\!$ and $\overrightarrow{AD\,\,}\!$ with
the common initial point $A$ are used. They are completed up to the parallelogram
$ABCD$. Then the diagonal of this parallelogram is taken for the geometric
realization of the sum: $\bold a+\bold b=\overrightarrow{AB\,\,}\!
+\overrightarrow{AD\,\,}\!=\overrightarrow{AC\,\,}\!$.
\myexercise{7.1} Prove the equality $\bold a+\bold b=\bold b+\bold a$ for the
case where $\bold a\parallel\bold b$. For this purpose consider the subcases
$$
\xalignat 2
&\text{\rm 1) \ }\bold a\upuparrows\bold b;
&&\text{\rm 2) \ }\bold a\updownarrows\bold b\text{ \ and \ }|\bold a|>|\bold b|;
\\
&\text{\rm 3) \ }\bold a\updownarrows\bold b\text{ \ and \ }|\bold a|=|\bold b|;
&&\text{\rm 4) \ }\bold a\updownarrows\bold b\text{ \ and \ }|\bold a|<|\bold b|.
\endxalignat
$$
\endproclaim
\parshape 11 4.8cm 5.2cm 4.8cm 5.2cm 4.8cm 5.2cm 4.8cm 5.2cm 4.8cm 5.2cm 
4.8cm 5.2cm 4.8cm 5.2cm 4.8cm 5.2cm 4.8cm 5.2cm 4.8cm 5.2cm 0cm 10cm 
     The next property in the theorem~\mythetheorem{7.1} is the associativity 
of the operation of vector addition. In order to prove this property
\vadjust{\vskip 5pt\hbox to 0pt{\kern 5pt\includegraphics{angemeng08.eps}\hss}\vskip -5pt}we choose some arbitrary initial point $A$ and 
construct the following geometric realizations of the vectors:
$\bold a=\overrightarrow{AB\,\,}\!$, $\bold b=\overrightarrow{BC\,\,}\!$, and
$\bold c=\overrightarrow{CD\,\,}$. Applying the triangle rule of vector addition
to the triangles $ABC$ and $ACD$, we get the relationships 
$$
\hskip -2em
\gathered
\bold a+\bold b=\overrightarrow{AB\,\,}\!+\overrightarrow{BC\,\,}\!=
\overrightarrow{AC\,\,}\!,\\
(\bold a+\bold b)+\bold c=\overrightarrow{AC\,\,}\!+
\overrightarrow{CD\,\,}\!=\overrightarrow{AD\,\,}\!
\endgathered
\mytag{7.2}
$$
(see Fig\.~7.2). Applying the same rule to the triangles $BCD$ and $ABD$, we
get the analogous relationships
$$
\hskip -2em
\gathered
\bold b+\bold c=\overrightarrow{BC\,\,}\!+\overrightarrow{CD\,\,}\!=
\overrightarrow{BD\,\,}\!,\\
\bold a+(\bold b+\bold c)=\overrightarrow{AB\,\,}\!+
\overrightarrow{BD\,\,}\!=\overrightarrow{AD\,\,}\!.
\endgathered
\mytag{7.3}
$$
The required relationship $(\bold a+\bold b)+\bold c=\bold a+(\bold b+\bold c)$
now is immediate from the formulas \mythetag{7.2} and \mythetag{7.3}.\par
     {\bf A remark}. The tetragon $ABCD$ in Fig\.~7.2 is not necessarily
planar. For this reason the line $CD$ is shown as if it goes under the line 
$AB$, while the line $BD$ is shown going over the line $AC$.\par  
     The feature of the null vector $\bold a+\bold 0=\bold a$ is immediate
from the triangle rule for vector addition. Indeed, if an initial point $A$
for the vector $\bold a$ is chosen and if its geometric realization 
$\overrightarrow{AB\,\,}\!$ is built, then the null vector $\bold 0$ is
presented by its geometric realization $\overrightarrow{BB\,\,}\!$. From 
the definition~\mythedefinition{5.1} we derive $\overrightarrow{AB\,\,}\!+
\overrightarrow{BB\,\,}\!=\overrightarrow{AB\,\,}\!$ which yields 
$\bold a+\bold 0=\bold a$.\par
     The existence of an opposite vector is also easy to prove. Assume that
the vector $\bold a$ is presented by its geometric realization 
$\overrightarrow{AB\,\,}\!$. Let's consider the opposite geometric vector 
$\overrightarrow{BA\,\,}\!$ and let's denote through $\bold a'$ the 
corresponding free vector. Then
$$
\bold a+\bold a'=\overrightarrow{AB\,\,}\!+\overrightarrow{BA\,\,}\!
=\overrightarrow{AA\,\,}\!=\bold 0.
$$\par
     The distributivity of multiplication over the vector addition follows
from the properties of the {\it homothety transformation} in the Euclidean
space $\Bbb E$ (see~\S\,11 of Chapter~\uppercase\expandafter{\romannumeral 6} 
in \mycite{6}). It is sometimes called the {\it similarity transformation},
which is not quite exact. Similarity transformations constitute a larger
class of trans\-formations that comprises homothety transformations as a 
subclass within it.\par
% \parshape 3 0cm 10cm 0cm 10cm 5cm 5cm
     Let $\bold a\nparallel\bold b$ and let the sum of vectors $\bold a
+\bold b$ is calculated according to the triangle rule as shown in Fig\.~7.3.
\vadjust{\vskip 5pt\hbox to 0pt{\kern -8pt\includegraphics{angemeng09.eps}\hss}\vskip 115pt}Assume that $k>0$. Let's construct the
homothety transformation $h_{k\sssize A}\!:\Bbb E\to\Bbb E$ with the center
at the point $A$ and with the homothety factor $k$. Let's denote through
$E$ the image of the point $B$ under the transformation $h_{k\sssize A}$ 
and let's denote through $F$ the image of the point $C$ under this
transformation:
$$
\xalignat 2
&\hskip -2em
E=h_{k\sssize A}(B),
&&F=h_{k\sssize A}(C).
\endxalignat
$$
Due to the properties of the homothety the line $EF$ is parallel to the line 
$BC$ and we have the following relationships:
$$
\xalignat 2
&\hskip -2em
\overrightarrow{EF\,\,}\!\!\upuparrows\overrightarrow{BC\,\,}\!,
&&|EF|=|k|\cdot|BC|,\\
&\hskip -2em
\overrightarrow{AE\,\,}\!\!\upuparrows\overrightarrow{AB\,\,}\!,
&&|AE|=|k|\cdot|AB|,
\mytag{7.4}\\
&\hskip -2em
\overrightarrow{AF\,\,}\!\!\upuparrows\overrightarrow{AC\,\,}\!,
&&|AF|=|k|\cdot|AC|.
\endxalignat 
$$
Comparing \mythetag{7.4} with \mythetag{6.1} and taking into account 
that we consider the case $k>0$, from \mythetag{7.4} we derive
$$
\xalignat 3
&\overrightarrow{AE\,\,}\!\!=k\cdot\overrightarrow{AB\,\,}\!,
&&\overrightarrow{EF\,\,}\!\!=k\cdot\overrightarrow{BC\,\,}\!,
&&\overrightarrow{AF\,\,}\!\!=k\cdot\overrightarrow{AC\,\,}\!.
\qquad
\mytag{7.5}
\endxalignat 
$$
The relationships \mythetag{7.5} are sufficient for to prove the distributivity
of the mutiplication of vectors by numbers over the operation of vector addition.
Indeed, from \mythetag{7.5} we obtain:
$$
\hskip -2em
\gathered
k\cdot (\bold a+\bold b)=k\cdot(\overrightarrow{AB\,\,}\!
+\overrightarrow{BC\,\,}\!)=k\cdot\overrightarrow{AC\,\,}\!
=\overrightarrow{AF\,\,}\!=\\
=\overrightarrow{AE\,\,}\!+\overrightarrow{EF\,\,}\!=
k\cdot\overrightarrow{AB\,\,}\!+k\cdot\overrightarrow{BC\,\,}\!
=k\cdot\bold a+k\cdot\bold b.
\endgathered
\mytag
{7.6}
$$\par
     The case where $\bold a\nparallel\bold b$ and $k<0$ is very similar to
the case just above. In this case Fig\.~7.3 is replaced by Fig\.~7.4. Instead of 
the relationships \mythetag{7.4} we have the relationships
$$
\xalignat 2
&\hskip -2em
\overrightarrow{EF\,\,}\!\!\updownarrows\overrightarrow{BC\,\,}\!,
&&|EF|=|k|\cdot|BC|,\\
&\hskip -2em
\overrightarrow{AE\,\,}\!\!\updownarrows\overrightarrow{AB\,\,}\!,
&&|AE|=|k|\cdot|AB|,
\mytag{7.7}\\
&\hskip -2em
\overrightarrow{AF\,\,}\!\!\updownarrows\overrightarrow{AC\,\,}\!,
&&|AF|=|k|\cdot|AC|.
\endxalignat 
$$
Taking into account $k<0$ from \mythetag{7.7} we derive \mythetag{7.5} and 
\mythetag{7.6}.\par
     In the case $k=0$ the relationship $k\cdot (\bold a+\bold b)
=k\cdot\bold a+k\cdot\bold b$ reduces to the equality $\bold 0=\bold 0
+\bold 0$. This equality is trivially valid.\par 
\myexercise{7.2} Prove that $k\cdot (\bold a+\bold b)=k\cdot\bold a
+k\cdot\bold b$ for the case $\bold a\parallel\bold b$. For this purpose 
consider the subcases
$$
\xalignat 2
&\text{\rm 1) \ }\bold a\upuparrows\bold b;
&&\text{\rm 2) \ }\bold a\updownarrows\bold b\text{ \ and \ }|\bold a|
>|\bold b|;
\\
&\text{\rm 3) \ }\bold a\updownarrows\bold b\text{ \ and \ }|\bold a|
=|\bold b|;
&&\text{\rm 4) \ }\bold a\updownarrows\bold b\text{ \ and \ }|\bold a|
<|\bold b|.
\endxalignat
$$
In each of these subcases consider two options: $k>0$ and $k<0$.
\endproclaim
     Let's proceed to proving the distributivity of multiplication of vectors
by numbers over the addition of numbers. The case $\bold a=\bold 0$ in the
equality $(k+q)\cdot\bold a=k\cdot\bold a+q\cdot\bold a$ is trivial. In this
case the equality $(k+q)\cdot\bold a=k\cdot\bold a+q\cdot\bold a$ reduces to
$\bold 0=\bold 0+\bold 0$.\par
     The cases $k=0$ and $q=0$ are also trivial. In these cases the equality
$(k+q)\cdot\bold a=k\cdot\bold a+q\cdot\bold a$ reduces to the equalities 
$q\cdot\bold a=\bold 0+q\cdot\bold a$ and $k\cdot\bold a
=k\cdot\bold a+\bold 0$ respectively.\par
\parshape 10 0cm 10cm 0cm 10cm 5cm 5cm 5cm 5cm 5cm 5cm 5cm 5cm 
5cm 5cm 5cm 5cm 5cm 5cm 0cm 10cm 
     Let's consider the case $\bold a\neq\bold 0$ and for the sake of certainty 
let's assume that $k>0$ and $q>0$. \vadjust{\vskip 5pt\hbox to 0pt{\kern 5pt
\includegraphics{angemeng10.eps}\hss}\vskip -5pt}Let's choose some
arbitrary point $A$ and let's build the geometric realizations of the vector
$k\cdot\bold a$ with the initial point $A$. Let $B$ be the terminal point 
of this geometric realization. Then $\overrightarrow{AB\,\,}\!
=k\cdot\bold a$. Similarly, we construct the geometric realization 
$\overrightarrow{BC\,\,}\!=q\cdot\bold a$. Due to $k>0$ and $q>0$ the vectors
$\overrightarrow{AB\,\,}\!$ and $\overrightarrow{BC\,\,}\!$ are codirected 
to the vector $\bold a$. These vectors lie on the line $AB$ (see Fig\.~7.5).
The sum of these two vectors
$$
\hskip -2em
\overrightarrow{AC\,\,}\!=\overrightarrow{AB\,\,}\!
+\overrightarrow{BC\,\,}\!
\mytag{7.8}
$$
lie on the same line and it is codirected to the vector $\bold a$. The length of
the vector $\overrightarrow{AC\,\,}\!$ is given by the formula
$$
\hskip -2em
|AC|=|AB|+|BC|=k\,|\bold a|+q\,|\bold a|=(k+q)\,|\bold a|.
\mytag{7.9}
$$
Due to $\overrightarrow{AC\,\,}\!\upuparrows\bold a$ and due to $k+q>0$ from
\mythetag{7.9} we derive
$$
\hskip -2em
\overrightarrow{AC\,\,}\!=(k+q)\cdot\bold a. 
\mytag{7.10}
$$
Let's substitute \mythetag{7.10} and \mythetag{7.8} and take into account 
the relationships $\overrightarrow{AB\,\,}\!=k\cdot\bold a$ and 
$\overrightarrow{BC\,\,}\!=q\cdot\bold a$ which follow from our constructions. 
As a result we get the required equality 
$(k+q)\cdot\bold a=k\cdot\bold a+q\cdot\bold a$.\par
\myexercise{7.3} Prove that $(k+q)\cdot\bold a=k\cdot\bold a+q\cdot\bold a$ 
for the case where $\bold a\neq\bold 0$, while $k$ and $q$ are two numbers 
of mutually opposite signs. For the case consider the subcases
$$
\xalignat 3
&\text{\rm 1) \ }|k|>|q|;
&&\text{\rm 2) \ }|k|=|q|;
&&\text{\rm 3) \ }|k|<|q|.
\endxalignat
$$
\endproclaim
     The associativity of the multiplication of vectors by numbers is expressed
by the equality $(k\,q)\cdot\bold a=k\cdot(q\cdot\bold a)$. If $\bold a=\bold 0$,
this equality is trivially fulfilled. It reduces to $\bold 0=\bold 0$. If $k=0$ 
or if $q=0$, it is also trivial. In this case it reduces to $\bold 0=\bold 0$.
\par
\parshape 10 0cm 10cm 0cm 10cm 5cm 5cm 5cm 5cm 5cm 5cm 5cm 5cm 
5cm 5cm 5cm 5cm 5cm 5cm 0cm 10cm 
     Let's consider the case $\bold a\neq\bold 0$, $k>0$, and $q>0$. 
\vadjust{\vskip 5pt\hbox to 0pt{\kern 5pt\includegraphics{angemeng11.eps}\hss}\vskip -5pt}Let's choose some arbitrary point $A$ in the space 
$\Bbb E$ and build the geometric realization of the vector $q\cdot\bold a$ with
the initial point $A$. Let $B$ be the terminal point of this geometric realization. 
Then $\overrightarrow{AB\,\,}\!=q\cdot\bold a$ (see Fig.~7.6). Due to $q>0$ the
vector $\overrightarrow{AB\,\,}\!$ is codirected with the vector $\bold a$. Let's
build the vector $\overrightarrow{AC\,\,}\!$ as the product 
$\overrightarrow{AC\,\,}\!=k\cdot\overrightarrow{AB\,\,}\!=
k\cdot(q\cdot\bold a)$ relying upon the definition~\mythedefinition{6.1}. Due to
$k>0$ the vector $\overrightarrow{AC\,\,}\!$ is also codirected with $\bold a$. 
The lengths of $\overrightarrow{AB\,\,}\!$ and $\overrightarrow{AC\,\,}\!$ are given
by the formulas
$$
\xalignat 2
&\hskip -2em
|AB|=q\,|\bold a|,
&&|AC|=k\,|AB|.
\mytag{7.11}
\endxalignat
$$
From the relationships \mythetag{7.11} we derive the equality 
$$
\hskip -2em
|AC|=k\,(q\,|\bold a|)=(k\,q)\,|\bold a|. 
\mytag{7.12}
$$
The equality \mythetag{7.12} combined with $\overrightarrow{AC\,\,}\!\upuparrows
\bold a$ and $k\,q>0$ yields $\overrightarrow{AC\,\,}\!=(k\,q)\cdot\bold a$. By
our construction $\overrightarrow{AC\,\,}\!=k\cdot\overrightarrow{AB\,\,}\!=
k\cdot(q\cdot\bold a)$. As a result now we immediately derive the required equality 
$(k\,q)\cdot\bold a=k\cdot(q\cdot\bold a)$.\par
\myexercise{7.4} Prove the equality $(k\,q)\cdot\bold a=k\cdot(q\cdot\bold a)$ 
in the case where $\bold a\neq\bold 0$, while the numbers $k$ and $q$ are of opposite
signs. For this case consider the following two subcases:
$$
\xalignat 2
&\text{\rm 1) \ }k>0\text{\ \ and \ }q<0;
&&\text{\rm 2) \ }k<0\text{\ \ and \ }q>0.
\endxalignat
$$
\endproclaim
     The last item 8 in the theorem~\mythetheorem{7.1} is trivial. It is immediate 
from the definition~\mythedefinition{6.1}.
\head
\SectionNum{8}{28} Vectorial expressions and their transformations.
\endhead
\rightheadtext{\S\,8. Vectorial expressions and \dots}
      The properties of the algebraic operations with vectors listed in the
theorem~\mythetheorem{7.1} are used in transforming vectorial expressions. Saying
a vectorial expression one usually assumes a formula such that it yields a vector
upon performing calculations according to this formula. In this section we consider
some examples of vectorial expressions and learn some methods of transforming these
expressions.\par
     Assume that a list of several vectors $\bold a_1,\,\ldots,\,\bold a_n$ is given. 
Then one can write their sum setting brackets in various ways:
$$
\hskip -2em
\gathered
(\bold a_1+\bold a_2)+(\bold a_3+(\bold a_4+\ldots+(\bold a_{n-1}
+\bold a_n)\ldots)),\\
(\ldots((((\bold a_1+\bold a_2)+\bold a_3)+\bold a_4)+\ldots+\bold a_{n-1})
+\bold a_n).
\endgathered
\mytag{8.1}
$$
There are many ways of setting brackets. The formulas \mythetag{8.1} show only
two of them. However, despite the abundance of the ways for brackets setting, due
to the associativity of the vector addition (see item 2 in the 
theorem~\mythetheorem{7.1}) all of the expressions like \mythetag{8.1} yield
the same result. For this reason the sum of the vectors $\bold a_1,\,\ldots,
\,\bold a_n$ can be written without brackets at all:
$$
\hskip -2em
\bold a_1+\bold a_2+\bold a_3+\bold a_4+\ldots+\bold a_{n-1}+\bold a_n.
\mytag{8.2}
$$
In order to make the formula \mythetag{8.2} more concise the summation sign is used. 
Then the formula \mythetag{8.2} looks like
$$
\hskip -2em
\bold a_1+\bold a_2+\ldots+\bold a_n=\sum^n_{i=1}\bold a_i.
\mytag{8.3}
$$
\par
     The variable $i$ in the formula \mythetag{8.3} plays the role of the cycling 
variable in the summation cycle. It is called a {\it summation index}. This variable
takes all integer values ranging from $i=1$ to $i=n$. The sum \mythetag{8.3} itself
does not depend on the variable $i$. The symbol $i$ in the formula \mythetag{8.3} 
can be replaced by any other symbol, e.\,g\. by the symbol $j$ or by the symbol $k$:
$$
\hskip -2em
\sum^n_{i=1}\bold a_i=\sum^n_{j=1}\bold a_j=\sum^n_{k=1}\bold a_k.
\mytag{8.4}
$$
The trick with changing (redesignating) a summation index used in \mythetag{8.4}
is often applied for transforming expressions with sums.\par
     The commutativity of the vector addition (see item 1 in the 
theorem~\mythetheorem{7.1}) means that we can change the order of summands in
sums of vectors. For instance, in the sum \mythetag{8.2} we can write
$$
\bold a_1+\bold a_2+\ldots+\bold a_n
=\bold a_n+\bold a_{n-1}+\ldots+\bold a_1.
$$
The most often application for the commutativity of the vector addition is 
changing the summation order in multiple sums. Assume that a collection of 
vectors $\bold a_{ij}$ is given which is indexed by two indices $i$ and $j$, 
where $i=1,\,\ldots,\,m$ and $j=1,\,\ldots,\,n$. Then we have the equality
$$
\hskip -2em
\sum^m_{i=1}\sum^n_{j=1}\bold a_{ij}=\sum^n_{j=1}\sum^m_{i=1}\bold a_{ij}
\mytag{8.5}
$$
that follows from the commutativity of the vector addition. In the same time
we have the equality
$$
\hskip -2em
\sum^m_{i=1}\sum^n_{j=1}\bold a_{ij}
=\sum^m_{j=1}\sum^n_{i=1}\bold a_{j\kern 0.4pt i} 
\mytag{8.6}
$$
which is obtained by redesignating indices. Both methods of transforming
multiple sums \mythetag{8.5} and \mythetag{8.6} are used in dealing with 
vectors.\par
     The third item in the theorem~\mythetheorem{7.1} describes the
property of the null vector. This property is often used in calculations. 
If the sum of a part of summands in \mythetag{8.3} is zero, e\.\,g\. if
the equality 
$$
\bold a_{k+1}+\ldots+\bold a_n=\sum^n_{i=k+1}\bold a_i=\bold 0
$$
is fulfilled, then the sum \mythetag{8.3} can be transformed as follows:
$$
\bold a_1+\ldots+\bold a_n=\sum^n_{i=1}\bold a_i=\sum^k_{i=1}\bold a_i
=\bold a_1+\ldots+\bold a_k.
$$\par
     The fourth item in the theorem~\mythetheorem{7.1} declares the existence
of an opposite vector $\bold a'$ for each vector $\bold a$. Due to this item we
can define the subtraction of vectors.
\mydefinition{8.1} The difference of two vectors $\bold a-\bold b$ is the sum
of the vector $\bold a$ with the vector $\bold b'$ opposite to the vector
$\bold b$. This fact is written as the equality 
$$
\hskip -2em
\bold a-\bold b=\bold a+\bold b'.
\mytag{8.7}
$$
\enddefinition
\myexercise{8.1} Using the definitions~\mythedefinition{6.1} and
\mythedefinition{8.1}, show that the opposite vector $\bold a'$ is produced
from the vector $\bold a$ by multiplying it by the number $-1$, i\.\,e\. 
$$
\hskip -2em
\bold a'=(-1)\cdot\bold a.
\mytag{8.8}
$$
Due to \mythetag{8.8} the vector $\bold a'$ opposite to $\bold a$ is denoted
through $-\bold a$ and we write $\bold a'=-\bold a=(-1)\cdot\bold a$.
\endproclaim
     The distributivity properties of the multiplication of vectors by numbers
(see items 5 and 6 in the theorem~\mythetheorem{7.1}) are used for expanding 
expressions and for collecting similar terms in them:
$$
\gather
\hskip -2em
\alpha\cdot\!\biggl(\,\sum^n_{i=1}\bold a_i\biggr)=
\sum^n_{i=1}\alpha\cdot\bold a_i,
\mytag{8.9}\\
\hskip -2em
\biggl(\,\sum^n_{i=1}\alpha_i\biggr)\!\cdot\bold a=
\sum^n_{i=1}\alpha_i\cdot\bold a.
\mytag{8.10}
\endgather
$$
Transformations like \mythetag{8.9} and \mythetag{8.10} can be found in various
calculations with vectors.\par
\myexercise{8.2} Using the relationships \mythetag{8.7} and \mythetag{8.8}, 
prove the following properties of the operation of subtraction:
$$
\xalignat 2
&\bold a-\bold a=\bold 0;
&&(\bold a+\bold b)-\bold c=\bold a+(\bold b-\bold c);\\
&(\bold a-\bold b)+\bold c=\bold a-(\bold b-\bold c);
&&(\bold a-\bold b)-\bold c=\bold a-(\bold b+\bold c);\\
&\alpha\cdot(\bold a-\bold b)=\alpha\cdot\bold a-\alpha\cdot\bold b;
&&(\alpha-\beta)\cdot\bold a=\alpha\cdot\bold a-\beta\cdot\bold a.
\endxalignat
$$
Here $\bold a$, $\bold b$, and\/ $\bold c$ are vectors, while $\alpha$ and\/ 
$\beta$ are numbers.
\endproclaim
     The associativity of the multiplication of vectors by numbers (see item
7 in the theorem~\mythetheorem{7.1}) is used expanding vectorial expressions. 
Here is an example of such usage:
$$
\beta\cdot\!\biggl(\,\sum^n_{i=1}\alpha_i\cdot\bold a_i\biggr)
=\sum^n_{i=1}\beta\cdot(\alpha_i\cdot\bold a_i)
=\sum^n_{i=1}(\beta\,\alpha_i)\cdot\bold a_i.
\quad
\mytag{8.11}
$$
In multiple sums this property is combined with the commutativity of the
regular multiplication of numbers by numbers:
$$
\sum^m_{i=1}\alpha_i\cdot\!\biggl(\sum^n_{j=1}\beta_j\cdot
\bold a_{ij}\biggr)
=\sum^n_{j=1}\beta_j\cdot\!\biggl(\sum^m_{i=1}\alpha_i\cdot
\bold a_{ij}\biggr).
$$\par
     {\bf A remark}. It is important to note that the associativity of the
multiplication of vectors by numbers is that very property because of which
one can omit the dot sign in writing a product of a number and a vector:
$$
\alpha\,\bold a=\alpha\cdot\bold a.
$$
Below I use both forms of writing for products of vectors by numbers intending
to more clarity, conciseness, and aesthetic beauty of formulas.\par
     The last item 8 of the theorem~\mythetheorem{7.1} expresses the property
of the numeric unity in the form of the relationship $1\cdot\bold a=\bold a$.
This property is used in collecting similar terms and in finding common factors.
Let's consider an example:
$$
\gather
\bold a+3\cdot\bold b+2\cdot\bold a+\bold b=\bold a+2\cdot\bold a
+3\cdot\bold b+\bold b=1\cdot\bold a+2\cdot\bold a\,+\\
+\,3\cdot\bold b+1\cdot\bold b=(1+2)\cdot\bold a+(3+1)\cdot\bold b
=3\cdot\bold a+4\cdot\bold b.
\endgather
$$ 
\myexercise{8.3} Using the relationship $1\cdot\bold a=\bold a$, prove that
the conditions $\alpha\cdot\bold a=\bold 0$ and $\alpha\neq 0$ imply 
$\bold a=\bold 0$.
\endproclaim
\head
\SectionNum{9}{32} Linear combinations. Triviality, non-triviality, 
and vanishing.
\endhead
\rightheadtext{\S\,9. Linear combinations.}
     Assume that some set of $n$ free vectors $\bold a_1,\,\ldots,
\,\bold a_n$ is given. One can call it a collection of $n$ vectors,
a system of $n$ vectors, or a family of $n$ vectors either.\par
     Using the operation of vector addition and the operation of 
multiplication of vectors by numbers, \pagebreak one can compose some 
vectorial expression of the vectors $\bold a_1,\,\ldots,\,\bold a_n$. 
It is quite likely that this expression will comprise sums of vectors 
taken with some numeric coefficients.
\mydefinition{9.1}\!An expression of the form $\alpha_1\,\bold a_1+\ldots+
\alpha_n\,\bold a_n$\linebreak composed of the vectors $\bold a_1,\,\ldots,
\,\bold a_n$ is called a {\it linear combination} of these vectors. 
The numbers $\alpha_1,\,\ldots,\,\alpha_n$ are called the {\it 
coefficients of a linear combination}. If
$$
\hskip -2em
\alpha_1\,\bold a_1+\ldots+
\alpha_n\,\bold a_n=\bold b,
\mytag{9.1}
$$
then the vector $\bold b$ is called the {\it value of a linear combination}.
\enddefinition
     In complicated vectorial expressions linear combinations of the vectors
$\bold a_1,\,\ldots,\,\bold a_n$  can be multiplied by numbers and can be added 
to other linear combinations which can also be multiplied by some numbers. Then 
these sums can be multiplied by numbers and again can be added to other 
subexpressions of the same sort. This process can be repeated several times. 
However, upon expanding, upon applying the formula \mythetag{8.11}, and upon 
collecting similar terms with the use of the formula \mythetag{8.10} all such
complicated vectorial expressions reduce to linear combinations of the vectors
$\bold a_1,\,\ldots,\,\bold a_n$. Let's formulate this fact as a theorem. 
\mytheorem{9.1} Each vectorial expression composed of vectors $\bold a_1,\,\ldots,
\,\bold a_n$ by means of the operations of addition and multiplication by numbers 
can be transformed to some linear combination of these vectors $\bold a_1,\,\ldots,
\,\bold a_n$.
\endproclaim
      The value of a linear combination does not depend on the order of summands 
in it. For this reason linear combinations differing only in order of summands
are assumed to be coinciding. For example, the expressions $\alpha_1\,\bold a_1+\ldots+
\alpha_n\,\bold a_n$ and $\alpha_n\,\bold a_n+\ldots+\alpha_1\,\bold a_1$ are assumed
to define the same linear combination. 
\mydefinition{9.2} A linear combination $\alpha_1\,\bold a_1+\ldots+\alpha_n\,\bold a_n$ 
composed of the vectors $\bold a_1,\,\ldots,\,\bold a_n$ is called {\it trivial} if 
all of its coefficients are zero, i\.\,e\. if $\alpha_1=\ldots=\alpha_n=0$.
\enddefinition
\mydefinition{9.3} 
\parshape 3 0cm 10cm 0cm 10cm 5cm 5cm 
A linear combination $\alpha_1\,\bold a_1+\ldots+\alpha_n\,\bold a_n$
\vadjust{\vskip 5pt\hbox to 0pt{\kern 5pt\includegraphics{angemeng12.eps}\hss}\vskip -5pt}composed of vectors $\bold a_1,\,\ldots,\,\bold a_n$
is called {\it vanishing} or {\it being equal to zero} if its value is equal to the 
null vector, i\.\,e\. if the vector $\bold b$ in \mythetag{9.1} is equal to zero.
\enddefinition
\parshape 6 5cm 5cm 5cm 5cm 5cm 5cm 5cm 5cm 5cm 5cm 0cm 10cm 
      Each trivial linear combination is equal to zero. However, the converse is 
not valid. Noe each vanishing linear combination is trivial. In Fig\.~9.1 we have 
a triangle $ABC$. Its sides are marked as vectors $\bold a_1=\overrightarrow{AB\,\,}\!$, 
$\bold a_2=\overrightarrow{BC\,\,}\!$, and $\bold a_3=\overrightarrow{CA\,\,}\!$. 
By construction the sum of these three vectors $\bold a_1$, $\bold a_2$, and $\bold a_3$
in Fig\.~9.1 is zero: 
$$
\hskip -2em
\bold a_1+\bold a_2+\bold a_3
=\overrightarrow{AB\,\,}\!+\overrightarrow{BC\,\,}\!
+\overrightarrow{CA\,\,}\!=\bold 0. 
\mytag{9.2}
$$
The equality \mythetag{9.2} can be written as follows:
$$
\hskip -2em
1\cdot\bold a_1+1\cdot\bold a_2+1\cdot\bold a_3=\bold 0. 
\mytag{9.3}
$$
It is easy to see that the linear combination in the left hand side of the equality
\mythetag{9.3} is not trivial (see Definition~\mythedefinition{9.2}), however, 
it is equal to zero according to the definition~\mythedefinition{9.3}.
\mydefinition{9.4} A linear combination $\alpha_1\,\bold a_1+\ldots+\alpha_n\,\bold a_n$
composed of the vectors $\bold a_1,\,\ldots,\,\bold a_n$ is called {\it non-trivial} 
if it is not trivial, i\.\,e\. at least one of its coefficients $\alpha_1,\,\ldots,\,
\alpha_n$ is not equal to zero.
\enddefinition
\head
\SectionNum{10}{34} Linear dependence and linear independence. 
\endhead
\rightheadtext{\S\,10. Linear dependence and \dots}
\mydefinition{10.1} A system of vectors $\bold a_1,\,\ldots,\,\bold a_n$ 
is called {\it linearly dependent} if there is a non-trivial linear combination
of these vectors which is equal to zero. 
\enddefinition
     The vectors $\bold a_1$, $\bold a_2$, $\bold a_3$ shown in Fig\.~9.1
is an example of a linearly dependent set of vectors.\par
     It is important to note that the linear dependence is a property of systems 
of vectors, it is not a property of linear combinations. Linear combinations in 
the definition~\mythedefinition{10.1} are only tools for revealing the linear 
dependence.\par
     It is also important to note that the linear dependence, once it is present
in a collection of vectors $\bold a_1,\,\ldots,\,\bold a_n$, does not depend on the 
order of vectors in this collection. This follows from the fact that the value of 
any linear combination and its triviality or non-triviality are not destroyed
if we transpose its summands. 
\mydefinition{10.2} A system of vectors $\bold a_1,\,\ldots,\,\bold a_n$ is called
{\it linearly independent}, if it is not linearly dependent in the sense of the
definition~\mythedefinition{10.1}, i\.\,e\. if there is no linear combination of
these vectors being non-trivial and being equal to zero simultaneously. 
\enddefinition
     One can prove the existence of a linear combination with the required 
properties in the definition~\mythedefinition{10.1} by finding an example of such 
a linear combination. However, proving the non-existence in the 
definition~\mythedefinition{10.2} is more difficult. For this reason the following
theorem is used. 
\mytheoremwithtitle{10.1}{\ ({\bf linear independence criterion})} A system of
vectors $\bold a_1,\,\ldots,\,\bold a_n$ is linearly independent if and only if
vanishing of a linear combination of these vectors implies its triviality. 
\endproclaim
\demo{Proof} The proof is based on a simple logical reasoning. Indeed, the
non-existence of a linear combination of the vectors $\bold a_1,\,\ldots,\,\bold a_n$, 
being non-trivial and vanishing simultaneously means that a linear combination of
these vectors is inevitably trivial whenever it is equal to zero. In other words
vanishing of a linear combination of these vectors implies triviality of this linear
combination. The theorem~\mythetheorem{10.1} is proved.\qed\enddemo
\mytheorem{10.2} A system of vectors $\bold a_1,\,\ldots,\,\bold a_n$ is linearly
independent if and only if non-triviality of a linear combination of these vectors
implies that it is not equal to zero. 
\endproclaim
     The theorem~\mythetheorem{10.2} is very similar to the theorem~\mythetheorem{10.1}.
However, it is less popular and is less often used. 
\myexercise{10.1} Prove the theorem~\mythetheorem{10.2} using the analogy with the 
theorem~\mythetheorem{10.1}. 
\endproclaim
\head
\SectionNum{11}{36} Properties of the linear dependence.
\endhead
\rightheadtext{\S\,11. Properties of the linear dependence.}
\mydefinition{11.1} The vector $\bold b$ is said {\it to be expressed as a linear
combination} of the vectors $\bold a_1,\,\ldots,\,\bold a_n$ if it is the value
of some linear combination composed of the vectors $\bold a_1,\,\ldots,\,\bold a_n$ 
(see \mythetag{9.1}). In this situation for the sake of brevity the vector $\bold b$
is sometimes said {\it to be linearly expressed} through the vectors $\bold a_1,\,
\ldots,\,\bold a_n$ or {\it to be expressed in a linear way\/} through $\bold a_1,
\,\ldots,\,\bold a_n$. 
\enddefinition
     There are five basic properties of the linear dependence of vectors. We 
formulate them as a theorem. 
\mytheorem{11.1} The relation of the linear dependence for a system of vectors 
possesses the following basic properties:
\roster
\rosteritemwd=0pt
\item"1)" a system of vectors comprising the null vector is linearly dependent;
\item"2)" a system of vectors comprising a linearly dependent subsystem is linearly 
          dependent itself;
\item"3)" if a system of vectors is linearly dependent, then at least one of these
vectors is expressed in a linear way through other vectors of this system;
\item"4)" if a system of vectors $\bold a_1,\,\ldots,\,\bold a_n$ is linearly
independent, while complementing it with one more vector $\bold a_{n+1}$ makes the
system linearly dependent, then the vector $\bold a_{n+1}$ is linearly expressed 
through the vectors $\bold a_1,\,\ldots,\,\bold a_n$;
\item"5)" if a vector $\bold b$ is linearly expressed through some $m$ vectors 
$\bold a_1,\,\ldots,\,\bold a_m$ and if each of the vectors $\bold a_1,\,\ldots,
\,\bold a_m$ is linearly expressed through some other $n$ vectors $\bold c_1,\,\ldots,
\,\bold c_n$, then the vector $\bold b$ is linearly expressed through the vectors 
$\bold c_1,\,\ldots,\,\bold c_n$.
\endroster
\endproclaim
     The properties 1)--5) in the theorem~\mythetheorem{11.1} are relatively
simple. Their proofs are purely algebraic, they do not require drawings. I do 
not prove them in this book since the reader can find their proofs in \S\,3 of 
Chapter~\uppercase\expandafter{\romannumeral 1} in the book \mycite{1}.\par
     Apart from the properties 1)--5) listed in the theorem~\mythetheorem{11.1},
there is one more property which is formulated separately. 
\mytheoremwithtitle{11.2}{\ ({\bf Steinitz})} If the vectors $\bold a_1,\,
\ldots,\,\bold a_n$ are linearly independent and if each of them is linearly 
expressed through some other vectors $\bold b_1,\,\ldots,\,\bold b_m$, then 
$m\geqslant n$.
\endproclaim
     The Steinitz theorem~\mythetheorem{11.2} is very important in studying
multidimensional spaces. We do not use it in studying the three-dimensional 
space $\Bbb E$ in this book. 
\head
\SectionNum{12}{37} Linear dependence for $n=1$.
\endhead
\rightheadtext{\S\,12. Linear dependence for $n=1$.}
     Let's consider the case of a system composed of a single vector $\bold a_1$
and apply the definition of the linear dependence~\mythedefinition{10.1} to this
system. The linear dependence of such a system of one vector $\bold a_1$ means
that there is a linear combination of this single vector which is non-trivial
and equal to zero at the same time:
$$
\hskip -2em
\alpha_1\,\bold a_1=\bold 0.
\mytag{12.1}
$$
Non-triviality of the linear combination in the left hand side of \mythetag{12.1} 
means that $\alpha_1\neq 0$. Due to $\alpha_1\neq 0$ from \mythetag{12.1} we derive 
$$
\hskip -2em
\bold a_1=\bold 0
\mytag{12.2}
$$
(see Exercise~\mytheexercise{8.3}). Thus, the linear dependence of a system composed
of one vector $\bold a_1$ yields $\bold a_1=\bold 0$.\par
     The converse proposition is also valid. Indeed, assume that the equality
\mythetag{12.2} is fulfilled. Let's write it as follows:
$$
\pagebreak
\hskip -2em
1\cdot\bold a_1=\bold 0.
\mytag{12.3}
$$
The left hand side of the equality \mythetag{12.3} is a non-trivial linear combination
for the system of one vector $\bold a_1$ which is equal to zero. Its existence means
that such a system of one vector is linearly dependent. We write this result as a theorem. 
\mytheorem{12.1} A system composed of a single vector $\bold a_1$ is linearly dependent 
if and only if this vector is zero, i.\,e\. $\bold a_1=\bold 0$. 
\endproclaim
\head
\SectionNum{13}{38} Linear dependence for $n=2$.\\
Collinearity of vectors.
\endhead
\rightheadtext{\S\,13. Linear dependence for $n=2$.}
     Let's consider a system composed of two vectors $\bold a_1$ and $\bold a_2$. 
Applying the definition of the linear dependence~\mythedefinition{10.1} to it, we 
get the existence of a linear combination of these vectors which is non-trivial and 
equal to zero simultaneously:
$$
\hskip -2em
\alpha_1\,\bold a_1+\alpha_2\,\bold a_2=\bold 0.
\mytag{13.1}
$$
The non-triviality of the linear combination in the left hand side of the equality  
\mythetag{13.1} means that $\alpha_1\neq 0$ or $\alpha_2\neq 0$. Since the linear
dependence is not sensitive to the order of vectors in a system, without loss of
generality we can assume that $\alpha_1\neq 0$. Then the equality \mythetag{13.1} can
be written as   
$$
\hskip -2em
\bold a_1=-\frac{\alpha_2}{\alpha_1}\,\bold a_2.
\mytag{13.2}
$$
Let's denote $\beta_2=-\alpha_2/\alpha_1$ and write the equality \mythetag{13.2} as 
$$
\hskip -2em
\bold a_1=\beta_2\,\bold a_2.
\mytag{13.3}
$$
Note that the relationship \mythetag{13.3} could also be derived by means of the
item 3 of the theorem~\mythetheorem{11.1}.\par
     According to \mythetag{13.3}, the vector $\bold a_1$ is produced from the
vector $\bold a_2$ by multiplying it by the number $\beta_2$. In multiplying a vector 
by a number it length is usually changed (see Formulas \mythetag{6.1}, \mythetag{6.2}, 
\mythetag{6.3}, and Figs\.~6.1, 6.2, and 6.3). As for its direction, it either is 
preserved or is changed for the opposite one. In both of these cases the vector 
$\beta_2\,\bold a_2$ is parallel to the vector $\bold a_2$. If $\beta_2=0$m the vector
$\beta_2\,\bold a_2$ appears to be the null vector. Such a vector does not have its own 
direction, the null vector is assumed to be parallel to any other vector by definition. 
As a result of the above considerations the equality \mythetag{13.3} yields
$$
\hskip -2em
\bold a_1\parallel\bold a_2.
\mytag{13.4}
$$\par
     In the case of vectors for to denote their parallelism a special term 
{\it collinearity} is used. 
\mydefinition{13.1} Two free vectors $\bold a_1$ and $\bold a_2$ are called
{\it collinear}, if their geometric realizations are parallel to some straight 
line common to both of them. 
\enddefinition
     As we have seen above, in the case of two vectors their linear dependence
implies the collinearity of these vectors. The converse proposition is also
valid. Assume that the relationship \mythetag{13.4} is fulfilled. If both vectors 
$\bold a_1$ and $\bold a_2$ are zero, then the equality \mythetag{13.3} is
fulfilled where we can choose $\beta_2=1$. If at least one the two vectors is
nonzero, then up to a possible renumerating these vectors we can assume that
$\bold a_2\neq\bold 0$. Having built geometric realizations $\bold a_2
=\overrightarrow{AB\,\,}\!$ and $\bold a_1=\overrightarrow{AC\,\,}\!$, one can
choose the coefficient $\beta_2$ on the base of the Figs\.~6.1, 6.2, or 6.3 and
on the base of the formulas \mythetag{6.1}, \mythetag{6.2}, \mythetag{6.3} so that
the equality \mythetag{13.3} is fulfilled in this case either. As for the equality 
\mythetag{13.3} itself, we write it as follows:
$$
\hskip -2em
1\cdot\bold a_1+(-\beta_2)\cdot\bold a_2=\bold 0.
\mytag{13.5}
$$
Since $1\neq 0$, the left hand side of the equality \mythetag{13.5} is a non-trivial
linear combination of the vectors $\bold a_1$ and $\bold a_2$ which is equal to zero. 
The existence of such a linear combination means that the vectors $\bold a_1$ and
$\bold a_2$ are linearly dependent. Thus, the converse proposition that the collinearity 
of two vectors implies their linear \pagebreak dependence is proved.\par
     Combining the direct and converse propositions proved above, one can formulate
them as a single theorem. 
\mytheorem{13.1} A system of two vectors $\bold a_1$ and $\bold a_2$ is
linearly dependent if and only if these vectors are collinear, i\.\,e\. 
$\bold a_1\parallel\bold a_2$.
\endproclaim
\head
\SectionNum{14}{40} Linear dependence for $n=3$.\\
Coplanartity of vectors.
\endhead
\rightheadtext{\S\,14. Linear dependence for $n=3$.}
\parshape 3 0cm 10cm 0cm 10cm 5cm 5cm 
     Let;s consider a system composed of three vectors $\bold a_1$, $\bold a_2$, 
and $\bold a_3$. Assume that it is linearly dependent. \vadjust{\vskip 5pt\hbox to 
0pt{\kern 5pt\includegraphics{angemeng13.eps}\hss}\vskip -5pt}Applying 
the item 3 of the theorem~\mythetheorem{11.1} to this system, we get that one of 
the three vectors is linearly expressed through the other two vectors. Taking into 
account the possibility of renumerating our vectors, we can assume that the vector 
$\bold a_1$ is expressed through the vectors $\bold a_2$ and $\bold a_3$:
$$
\bold a_1=\beta_2\,\bold a_2+\beta_3\,\bold a_3.\quad
\mytag{14.1}
$$\par
% \parshape 2 5cm 5cm 0cm 10cm 
Let $A$ be some arbitrary point of the space $\Bbb E$. Let's build the geometric 
realizations $\bold a_2=\overrightarrow{AB\,\,}\!$ and $\beta_2\,\bold a_2
=\overrightarrow{AC\,\,}\!$. Then at the point $C$ we build the geometric
realizations of the vectors $\bold a_3=\overrightarrow{CD\,\,}\!$ and $\beta_3
\,\bold a_3=\overrightarrow{CE\,\,}\!$. The vectors $\overrightarrow{AC\,\,}\!$ 
and $\overrightarrow{CE\,\,}\!$ constitute two sides of the triangle $ACE$ (see
Fig\.~14.1). Then the sum of the vectors \mythetag{14.1} is presented by the
third side $\bold a_1=\overrightarrow{AE\,\,}\!$.\par
     The triangle $ACE$ is a planar form. The geometric realizations of the 
vectors $\bold a_1$, $\bold a_2$, and $\bold a_3$ lie on the sides of this 
triangle. Therefore they lie on the plane $ACE$. Instead of 
$\bold a_1=\overrightarrow{AE\,\,}\!$, 
$\bold a_2=\overrightarrow{AB\,\,}\!$, and 
$\bold a_3=\overrightarrow{CD\,\,}\!$ by means of parallel translations we can 
build some other geometric realizations of these three vectors. These geometric
realizations do not lie on the plane $ACE$, but they keep parallelism to this plane.
\mydefinition{14.1} Three free vectors $\bold a_1$, $\bold a_2$, and $\bold a_3$ 
are called {\it coplanar} if their geometric realizations are parallel to some
plane common to all three of them. 
\enddefinition
\mylemma{14.1} The linear dependence of three vecctors $\bold a_1$, $\bold a_2$,
$\bold a_3$ implies their coplanarity. 
\endproclaim
\myexercise{14.1} The above considerations yield a proof for the 
lemma~\mythelemma{14.1} on the base of the formula \mythetag{14.1} in the case 
where
$$
\xalignat 3
&\hskip -2em
\bold a_2\neq\bold 0,
&&\bold a_3\neq\bold 0,
&&\bold a_2\nparallel\bold a_3.
\quad
\mytag{14.2}
\endxalignat 
$$
Consider special cases where one or several conditions \mythetag{14.2} are not
fulfilled and derive the lemma~\mythelemma{14.1} from the formula \mythetag{14.1} 
in those special cases.
\endproclaim
\mylemma{14.2} The coplanarity of three vectors $\bold a_1$, $\bold a_2$, 
$\bold a_3$ imply their linear dependence. 
\endproclaim
\demo{Proof} \parshape 3 0cm 10cm 0cm 10cm 5cm 5cm 
If $\bold a_2=\bold 0$ or $\bold a_3=\bold 0$, then the propositions of the 
lemma~\mythelemma{14.2} follows from the first item of the theorem~\mythetheorem{11.1}. 
\vadjust{\vskip 5pt\hbox to 0pt{\kern 5pt
\includegraphics{angemeng14.eps}\hss}\vskip -5pt}If
$\bold a_2\neq\bold 0$, $\bold a_3\neq\bold 0$, $\bold a_2\parallel\bold a_3$, 
then the proposition of the lemma~\mythelemma{14.2} follows from the 
theorem~\mythetheorem{13.1} and from the item 2 of the 
theorem~\mythetheorem{11.1}. Therefore, in order to complete the proof of the 
lemma~\mythelemma{14.2} we should consider the case where all of the three conditions
\mythetag{14.2} are fulfilled.\par
\parshape 3 5cm 5cm 5cm 5cm 0cm 10cm 
     Let $A$ be some arbitrary point of the space $\Bbb E$. At this point we build
the geometric realizations of the vectors $\bold a_1=\overrightarrow{AD\,\,}\!$, 
$\bold a_2=\overrightarrow{AC\,\,}\!$, and 
$\bold a_3=\overrightarrow{AB\,\,}\!$ (see Fig\.~14.2). Due to the coplanarity of 
the vectors $\bold a_1$, $\bold a_2$, and $\bold a_3$ their geometric realizations
$\overrightarrow{AB\,\,}\!$, $\overrightarrow{AC\,\,}\!$, and  
$\overrightarrow{AD\,\,}\!$ lie on a plane. Let's denote this plane through $\alpha$. 
Through the point $D$ we draw a line parallel to the vector $\bold a_3\neq\bold 0$.
Such a line is unique and it lies on the plane $\alpha$. This line intersects
the line comprising the vector $\bold a_2=\overrightarrow{AC\,\,}\!$ at some unique 
point $E$ since $\bold a_2\neq\bold 0$ and $\bold a_2\nparallel\bold a_3$. 
Considering the points $A$, $E$, and $D$ in Fig\.~14.2, we derive the equality
$$
\hskip -2em
\bold a_1=\overrightarrow{AD\,\,}\!=\overrightarrow{AE\,\,}\!
+\overrightarrow{ED\,\,}\!.
\mytag{14.3}
$$\par
     The vector $\overrightarrow{AE\,\,}\!$ is collinear to the vector
$\overrightarrow{AC\,\,}\!=\bold a_2\neq\bold 0$ since these vectors lie on the
same line. For this reason there is a number $\beta_2$ such that
$\overrightarrow{AE\,\,}\!=\beta_2\,\bold a_2$. The vector $\overrightarrow{ED\,\,}\!$ 
is collinear to the vector $\overrightarrow{AB\,\,}\!=\bold a_3\neq\bold 0$ since
these vectors lie on parallel lines. Hence $\overrightarrow{ED\,\,}\!
=\beta_3\,\bold a_3$ for some number $\beta_3$. Upon substituting
$$
\xalignat 2
&\overrightarrow{AE\,\,}\!=\beta_2
\,\bold a_2,
&&\overrightarrow{ED\,\,}\!=\beta_3\,\bold a_3
\endxalignat
$$ 
into the equality \mythetag{14.3} this equality takes the form of 
\mythetag{14.1}.\par
      The last step in proving the lemma~\mythelemma{14.2} consists in writing
the equality \mythetag{14.1} in the following form: 
$$
\hskip -2em
1\cdot\bold a_1+(-\beta_2)\cdot\bold a_2+(-\beta_3)\cdot\bold a_3
=\bold 0.
\mytag{14.4}
$$
Since $1\neq 0$, the left hand side of the equality \mythetag{14.4} is a
non-trivial linear combination of the vectors $\bold a_1$, $\bold a_2$, 
$\bold a_3$ which is equal to zero. The existence of such a linear combination
means that the vectors $\bold a_1$, $\bold a_2$, and $\bold a_3$ are linearly
dependent.\qed\enddemo
\noindent     
The following theorem is derived from the lemmas~\mythelemma{14.1} and
\mythelemma{14.2}.
\mytheorem{14.1} A system of three vectors $\bold a_1$, $\bold a_2$, $\bold a_3$
is linearly dependent if and only if these vectors are coplanar.
\endproclaim
\head
\SectionNum{15}{42} Linear dependence for $n\geqslant 4$.
\endhead
\rightheadtext{\S\,15. Linear dependence for $n\geqslant 4$.}
\mytheorem{15.1} Any system consisting of four vectors 
\,$\bold a_1$, $\bold a_2$, $\bold a_3$, $\bold a_4$ in the space\/ $\Bbb E$
is linearly dependent.
\endproclaim
\mytheorem{15.2} Any system consisting of more than four vectors in the space\/ 
$\Bbb E$ is linearly dependent.
\endproclaim
     The theorem~\mythetheorem{15.2} follows from the theorem~\mythetheorem{15.1}
due to the item 3 of the theorem~\mythetheorem{11.1}. Therefore it is sufficient
to prove the theorem~\mythetheorem{15.1}. The theorem~\mythetheorem{15.1} itself
expresses a property of the three-dimensional space $\Bbb E$.  
\demo{Proof of the theorem~\mythetheorem{15.1}} Let's choose the subsystem
composed by three vectors $\bold a_1$, $\bold a_2$, $\bold a_3$ within the 
system of four vectors $\bold a_1$, $\bold a_2$, $\bold a_3$, $\bold a_4$. If these
three vectors \vadjust{\vskip 5pt\hbox to 0pt{\kern 30pt
\includegraphics{angemeng15.eps}\hss}\vskip 175pt}are linearly dependent,
then in order to prove the linear dependence of the vectors $\bold a_1$, $\bold a_2$, 
$\bold a_3$, $\bold a_4$ it is sufficient to apply the item 3 of the 
theorem~\mythetheorem{11.1}. Therefore in what fallows we consider the case where
the vectors $\bold a_1$, $\bold a_2$, $\bold a_3$ are linearly independent.\par
     From the linear independence of the vectors $\bold a_1$, $\bold a_2$, $\bold a_3$,
according to the theorem~\mythetheorem{14.1}, we derive their non-coplanarity. 
Moreove, from the linear independence of $\bold a_1$, $\bold a_2$, $\bold a_3$ due
to the item 3 of the theorem~\mythetheorem{11.1} we derive the linear independence of
any smaller subsystem within the system of these three vectors. In particular, the 
vectors $\bold a_1$, $\bold a_2$, $\bold a_3$ are nonzero and the vectors $\bold a_1$ 
and $\bold a_1$ are not collinear (see Theorems~\mythetheorem{12.1} and 
\mythetheorem{13.1}), i\.\,e\.
$$
\xalignat 4
&\bold a_1\neq\bold 0, 
&&\bold a_2\neq\bold 0, 
&&\bold a_3\neq\bold 0,
&&\bold a_1\nparallel\bold a_2.
\ \qquad
\mytag{15.1}
\endxalignat
$$\par
     Let $A$ be some arbitrary point of the space $\Bbb E$. Let's build the 
geometric realizations $\bold a_1=\overrightarrow{AB\,\,}\!$, 
$\bold a_2=\overrightarrow{AC\,\,}\!$, 
$\bold a_3=\overrightarrow{AD\,\,}\!$, $\bold a_4=\overrightarrow{AE\,\,}\!$
with the initial point $A$. Due to the condition $\bold a_1\nparallel\bold a_2$ 
in \mythetag{15.1} the vectors $\overrightarrow{AB\,\,}\!$ and
$\overrightarrow{AC\,\,}\!$ define a plane (see Fig\.~15.1). Let's denothe
this plane through $\alpha$. The vector $\overrightarrow{AD\,\,}\!$ does not
lie on the plane $\alpha$ and it is not parallel to this plane since the vectors
$\bold a_1$, $\bold a_2$, $\bold a_3$ are not coplanar.\par
     Let's draw a line passing through the terminal point of the vector 
$\bold a_4=\overrightarrow{AE\,\,}\!$ and being parallel to the vector $\bold a_3$. 
Since $\bold a_3\nparallel\alpha$, this line crosses the plane $\alpha$ at some unique
point $F$ and we have 
$$
\hskip -2em
\bold a_4=\overrightarrow{AE\,\,}\!=\overrightarrow{AF\,\,}\!
+\overrightarrow{FE\,\,}\!.
\mytag{15.2}
$$
Now let's draw a line passing through the point $F$ and being parallel
to the vector $\bold a_2$. Such a line lies on the plane $\alpha$. Due to
$\bold a_1\nparallel\bold a_2$ this line intersects the line $AB$ at some 
unique point $G$. Hence we have the following equality:
$$
\hskip -2em
\overrightarrow{AF\,\,}\!=\overrightarrow{AG\,\,}\!
+\overrightarrow{GF\,\,}\!.
\mytag{15.3}
$$\par
      Note that the vector $\overrightarrow{AG\,\,}\!$ lies on the same line
as the vector $\bold a_1=\overrightarrow{AB\,\,}\!$. From \mythetag{15.1} we get
$\bold a_1\neq\bold 0$. Hence there is a number $\beta_1$ such that 
$\overrightarrow{AG\,\,}\!=\beta_1\,\bold a_1$. Similarly, from  
$\overrightarrow{GF\,\,}\!\parallel\bold a_2$ and $\bold a_2\neq\bold 0$ we derive 
$\overrightarrow{GF\,\,}\!=\beta_2\,\bold a_2$ for some number 
$\beta_2$ and from $\overrightarrow{FE\,\,}\!\parallel\bold a_3$ and $a_3\neq\bold 0$ 
we derive that $\overrightarrow{FE\,\,}\!=\beta_3\,\bold a_3$ for some number
$\beta_3$. The rest is to substitute the obtained expressions for 
$\overrightarrow{AG\,\,}\!$, $\overrightarrow{GF\,\,}\!$, and $\overrightarrow{FE\,\,}\!$ 
into the formulas \mythetag{15.3} and \mythetag{15.2}. This yields
$$
\hskip -2em
\bold a_4=\beta_1\,\bold a_1+\beta_2\,\bold a_2
+\beta_3\,\bold a_3.
\mytag{15.4}
$$
The equality \mythetag{15.4} can be rewritten as follows:
$$
1\cdot\bold a_4+(-\beta_1)\cdot\bold a_1+(-\beta_2)\cdot\bold a_2
+(-\beta_3)\cdot\bold a_3=\bold 0.\quad
\mytag{15.5}
$$
Since $1\neq 0$, the left hand side of the equality \mythetag{15.5} is a non-trivial
linear combination of the vectors $\bold a_1$, $\bold a_2$, $\bold a_3$, $\bold a_4$ 
which is equal to zero. The existence of such a linear combination means that the
vectors $\bold a_1$, $\bold a_2$, $\bold a_3$, $\bold a_3$ are linearly dependent. 
The theorem~\mythetheorem{15.1} is proved.\qed\enddemo
\head
\SectionNum{16}{45} Bases on a line.
\endhead
\rightheadtext{\S\,16. Bases on a line.}
     Let $a$ be some line in the space $\Bbb E$. Let's consider free vectors parallel 
to the line $a$. They have geometric realizations lying on the line $a$. Restricting 
the freedom of moving such vectors, i\.\,e. forbidding geometric realizations outside
the line $a$, we obtain partially free vectors lying on the line $a$. 
\mydefinition{16.1} A system consisting of one non-zero vector $\bold e\neq\bold 0$
lying on a line $a$ is called a {\it basis} on this line. The vector $\bold e$ is
called the {\it basis vector} of this basis. 
\enddefinition
\parshape 7 5cm 5cm 5cm 5cm 5cm 5cm 5cm 5cm 5cm 5cm 5cm 5cm 0cm 10cm 
     Let $\bold e$ be the basis vector of some basis on the line $a$ and let $\bold x$ 
be some other vector lying on this line (see Fig\.~16.1). \vadjust{\vskip 5pt\hbox to 
0pt{\kern 5pt\includegraphics{angemeng16.eps}\hss}\vskip -5pt}Then 
$\bold x\parallel\bold e$ and hence there is a number $x$ such that the vector 
$\bold x$ is expressed through $\bold a$ by means of the formula
$$
\hskip -2em
\bold x=x\,\bold e.
\mytag{16.1}
$$
The number $x$ in the formula \mythetag{16.1} is called the {\it coordinate\/} of the
vector $\bold x$ in the basis $\bold e$, while the formula \mythetag{16.1} itself 
is called the {\it expansion\/} of the vector $\bold x$ in this basis.\par
      When writing the coordinates of vectors extracted from their expansions 
\mythetag{16.1} \pagebreak in a basis these coordinates are usually surrounded 
with double vertical lines 
$$
\hskip -2em
\bold x\mapsto\Vert x\Vert.
\mytag{16.2}
$$
Then these coordinated turn to matrices (see \mycite{7}). The mapping \mythetag{16.2} 
implements the basic idea of analytical geometry. This idea consists in replacing 
geometric objects by their numeric presentations. Bases in this case are tools for
such a transformation.\par
\head
\SectionNum{17}{46} Bases on a plane.
\endhead
\rightheadtext{\S\,17. Bases on a plane.}
     Let $\alpha$ be some plane in the space $\Bbb E$. Let's consider free vectors
parallel to the plane $\alpha$. They have geometric realizations lying on the plane
$\alpha$. Restricting the freedom of moving such vectors, i\.\,e. forbidding geometric 
realizations outside the plane $\alpha$, we obtain partially free vectors lying on 
the plane $\alpha$.
\mydefinition{17.1} An ordered pair of two non-collinear vectors $\bold e_1,\,\bold e_2$
lying on a plane $\alpha$ is called a {\it basis\/} on this plane. The vectors 
$\bold e_1$ and $\bold e_2$ are called the {\it basis vectors\/} of this basis. 
\enddefinition
\parshape 4 0cm 10cm 0cm 10cm 5cm 5cm 5cm 5cm
     In the definition~\mythedefinition{17.1} the term ``ordered system of vectors''
is used. \vadjust{\vskip 5pt\hbox to 0pt{\kern 0pt
\includegraphics{angemeng17.eps}\hss}\vskip -5pt}This term means a system 
of vectors in which some ordering of vectors is 
fixed: $\bold e_1$ is the first vector, $\bold e_2$ is the second vector. If we exchange
the vectors $\bold e_1$ and $\bold e_2$ and take $\bold e_2$ for the first vector, 
while $\bold e_1$ for the second vector, that would be another basis 
$\tilde\bold e_1,\,\tilde\bold e_2$ different from the basis $\bold e_1,\,\bold e_2$:
$$
\xalignat 2
&\tilde\bold e_1=\bold e_2,
&&\tilde\bold e_2=\bold e_1.
\endxalignat
$$\par
\parshape 3 5cm 5cm 5cm 5cm 0cm 10cm
     Let $\bold e_1,\,\bold e_2$ be a basis on a plane $\alpha$ and let $\bold x$ be
some vector lying on this place. Let's choose some arbitrary point
$O\in\alpha$ and let's build the geometric realizations of the three vectors 
$\bold e_1$, $\bold e_2$, and $\bold x$ with the initial point $O$:
$$
\xalignat 3
&\hskip -2em
\bold e_1=\overrightarrow{OA\,\,}\!,
&&\bold e_2=\overrightarrow{OB\,\,}\!,
&&\bold x=\overrightarrow{OC\,\,}\!.
\quad
\mytag{17.1}
\endxalignat
$$
Due to our choice of the vectors $\bold e_1$, $\bold e_2$, $\bold x$ and due to
$O\in\alpha$ the geometric realizations \mythetag{17.1} lie on the plane $\alpha$. 
Let's draw a line passing through the terminal point of the vector $\bold x$, i\.\,e\.
through the point $C$, and being parallel to the vector 
$\bold e_2=\overrightarrow{OB\,\,}\!$. Due to non-collinearity of the vectors
$\bold e_1\nparallel\bold e_2$ such a line intersects the line comprising the vector
$\bold e_1=\overrightarrow{OA\,\,}\!$ at some unique point $D$ (see Fig\.~17.1). 
This yields the equality
$$
\hskip -2em
\overrightarrow{OC\,\,}\!=\overrightarrow{OD\,\,}\!
+\overrightarrow{DC\,\,}\!.
\mytag{17.2}
$$\par
      The vector $\overrightarrow{OD\,\,}\!$ in \mythetag{17.2} is collinear with
the vector $\bold e_1=\overrightarrow{OA\,\,}\!$, while the vector
$\overrightarrow{DC\,\,}\!$ is collinear with the vector 
$\bold e_2=\overrightarrow{OB\,\,}\!$. For this reason there are 
two numbers $x_1$ and $x_2$ such that
$$
\xalignat 2
&\hskip -2em
\overrightarrow{OD\,\,}\!=x_1\,\bold e_1,
&&\overrightarrow{DC\,\,}\!=x_2\,\bold e_2.
\mytag{17.3}
\endxalignat
$$
Upon substituting \mythetag{17.3} into \mythetag{17.2} and taking into account 
the formulas \mythetag{17.1} we get the equality 
$$
\hskip -2em
\bold x=x_1\,\bold e_1+x_2\,\bold e_2.
\mytag{17.4}
$$\par
     The formula \mythetag{17.4} is analogous to the formula \mythetag{16.1}. 
The numbers $x_1$ and $x_2$ are called the {\it coordinates\/} of the vector
$\bold x$ in the basis $\bold e_1,\,\bold e_2$, while the formula 
\mythetag{17.4} itself is called the {\it expansion\/} of the vector $\bold x$ 
in this basis. When writing the coordinates of vectors they are usually arranged into
columns and surrounded with double vertical lines 
$$
\pagebreak
\hskip -2em
\bold x\mapsto\Vmatrix x_1\\\vspace{1ex}x_2\endVmatrix.
\mytag{17.5}
$$
The column of two numbers $x_1$ and $x_2$ in \mythetag{17.5} is called the
{\it coordinate column\/} of the vector $\bold x$. 
\head
\SectionNum{18}{48} Bases in the space.
\endhead
\rightheadtext{\S\,18. Bases in the space.}
\mydefinition{18.1} An ordered system of three non-coplanar vectors $\bold e_1,
\,\bold e_2,\,\bold e_3$ is called a {\it basis in the space\/} $\Bbb E$. 
\enddefinition
\parshape 11 0cm 10cm 0cm 10cm 0cm 10cm 0cm 10cm 0cm 10cm 0cm 10cm 0cm 10cm 
0cm 10cm 0cm 10cm 0cm 10cm 5cm 5cm 
     Let $\bold e_1,\,\bold e_2,\,\bold e_3$ be a basis in the space $\Bbb E$ and
let $\bold x$ be some vector. Let's choose some arbitrary point $O$ and build the
geometric realizations of all of the four vectors $\bold e_1$, $\bold e_2$, 
$\bold e_3$, and $\bold x$ with the common initial point $O$:
$$
\xalignat 2
&\hskip -2em
\bold e_1=\overrightarrow{OA\,\,}\!,
&&\bold e_2=\overrightarrow{OB\,\,}\!,\\
\vspace{-1.5ex}
\mytag{18.1}\\
\vspace{-1.5ex}
&\hskip -2em
\bold e_3=\overrightarrow{OC\,\,}\!,
&&\bold x=\overrightarrow{OD\,\,}\!.
\endxalignat
$$
The vectors $\bold e_1$ and $\bold e_2$ \vadjust{\vskip 5pt\hbox to 
0pt{\kern 0pt\includegraphics{angemeng18.eps}\hss}
\vskip -5pt}are not collinear since otherwise the whole 
system of three vectors $\bold e_1,\,\bold e_2,\,\bold e_3$ would be coplanar. 
For this reason the vectors $\bold e_1=\overrightarrow{OA\,\,}\!$ and 
$\bold e_2=\overrightarrow{OB\,\,}\!$ define a plane (this plane is denoted 
through $\alpha$ in Fig\~18.1) and they lie on this plane. The vector 
$\bold e_3=\overrightarrow{OC\,\,}\!$ does not lie on the plane $\alpha$ 
and it is not parallel to this plane (see Fig\.~18.1).\par
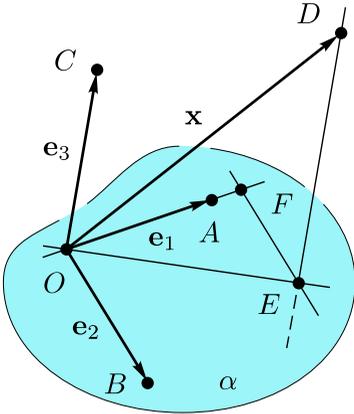
\parshape 1 5cm 5cm 
     Let's draw a line passing thro\-ugh the terminal point of the vector
$\bold x=\overrightarrow{OD\,\,}\!$ and being parallel to the vector 
$\bold e_3$. Such a line is not parallel to the plane $\alpha$ since 
$\bold e_3\nparallel\alpha$. It crosses the plane $\alpha$ at some unique point
$E$. As a result we get the equality
$$
\overrightarrow{OD\,\,}\!=\overrightarrow{OE\,\,}\!
+\overrightarrow{ED\,\,}\!.\quad
\mytag{18.2}
$$
\par
     Now let's draw a line passing through the point $E$ and being parallel to
the vector $\bold e_2$. The vectors $\bold e_1$ and $\bold e_2$ in \mythetag{18.1} 
are not collinear. For this reason such a line crosses the line comprising the vector
$\bold e_1$ at some unique point $F$. Considering the sequence of points $O$, $F$, $E$,
we derive 
$$
\hskip -2em
\overrightarrow{OE\,\,}\!=\overrightarrow{OF\,\,}\!
+\overrightarrow{FE\,\,}\!.
\mytag{18.3}
$$
Combining the equalities \mythetag{18.3} and \mythetag{18.2}, we obtain
$$
\hskip -2em
\overrightarrow{OD\,\,}\!=\overrightarrow{OF\,\,}\!
+\overrightarrow{FE\,\,}\!+\overrightarrow{ED\,\,}\!.
\mytag{18.4}
$$\par
     Note that, according to our geometric construction, the following collinearity 
conditions are fulfilled:
$$
\xalignat 3
&\hskip -2em
\overrightarrow{OF\,\,}\!\parallel\bold e_1,
&&\overrightarrow{FE\,\,}\!\parallel\bold e_2,
&&\overrightarrow{ED\,\,}\!\parallel\bold e_3.
\quad
\mytag{18.5}
\endxalignat
$$
From the collinearity condition $\overrightarrow{OF\,\,}\!\parallel\bold e_1$
in \mythetag{18.5} we derive the existence of a number $x_1$ such that 
$\overrightarrow{OF\,\,}\!=x_1\,\bold e_1$. From the other two collinearity
conditions in \mythetag{18.5} we derive the existence of two other numbers 
$x_2$ and $x_3$ such that $\overrightarrow{FE\,\,}\!=x_2\,\bold e_2$ and
$\overrightarrow{ED\,\,}\!=x_3\,\bold e_3$. As a result the equality 
\mythetag{18.4} is written as
$$
\hskip -2em
\bold x=x_1\,\bold e_1+x_2\,\bold e_2
+x_3\,\bold e_3.
\mytag{18.6}
$$\par
     The formula \mythetag{18.6} is analogous to the formulas \mythetag{16.1}
and \mythetag{17.4}. The numbers $x_1$, $x_2$, $x_3$ are called the 
{\it coordinates\/} of the vector $\bold x$ in the basis
$\bold e_1,\,\bold e_2,\,\bold e_3$, while the formula \mythetag{18.6} itself
is called the {\it expansion\/} of the vector $\bold x$ in this basis. In 
writing the coordinates of a vector they are usually arranged into a column
surrounded with two double vertical lines:
$$
\pagebreak
\hskip -2em
\bold x\mapsto\Vmatrix x_1\\\vspace{0.3ex}x_2\\\vspace{0.3ex}x_3
\endVmatrix.
\mytag{18.7}
$$
The column of the numbers $x_1$, $x_2$, $x_3$ in \mythetag{18.7} is called the
{\it coordinate column\/} of a vector $\bold x$.
\head
\SectionNum{19}{50} Uniqueness of the expansion of a vector 
in a basis.
\endhead
\rightheadtext{\S\,19. Uniqueness of the expansion \dots}
     Let $\bold e_1,\,\bold e_2,\,\bold e_3$ be some basis in the space $\Bbb E$. 
The geometric construction shown in Fig\.~18.1 can be applied to an arbitrary 
vector $\bold x$. It yields an expansion of \mythetag{18.6} of this vector in 
the basis $\bold e_1,\,\bold e_2,\,\bold e_3$. However, this construction is not 
a unique way foe expanding a vector in a basis. For example, instead of the plane
$OAB$ the plane $OAC$ can be taken for $\alpha$, while the line can be directed
parallel to the vector $\bold e_2$. The the line $EF$ is directed parallel to the
vector $\bold e_3$ and the point $F$ is obtained in its intersection with the
line $OA$ comprising the geometric realization of the vector $\bold e_1$. Such a
construction potentially could yield some other expansion of a vector $\bold x$ 
in the basis $\bold e_1,\,\bold e_2,\,\bold e_3$, i\.\,e\. an expansion with the
coefficients different from those of \mythetag{18.6}. The fact that actually this 
does not happen should certainly be proved.
\mytheorem{19.1} The expansion of an arbitrary vector $\bold x$ in a given basis 
$\bold e_1,\,\bold e_2,\,\bold e_3$ is unique. 
\endproclaim
\demo{Proof} The proof is by contradiction. Assume that the expansion 
\mythetag{18.6} is not unique and we have some other expansion of the vector
$\bold x$ in the basis $\bold e_1,\,\bold e_2,\,\bold e_3$:
$$
\hskip -2em
\bold x=\tilde x_1\,\bold e_1+\tilde x_2\,\bold e_2
+\tilde x_3\,\bold e_3.
\mytag{19.1}
$$
Let's subtract th expansion \mythetag{18.6} from the expansion \mythetag{19.1}. 
Upon collecting similar term we get 
$$
(\tilde x_1-x_1)\,\bold e_1+(\tilde x_2-x_2)\,\bold e_2
+(\tilde x_3-x_3)\,\bold e_3=\bold 0.
\quad
\mytag{19.2}
$$\par
     According to the definition~\mythedefinition{18.1}, the basis $\bold e_1,
\,\bold e_2,\,\bold e_3$ is a triple of non-coplanar vectors. Due to the
theorem~\mythetheorem{14.1} such a triple of vectors is linearly independent.
Looking at \mythetag{19.2}, we see that there we have a linear combination 
of the basis vectors $\bold e_1,\,\bold e_2,\,\bold e_3$ which is equal to
zero. Applying the theorem~\mythetheorem{10.1}, we conclude that this linear
combination is trivial, i\.\,e\.
$$
\xalignat 3
&\tilde x_1-x_1=0,
&&\tilde x_2-x_2=0,
&&\tilde x_3-x_3=0.
\quad\qquad
\mytag{19.3}
\endxalignat
$$
The equalities \mythetag{19.3} mean that the coefficients in the expansions
\mythetag{19.1} and \mythetag{18.6} do coincide, which contradicts the assumption 
that these expansions are different. The contradiction obtained proves the 
theorem~\mythetheorem{19.1}.\qed\enddemo
\myexercise{19.1} By analogy to the theorem~\mythetheorem{19.1} formulate 
and prove uniqueness theorems for expansions of vectors in bases on a plane
and in bases on a line. 
\endproclaim
\head
\SectionNum{20}{51} Index setting convention.
\endhead
\rightheadtext{\S\,20. Index setting convention.}
     The theorem~\mythetheorem{19.1} on the uniqueness of an expansion of
vector in a basis means that the mapping \mythetag{18.7}, which associates vectors
with their coordinates in some fixed basis, is a bijective mapping. This makes bases
an important tool for quantitative description of geometric objects. This tool was
improved substantially when a special index setting convention was admitted. This
convention is known as Einstein's tensorial notation.
\mydefinition{20.1} The {\it index setting convention}, which is also known as 
{\it Einstein's tensorial notation}, is a set of rules for placing indices in writing 
components of numeric arrays representing various geometric objects upon choosing 
some basis or some coordinate system. 
\enddefinition
     Einstein's tensorial notation is not a closed set of rules. When new types of
geometric objects are designed in science, new rules are added. For this reason 
below I formulate the index setting rules as they are needed.
\mydefinition{20.2} Basis vectors in a basis are enumerated by lower indices, while 
the coordinates of vectors expanded in a basis are enumerated by upper indices. 
\enddefinition
      The rule formulated in the definition~\mythedefinition{20.2} belongs to
Einstein's tensorial notation. According to this rule the formula \mythetag{18.6}
should be rewritten as 
$$
\hskip -2em
\bold x=x^1\,\bold e_1+x^2\,\bold e_2
+x^3\,\bold e_3=\sum^3_{i=1}x^i\,\bold e_i,
\mytag{20.1}
$$  
whire the mapping \mythetag{18.7} should be written in the following way:
$$
\hskip -2em
\bold x\mapsto\Vmatrix x^1\\\vspace{0.4ex}x^2\\\vspace{0.4ex}x^3
\endVmatrix.
\mytag{20.2}
$$
\myexercise{20.1} The rule from the definition~\mythedefinition{20.2} is applied
for bases on a line and for bases on a plane. Relying on this rule rewrite the
formulas \mythetag{16.1}, \mythetag{16.2}, \mythetag{17.4}, and \mythetag{17.5}. 
\endproclaim
\head
\SectionNum{21}{52} Usage of the coordinates of vectors.
\endhead
\rightheadtext{\S\,21. Usage of the coordinates of vectors.}
     Vectors can be used in solving various geometric problems, where the basic
algebraic operations with them are performed. These are the operation of vector
addition and the operation of multiplication of vectors by numbers. The usage of
bases and coordinates of vectors in bases relies on the following theorem.
\mytheorem{21.1} Let some basis $\bold e_1,\,\bold e_2,\,\bold e_3$ in
the space $\Bbb E$ be chosen and fixed. In this situation when adding vectors
their coordinates are added, while when multiplying a vector by a number its 
coordinates are multiplied by this number, i\.\,e\. if
$$
\xalignat 2
&\hskip -2em
\bold x\mapsto\Vmatrix x^1\\\vspace{0.4ex}x^2\\\vspace{0.4ex}x^3
\endVmatrix,
&&\bold y\mapsto\Vmatrix y^1\\\vspace{0.4ex}y^2\\\vspace{0.4ex}y^3
\endVmatrix,
\mytag{21.1}
\endxalignat
$$
then for $\bold x+\bold y$ and $\alpha\cdot\bold x$ we have the relationships
$$
\xalignat 2
&\bold x+\bold y\mapsto\Vmatrix x^1+y^1\\\vspace{0.4ex}x^2+y^2\\
\vspace{0.4ex}x^3+y^3\endVmatrix,
&&\alpha\cdot\bold x\mapsto\Vmatrix \alpha\,x^1\\\vspace{0.4ex}
\alpha\,x^2\\\vspace{0.4ex}\alpha\,x^3\endVmatrix.
\qquad
\mytag{21.2}
\endxalignat
$$
\endproclaim
\myexercise{21.1} Prove the theorem~\mythetheorem{21.1}, using the formulas 
\mythetag{21.1} and \mythetag{21.2} for this purpose, and prove the 
theorem~\mythetheorem{19.1}.
\endproclaim
\myexercise{21.2} Formulate and prove theorems analogous to the 
theorem~\mythetheorem{21.1} in the case of bases on a line and on a plane.
\endproclaim
\head
\SectionNum{22}{53} Changing a basis. Transition formulas\\
and transition matrices.
\endhead
\rightheadtext{\S\,22. Changing a basis.}
     When solving geometric problems sometimes one needs to change a basis 
replacing it with another basis. Let $\bold e_1,\,\bold e_2,\,\bold e_3$ 
and $\tilde\bold e_1,\,\tilde\bold e_2,\,\tilde\bold e_3$ be two bases in the
space $\Bbb E$. If the basis $\bold e_1,\,\bold e_2,\,\bold e_3$ is replaced
by the basis $\tilde\bold e_1,\,\tilde\bold e_2,\,\tilde\bold e_3$, then
$\bold e_1,\,\bold e_2,\,\bold e_3$ is usually called the {\it old basis}, 
while $\tilde\bold e_1,\,\tilde\bold e_2,\,\tilde\bold e_3$ is called the 
{\it new basis}. The procedure of changing an old basis for a new one can be 
understood as a {\it transition\/} from an old basis to a new basis or, in 
other words, as a {\it direct transition}. Conversely, changing a new basis
for an old one is understood as an {\it inverse transition}.\par
     Let $\bold e_1,\,\bold e_2,\,\bold e_3$ and $\tilde\bold e_1,
\,\tilde\bold e_2,\,\tilde\bold e_3$ be two bases in the space $\Bbb E$, where
$\bold e_1,\,\bold e_2,\,\bold e_3$ is an old basis and $\tilde\bold e_1,
\,\tilde\bold e_2,\,\tilde\bold e_3$ is a new basis. In the direct transition
procedure vectors of a new basis are expanded in an old basis, i\.\,e\. we
have
$$
\align
\hskip -2em
\tilde\bold e_1&=S^1_1\,\bold e_1+S^2_1\,\bold e_2
+S^3_1\,\bold e_3,\\
\hskip -2em
\tilde\bold e_2&=S^1_2\,\bold e_1+S^2_2\,\bold e_2
+S^3_2\,\bold e_3,
\mytag{22.1}\\
\hskip -2em
\tilde\bold e_3&=S^1_3\,\bold e_1+S^2_3\,\bold e_2
+S^3_3\,\bold e_3.
\endalign
$$
The formulas \mythetag{22.1} are called the {\it direct transition formulas}.
The numeric coefficients $S^1_1$, $S^2_1$, $S^3_1$ in \mythetag{22.1} are the
coordinates of the vector $\tilde\bold e_1$ expanded in the old basis. According
to the definition~\mythedefinition{20.2}, they are enumerated by an upper index. 
The lower index $1$ of them is the number of the vector $\tilde\bold e_1$ of 
which they are the coordinates. It is used in order to distinguish the coordinates
of the vector $\tilde\bold e_1$ from the coordinates of $\tilde \bold e_2$ and 
$\tilde\bold e_3$ in the second and in the third formulas \mythetag{22.1}.\par
     Let's apply the first mapping \mythetag{20.2} to the transition formulas
and write the coordinates of the vectors $\tilde\bold e_1,\,\tilde\bold e_2,
\,\tilde\bold e_3$ as columns:
$$
\xalignat 3
&\tilde\bold e_1\mapsto
\Vmatrix S^1_1\\\vspace{1ex}S^2_1\\\vspace{1ex}S^3_1\endVmatrix,
&&\tilde\bold e_2\mapsto
\Vmatrix S^1_2\\\vspace{1ex}S^2_2\\\vspace{1ex}S^3_2\endVmatrix,
&&\tilde\bold e_3\mapsto
\Vmatrix S^1_3\\\vspace{1ex}S^2_3\\\vspace{1ex}S^3_3\endVmatrix.
\ \qquad
\mytag{22.2}
\endxalignat
$$
The columns \mythetag{22.2} are usually glued into a single matrix. Such a matrix
is naturally denoted through $S$:
$$
\hskip -2em
S=\Vmatrix 
S^1_1 & S^1_2 & S^1_3\\\vspace{1ex}
S^2_1 & S^2_2 & S^2_3\\\vspace{1ex}
S^3_1 & S^3_2 & S^3_3
\endVmatrix.
\mytag{22.3}
$$
\mydefinition{22.1} The matrix \mythetag{22.3} whose components are determined
by the direct transition formulas \mythetag{22.1} is called the {\it direct
transition matrix}.
\enddefinition
     Note that the components of the direct transition matrix $S^{\kern 0.5pt i}_j$ 
are enumerated by two indices one of which is an upper index, while the other is a
lower index. These indices define the position of the element $S^{\kern 0.5pt i}_j$ 
in the matrix $S$: the upper index $i$ is a row number, while the lower index $j$ 
is a column number. This notation is a part of a general rule. 
\mydefinition{22.2} If elements of a double index array are enumerated by indices 
on different levels, then in composing a matrix of these elements the upper index is
used as a row number, while the lower index is used a column number.
\enddefinition
\mydefinition{22.3} If elements of a double index array are enumerated by indices
on the same level, then in composing a matrix of these elements the first index is
used as a row number, while the second index is used a column number.
\enddefinition
     The definitions~\mythedefinition{22.2} and \mythedefinition{22.3} can be 
considered as a part of the index setting convention from the 
definition~\mythedefinition{20.1}, though formally they are not since they do
not define the positions of array indices, but the way to visualize the array 
as a matrix.\par
     The direct transition formulas \mythetag{22.1} can be written in a concise
form using the summation sign for this purpose:
$$
\hskip -2em
\tilde\bold e_j=\sum^3_{i=1}S^{\kern 0.5pt i}_j\,\bold e_i
\text{, \ where \ }j=1,\,2,\,3.
\mytag{22.4}
$$
There is another way to write the formulas \mythetag{22.1} concisely. It is
based on the matrix multiplication (see \mycite{7}): 
$$
\Vmatrix 
\tilde\bold e_1 & \tilde\bold e_2 & \tilde\bold e_3
\endVmatrix
=\Vmatrix 
\bold e_1 & \bold e_2 & \bold e_3
\endVmatrix
\cdot
\Vmatrix 
S^1_1 & S^1_2 & S^1_3\\\vspace{1ex}
S^2_1 & S^2_2 & S^2_3\\\vspace{1ex}
S^3_1 & S^3_2 & S^3_3
\endVmatrix.
\quad
\mytag{22.5}
$$
Note that the basis vectors $\tilde\bold e_1,\,\tilde\bold e_2,
\,\tilde\bold e_3$ and $\bold e_1,\,\bold e_2,\,\bold e_3$ in the formula
\mythetag{22.5} are written in rows. This fact is an instance of a general 
rule.
\mydefinition{22.4} If elements of a single index array are enumerated by 
lower indices, then in matrix presentation they are written in a row,
i\.\,e\. they constitute a matrix whose height is equal to unity.
\enddefinition
\mydefinition{22.5} If elements of a single index array are enumerated by 
upper indices, then in matrix presentation they are written in a column,
i\.\,e\. they constitute a matrix whose width is equal to unity.
\enddefinition
     Note that writing the components of a vector $\bold x$ as a column
in the formula \mythetag{20.2} is concordant with the rule from the
definition~\mythedefinition{22.5}, while the formula \mythetag{18.7} violates
two rules at once --- the rule from the definition~\mythedefinition{20.2} and
the rule from the definition~\mythedefinition{22.5}.\par
     Now let's consider the inverse transition from the new basis 
$\tilde\bold e_1,\,\tilde\bold e_2,\,\tilde\bold e_3$ to the old basis $\bold e_1,
\,\bold e_2,\,\bold e_3$. In the inverse transition procedure the vectors 
of an old basis are expanded in a new basis:
$$
\align
\hskip -2em
\bold e_1&=T^1_1\,\tilde\bold e_1+T^2_1\,\tilde\bold e_2
+T^3_1\,\tilde\bold e_3,\\
\hskip -2em
\bold e_2&=T^1_2\,\tilde\bold e_1+T^2_2\,\tilde\bold e_2
+T^3_2\,\tilde\bold e_3,
\mytag{22.6}\\
\hskip -2em
\bold e_3&=T^1_3\,\tilde\bold e_1+T^2_3\,\tilde\bold e_2
+T^3_3\,\tilde\bold e_3.
\endalign
$$
The formulas \mythetag{22.6} are called the {\it inverse transition formulas}.
The numeric coefficients in the formulas \mythetag{22.6} are coordinates of the
vectors $\bold e_1$, $\bold e_2$, $\bold e_3$ in their expansions in the new
basis. These coefficients are arranged into columns:
$$
\xalignat 3
&\bold e_1\mapsto
\Vmatrix T^1_1\\\vspace{1ex}T^2_1\\\vspace{1ex}T^3_1\endVmatrix,
&&\bold e_2\mapsto
\Vmatrix T^1_2\\\vspace{1ex}T^2_2\\\vspace{1ex}T^3_2\endVmatrix,
&&\bold e_3\mapsto
\Vmatrix T^1_3\\\vspace{1ex}T^2_3\\\vspace{1ex}T^3_3\endVmatrix.
\ \qquad
\mytag{22.7}
\endxalignat
$$
Then the columns \mythetag{22.7} are united into a matrix:
$$
\hskip -2em
T=\Vmatrix 
T^1_1 & T^1_2 & T^1_3\\\vspace{1ex}
T^2_1 & T^2_2 & T^2_3\\\vspace{1ex}
T^3_1 & T^3_2 & T^3_3
\endVmatrix.
\mytag{22.8}
$$
\mydefinition{22.6} The matrix \mythetag{22.8} whose components are
determined by the inverse transition formulas \mythetag{22.6} is called
the {\it inverse transition matrix}.
\enddefinition
     The inverse transition formulas \mythetag{22.6} have a concise form, 
analogous to the formula \mythetag{22.4}:
$$
\hskip -2em
\bold e_j=\sum^3_{i=1}T^{\kern 0.5pt i}_j\,\tilde\bold e_i
\text{, \ where \ }j=1,\,2,\,3.
\mytag{22.9}
$$
There is also a matrix form of the formulas \mythetag{22.6}:
$$
\Vmatrix 
\bold e_1 & \bold e_2 & \bold e_3
\endVmatrix
=\Vmatrix 
\tilde\bold e_1 & \tilde\bold e_2 & \tilde\bold e_3
\endVmatrix
\cdot
\Vmatrix 
T^1_1 & T^1_2 & T^1_3\\\vspace{1ex}
T^2_1 & T^2_2 & T^2_3\\\vspace{1ex}
T^3_1 & T^3_2 & T^3_3
\endVmatrix.
\quad
\mytag{22.10}
$$
The formula \mythetag{22.10} is analogous to the formula \mythetag{22.5}.
\myexercise{20.1} By analogy to \mythetag{22.1} and \mythetag{22.6} write
the transition formulas for bases on a plane and for bases on a line (see
Definitions~\mythedefinition{16.1} and \mythedefinition{17.1}). Write also
the concise and matrix versions of these formulas. 
\endproclaim
\head
\SectionNum{23}{57} Some information on transition matrices.
\endhead
\rightheadtext{\S\,23. Some information \dots}
\mytheorem{23.1} The matrices $S$ and $T$ whose components are determined
by the transition formulas \mythetag{22.1} and \mythetag{22.6} are inverse
to each other, i.\,e\. $T=S^{-1}$ and $S=T^{-1}$. 
\endproclaim
    I do not prove the theorem~\mythetheorem{23.1} in this book. The reader 
can find this theorem and its proof in \mycite{1}.\par
    The relationships $T=S^{-1}$ and $S=T^{-1}$ from the 
theorem~\mythetheorem{23.1} mean that the product of $S$ by $T$ and the 
product of $T$ by $S$ both are equal to the unit matrix (see \mycite{7}): 
$$
\xalignat 2
&\hskip -2em
S\cdot T=1
&&T\cdot S=1.
\mytag{23.1}
\endxalignat
$$
Let's recall that the unit matrix is a square $n\times n$ matrix that has 
ones on the main diagonal and zeros in all other positions. Such a matrix 
is often denoted through the same symbol $1$ as the numeric unity. Therefore
we can write
$$
\hskip -2em
1=\Vmatrix 1 & 0 & 0\\
\vspace{0.5ex}
0 & 1 & 0\\
\vspace{0.5ex}
0 & 0 & 1
\endVmatrix.
\mytag{23.2}
$$
In order to denote the components of the unit matrix \mythetag{23.2} the symbol
$\delta$ is used. The indices enumerating rows and columns can be placed either 
on the same level or on different levels:
$$
\delta^{\kern 0.2pt ij}=\delta^{\kern 0.2pt i}_j=\delta_{ij}
=\cases
1&\text{for \ }i=j,\\
0&\text{for \ }i\neq j.
\endcases
\mytag{23.3}
$$
\mydefinition{23.1} The double index numeric array $\delta$ determined by the
formula \mythetag{23.3} is called the {\it Kronecker symbol\kern 2pt\footnotemark} 
or the {\it Kronecker delta}. 
\enddefinition
\footnotetext{\ Don't mix with the Kronecker symbol used in number theory
(see \mycite{8}).}
\adjustfootnotemark{-1}
     The positions of indices in the Kronecker symbol are determined by a context
where it is used. For example, the relationships \mythetag{23.1} can be written
in components. In this particular case the indices of the Kronecker symbol are 
placed on different levels:
$$
\xalignat 2
&\hskip -2em
\sum^3_{j=1}S^{\kern 0.5pt i}_j\,T^j_k=\delta^{\kern 0.2pt i}_k,
&&\sum^3_{j=1}T^{\kern 0.5pt i}_j\,S^j_k=\delta^{\kern 0.2pt i}_k.
\mytag{23.4}
\endxalignat
$$
Such a placement of the indices in the Kronecker symbol in \mythetag{23.4} is
inherited from the transition matrices $S$ and $T$.\par
     Noter that the transition matrices $S$ and $T$ are square matrices. For such
matrices the concept of the {\it determinant} is introduced (see \mycite{7}). This
is a number calculated through the components of a matrix according to some special
formulas. In the case of the unit matrix \mythetag{23.2} these formulas yield
$$
\hskip -2em
\det 1=1.
\mytag{23.5}
$$
The following fact is also well known. Its proof can be found in the book \mycite{7}.
\mytheorem{23.2} The determinant of a product of matrices is equal to the product of
their \pagebreak determinants. 
\endproclaim
     Let's apply the theorem~\mythetheorem{23.2} and the formula \mythetag{23.5}
to the products of the matrices $S$ and $T$ in \mythetag{23.1}. This yields
$$
\hskip -2em
\det S\cdot\det T=1.
\mytag{23.6}
$$
\mydefinition{23.2} A matrix with zero determinant is called {\it degenerate}. If the
determinant of a matrix is nonzero, such a matrix is called {\it non-degenerate}.
\enddefinition
     From the formula \mythetag{23.6} and the definition~\mythedefinition{23.2} we
immediately derive the following theorem. 
\mytheorem{23.3} For any two bases in the space $\Bbb E$ the corresponding transition 
matrices $S$ and $T$ are non-degenerate and the product of their determinants is equal
to the unity. 
\endproclaim
\mytheorem{23.4} Each non-degenerate $3\times 3$ matrix $S$ is a transition matrix
relating some basis $\bold e_1,\,\bold e_2,\,\bold e_3$ in the space $\Bbb E$ with
some other basis $\tilde\bold e_1,\,\tilde\bold e_2,\,\tilde\bold e_3$ in this space.
\endproclaim
     The theorem~\mythetheorem{23.4} is a strengthened version of the
theorem~\mythetheorem{23.3}. Its proof can be found in the book \mycite{1}.
\myexercise{23.1} Formulate theorems analogous to the theorems~\mythetheorem{23.1}, 
\mythetheorem{23.2}, and \mythetheorem{23.4} in the case of bases on a plane and in the
case of bases on a line.
\endproclaim
\head
\SectionNum{24}{59} Index setting in sums.
\endhead
\rightheadtext{\S\,24. Index setting in sums.}
    As we have already seen, in dealing with coordinates of vectors formulas with
sums arise. It is convenient to write these sums in a concise form using the summation
sign. Writing sums in this way one should follow some rules, which are listed below
as definitions.
\mydefinition{24.1} Each summation sign in a formula has its scope. This scope 
begins immediately after the summation sign to the right of it and ranges up to 
some delimiter:
\roster
\rosteritemwd=0pt
\item"1)" the end of the formula;
\item"2)" the equality sign;
\item"3)" the plus sign {\tencyr\char '074}$+${\tencyr\char '076} or the minus 
sign {\tencyr\char '074}$-${\tencyr\char '076} not enclosed into brackets
opened after a summation sign in question;
\item"4)" the closing bracket whose opening bracket precedes the summation sign
in question. 
\endroster
\enddefinition
     Let's recall that a summation sign is present in a formula and if some
variable is used as a cycling variable in this summation sign, such a variable
is called a {\it summation index\/} (see Formula \mythetag{8.3} and the comment 
to it). 
\mydefinition{24.2} Each summation index can be used only within the scope
of the corresponding summation sign. 
\enddefinition
     Apart from simple sums, multiple sums can be used in formulas. They obey the
following rule.
\mydefinition{24.3} A variable cannot be used as a summation index in more than 
one summation signs of a multiple sum. 
\enddefinition
\mydefinition{24.4} Variables which are not summation indices are called
{\it free variables}.
\enddefinition
     Summation indices as well as free variables can be used as indices  
enumerating basis vectors and array components. The following terminology
goes along with this usage.
\mydefinition{24.5} A free variable which is used as an index is called
a {\it free index}.
\enddefinition
     In the definitions~\mythedefinition{24.1}, \mythedefinition{24.2} 
and \mythedefinition{24.3} the commonly admitted rules are listed. Apart from them
there are more special rules which are used within the framework of Einstein's
tensorial notation (see Definition~\mythedefinition{20.1}). 
\mydefinition{24.6} If an expression is a simple sum or a multiple sum and if
each summand of it does not comprise other sums, then each free index should have
exactly one entry in this expression, while each summation index should enter 
twice --- once\linebreak as an upper index and once as a lower index.
\enddefinition
\mydefinition{24.7} The expression built according to the 
definition~\mythedefinition{24.6} can be used for composing sums with
numeric coefficients. Then all summands in such sums should have the same 
set of free indices and each free index should be on the same level (upper
or lower) in all summands. Regardless to the number of summands, in counting
the number of entries to the whole sum each free index is assumed to be 
entering only once. The level of a free index in the sum (upper or lower) is 
determined by its level in each summand. 
\enddefinition
      Lets consider some expressions as examples:
$$
\gather
\hskip -2em
\sum^3_{i=1}a^i\,b_i,
\kern 2em
\sum^3_{i=1}a^i\,g_{ij}, 
\kern 2em
\sum^3_{i=1}\sum^3_{k=1}a^i\,b^k\,g_{ik}, 
\mytag{24.1}\\
\hskip -2em
\sum^3_{k=1}2\,A^i_k\,b^m\,u^k\,v_q+\sum^3_{j=1}3\,C^{ijm}_{jq}
-\sum^3_{r=1}\sum^3_{s=1}C^r_{qs}
\,B^{ism}\,v_r.\quad
\mytag{24.2}
\endgather
$$
\myexercise{24.1} Verify that each expression in \mythetag{24.1} satisfies the
definition~\mythedefinition{24.6}, while the expression \mythetag{24.2} satisfies 
the de\-finition~\mythedefinition{24.7}. 
\endproclaim
\mydefinition{24.8} Sums composed according to the definition~\mythedefinition{24.7}
can enter as subexpressions into simple and multiple sums which  will be external sums 
with respect to them. Then some of their free indices or all of their free indices 
can turn into summation indices. Those of free indices that remain free are included 
into the list of free indices of the whole expression.\par
     In counting the number of entries of an index in a sum included into an external
simple or multiple sum the rule from the definition~\mythedefinition{24.7} is applied. 
Taking into account this rue, each free index of the ultimate expression should enter
it exactly once, while each summation index should enter it exactly twice --- once as
an upper index and once as a lower index. 
     In counting the number of entries of an index in a sum included into an external
simple or multiple sum the rule from the definition~\mythedefinition{24.7} is applied. 
Taking into account this rue, each free index of the ultimate expression should enter
it exactly once, while each summation index should enter it exactly twice --- once as
an upper index and once as a lower index. 
     In counting the number of entries of an index in a sum included into an external
simple or multiple sum the rule from the definition~\mythedefinition{24.7} is applied. 
Taking into account this rue, each free index of the ultimate expression should enter
it exactly once, while each summation index should enter it exactly twice --- once as
an upper index and once as a lower index. 
     In counting the number of entries of an index in a sum included into an external
simple or multiple sum the rule from the definition~\mythedefinition{24.7} is applied. 
Taking into account this rue, each free index of the ultimate expression should enter
it exactly once, while each summation index should enter it exactly twice --- once as
an upper index and once as a lower index. 
\enddefinition
     As an example we consider the following expression which comprises inner and outer 
sums:
$$
\hskip -2em
\sum^3_{k=1}A^i_k\biggl(\!B^k_q+\sum^3_{i=1}C^k_{iq}\,u^i\!\biggr).
\mytag{24.3}
$$
\myexercise{24.2} Make sure that the expression \mythetag{24.3} satisfies the
definition~\mythedefinition{24.8}. 
\endproclaim
\myexercise{24.3} Open the brackets in \mythetag{24.3} and verify that the resulting 
expression satisfies the definition~\mythedefinition{24.7}. 
\endproclaim
     The expressions built according to the definitions~\mythedefinition{24.6}, 
\mythedefinition{24.7}, and \mythedefinition{24.8} can be used for composing
equalities. In composing equalities the following rule is applied.  
\mydefinition{24.9} Both sides of an equality should have the same set of free indices
and each free index should have the same position (upper or lower) in both sides of an 
equality. 
\enddefinition
\head
\SectionNum{25}{63} Transformation of the coordinates of vectors
under a change of a basis.
\endhead
\rightheadtext{\S\,25. Transformation of the coordinates \dots}
     Let $\bold e_1,\,\bold e_2,\,\bold e_3$ and $\tilde\bold e_1,
\,\tilde\bold e_2,\,\tilde\bold e_3$ be two bases in the space $\Bbb E$ and let
$\bold e_1,\,\bold e_2,\,\bold e_3$ be changed for the basis $\tilde\bold e_1,
\,\tilde\bold e_2,\,\tilde\bold e_3$. As we already mentioned, in this case the basis
$\bold e_1,\,\bold e_2,\,\bold e_3$ is called an old basis, while $\tilde\bold e_1,
\,\tilde\bold e_2,\,\tilde\bold e_3$ is called a new basis.\par
     Let's consider some arbitrary vector $\bold x$ in the space $\Bbb E$. Expanding 
this vector in the old basis $\bold e_1,\,\bold e_2,\,\bold e_3$ and in the new
basis $\tilde\bold e_1,\,\tilde\bold e_2,\,\tilde\bold e_3$, we get two sets of its
coordinates:
$$
\xalignat 2
&\hskip -2em
\bold x\mapsto\Vmatrix x^1\\\vspace{0.4ex}x^2\\\vspace{0.4ex}x^3
\endVmatrix,
&&\bold x\mapsto\Vmatrix\tilde x^1\\\vspace{0.4ex}\tilde x^2\\
\vspace{0.4ex}\tilde x^3
\endVmatrix.
\mytag{25.1}
\endxalignat
$$
Both mappings \mythetag{25.1} are bijective. For this reason there is a bijective
correspondence between two sets of numbers $x^1,\,x^2,\,x^3$ and $\tilde x^1,
\,\tilde x^2,\,\tilde x^3$. In order to get explicit formulas expressing the coordinates 
of the vector $\bold x$ in the new basis through  its coordinates in the old basis 
we use the expansion \mythetag{20.1}: 
$$
\hskip -2em
\bold x=\sum^3_{j=1}x^{\kern 0.4pt j}\,\bold e_j.
\mytag{25.2}
$$  
Let's apply the inverse transition formula \mythetag{22.9} in order to express
the vector $\bold e_j$ in \mythetag{25.2} through the vectors of the new basis 
$\tilde\bold e_1,\,\tilde\bold e_2,\,\tilde\bold e_3$. Upon substituting 
\mythetag{22.9} into \mythetag{25.2}, we get 
$$
\hskip -2em
\gathered
\bold x=\sum^3_{j=1}x^{\kern 0.4pt j}\biggl(\,\sum^3_{i=1}
T^{\kern 0.5pt i}_j\,\tilde\bold e_i\!\biggr)=\\
=\sum^3_{j=1}\sum^3_{i=1}x^{\kern 0.4pt j}\,T^{\kern 0.5pt i}_j
\,\tilde\bold e_i=\sum^3_{i=1}\biggl(\,\sum^3_{j=1}
T^{\kern 0.5pt i}_j\,x^{\kern 0.4pt j}\!\biggr)\tilde\bold e_i.
\endgathered
\mytag{25.3}
$$  
The formula \mythetag{25.3} expresses the vector $\bold x$ as a linear combination
of the basis vectors $\tilde\bold e_1,\,\tilde\bold e_2,
\,\tilde\bold e_3$, i\.\,e\. it is an expansion of the vector $\bold x$ in the new
basis. Due to the uniqueness of the expansion of a vector in a basis (see 
Theorem~\mythetheorem{19.1}) the coefficients of such an expansion should coincide 
with the coordinates of the vector $\bold x$ in the new basis $\tilde\bold e_1,
\,\tilde\bold e_2,\,\tilde\bold e_3$:
$$
\hskip -2em
\tilde x^i=\sum^3_{j=1}
T^{\kern 0.5pt i}_j\,x^{\kern 0.4pt j}\text{, \ there \ }
i=1,\,2,\,3.
\mytag{25.4}
$$\par
     The formulas \mythetag{25.4} expressing the coordinates of an arbitrary vector
$\bold x$ in a new basis through its coordinates in an old basis are called the 
{\it direct transformation formulas}. Accordingly, the formulas expressing the 
coordinates of an arbitrary vector $\bold x$ in an old basis through its coordinates
in a new basis are called the {\it inverse transformation formulas}. The latter 
formulas need not be derived separately. It is sufficient to move the tilde sign from
the left hand side of the formulas \mythetag{25.4} to their right hand side and replace
$T^{\kern 0.5pt i}_j$ with $S^{\kern 0.5pt i}_j$. This yields
$$
\hskip -2em
x^i=\sum^3_{j=1}
S^{\kern 0.5pt i}_j\,\tilde x^{\kern 0.4pt j}\text{, \ where \ }
i=1,\,2,\,3.
\mytag{25.5}
$$\par
     The direct transformation formulas \mythetag{25.4} have the expanded form where
the summation is performed explicitly:
$$
\align
\hskip -2em
\tilde x^1&=T^1_1\,x^1+T^1_2\,x^2+T^1_3\,x^3,\\
\hskip -2em
\tilde x^2&=T^2_1\,x^1+T^2_2\,x^2+T^2_3\,x^3,
\mytag{25.6}\\
\hskip -2em
\tilde x^3&=T^3_1\,x^1+T^3_2\,x^2+T^3_3\,x^3.
\endalign
$$
The same is true for the inverse transformation formulas \mythetag{25.5}:
$$
\align
\hskip -2em
x^1&=S^1_1\,\tilde x^1+S^1_2\,\tilde x^2+S^1_3\,\tilde x^3,\\
\hskip -2em
x^2&=S^2_1\,\tilde x^1+S^2_2\,\tilde x^2+S^2_3\,\tilde x^3,
\mytag{25.7}\\
\hskip -2em
x^3&=S^3_1\,\tilde x^1+S^3_2\,\tilde x^2+S^3_3\,\tilde x^3.
\endalign
$$
Along with \mythetag{25.6}, there is the matrix form of the formulas 
\mythetag{25.4}:
$$
\hskip -2em
\Vmatrix\tilde x^1\\\vspace{0.4ex}\tilde x^2\\
\vspace{0.4ex}\tilde x^3
\endVmatrix
=\Vmatrix 
T^1_1 & T^1_2 & T^1_3\\\vspace{1ex}
T^2_1 & T^2_2 & T^2_3\\\vspace{1ex}
T^3_1 & T^3_2 & T^3_3
\endVmatrix\cdot\Vmatrix x^1\\\vspace{0.4ex}x^2\\\vspace{0.4ex}x^3
\endVmatrix
\mytag{25.8}
$$
Similarly, the inverse transformation formulas \mythetag{25.5}, along with
the expanded form \mythetag{25.7}, have the matrix form either
$$
\hskip -2em
\Vmatrix x^1\\\vspace{0.4ex}x^2\\\vspace{0.4ex}x^3
\endVmatrix
=\Vmatrix 
S^1_1 & S^1_2 & S^1_3\\\vspace{1ex}
S^2_1 & S^2_2 & S^2_3\\\vspace{1ex}
S^3_1 & S^3_2 & S^3_3
\endVmatrix\cdot\Vmatrix\tilde x^1\\\vspace{0.4ex}\tilde x^2\\
\vspace{0.4ex}\tilde x^3
\endVmatrix
\mytag{25.9}
$$
\myexercise{25.1} Write the analogs of the transformation formulas \mythetag{25.4},
\mythetag{25.5}, \mythetag{25.6}, \mythetag{25.7}, \mythetag{25.8}, \mythetag{25.9} 
for vectors on a plane and for vectors on a line expanded in corresponding bases on 
a plane and on a line. 
\endproclaim
\head
\SectionNum{26}{65} Scalar product.
\endhead
\rightheadtext{\S\,26. Scalar product.}
\parshape 3 0cm 10cm 0cm 10cm 4cm 6cm 
     Let $\bold a$ and $\bold b$ be two nonzero free vectors. Let's build their
geometric realizations $\bold a=\overrightarrow{OA\,\,}\!$ and 
$\bold b=\overrightarrow{OB\,\,}\!$ at some arbitrary point $O$. 
The smaller of two angles formed by the rays $[OA)$ 
and $[OB)$ at the point $O$ is called the {\it angle between vectors\/}
$\overrightarrow{OA\,\,}\!$ and $\overrightarrow{OB\,\,}\!$. In Fig\.~26.1 this
angle is denoted through $\varphi$. The value of the angle $\varphi$ ranges from
$0$ to $\pi$:
$$
0\leqslant\varphi\leqslant\pi.
$$\par
\parshape 2 4cm 6cm 0cm 10cm 
     The lengths of the vectors $\overrightarrow{OA\,\,}\!$ and 
$\overrightarrow{OB\,\,}\!$ \vadjust{\vskip 5pt\hbox to 0pt{\kern 
0pt\includegraphics{angemeng19.eps}\hss}\vskip -5pt}do not depend on the choice of a point $O$ (see
Definitions~\mythedefinition{3.1} and \mythedefinition{4.2}). The same is true
for the angle between them. Therefore we can deal with the lengths of the free 
vectors $\bold a$ and $\bold b$ and with the angle between them:
$$
\xalignat 3
&|\bold a|=|OA|, 
&&|\bold b|=|OB|,
&&\widehat{\bold a\bold b}=\widehat{AOB}=\varphi.
\endxalignat
$$
In the case where $\bold a=\bold 0$ or $\bold b=\bold 0$ the lengths of the 
vectors $\bold a$ and  $\bold b$ are defined, but the angle between 
these vectors is not defined.
\mydefinition{26.1} The {\it scalar product\/} of two nonzero vectors
$\bold a$ and $\bold b$ is a number equal to the product of their lengths
and the cosine of the angle between them:
$$
\hskip -2em
(\bold a,\bold b)=|\bold a|\,|\bold b|\,\cos\varphi.
\mytag{26.1}
$$
In the case where $\bold a=\bold 0$ or $\bold b=\bold 0$ the scalar product 
$(\bold a,\bold b)$ is assumed to be equal to zero by definition.
\enddefinition
     A comma is the multiplication sign in the writing the scalar product, 
not by itself, but together with round brackets surrounding the whole expression.
These brackets are natural delimiters for multiplicands: the first
multiplicand is an expression between the opening bracket and the comma, 
while the second multiplicand is an expression between the comma and the
closing bracket. Therefore in complicated expressions no auxiliary delimiters 
are required. For example, in the formula
$$
(\bold a+\bold b,\bold c+\bold d)
$$
the sums $\bold a+\bold b$ and $\bold c+\bold d$ are calculated first, then the
scalar multiplication is performed.\par
     {\bf A remark}. Often the scalar product is written as $\bold a\cdot\bold b$. 
Even the special term {\tencyr\char '074}dot product{\tencyr\char '076} is used. 
However, to my mind, this notation is not good. It is misleading since the dot
sign is used for denoting the product of a vector and a number and for denoting 
the product of two numbers. 
\head
\SectionNum{27}{67} Orthogonal projection onto a line.
\endhead
\rightheadtext{\S\,27. Orthogonal projection onto a line.}
\parshape 13 0cm 10cm 0cm 10cm 5cm 5cm 5cm 5cm 5cm 5cm 5cm 5cm 5cm 5cm 
5cm 5cm 5cm 5cm 5cm 5cm 5cm 5cm 5cm 5cm 0cm 10cm 
     Let $\bold a$ and $\bold b$ be two free vectors such that $\bold a\neq\bold 0$. 
Let's build their geometric realizations $\bold a
=\overrightarrow{OA\,\,}\!$ and $\bold b=\overrightarrow{OB\,\,}\!$ at some 
arbitrary point $O$. \vadjust{\vskip 5pt\hbox to 0pt{\kern 0pt
\includegraphics{angemeng20.eps}\hss}\vskip -5pt}The nonzero vector
$\overrightarrow{OA\,\,}\!$ defines a line. Let's drop the perpendicular from the
terminal point of the vector $\overrightarrow{OB\,\,}\!$, i\.\,e\. from the point
$B$, to this line and let's denote through $C$ the base of this perpendicular
(see Fig\.~27.1). In the special case where $\bold b\parallel\bold a$ and where
the point $B$ lies on the line $OA$ we choose the point $C$ coinciding with the
point $B$.\par
    The point $C$ determines two vectors $\overrightarrow{OC\,\,}\!$ and 
$\overrightarrow{CB\,\,}\!$. The vector $\overrightarrow{OC\,\,}\!$ is collinear to
the vector $\bold a$, while the vector $\overrightarrow{CB\,\,}\!$ is perpendicular 
to it. By means of parallel translations one can replicate the vectors 
$\overrightarrow{OC\,\,}\!$ and $\overrightarrow{CB\,\,}\!$ up to free vectors 
$\bold b_{\sssize\parallel}$ and $\bold b_{\sssize\perp}$ respectively. Note that 
the point $C$ is uniquely determined by the point $B$ and by the line $OA$ (see
Theorem 6.5 in Chapter~\uppercase\expandafter{\romannumeral 3} of the book 
\mycite{7}). For this reason the vectors $\bold b_{\sssize\parallel}$ and
$\bold b_{\sssize\perp}$ do not depend on the choice of a point $O$ and we can 
formulate the following theorem. 
\mytheorem{27.1} For any nonzero vector $\bold a\neq\bold 0$ and for any
vector $\bold b$ there two unique vectors $\bold b_{\sssize\parallel}$ and
$\bold b_{\sssize\perp}$ such that the vector $\bold b_{\sssize\parallel}$ is
collinear to $\bold a$, the vector $\bold b_{\sssize\perp}$ is perpendicular to
$\bold a$, and they both satisfy the equality being the expansion of the vector
$\bold b$:
$$
\hskip -2em
\bold b=\bold b_{\sssize\parallel}+\bold b_{\sssize\perp}.
\mytag{27.1}
$$
\endproclaim
     One should recall the special case where the point $C$ coincides with the 
point $B$. In this case $\bold b_{\sssize\perp}=\bold 0$ and we cannot verify
visually the orthogonality of the vectors $\bold b_{\sssize\perp}$ and $\bold a$. 
In order to extend the theorem~\mythetheorem{27.1} to this special case the 
following definition is introduced.
\mydefinition{27.1} All null vectors are assumed to be perpendicular
to each other and each null vector is assumed to be perpendicular to
any nonzero vector.
\enddefinition
     Like the definition~\mythedefinition{3.2}, the definition~\mythedefinition{27.1} 
is formulated for geometric vectors. Upon passing to free vectors it is convenient to
unite the definition \mythedefinition{3.2} with the definition~\mythedefinition{27.1} 
and then formulate the following definition.
\mydefinition{27.2} A free null vector $\bold 0$ of any physical nature is codirected
to itself and to any other vector. A free null vector $\bold 0$ of any physical nature
is perpendicular to itself and to any other vector.
\enddefinition
     When taking into account the definition~\mythedefinition{27.2}, the 
theorem~\mythetheorem{27.1} is proved by the constructions preceding it, while the
expansion \mythetag{27.1} follows from the evident equality
$$
\overrightarrow{OB\,\,}\!=\overrightarrow{OC\,\,}\!
+\overrightarrow{CB\,\,}\!.
$$\par
     Assume that a vector $\bold a\neq\bold 0$ is fixed. In this situation the
theorem \mythetheorem{27.1} provides a mapping $\pi_{\bold a}$ that associates
each vector $\bold b$ with its parallel component $\bold b_{\sssize\parallel}$.
\mydefinition{27.3} The mapping $\pi_{\bold a}$ that associates each free vector 
$\bold b$ with its parallel component $\bold b_{\sssize\parallel}$ in the
expansion \mythetag{27.1} is called the {\it orthogonal projection onto a line
given by the vector\linebreak $\bold a\neq\bold 0$} or, more exactly, the 
{\it orthogonal projection onto the direction of the vector $\bold a\neq\bold 0$}.
\enddefinition
      The orthogonal projection $\pi_{\bold a}$ is closely related to the scalar 
product of vectors. This relation is established by the following theorem. 
\mytheorem{27.2} For each nonzero vector $\bold a\neq\bold 0$ \pagebreak
and for any vector $\bold b$ the vector $\pi_{\bold a}(\bold b)$ is calculated
by means of the formula
$$
\hskip -2em
\pi_{\bold a}(\bold b)=\frac{(\bold b,\bold a)}{|\bold a|^2}\,
\bold a.
\mytag{27.2}
$$
\endproclaim
\demo{Proof} If $\bold b=\bold 0$ both sides of the equality \mythetag{27.2} 
do vanish and it is trivially fulfilled. Therefore we can assume that 
$\bold b\neq\bold 0$.\par
     It is easy to see that the vectors in two sides of the equality \mythetag{27.2} 
are collinear. For the beginning let's prove that the lengths of these two vectors 
are equal to each other. The length of the vector $\pi_{\bold a}(\bold b)$ is calculated 
according to Fig\.~27.1:
$$
\hskip -2em
|\pi_{\bold a}(\bold b)|=|\bold b_{\sssize\parallel}|
=|\bold b|\,|\cos\varphi|.
\mytag{27.3}
$$
The length of the vector in the right hand side of the formula \mythetag{27.2} is
determined by the formula itself:
$$
\biggl|\frac{(\bold b,\bold a)}{|\bold a|^2}\,
\bold a\kern 0.7pt\biggr|=\frac{|(\bold b,\bold a)|}{|\bold a|^2}\,
|\bold a|=\frac{|\bold b|\,|\bold a|\,|\cos\varphi|}{|\bold a|}
=|\bold b|\,|\cos\varphi|.
\quad
\mytag{27.4}
$$
Comparing the results of \mythetag{27.3} and \mythetag{27.4}, we conclude that
$$
\hskip -2em
|\pi_{\bold a}(\bold b)|=\biggl|\frac{(\bold b,\bold a)}{|\bold a|^2}\,
\bold a\kern 0.7pt\biggr|.
\mytag{27.5}
$$\par
      Due to \mythetag{27.5} in order to prove the equality \mythetag{27.2} it is
sufficient to prove the codirectedness of vectors 
$$
\hskip -2em
\pi_{\bold a}(\bold b)\upuparrows\frac{(\bold b,\bold a)}{|\bold a|^2}\,
\bold a.
\mytag{27.6}
$$
Since $\pi_{\bold a}(\bold b)=\bold b_{\sssize\parallel}$, again applying 
Fig\.~27.1, we consider the following three possible cases:
$$
\pagebreak
\xalignat 3
&0\leqslant\varphi<\pi/2,
&&\varphi=\pi/2,
&&\pi/2<\varphi\leqslant\pi.
\endxalignat
$$
In the first case both vectors \mythetag{27.6} are codirected with the vector
$\bold a\neq\bold 0$. Hence they are codirected with each other.\par
     In the second case both vectors \mythetag{27.6} are equal to zero. They are
codirected according to the definition~\mythedefinition{3.2}.\par
     In the third case both vectors \mythetag{27.6} are opposite to the vector
$\bold a\neq\bold 0$. Therefore they are again codirected with each other. 
The relationship \mythetag{27.6} and the theorem~\mythetheorem{27.2} in whole 
are proved.\qed\enddemo
\mydefinition{27.4} A mapping $f$ acting from the set of all free vectors to the
set of all free vectors is called a {\it linear mapping} if it possesses the 
following two properties:
\roster
\rosteritemwd=5pt
\item"1)" $f(\bold a+\bold b)=f(\bold a)+f(\bold b)$;
\item"2)" $f(\alpha\,\bold a)=\alpha\,f(\bold a)$.
\endroster
The properties 1) and 2), which should be fulfilled for any two vectors $\bold a$ 
and $\bold b$ and for any number $\alpha$, constitute a property which is called
the {\it linearity}.
\enddefinition
\mytheorem{27.3} For any nonzero vector $\bold a\neq\bold 0$ the orthogonal
projection $\pi_{\bold a}$ onto a line given by the vector $\bold a$ is a linear
mapping.
\endproclaim 
\parshape 2 0cm 10cm 5cm 5cm 
     In order to prove the theorem~\mythetheorem{27.3} we need the following 
auxiliary lemma. 
\mylemma{27.1}\parshape 1 5cm 5cm 
For any nonze\-ro vector $\bold a\neq\bold 0$ the sum of two 
vectors collinear to $\bold a$ is a vector collinear to $\bold a$ and the sum
of two vectors perpendicular to $\bold a$ is a vector perpendicular to 
$\bold a$.
\endproclaim 
\demo{Proof of the lemma \mythelemma{27.1}}\linebreak
\parshape 1 5cm 5cm 
The first proposition of the lem\-ma
is obvious. It follows immediately from the definition~\mythedefinition{5.1}.\par
\parshape 4 5cm 5cm 5cm 5cm 5cm 5cm 0cm 10cm  
     Let's prove the second proposition. \vadjust{\vskip 5pt\hbox to 
0pt{\kern 0pt\includegraphics{angemeng21.eps}\hss}
\vskip -5pt}Let $\bold b$ and $\bold c$ be two vectors, 
such that $\bold b\perp\bold a$ and 
$\bold c\perp\bold a$. In the cases $\bold b=\bold 0$, $\bold c=\bold 0$, 
$\bold b+\bold c=\bold 0$, and in the case $\bold b\parallel\bold c$ the second
proposition of the lemma is also obvious. In these cases geometric realizations
of all the three vectors $\bold b$, $\bold c$, and $\bold b+\bold c$ can be chosen
such that they lie on the same line. Such a line is perpendicular to the vector 
$\bold a$.\par
     Let's consider the case where $\bold b\nparallel\bold c$. Let's build a
geometric realization of the vector $\bold b=\overrightarrow{OB\,\,}\!$ with
initial point at some arbitrary point $O$. Then we lay the vector $\bold c
=\overrightarrow{BC\,\,}\!$ at the terminal point $B$ of the vector
$\overrightarrow{OB\,\,}\!$. Since $\bold b\nparallel\bold c$, the points $O$, 
$B$, and $C$ do not lie on a single straight line altogether. Hence they 
determine a plane. We denote this plane through $\alpha$. The sum of vectors
$\overrightarrow{OC\,\,}\!=\overrightarrow{OB\,\,}\!+\overrightarrow{BC\,\,}\!$ 
lies on this plane. It is a geometric realization for the vector $\bold b+\bold c$,
i.\.\,e\. $\overrightarrow{OC\,\,}\!=\bold b+\bold c$.\par
     We build two geometric realizations $\overrightarrow{OA\,\,}\!$ and 
$\overrightarrow{BD\,\,}\!$ for the vector $\bold a$. These are two different 
geometric vectors. From the equality $\overrightarrow{OA\,\,}\!
=\overrightarrow{BD\,\,}\!$ we derive that the lines $OA$ and $BD$ are
parallel. Due to $\bold b\perp\bold a$ and $\bold c\perp\bold a$ the line
$BD$ is perpendicular to the pair of crosswise intersecting lines $OB$ and $BC$ 
lying on the plane $\alpha$. Hence it is perpendicular to this plane. From
$BD\perp\alpha$ and $BD\parallel OA$ we derive $OA\perp\alpha$, while from 
$OA\perp\alpha$ we derive $\overrightarrow{OA\,\,}\!\perp\overrightarrow{OC\,\,}\!$. 
Hence the sum of vectors $\bold b+\bold c$ is perpendicular to the vector $\bold a$.
The lemma~\mythelemma{27.1} is proved. 
\qed\enddemo
\demo{Proof of the theorem \mythetheorem{27.3}} According to the 
definition \mythedefinition{27.4}, in order to prove the theorem
we need to verify two linearity conditions for the mapping $\pi_{\bold a}$. 
The first of these conditions in our particular case is written as the 
equality 
$$
\hskip -2em
\pi_{\bold a}(\bold b+\bold c)=\pi_{\bold a}(\bold b)
+\pi_{\bold a}(\bold c).
\mytag{27.7}
$$
Let's denote $\bold d=\bold b+\bold c$. According to the theorem~\mythetheorem{27.1},
there are expansions of the form \mythetag{27.1} for the vectors $\bold b$, 
$\bold c$, and $\bold d$:
$$
\allowdisplaybreaks
\align
\hskip -2em
\bold b&=\bold b_{\sssize\parallel}+\bold b_{\sssize\perp},
\mytag{27.8}\\
\hskip -2em
\bold c&=\bold c_{\sssize\parallel}+\bold c_{\sssize\perp},
\mytag{27.9}\\
\hskip -2em
\bold d&=\bold d_{\sssize\parallel}+\bold d_{\sssize\perp}.
\mytag{27.10}
\endalign
$$
According the same theorem~\mythetheorem{27.1}, the components of the
expansions \mythetag{27.8}, \mythetag{27.9}, \mythetag{27.10} are uniquely
fixed by the conditions
$$
\xalignat 2
&\hskip -2em
\bold b_{\sssize\parallel}\parallel\bold a,
&&\bold b_{\sssize\perp}\perp\bold a,
\mytag{27.11}\\
&\hskip -2em
\bold c_{\sssize\parallel}\parallel\bold a,
&&\bold c_{\sssize\perp}\perp\bold a,
\mytag{27.12}\\
&\hskip -2em
\bold d_{\sssize\parallel}\parallel\bold a,
&&\bold d_{\sssize\perp}\perp\bold a.
\mytag{27.13}
\endxalignat 
$$
Adding the equalities \mythetag{27.8} and \mythetag{27.9}, we get
$$
\hskip -2em
\bold d=\bold b+\bold c=
(\bold b_{\sssize\parallel}+\bold c_{\sssize\parallel})
+(\bold b_{\sssize\perp}+\bold c_{\sssize\perp}).
\mytag{27.14}
$$
Due to \mythetag{27.11} and \mythetag{27.12} we can apply the 
lemma~\mythelemma{27.1} to the components of the expansion
\mythetag{27.14}. This yields 
$$
\xalignat 2
&\hskip -2em
(\bold b_{\sssize\parallel}+\bold c_{\sssize\parallel})
\parallel\bold a,
&&(\bold b_{\sssize\perp}+\bold c_{\sssize\perp})
\perp\bold a.
\mytag{27.15}
\endxalignat 
$$
The rest is to compare \mythetag{27.14} with \mythetag{27.10} and 
\mythetag{27.15} with \mythetag{27.13}. From this comparison, applying
the theorem~\mythetheorem{27.1}, we derive the following relationships:
$$
\xalignat 2
&\hskip -2em
\bold d_{\sssize\parallel}=\bold b_{\sssize\parallel}
+\bold c_{\sssize\parallel},
&&\bold d_{\sssize\perp}=\bold b_{\sssize\perp}
+\bold c_{\sssize\perp}.
\mytag{27.16}
\endxalignat 
$$
According to the definition~\mythedefinition{27.3}, the first of the above
relationships \mythetag{27.16} is equivalent to the equality \mythetag{27.7}
which was to be verified.\par
     Let's proceed to proving the second linearity condition for the
mapping $\pi_{\bold a}$. It is written as follows:
$$
\hskip -2em
\pi_{\bold a}(\alpha\,\bold b)=\alpha\,\pi_{\bold a}(\bold b).
\mytag{27.17}
$$
Let's denote $\bold e=\alpha\,\bold b$ and then, applying the 
theorem~\mythetheorem{27.1}, write 
$$
\align
\hskip -2em
\bold b&=\bold b_{\sssize\parallel}+\bold b_{\sssize\perp},
\mytag{27.18}\\
\hskip -2em
\bold e&=\bold e_{\sssize\parallel}+\bold e_{\sssize\perp}.
\mytag{27.19}
\endalign
$$
According to the theorem~\mythetheorem{27.1}, the components of the expansions
\mythetag{27.18} and \mythetag{27.19} are uniquely fixed by the conditions
$$
\xalignat 2
&\hskip -2em
\bold b_{\sssize\parallel}\parallel\bold a,
&&\bold b_{\sssize\perp}\perp\bold a,
\mytag{27.20}\\
&\hskip -2em
\bold e_{\sssize\parallel}\parallel\bold a,
&&\bold e_{\sssize\perp}\perp\bold a.
\mytag{27.21}
\endxalignat 
$$
Let's multiply both sides of \mythetag{27.18} by $\alpha$. Then we get
$$
\hskip -2em
\bold e=\alpha\bold b=\alpha\,\bold b_{\sssize\parallel}
+\alpha\bold b_{\sssize\perp}.
\mytag{27.22}
$$
Multiplying a vector by the number $\alpha$, we get a vector collinear to the
initial vector. For this reason from \mythetag{27.20} we derive 
$$
\xalignat 2
&\hskip -2em
(\alpha\,\bold b_{\sssize\parallel})\parallel\bold a,
&&(\alpha\,\bold b_{\sssize\perp})\perp\bold a,
\mytag{27.23}
\endxalignat 
$$
Let's compare \mythetag{27.22} with \mythetag{27.19} and \mythetag{27.23} with
\mythetag{27.21}. Then, applying the theorem~\mythetheorem{27.1}, we obtain
$$
\xalignat 2
&\hskip -2em
\bold e_{\sssize\parallel}=\alpha\,\bold b_{\sssize\parallel},
&&\bold e_{\sssize\perp}=\alpha\,\bold b_{\sssize\perp}.
\mytag{27.24}
\endxalignat 
$$
According to the definition~\mythedefinition{27.3}, the first of the equalities
\mythetag{27.24} is equivalent to the required equality \mythetag{27.17}.
The theorem~\mythetheorem{27.3} is proved.
\qed\enddemo
\head
\SectionNum{28}{73} Properties of the scalar product.
\endhead
\rightheadtext{\S\,28. Properties of the scalar product.}
\mytheorem{28.1} The scalar product of vectors possesses the following four
properties which are fulfilled for any three vectors $\bold a$, $\bold b$, 
$\bold c$ and for any number $\alpha$:
\roster
\rosteritemwd=5pt
\item"1)" $(\bold a,\bold b)=(\bold b,\bold a)$;
\item"2)" $(\bold a+\bold b,\bold c)=(\bold a,\bold c)
+(\bold b,\bold c)$;
\item"3)" $(\alpha\,\bold a,\bold c)=\alpha\,(\bold a,\bold c)$;
\item"4)" $(\bold a,\bold a)\geqslant 0$ and $(\bold a,\bold a)=0$ 
implies $\bold a=\bold 0$.
\endroster
\endproclaim
\mydefinition{28.1} The property 1) in the theorem~\mythetheorem{28.1}
is called the property of {\it symmetry}; the properties 2) and 3) are called 
the properties of {\it linearity with respect to the first multiplicand};
the property 4) is called the property of {\it positivity\/}.
\enddefinition
\demo{Proof of the theorem~\mythetheorem{28.1}} The property of symmetry 1) 
is immediate from the definition~\mythedefinition{26.1} and the formula
\mythetag{26.1} in this definition. Indeed, if one of the vectors $\bold a$ 
or $\bold b$ is equal to zero, then both sides of the equality 
$(\bold a,\bold b)=(\bold b,\bold a)$ equal to zero. Hence the equality is
fulfilled in this case.\par
     In the case of the nonzero vectors $\bold a$ and $\bold b$ the angle
$\varphi$ is determined by the pair of vectors $\bold a$ and $\bold b$ according
to Fig\.~26.1, it does not depend on the order of vectors in this pair. Therefore 
the equality $(\bold a,\bold b)=(\bold b,\bold a)$ in this case is reduced to 
the equality
$$
|\bold a|\,|\bold b|\,\cos\varphi=|\bold b|\,|\bold a|\,\cos\varphi.
$$
It is obviously fulfilled since $|\bold a|$ and $|\bold b|$ are numbers complemented
by some measure units depending on the physical nature of the vectors $\bold a$ and 
$\bold b$.\par
      Let's consider the properties of linearity 2) and 3). If $\bold c=\bold 0$,
then both sides of the equalities $(\bold a+\bold b,\bold c)=(\bold a,\bold c)
+(\bold b,\bold c)$ and $(\alpha\,\bold a,\bold b)=\alpha\,(\bold a,\bold b)$ are 
equal to zero. Hence these equalities are fulfilled in this case.\par
     If $\bold c\neq\bold 0$, then we apply the theorems~\mythetheorem{27.2} and
\mythetheorem{27.3}. From the theorems~\mythetheorem{27.2} we derive the following
equalities:
$$
\gather
\pi_{\bold c}(\bold a+\bold b)-\pi_{\bold c}(\bold a)
-\pi_{\bold c}(\bold b)=\frac{(\bold a+\bold b,\bold c)
-(\bold a,\bold c)-(\bold b,\bold c)}{|\bold c|^2}\,
\bold c,\\
\vspace{2ex}
\pi_{\bold c}(\alpha\,\bold a)-\alpha\,\pi_{\bold c}(\bold a)
=\frac{(\alpha\,\bold a,\bold c)
-\alpha\,(\bold a,\bold c)}{|\bold c|^2}\,\bold c.
\endgather
$$ 
Due to the theorem~\mythetheorem{27.3} the mapping $\pi_{\bold c}$ is a linear
mapping (see Definition~\mythedefinition{27.4}). Therefore the left hand sides 
of the above equalities are zero. Now due to $\bold c\neq\bold 0$ we conclude
that the numerators of the fractions in their right hand sides are also zero.
This fact proves the properties 2) and 3) from the theorem~\mythetheorem{28.1}
are valid in the case $\bold c\neq\bold 0$.\par
     According to the definition~\mythedefinition{26.1} the scalar product 
$(\bold a,\bold a)$ is equal to zero for $\bold a=\bold 0$. Otherwise, if
$\bold a\neq\bold 0$, the formula \mythetag{26.1} is applied where we should set
$\bold b=\bold a$. This yields $\varphi=0$ and 
$$
(\bold a,\bold a)=|\bold a|^2>0.
$$ 
This inequality proves the property 4) and completes the proof of the 
theorem~\mythetheorem{28.1} in whole. 
\qed\enddemo
\mytheorem{28.2} Apart from the properties 1)--4), the scalar product
of vectors possesses the following two properties fulfilled for any three
vectors $\bold a$, $\bold b$, $\bold c$ and for any number $\alpha$:
\roster
\rosteritemwd=5pt
\item"5)" $(\bold c,\bold a+\bold b)=(\bold c,\bold a)
+(\bold a,\bold b)$;
\item"6)" $(\bold c,\alpha\,\bold a)=\alpha\,(\bold c,\bold a)$.
\endroster
\endproclaim
\mydefinition{28.2} The properties 5) and 6) in the theorem~\mythetheorem{28.2}
are called the properties of  {\it linearity with respect to the second 
multiplicand}.
\enddefinition
     The properties 5) and 6) are easily derived from the properties 2) and 3) 
by applying the property 1). Indeed, we have
$$
\gather
(\bold c,\bold a+\bold b)=(\bold a+\bold b,\bold c)
=(\bold a,\bold c)+(\bold b,\bold c)
=(\bold c,\bold a)+(\bold c,\bold b),\\
(\bold c,\alpha\,\bold a)=(\alpha\,\bold a,\bold c)
=\alpha\,(\bold a,\bold c)=\alpha\,(\bold c,\bold a).
\endgather
$$
These calculations prove the theorem~\mythetheorem{28.2}. 
\head
\SectionNum{29}{75} Calculation of the scalar product through the 
coordinates of vectors in a skew-angular basis.
\endhead
\rightheadtext{\S\,29. \dots in a skew-angular basis.}
     Let $\bold e_1,\,\bold e_2,\,\bold e_3$ be some arbitrary basis in the
space $\Bbb E$. According to the definition~\mythedefinition{18.1}, this is
an ordered triple of non-coplanar vectors. The arbitrariness of a basis means
that no auxiliary restrictions are imposed onto the vectors $\bold e_1,\,\bold e_2,
\,\bold e_3$, except for non-coplanarity. In particular, this means that the
angles between the vectors $\bold e_1,\,\bold e_2,\,\bold e_3$ in an arbitrary
basis should not be right angles. For this reason such a basis is called a
{\it skew-angular basis\/} and abbreviated as SAB. 
\mydefinition{29.1} In this book a {\it skew-angular basis\/} (SAB) is 
understood as an {\it arbitrary basis}.
\enddefinition
     Thus, let $\bold e_1,\,\bold e_2,\,\bold e_3$ be some skew-angular
basis in the space $\Bbb E$ and let $\bold a$ and $\bold b$ be two free vectors
given by its coordinates in this basis. We write this fact as follows:
$$
\xalignat 2
&\hskip -2em
\bold a=\Vmatrix a^1\\\vspace{0.4ex}a^2\\\vspace{0.4ex}a^3
\endVmatrix,
&&\bold b=\Vmatrix b^1\\\vspace{0.4ex}b^2\\\vspace{0.4ex}b^3
\endVmatrix
\mytag{29.1}
\endxalignat
$$
Unlike \mythetag{21.1}, instead of the arrow sign in \mythetag{29.1} we use 
the equality sign. Doing this, we emphasize the fact that once a basis is fixed,
vectors are uniquely identified with their coordinates.\par
     The conditional writing \mythetag{29.1} means that the vectors $\bold a$ 
and $\bold b$ are presented by the following expansions:
$$
\xalignat 2
&\hskip -2em
\bold a=\sum^3_{i=1}a^i\,\bold e_i,
&&\bold b=\sum^3_{j=1}b^{\kern 0.6pt j}\,\bold e_j.
\mytag{29.2}
\endxalignat
$$
Substituting \mythetag{29.2} into the scalar product $(\bold a,\bold b)$, we get
$$
\hskip -2em
(\bold a,\bold b)=\biggl(\,\sum^3_{i=1}a^i\,\bold e_i,
\,\sum^3_{j=1}b^{\kern 0.6pt j}\,\bold e_j\biggr).
\mytag{29.3}
$$
In order to transform the formulas \mythetag{29.3} we apply the properties 2) 
and 5) of the scalar product from the theorems~\mythetheorem{28.1} and
\mythetheorem{28.2}. Due to these properties we can take the summation signs
over $i$ and $j$ out of the brackets of the scalar product:
$$
\hskip -2em
(\bold a,\bold b)=\shave{\sum^3_{i=1}\sum^3_{j=1}}\,(a^i\,\bold e_i,
b^{\kern 0.6pt j}\,\bold e_j).
\mytag{29.4}
$$
Then we apply the properties 3) and 6) from the theorems~\mythetheorem{28.1} 
and \mythetheorem{28.2}. Due to these properties we can take the numeric factors
$a^i$ and $b^{\kern 0.6pt j}$ out of the brackets of the scalar product in
\mythetag{29.4}:
$$
\hskip -2em
(\bold a,\bold b)=\sum^3_{i=1}\sum^3_{j=1}a^i\,b^{\kern 0.6pt j}
\,(\bold e_i,\bold e_j).
\mytag{29.5}
$$\par
      The quantities $(\bold e_i,\bold e_j)$ in the formula \mythetag{29.5} depend
on a basis $\bold e_1,\,\bold e_2,\,\bold e_3$, namely on the lengths of the basis
vectors and on the angles between them. They do not depend on the vectors $\bold a$ 
and $\bold b$. The quantities $(\bold e_i,\bold e_j)$ constitute an array of nine
numbers
$$
\hskip -2em
g_{ij}=(\bold e_i,\bold e_j)
\mytag{29.6}
$$
enumerated by two lower indices. The components of the array
\mythetag{29.6} are usually arranged into a square matrix:
$$
\hskip -2em
G=\Vmatrix 
g_{11} & g_{12} & g_{13}\\
\vspace{1ex}
g_{21} & g_{22} & g_{23}\\
\vspace{1ex}
g_{31} & g_{32} & g_{33}
\endVmatrix
\mytag{29.7}
$$
\mydefinition{29.2} The matrix \mythetag{29.7} with the components \mythetag{29.6} 
is called the Gram matrix of a basis $\bold e_1,\,\bold e_2,\,\bold e_3$.
\enddefinition
     Taking into account the notations \mythetag{29.6}, we write the formula 
\mythetag{29.5} in the following way:
$$
\hskip -2em
(\bold a,\bold b)=\sum^3_{i=1}\sum^3_{j=1}a^i\,b^{\kern 0.6pt j}
\,g_{ij}.
\mytag{29.8}
$$
\mydefinition{29.3} The formula \mythetag{29.8} is called the {\it formula
for calculating the scalar product through the coordinates of vectors in a 
skew-angular basis}.
\enddefinition
     The formula \mythetag{29.8} can be written in the matrix form
$$
(\bold a,\bold b)=
\Vmatrix a^1 & a^2 & a^3
\endVmatrix\cdot\Vmatrix 
g_{11} & g_{12} & g_{13}\\
\vspace{1ex}
g_{21} & g_{22} & g_{23}\\
\vspace{1ex}
g_{31} & g_{32} & g_{33}
\endVmatrix\cdot
\Vmatrix b^1\\\vspace{0.4ex}b^2\\\vspace{0.4ex}b^3
\endVmatrix\quad
\mytag{29.9}
$$
Note that the coordinate column of the vector $\bold b$ in the formula 
\mythetag{29.9} is used as it is, while the coordinate column of the vector 
$\bold a$ is transformed into a row. Such a transformation is known as 
{\it matrix transposing} (see \mycite{7}). 
\mydefinition{29.4} A transformation of a rectangular matrix under which
the element in the intersection of $i$-th row and $j$-th column is taken
to the intersection of $j$-th row and $i$-th column is called the
{\it matrix transposing}. It is denoted by means of the sign $\top$. 
In the {\TeX} and {La\TeX} computer packages this sign is coded by the
operator {\tt $\backslash$top}. 
\enddefinition
     The operation of matrix transposing can be understood as the mirror 
reflection with respect to the main diagonal of a matrix\par
     Taking into account the notations \mythetag{29.1}, \mythetag{29.7}, and
the de\-finition~\mythedefinition{29.4}, we can write the matrix formula
\mythetag{29.9} as follows:
$$
\hskip -2em
(\bold a,\bold b)=\bold a^{\sssize\top}\kern -2pt\cdot G\cdot\bold b.
\mytag{29.10}
$$
In the right hand side of the formula \mythetag{29.10} the vectors $\bold a$ 
and $\bold b$ are presented by their coordinate columns, while the transformation
of one of them into a row is written through the matrix transposing.
\myexercise{29.1} Show that for an arbitrary rectangular matrix $A$ the equality
$(A^{\sssize\top})^{\sssize\top}=A$ is fulfilled. 
\endproclaim
\myexercise{29.2} Show that for the product of two matrices $A$ and $B$ the 
equality \pagebreak $(A\kern 1pt \cdot B)^{\sssize\top}=B^{\sssize\top}
\cdot A^{\sssize\top}$ is fulfilled. 
\endproclaim
\myexercise{29.3} Define the Gram matrices for bases on a line and for bases on 
a plane. Write analogs of the formulas \mythetag{29.8}, \mythetag{29.9}, and 
\mythetag{29.10} for the scalar product of vectors lying on a line and on a plane.
\endproclaim
\head
\SectionNum{30}{79} Symmetry of the Gram matrix.
\endhead
\rightheadtext{\S\,30. Symmetry of the Gram matrix.}
\mydefinition{30.1} A square matrix $A$ is called {\it symmetric}, if it is 
preserved under transposing, i\.\,e\. if the following equality is fulfilled:
$A^{\sssize\top}=A$. 
\enddefinition
     Gram matrices possesses many important properties. One of these properties
is their symmetry. 
\mytheorem{30.1} The Gram matrix $G$ of any basis\/ $\bold e_1,\,\bold e_2,
\,\bold e_3$ in the space\/ $\Bbb E$ is symmetric. 
\endproclaim
\demo{Proof} According to the definition~\mythedefinition{30.1} the symmetry of $G$ 
is expressed by the formula $G^{\sssize\top}=G$. According to the 
definition~\mythedefinition{29.4}, the equality $G^{\sssize\top}=G$ is equivalent
to the relationship
$$
\hskip -2em
g_{ij}=g_{j\kern 0.5pt i}
\mytag{30.1}
$$
for the components of the matrix $G$. As for the relationship \mythetag{30.1}, upon
applying \mythetag{29.6}, it reduces to the equality
$$
(\bold e_i,\bold e_j)=(\bold e_j,\bold e_i)
$$
which is fulfilled due to the symmetry of the scalar product (see 
Theorem~\mythetheorem{28.1} and Definition~\mythedefinition{28.1}).
\qed\enddemo
     Note that the coordinate columns of the vectors $\bold a$ and $\bold b$
enter the right hand side of the formula \mythetag{29.10} in somewhat unequal way
--- one of them is transposed, the other is not transposed. The symmetry of the
matrix $G$ eliminates this difference. Redesignating the indices $i$ and $j$ in the
double sum \mythetag{29.8} and taking into account the relationship \mythetag{30.1} 
for the components of the Gram matrix, we get 
$$
\hskip -2em
(\bold a,\bold b)=\sum^3_{j=1}\sum^3_{i=1}a^{\kern 0.2pt j}
\,b^{\kern 0.4pt i}\,g_{j\kern 0.5pt i}=\sum^3_{i=1}\sum^3_{j=1}
b^{\kern 0.4pt i}\,a^{\kern 0.2pt j}\,g_{ij}.
\mytag{30.2}
$$
In the matrix form the formula \mythetag{30.2} is written as follows:
$$
\hskip -2em
(\bold a,\bold b)=\bold b^{\sssize\top}\kern -2pt\cdot G\cdot\bold a.
\mytag{30.3}
$$
The formula \mythetag{30.3} is analogous to the formula \mythetag{29.10}, but in this
formula the coordinate column of the vector $\bold b$ is transposed, while the coordinate
column of the vector $\bold a$ is not transposed. 
\myexercise{30.1} Formulate and prove a theorem analogous to the theorem~\mythetheorem{30.1} 
for bases on a plane. Is it necessary to formulate such a theorem for bases on a line. 
\endproclaim
\head
\SectionNum{31}{80} Orthonormal basis.
\endhead
\rightheadtext{\S\,31. Orthonormal basis.}
\mydefinition{31.1} A basis on a straight line consisting of a nonzero vector $\bold e$ 
is called an {\it orthonormal basis}, if $\bold e$ is a unit vector, i\.\,e\. if
$|\bold e|=1$.
\enddefinition
\mydefinition{31.2} A basis on a plane, consisting of two non-collinear vectors
$\bold e_1,\,\bold e_2$, is called an {\it orthonormal basis}, if the vectors $\bold e_1$ 
and $\bold e_2$ are two vectors of the unit lengths perpendicular to each other. 
\enddefinition
\mydefinition{31.3} A basis in the space $\Bbb E$ consisting of three non-coplanar
vectors $\bold e_1,\,\bold e_2,\,\bold e_3$ is called an {\it orthonormal basis} if the
vectors $\bold e_1$, $\bold e_2$, $\bold e_3$ are three vectors of the unit lengths
perpendicular to each other. 
\enddefinition
     In order to denote an orthonormal basis in each of the there cases listed above 
we use the abbreviation ONB. According to the definition~\mythedefinition{29.1} the
orthonormal basis is not opposed to a skew-angular basis SAB, it is a special case 
of such a basis.\par 
     Note that the unit lengths of the basis vectors of an orthonormal basis in the
definitions~\mythedefinition{31.1}, \mythedefinition{31.2}, and 
\mythedefinition{31.2} mean that their lengths are not one centimeter, not one meter,
not one kilometer, but the pure numeric unity. For this reason all geometric 
realizations of such vectors are conditionally geometric (see \S\,2). Like velocity 
vectors, acceleration vectors, and many other physical quantities, basis vectors of 
an orthonormal basis can be drawn only upon choosing some scaling factor. Such a factor
in this particular case is needed for to transform the numeric unity into a unit of 
length.\par
\head
\SectionNum{32}{81} Gram matrix of an orthonormal basis.
\endhead
\rightheadtext{\S\,32. \dots of an orthonormal basis.}
     Let $\bold e_1,\,\bold e_2,\,\bold e_3$ be some orthonormal basis in the space
$\Bbb E$. According to the definition~\mythedefinition{31.3} the vectors $\bold e_1$, 
$\bold e_2$, $\bold e_3$ satisfy the following relationships:
$$
\xalignat 3
&\hskip -2em
|\bold e_1|=1, &&|\bold e_2|=1, &&|\bold e_3|=1,\quad\\
\vspace{-1.5ex}
\mytag{32.1}\\
\vspace{-1.5ex}
&\hskip -2em
\bold e_1\perp\bold e_2, &&\bold e_2\perp\bold e_3, 
&&\bold e_3\perp\bold e_1.\quad 
\endxalignat
$$
Applying \mythetag{32.1} and \mythetag{29.6}, we find the components of the Gram matrix
for the orthonormal basis $\bold e_1,\,\bold e_2,\,\bold e_3$:
$$
\hskip -2em
g_{ij}=\cases
1&\text{for \ }i=j,\\
0&\text{for \ }i\neq j.
\endcases
\mytag{32.2}
$$
From \mythetag{32.2} we immediately derive the following theorem.
\mytheorem{32.1} The Gram matrix \mythetag{29.7} of any orthonormal basis is a unit
matrix:
$$
G=\Vmatrix 1 & 0 & 0\\
\vspace{0.5ex}
0 & 1 & 0\\
\vspace{0.5ex}
0 & 0 & 1
\endVmatrix=1.
\mytag{32.3}
$$ 
\endproclaim
     Let's recall that the components of a unit matrix constitute a numeric array
$\delta$ which is called the Kronecker symbol (see Definition~\mythedefinition{23.1}). 
Therefore, taking into account \mythetag{23.3}, the equality \mythetag{32.2} can be
written as:
$$
\hskip -2em
g_{ij}=\delta_{ij}.
\mytag{32.4}
$$
The Kronecker symbol in \mythetag{32.4} inherits the lower position of indices 
from $g_{ij}$. Therefore it is different from the Kronecker symbol in \mythetag{23.4}. 
Despite being absolutely identical, the components of the unit matrix in \mythetag{32.3} 
and in \mythetag{23.1} are of absolutely different nature. Having negotiated to use
indices on upper and lower levels (see Definition~\mythedefinition{20.1}), now we are 
able to reflect this difference in denoting these components.
\head
\SectionNum{33}{82} Calculation of the scalar product through the coordinates 
of vectors in an orthonormal basis.
\endhead
\rightheadtext{\S\,33. \dots in an orthonormal basis.}
      According to the definition~\mythedefinition{29.1} the term {\it skew-angular
basis} is used as a synonym of an arbitrary basis. For this reason an orthonormal
basis is a special case of a skew-angular basis and we can use the formula 
\mythetag{29.8}, taking into account \mythetag{32.4}:
$$
\hskip -2em
(\bold a,\bold b)=\sum^3_{i=1}\sum^3_{j=1}a^i\,b^{\kern 0.4pt i}
\,\delta_{ij}.
\mytag{33.1}
$$
In calculating the sum over $j$ in \mythetag{33.1}, it is the inner sum here, the
index $j$ runs over three values and only for one of these three values, where $j=i$,
the Kronecker symbol $\delta_{ij}$ is nonzero. For this reason wee can retain only
one summand of the inner sum over $j$ in \mythetag{33.1}, omitting other two summands:
$$
\hskip -2em
(\bold a,\bold b)=\sum^3_{i=1}a^i\,b^{\kern 0.6pt i}
\,\delta_{ii}.
\mytag{33.2}
$$
We know that $\delta_{ii}=1$. Therefore the formula \mythetag{33.2} turns to
$$
\hskip -2em
(\bold a,\bold b)=\sum^3_{i=1}a^i\,b^{\kern 0.6pt i}.
\mytag{33.3}
$$
\mydefinition{33.1} The formula \mythetag{33.3} is called the {\it formula for
calculating the scalar product through the coordinates of vectors in an orthonormal
basis}.
\enddefinition
     Note that the sums in the formula \mythetag{29.8} satisfy the index setting rule 
from the definition~\mythedefinition{24.8}, while the sum in the formula \mythetag{33.3} 
breaks this rule. In this formula the summation index has two entries and both of them 
are in the upper positions. This is a peculiarity of an orthonormal basis. It is more
symmetric as compared to a general skew-angular basis and this symmetry hides some
rules that reveal in a general non-symmetric bases.\par
     The formula \mythetag{33.3} has the following matrix form:
$$
\hskip -2em
(\bold a,\bold b)=
\Vmatrix a^1 & a^2 & a^3
\endVmatrix\cdot
\Vmatrix b^1\\\vspace{0.4ex}b^2\\\vspace{0.4ex}b^3
\endVmatrix.
\mytag{33.4}
$$
Taking into account the notations \mythetag{29.1} and taking into account the
definition~\mythedefinition{29.4}, the formula \mythetag{33.4} can be abbreviated 
to
$$
\hskip -2em
(\bold a,\bold b)=\bold a^{\sssize\top}\kern -2pt\cdot\bold b.
\mytag{33.5}
$$
The formula \mythetag{33.4} can be derived from the formula \mythetag{29.9},
while the formula \mythetag{33.5} can be derived from the formula \mythetag{29.10}.
\head
\SectionNum{34}{83} Right and left triples of vectors.
The concept of orientation.
\endhead
\rightheadtext{\S\,34. Right and left triples of vectors.}
\mydefinition{34.1} An {\it ordered triple of vectors\/} is a list of three vectors 
for which the order of listing vectors is fixed. 
\enddefinition
\mydefinition{34.2}\parshape 3 0cm 10cm 0cm 10cm 5cm 5cm \ An ordered triple of 
non-coplanar vectors $\bold a_1,\,\bold a_2,\,\bold a_3$ is called a {\it right 
triple} if, when observing from the end of the third vector, the shortest rotation 
from the first vector toward the second vector is seen as a counterclockwise 
rotation.
\enddefinition
\parshape 7 5cm 5cm 5cm 5cm 5cm 5cm 5cm 5cm 5cm 5cm 5cm 5cm 0cm 10cm 
     In the definition~\mythedefinition{34.2} we implicitly assume that the
geometric realizations of the vectors $\bold a_1$, $\bold a_2,\,\bold a_3$ with 
some common initial point are considered \vadjust{\vskip 5pt\hbox to 
0pt{\kern 0pt\includegraphics{angemeng22.eps}\hss}\vskip -5pt}as
it is shown in Fig\.~34.1. 
\mydefinition{34.3} An ordered triple of non-coplanar vectors $\bold a_1,
\,\bold a_2,\,\bold a_3$ is called a {\it left triple} if, when observing from 
the end of the third vector, the shortest rotation from the first vector toward 
the second vector is seen as a clockwise rotation.
\enddefinition
     A given rotation about a given axis when observing from a given position 
could be either a clockwise rotation or a counterclockwise rotation. No other
options are available. For this reason each ordered triple of non-coplanar vectors 
is either left or right. No other triples sorted by this criterion are available. 
\mydefinition{34.4} The property of ordered triples of non-coplanar vectors to be
left or right is called their {\it orientation}.
\enddefinition
\head
\SectionNum{35}{84} Vector product.
\endhead
\rightheadtext{\S\,35. Vector product.}
     Let $\bold a$ and $\bold b$ be two non-collinear free vectors. Let's lay
their geometric realizations $\bold a=\overrightarrow{OA\,\,}\!$ and
$\bold b=\overrightarrow{OB\,\,}\!$ at some arbitrary point $O$. In this case
the vectors $\bold a$ and $\bold b$ define a plane $AOB$ and lie on this plane. 
The angle $\varphi$ between the vectors $\bold a$ and $\bold b$ is determined
according to Fig\.~26.1. Due to $\bold a\nparallel\bold b$ this angle ranges
in the interval $0<\varphi<\pi$ and hence $\sin\varphi\neq 0$.
\par
\parshape 4 0cm 10cm 0cm 10cm 0cm 10cm 5cm 5cm 
% 0cm 10cm 0cm 10cm 0cm 10cm
% 0cm 10cm 0cm 10cm 5cm 5cm 
     Let's draw a line through the point $O$ perpendicular to the plane $AOB$
and denote this line through $c$. The line $c$ is perpendicular to the vectors 
$\bold a=\overrightarrow{OA\,\,}\!$ and $\bold b=\overrightarrow{OB\,\,}\!$:
$$
\hskip -2em
\aligned
&c\perp\bold a,\\
&c\perp\bold b.
\endaligned
\mytag{35.1}
$$ 
It is clear that the conditions \mythetag{35.1} fix a unique line $c$ passing 
through the point $O$ \vadjust{\vskip 5pt\hbox to 
0pt{\kern 0pt\includegraphics{angemeng23.eps}\hss}
\vskip -5pt}(see Theorems 1.1 and 1.3 in Chapter 
\uppercase\expandafter{\romannumeral 4} of the book \mycite{6}). 
\par
\parshape 4 5cm 5cm 5cm 5cm 5cm 5cm 0cm 10cm 
There are two directions on the line $c$. In Fig\~35.1 they are given by the
vectors $\bold c$ and $\tilde\bold c$. The vectors $\bold a,\,\bold b,\,\bold c$ 
constitute a right triple, while $\bold a,\,\bold b,\,\tilde\bold c$ is a left
triple. So, specifying the orientation chooses one of two possible directions 
on the line $c$.
\mydefinition{35.1} The {\it vector product\/} of two non-collinear vectors
$\bold a$ and $\bold b$ is a vector $\bold c=[\bold a,\bold b]$ which is
determined by the following three conditions:
\roster
\rosteritemwd=5pt
\item"1)" $\bold c\perp\bold a$ and $\bold c\perp\bold b$;
\item"2)" the vectors $\bold a,\,\bold b,\,\bold c$ form a right triple;
\item"3)" $|\bold c|=|\bold a|\,|\bold b|\,\sin\varphi$.
\endroster
In the case of collinear vectors $\bold a$ and $\bold b$ their vector product 
$[\bold a,\bold b]$ is taken to be zero by definition.
\enddefinition
     A comma is the multiplication sign in the writing the vector product,
not by itself, but together with square brackets surrounding the whole
expression. These brackets are natural delimiters for multiplicands: the first 
multiplicand is an expression between the opening bracket and the comma, while 
the second multiplicand is an expression between the comma and the closing 
bracket. Therefore in complicated expressions no auxiliary delimiters are required. 
For example, in the formula
$$
[\bold a+\bold b,\bold c+\bold d]
$$
the sums $\bold a+\bold b$ and $\bold c+\bold d$ are calculated first, then the
vector multiplication is performed.\par
     {\bf A remark}. Often the vector product is written as $\bold a\times\bold b$. 
Even the special term {\tencyr\char '074}cross product{\tencyr\char '076} is used. 
However, to my mind, this notation is not good. It is misleading since the cross
sign is sometimes used for denoting the product of numbers when a large formula
is split into several lines.\par
     {\bf A remark}. The physical nature of the vector product $[\bold a,\bold b]$ 
often differs from the nature of its multiplicands $\bold a$ and $\bold b$. Even
if the lengths of the vectors $\bold a$ and $\bold b$ are measured in length units, 
the length of their product $[\bold a,\bold b]$ is measured in units of area. 
\myexercise{35.1} Show that the vector product $\bold c=[\bold a,\bold b]$ of two
free vectors $\bold a$ and $\bold b$ is a free vector and, being a free vector,
it does not depend on where the point $O$ in Fig\.~35.1 is placed.
\endproclaim
\head
\SectionNum{36}{86} Orthogonal projection onto a plane.
\endhead
\rightheadtext{\S\,36. Orthogonal projection onto a plane.}
     Let $\bold a\neq\bold 0$ be some nonzero free vector. According to the
theorem~\mythetheorem{27.1}, each free vector $\bold b$ has the expansion
$$
\hskip -2em
\bold b=\bold b_{\sssize\parallel}+\bold b_{\sssize\perp}
\mytag{36.1}
$$
relative to the vector $\bold a$, where the vector $\bold b_{\sssize\parallel}$ 
is collinear to the vector $\bold a$, while the vector $\bold b_{\sssize\perp}$ 
is perpendicular to the vector $\bold a$. Recall that through $\pi_{\bold a}$ 
we denoted a mapping that associates each vector $\bold b$ with its component 
$\bold b_{\sssize\parallel}$ in the expansion \mythetag{36.1}. Such a mapping
was called the orthogonal projection onto the direction of the vector
$\bold a\neq\bold 0$ (see Definition~\mythedefinition{27.3}).
\mydefinition{36.1} The mapping $\pi_{\!{\sssize\perp}\bold a}$ that associates
each free vector $\bold b$ with its perpendicular component $\bold b_{\sssize\perp}$ 
in the expansion \mythetag{36.1} is called the {\it orthogonal projection onto
a plane perpendicular to the vector $\bold a\neq\bold 0$} or, more exactly, the
{\it orthogonal projection onto the orthogonal complement of the vector
$\bold a\neq\bold 0$}.
\enddefinition
\mydefinition{36.2} The {\it orthogonal complement\/} \pagebreak of a free vector 
$\bold a$ is the collection of all free vectors $\bold x$ perpendicular to 
$\bold a$:
$$
\hskip -2em
\alpha=\{\bold x\!:\,\bold x\perp\bold a\}.
\mytag{36.2}
$$
\enddefinition
\parshape 3 0cm 10cm 0cm 10cm 5cm 5cm 
     The orthogonal complement \mythetag{36.2} of a nonzero vector $\bold a\neq\bold 0$ 
can be visualized as a plane if we choose one of its geometric realizations $\bold a
=\overrightarrow{OA\,\,}\!$. \vadjust{\vskip 5pt\hbox to 
0pt{\kern 0pt\includegraphics{angemeng24.eps}\hss}
\vskip -5pt}Indeed, let's lay various vectors perpendicular to $\bold a$ at the point 
$O$. The ending points of such vectors fill the plane $\alpha$ shown in Fig\.~36.1.
\par
\parshape 6 5cm 5cm 5cm 5cm 5cm 5cm 5cm 5cm 5cm 5cm 0cm 10cm 
     The properties of the orthogonal projections onto a line $\pi_{\bold a}$ from the
definition~\mythedefinition{27.1} and the orthogonal projections onto a\linebreak plane 
$\pi_{\!{\sssize\perp}\bold a}$ from the definition\linebreak \mythedefinition{36.1} 
are very similar. Indeed, we have the following theorem. 
\mytheorem{36.1} For any nonzero vector $\bold a\neq\bold 0$ the orthogonal projection
$\pi_{\!{\sssize\perp}\bold a}$ onto a plane perpendicular to the vector $\bold a$ is a 
linear mapping.
\endproclaim 
\demo{Proof} In order to prove the theorem~\mythetheorem{36.1} we write the relationship 
\mythetag{36.1} as follows:
$$
\hskip -2em
\bold b=\pi_{\bold a}(\bold b)+\pi_{\!{\sssize\perp}
\bold a}(\bold b).
\mytag{36.3}
$$
The relationship \mythetag{36.3} is an identity, it is fulfilled for any vector $\bold b$. 
First we replace the vector $\bold b$ by $\bold b+\bold c$ in \mythetag{36.3}, then we
replace $\bold b$ by $\alpha\,\bold b$ in \mythetag{36.3}. As a result we get two 
relationships
$$
\align
&\hskip -2em
\pi_{\!{\sssize\perp}\bold a}(\bold b+\bold c)
=\bold b+\bold c-\pi_{\bold a}(\bold b+\bold c),
\mytag{36.4}\\
&\hskip -2em
\pi_{\!{\sssize\perp}\bold a}(\alpha\,\bold b)
=\alpha\,\bold b-\pi_{\bold a}(\alpha\,\bold b).
\mytag{36.5}
\endalign
$$\par
    Due to the theorem~\mythetheorem{27.3} the mapping $\pi_{\bold a}$ is a linear
mapping. For this reason the relationships \mythetag{36.4} and \mythetag{36.5} can 
be transformed into the following two relationships:
$$
\align
&\hskip -2em
\pi_{\!{\sssize\perp}\bold a}(\bold b+\bold c)
=\bold b-\pi_{\bold a}(\bold b)+\bold c-\pi_{\bold a}(\bold c),
\mytag{36.6}\\
&\hskip -2em
\pi_{\!{\sssize\perp}\bold a}(\alpha\,\bold b)
=\alpha\,(\bold b-\pi_{\bold a}(\bold b)).
\mytag{36.7}
\endalign
$$
The rest is to apply the identity \mythetag{36.3} to the relationships 
\mythetag{36.6} and \mythetag{36.7}. As a result we get
$$
\align
&\hskip -2em
\pi_{\!{\sssize\perp}\bold a}(\bold b+\bold c)
=\pi_{\!{\sssize\perp}\bold a}(\bold b)+
\pi_{\!{\sssize\perp}\bold a}(\bold c),
\mytag{36.8}\\
&\hskip -2em
\pi_{\!{\sssize\perp}\bold a}(\alpha\,\bold b)
=\alpha\,\pi_{\!{\sssize\perp}\bold a}(\bold b).
\mytag{36.9}
\endalign
$$
The relationships \mythetag{36.8} and \mythetag{36.9} are exactly the linearity
conditions from the definition~\mythedefinition{27.4} written for 
the mapping $\pi_{\!{\sssize\perp}\bold a}$. The theorem~\mythetheorem{36.1} is
proved. 
\qed\enddemo
\head
\SectionNum{37}{88} Rotation about an axis.
\endhead
\rightheadtext{\S\,37. Rotation about an axis.}
\parshape 7 0cm 10cm 0cm 10cm 0cm 10cm 0cm 10cm 0cm 10cm 0cm 10cm 
5cm 5cm
     Let $\bold a\neq\bold 0$ be some nonzero free vector and let $\bold b$ be some
arbitrary free vector. Let's lay the vector $\bold b=\overrightarrow{BO\,\,}\!$ at
some arbitrary point $B$. Then we lay the vector $\bold a=\overrightarrow{OA\,\,}\!$
at the terminal point of the vector $\overrightarrow{BO\,\,}\!$. The vector 
$\bold a=\overrightarrow{OA\,\,}\!$ is nonzero. For this reason it defines a line
$OA$. We take this line for the rotation axis. Let's denote through 
$\theta_{\bold a}^\varphi$ the rotation of the space $\Bbb E$ about the axis $OA$ 
by the angle $\varphi$ (see Fig\.~37.1). The vector $\bold a=\overrightarrow{OA\,\,}\!$ 
fixes one of two directions on the rotation axis. At the same time this vector fixes 
the positive direction of rotation about the axis $OA$.
\mydefinition{37.1} \parshape 7 5cm 5cm 5cm 5cm 5cm 5cm 5cm 5cm 5cm 5cm 5cm 5cm 
0cm 10cm 
The rotation about an axis $OA$ with with fixed direction
\vadjust{\vskip 5pt\hbox to 
0pt{\kern 0pt\includegraphics{angemeng25.eps}\hss}
\vskip -5pt}$\bold a=\overrightarrow{OA\,\,}\!$ on it is called a positive rotation
if, being observed from the terminal point of the vector $\overrightarrow{OA\,\,}\!$,
i\.\,e\. when looking from the point $A$ toward the point $O$, it occurs in the 
counterclockwise direction. 
\enddefinition
     Taking into account the definition~\mythedefinition{37.1}, we can consider the
rotation angle $\varphi$ as a signed quantity. If $\varphi>0$, the rotation 
$\theta_{\bold a}^\varphi$ occurs in the positive direction with respect to the 
vector $\bold a$, if $\varphi<0$, it occurs in the negative direction.\par
     Let's apply the rotation mapping $\theta_{\bold a}^\varphi$ to the vectors 
$\bold a=\overrightarrow{OA\,\,}\!$ and $\bold b=\overrightarrow{BO\,\,}\!$ in
Fig\.~37.1. The points $A$ and $O$ are on the rotation axis. For this reason 
under the rotation $\theta_{\bold a}^\varphi$ the points $A$ and $O$ stay at their 
places and the vector $\bold a=\overrightarrow{OA\,\,}\!$ does not change. As for
the vector $\bold b=\overrightarrow{BO\,\,}\!$, it is mapped onto another vector
$\overrightarrow{\tilde BO\,\,}\!$. Now, applying parallel translations, we can 
replicate the vector $\overrightarrow{\tilde BO\,\,}\!$ up to a free vector 
$\tbb=\overrightarrow{\tilde BO\,\,}\!$ (see Definitions~\mythedefinition{4.1} and
\mythedefinition{4.2}). The vector $\tbb$ is said to be produced from the vector
$\bold b$ by applying the mapping $\theta_{\bold a}^\varphi$ and is written as
$$
\hskip -2em
\tbb=\theta_{\bold a}^\varphi(\bold b). 
\mytag{37.1}
$$
\mylemma{37.1} The free vector $\tbb=\theta_{\bold a}^\varphi(\bold b)$ 
in \mythetag{37.1} produced from a free vector $\bold b$ by means of the
rotation mapping $\theta_{\bold a}^\varphi$ does not depend on the choice of a
geometric realization of the vector $\bold a$ defining the rotation axis and 
on a geometric realization of the vector $\bold b$ itself.
\endproclaim
\mydefinition{37.2} The mapping $\theta_{\bold a}^\varphi$ acting upon
free vectors of the space $\Bbb E$ and taking them to other free vectors in
$\Bbb E$ is called the rotation by the angle $\varphi$ about the vector $a$.
\enddefinition
\myexercise{37.1} Rotations and parallel translations belong to the class of
mappings preserving lengths of segments and measures of angles. They take
each segment to a congruent segment and each angle to a congruent angle (see
\mycite{6}). Let $p$ be some parallel translation, let $p^{-1}$ be its inverse
parallel translation, and let $\theta$ be a rotation by some angle about some
axis. Prove that the composite mapping $\tilde\theta=p\compos\theta\compos 
p^{-1}$ is the rotation by the same angle about the axis produced from the axis
of $\theta$ by applying the parallel translation $p$ to it.  
\endproclaim
\myexercise{37.2} Apply the result of the exercise~\mytheexercise{37.1} for
proving the lemma~\mythelemma{37.1}.
\endproclaim
\mytheorem{37.1} For any nonzero free vector $\bold a\neq\bold 0$ and for 
any angle $\varphi$ the rotation $\theta_{\bold a}^\varphi$ by the angle 
$\varphi$ about the vector $\bold a$ is a linear mapping of free vectors. 
\endproclaim
\demo{Proof} In order to prove the theorem we need to inspect the 
conditions 1) and 2) from the definition~\mythedefinition{27.4} for the
mapping $\theta_{\bold a}^\varphi$. Let's begin with the first of these
conditions. Assume that $\bold b$ and $\bold c$ are two free vectors. Let's
build their geometric realizations $\bold b=\overrightarrow{BC\,\,}\!$ and
$\bold c=\overrightarrow{CO\,\,}\!$. Then the vector $\overrightarrow{BO\,\,}\!$ 
is the geometric realization for the sum of vectors $\bold b+\bold c$.\par
      Now let's choose a geometric realization for the vector $\bold a$. 
It determines the rotation axis. According to the lemma~\mythelemma{37.1}
the actual place of such a geometric realization does not matter for the
ultimate definition of the mapping $\theta_{\bold a}^\varphi$ as applied
to free vectors. But for the sake of certainty we choose 
$\bold a=\overrightarrow{OA\,\,}\!$.\par
      Let's apply the rotation by the angle $\varphi$ about the axis $OA$ 
to the points $B$, $C$, and $O$. The point $O$ is on the rotation axis. 
Therefore it is not moved. The points $B$ and $C$ are moved to the points 
$\tilde B$ and $\tilde C$ respectively, while the triangle $BCO$ is moved
to the triangle $\tilde B\tilde CO$. As a result we get the following
relationships:
$$
\align
&\hskip -2em
\theta_{\bold a}^\varphi\bigl(\,\overrightarrow{BC\,\,}\bigr)
=\overrightarrow{\tilde B\tilde C\,\,}\!,\\
&\hskip -2em
\theta_{\bold a}^\varphi\bigl(\,\overrightarrow{CO\,\,}\bigr)
=\overrightarrow{\tilde CO\,\,}\!,
\mytag{37.2}\\
&\hskip -2em
\theta_{\bold a}^\varphi\bigl(\,\overrightarrow{BO\,\,}\bigr)
=\overrightarrow{\tilde BO\,\,}\!.
\endalign
$$
But the vectors $\overrightarrow{\tilde B\tilde C\,\,}\!$ and
$\overrightarrow{\tilde CO\,\,}\!$ in \mythetag{37.2} are geometric
realizations for the vectors $\tbb=\theta_{\bold a}^\varphi(\bold b)$
and $\tilde\bold c=\theta_{\bold a}^\varphi(\bold c)$, while
$\overrightarrow{\tilde BO\,\,}\!$ is a geometric realization for the
vector $\theta_{\bold a}^\varphi(\bold b+\bold c)$. Hence we have
$$
\theta_{\bold a}^\varphi(\bold b+\bold c)
=\overrightarrow{\tilde BO\,\,}\!
=\overrightarrow{\tilde B\tilde C\,\,}\!
+\overrightarrow{\tilde CO\,\,}\!
=\theta_{\bold a}^\varphi(\bold b)
+\theta_{\bold a}^\varphi(\bold c). 
\mytag{37.3}
$$
The chain of equalities \mythetag{37.3} proves the first linearity
condition from the definition~\mythedefinition{27.4} as applied to
$\theta_{\bold a}^\varphi$.\par
     Let's proceed to proving the second linearity condition. It is more
simple than the first one. Multiplying a vector $\bold b$ by a number
$\alpha$, we make the length of its geometric realizations $|\alpha|$ 
times as greater. If $\alpha>0$, geometric realizations of the vector
$\alpha\,\bold b$ are codirected to geometric realizations of $\bold b$.
If $\alpha<0$, they are opposite to geometric realizations of $\bold b$.
And if $\alpha=0$, geometric realizations of the vector $\alpha\,\bold b$
do vanish. Let's apply the rotation by the angle $\varphi$ about some
geometric realization of the vector $\bold a\neq\bold 0$ some geometric
realizations of the vectors $\bold b$ and $\alpha\,\bold b$. Such a 
mapping preserves the lengths of vectors. Hence it preserves all the 
relations of their lengths. Moreover it maps straight lines to straight 
lines and preserves the order of points on that straight lines. Hence
codirected vectors are mapped to codirected ones and opposite vectors
to opposite ones respectively. As a result
$$
\hskip -2em
\theta_{\bold a}^\varphi(\alpha\,\bold b)
=\alpha\,\theta_{\bold a}^\varphi(\bold b).
\mytag{37.4}
$$
The relationship \mythetag{37.4} completes the proof of the linearity for
the mapping $\theta_{\bold a}^\varphi$ as applied to free vectors.
\qed\enddemo
\head
\SectionNum{38}{91} The relation of the vector product
with projections and rotations.
\endhead
\rightheadtext{\S\,38. The relation of the vector product \dots}
\parshape 7 0cm 10cm 0cm 10cm 0cm 10cm 0cm 10cm 0cm 10cm 0cm 10cm 5cm 5cm 
    Let's consider two non-collinear vectors $\bold a\nparallel
\bold b$ and their vector product $\bold c=[\bold a,
\bold b]$. The length of the vector $\bold c$ is determined by the lengths
of $\bold a$ and $\bold b$ and by the angle $\varphi$ between them:
$$
\pagebreak
|\bold c|=|\bold a|\,|\bold b|\,\sin\varphi.
\mytag{38.1}
$$ 
The vector $\bold c$ lies on the plane $\alpha$ \vadjust{\vskip 5pt\hbox 
to 0pt{\kern 0pt\includegraphics{angemeng26.eps}\hss}
\vskip -5pt}perpendicular to the vector $\bold b$ (see Fig\.~38.1). Let's 
denote through $\bold a_{\sssize\perp}$ the orthogonal projection of the 
vector $\bold a$ onto the plane $\alpha$, i\,\,e\. we set
$$
\bold a_{\sssize\perp}
=\pi_{\!{\sssize\perp}\bold b}(\bold a).
\mytag{38.2}
$$\par
\parshape 4 5cm 5cm 5cm 5cm 5cm 5cm 0cm 10cm
     The length of the above vector \mythetag{38.2} is determined by the 
formula $|\bold a_{\sssize\perp}|=|\bold a|\,\sin\varphi$. Comparing this 
formula with \mythetag{38.1} and taking into account Fig\.~38.1, we conclude 
that in order to superpose the vector $\bold a_{\sssize\perp}$ with the vector 
$\bold c$ one should first rotate it counterclockwise by the right angle about 
the vector $\bold b$ and then multiply by the negative number $-|\bold b|$. 
This yields the formula
$$
\hskip -2em
[\bold a,\bold b]=-|\bold b|\cdot
\theta_{\bold b}^{\pi\kern -1pt\sssize{/2}}
\bigl(\pi_{\!{\sssize\perp}\bold b}(\bold a)\bigr).
\mytag{38.3}
$$\par
     The formula \mythetag{38.3} sets the relation of the vector product with
the two mappings $\pi_{\!{\sssize\perp}\bold b}$ and 
$\theta_{\bold b}^{\pi\kern -1pt\sssize{/2}}$. One of them is the projection 
onto the orthogonal complement of the vector $\bold b$, while the other is the 
rotation by the angle $\pi/2$ about the vector $\bold b$. The formula 
\mythetag{38.3} is applicable provided the vector $\bold b$ is nonzero:
$$
\hskip -2em
\bold b\neq\bold 0,
\mytag{38.4}
$$ 
while the condition $\bold a\nparallel\bold b$ can be broken. If $\bold a
\parallel\bold b$ both sides of the formula \mythetag{38.4} vanish, but the
formula itself remains valid.
\head
\SectionNum{39}{92} Properties of the vector product.
\endhead
\rightheadtext{\S\,39. Properties of the vector product.}
\mytheorem{39.1} The vector product of vectors possesses the following four 
properties \pagebreak fulfilled for any three vectors $\bold a$, $\bold b$, 
$\bold c$ and for any number $\alpha$:
\roster
\rosteritemwd=5pt
\item"1)" $[\bold a,\bold b]=-[\bold b,\bold a]$;
\item"2)" $[\bold a+\bold b,\bold c]=[\bold a,\bold c]
+[\bold b,\bold c]$;
\item"3)" $[\alpha\,\bold a,\bold c]=\alpha\,[\bold a,\bold c]$;
\item"4)" $[\bold a,\bold b]=0$ if and only if the vectors $\bold a$ and
$\bold b$\newline are collinear, i\.\,e\. if $\bold a\parallel\bold b$.
\endroster
\endproclaim
\mydefinition{39.1} The property 1) in the theorem~\mythetheorem{39.1} is called
{\it anticommutativity}; the properties 2) and 3) are called the properties of
{\it linearity with respect to the first multiplicand}; the property 4) is called
the {\it vanishing condition}.
\enddefinition
\demo{Proof of the theorem~\mythetheorem{39.1}} The property of anticommutativity
1) is derived immediately from the definition~\mythedefinition{35.1}. Let
$\bold a\nparallel\bold b$. Exchanging the vectors $\bold a$ and $\bold b$, we do 
not violate the first and the third conditions for the triple of vectors $\bold a$, 
$\bold b$, $\bold c$ in the definition~\mythedefinition{35.1}, provided they are 
initially fulfilled. As for the direction of rotation in Fig\.~35.1, it changes
for the opposite one. Therefore, if the triple $\bold a,\,\bold b,\,\bold c$ is
right, the triple $\bold b,\,\bold a,\,\bold c$ is left. In order to get a right 
triple the vectors $\bold b$ and $\bold a$ should be complemented with the vector 
$-\bold c$. This yield the equality
$$
\hskip -2em
[\bold a,\bold b]=-[\bold b,\bold a]
\mytag{39.1}
$$
for the case $\bold a\nparallel\bold b$. If $\bold a\parallel\bold b$, both sides 
of the equality \mythetag{39.1} do vanish. So the equality remains valid in this 
case too.\par
     Let $\bold c\neq\bold 0$. The properties of linearity 1) and 2) for this case
are derived with the use of the formula \mythetag{38.3} and the 
theorems~\mythetheorem{36.1} and \mythetheorem{37.1}. Let's write the formula
\mythetag{38.3} as
$$
\hskip -2em
[\bold a,\bold c]=-|\bold c|\cdot
\theta_{\bold c}^{\pi\kern -1pt\sssize{/2}}
\bigl(\pi_{\!{\sssize\perp}\bold c}(\bold a)\bigr).
\mytag{39.2}
$$
Then we change the vector $\bold a$ in \mythetag{39.2} for the sum of vectors $\bold a
+\bold b$ and apply the theorems~\mythetheorem{36.1} and \mythetheorem{37.1}. This
yields
$$
\allowdisplaybreaks
\gather
[\bold a+\bold b,\bold c]=-|\bold c|\cdot
\theta_{\bold c}^{\pi\kern -1pt\sssize{/2}}
\bigl(\pi_{\!{\sssize\perp}\bold c}(\bold a+\bold b)\bigr)
=-|\bold c|\,\cdot\\
\cdot\,\theta_{\bold c}^{\pi\kern -1pt\sssize{/2}}
\bigl(\pi_{\!{\sssize\perp}\bold c}(\bold a)
+\pi_{\!{\sssize\perp}\bold c}(\bold b)\bigr)
=-|\bold c|\cdot
\theta_{\bold c}^{\pi\kern -1pt\sssize{/2}}
\bigl(\pi_{\!{\sssize\perp}\bold c}(\bold a)\bigr)
\,-\\
-\,|\bold c|\cdot
\theta_{\bold c}^{\pi\kern -1pt\sssize{/2}}
\bigl(\pi_{\!{\sssize\perp}\bold c}(\bold b)\bigr)
=[\bold a,\bold c]+[\bold b,\bold c].
\endgather
$$
Now we change the vector $\bold a$ in \mythetag{39.2} for the product $\alpha\,\bold a$
and then apply the theorems~\mythetheorem{36.1} and \mythetheorem{37.1} again:
$$
\gather
[\alpha\,\bold a,\bold c]=-|\bold c|\cdot
\theta_{\bold c}^{\pi\kern -1pt\sssize{/2}}
\bigl(\pi_{\!{\sssize\perp}\bold c}(\alpha\,\bold a)\bigr)
=-|\bold c|\cdot\theta_{\bold c}^{\pi\kern -1pt\sssize{/2}}
\bigl(\alpha\,\pi_{\!{\sssize\perp}\bold c}(\bold a)\bigr)=\\
=-\alpha\,|\bold c|\cdot
\theta_{\bold c}^{\pi\kern -1pt\sssize{/2}}
\bigl(\pi_{\!{\sssize\perp}\bold c}(\bold a)\bigr)
=\alpha\,[\bold a,\bold c].
\endgather
$$\par
     The calculations which are performed above prove the equalities 2) and
3) in the theorem~\mythetheorem{39.1} for the case $\bold c\neq\bold 0$. 
If $\bold c=\bold 0$, both sides of these equalities do vanish and they appear
to be trivially fulfilled.\par
     Let's proceed to proving the fourth item in the 
theorem~\mythetheorem{39.1}. For $\bold a\parallel\bold b$ the vector product
$[\bold a,\bold b]$ vanishes by the definition~\mythedefinition{35.1}. Let 
$\bold a\nparallel\bold b$. In this case both vectors $\bold a$ and $\bold b$ 
are nonzero, while the angle $\varphi$ between them differs from $0$ and $\pi$. 
For this reason $\sin\varphi\neq 0$. Summarizing these restrictions and applying
the item 3) of the definition~\mythedefinition{35.1}, we get
$$
|[\bold a,\bold b]|=|\bold a|\,|\bold b|\,\sin\varphi
\neq 0,
$$
i\.\,e\. for $\bold a\nparallel\bold b$ the vector product $[\bold a,\bold b]$
cannot vanish. The proof of the theorem~\mythetheorem{39.1} is over.\qed\enddemo	
\mytheorem{39.2} Apart from the properties 1)--4), the vector product possesses 
the following two properties which are fulfilled for any vectors $\bold a$, 
$\bold b$, $\bold c$ and for any number $\alpha$:
\roster
\rosteritemwd=5pt
\item"5)" $[\bold c,\bold a+\bold b]=[\bold c,\bold a]
+[\bold a,\bold b]$;
\item"6)" $[\bold c,\alpha\,\bold a]=\alpha\,[\bold c,\bold a]$.
\endroster
\endproclaim
\mydefinition{39.2} The properties 5) and 6) in the theorem~\mythetheorem{39.2}
are called the properties of {\it linearity with respect to the second multiplicand}.
\enddefinition
     The properties 5) and 6) are easily derived from the properties 2) and 3) by
applying the property 1). Indeed, we have
$$
\gather
[\bold c,\bold a+\bold b]=-[\bold a+\bold b,\bold c]
=-[\bold a,\bold c]-[\bold b,\bold c]
=[\bold c,\bold a]+[\bold c,\bold b],\\
[\bold c,\alpha\,\bold a]=-[\alpha\,\bold a,\bold c]
=-\alpha\,[\bold a,\bold c]=\alpha\,[\bold c,\bold a].
\endgather
$$
These calculations prove the theorem~\mythetheorem{39.2}. 
\head
\SectionNum{40}{95} Structural constants of the vector product.
\endhead
\rightheadtext{\S\,40. Structural constants \dots}
     Let $\bold e_1,\,\bold e_2,\,\bold e_3$ be some arbitrary basis in the space
$\Bbb E$. Let's take two vectors $\bold e_i$ and $\bold e_j$ of this basis and 
consider their vector product $[\bold e_i,\bold e_j]$. The vector $[\bold e_i,
\bold e_j]$ can be expanded in the basis $\bold e_1,\,\bold e_2,\,\bold e_3$. 
Such an expansion is usually written as follows:
$$
\hskip -2em
[\bold e_i,\bold e_j]=C^1_{ij}\,\bold e_1+
C^2_{ij}\,\bold e_2+C^3_{ij}\,\bold e_3.
\mytag{40.1}
$$
The expansion \mythetag{40.1} contains three coefficients $C^1_{ij}$, $C^2_{ij}$ 
and $C^3_{ij}$. However, the indices $i$ and $j$ in it run independently over
three values $1,\,2,\,3$. For this reason, actually, the formula \mythetag{40.1} 
represent nine expansions, the total number of coefficients in it is equal to
twenty seven.\par 
     The formula \mythetag{40.1} can be abbreviated in the following way:
$$
\hskip -2em
[\bold e_i,\bold e_j]=\sum^3_{k=1}C^k_{ij}\,\bold e_k.
\mytag{40.2}
$$
Let's apply the theorem~\mythetheorem{19.1} on the uniqueness of the expansion 
of a vector in a basis to the expansions of $[\bold e_i,\bold e_j]$ in 
\mythetag{40.1} or in \mythetag{40.2}. As a result we can formulate the following
theorem. 
\mytheorem{40.1} Each basis $\bold e_1,\,\bold e_2,\,\bold e_3$ in the space 
$\Bbb E$ is associated with a collection of twenty seven constants $C^k_{ij}$ 
which are determined uniquely by this basis through the expansions  \mythetag{40.2}.
\endproclaim
\mydefinition{40.1} The constants $C^k_{ij}$, \pagebreak which are uniquely determined
by a basis $\bold e_1,\,\bold e_2,\,\bold e_3$ through the expansions \mythetag{40.2}, 
are called the {\it structural constants of the vector product\/} in this basis.
\enddefinition
     The structural constants of the vector product are similar to the components of
the Gram matrix for a basis $\bold e_1,\,\bold e_2,\,\bold e_3$ (see 
Definition~\mythedefinition{29.2}). But they are more numerous and form a three index
array with two lower indices and one upper index. For this reason they cannot be placed
into a matrix.\par
\head
\SectionNum{41}{96} Calculation of the vector product through 
the coordinates of vectors in a skew-angular basis.
\endhead
\rightheadtext{\S\,41. \dots in a skew-angular basis.}
     Let $\bold e_1,\,\bold e_2,\,\bold e_3$ be some skew-angular basis. According to
the definition~\mythedefinition{29.1} the term {\it skew-angular basis} in this book 
is used as a synonym of an arbitrary basis. Let's choose some arbitrary vectors
$\bold a$ and $\bold b$ in the space $\Bbb E$ and consider their expansions 
$$
\xalignat 2
&\hskip -2em
\bold a=\sum^3_{i=1}a^i\,\bold e_i,
&&\bold b=\sum^3_{j=1}b^{\kern 0.6pt j}\,\bold e_j
\mytag{41.1}
\endxalignat
$$
in the basis $\bold e_1,\,\bold e_2,\,\bold e_3$. Substituting \mythetag{41.1} into
the vector product $[\bold a,\bold b]$, we get the following formula:
$$
\hskip -2em
[\bold a,\bold b]=\biggl[\,\sum^3_{i=1}a^i\,\bold e_i,
\,\sum^3_{j=1}b^{\kern 0.6pt j}\,\bold e_j\biggr].
\mytag{41.2}
$$
In order to transform the formula \mythetag{41.2} we apply the properties 2) and 5)
of the vector product (see Theorems~\mythetheorem{39.1} and \mythetheorem{39.2}).
Due to these properties we can bring the summation signs over $i$ and $j$ outside
the brackets of the vector product:
$$
\hskip -2em
[\bold a,\bold b]=\sum^3_{i=1}\sum^3_{j=1}[a^i\,\bold e_i,
b^{\kern 0.6pt j}\,\bold e_j].
\mytag{41.3}
$$
Now let's apply the properties 3) and 6) from the theorems~\mythetheorem{39.1} and
\mythetheorem{39.2}. \pagebreak Due to these properties we can bring the numeric 
factors $a^i$ and $b^{\kern 0.6pt j}$ outside the brackets of the vector product 
\nolinebreak\mythetag{41.3}:
$$
\hskip -2em
[\bold a,\bold b]=\sum^3_{i=1}\sum^3_{j=1}a^i\,b^{\kern 0.6pt j}
\,[\bold e_i,\bold e_j].
\mytag{41.4}
$$\par
     The vector products $[\bold e_i,\bold e_j]$ in the formula \mythetag{41.4}
can be replaced by their expansions \mythetag{40.2}. Upon substituting 
\mythetag{40.2} into \mythetag{41.4} the formula \mythetag{41.4} is written
as follows:
$$
\hskip -2em
[\bold a,\bold b]=\sum^3_{i=1}\sum^3_{j=1}\sum^3_{k=1}a^i
\,b^{\kern 0.6pt j}\,C^k_{ij}\,\bold e_k.
\mytag{41.5}
$$
\mydefinition{41.1} The formula \mythetag{41.5} is called the {\it formula for
calculating the vector product through the coordinates of vectors in a 
skew-angular basis}.
\enddefinition
\head
\SectionNum{42}{97} Structural constants of the vector product
in an orthonormal basis.
\endhead
\rightheadtext{\S\,42. Structural constants \dots}
\parshape 6 0cm 10cm 0cm 10cm 0cm 10cm 0cm 10cm 0cm 10cm 
5cm 5cm 
     Let's recall that an orthonormal basis (ONB) in the space $\Bbb E$ is a basis
composed by three unit vectors perpendicular to each other (see
Definition~\mythedefinition{31.3}). By their orientation, triples of non-coplanar 
vectors in the space $\Bbb E$ are subdivided into right and left triples (see Definition~\mythedefinition{34.4}). Therefore all bases in the space $\Bbb E$ are
subdivided into right bases and left bases, which applies to orthonormal bases
as well.\par
\parshape 8 5cm 5cm 5cm 5cm 5cm 5cm 5cm 5cm 5cm 5cm 5cm 5cm 
5cm 5cm 0cm 10cm
     Let's consider some right orthonormal basis $\bold e_1,\,\bold e_2,\,\bold e_3$. 
It is shown in Fig\.~42.1. \vadjust{\vskip 5pt\hbox to 
0pt{\kern 0pt\includegraphics{angemeng27.eps}\hss}
\vskip -5pt}Using the definition \mythedefinition{35.1}, one can calculate
various pairwise vector products of the vectors composing this basis. Since the
geometry of a right ONB is rather simple, we can perform these calculations up
to an explicit result and compose the multiplication table for $\bold e_1,
\,\bold e_2,\,\bold e_3$:
$$
\xalignat 3
&[\bold e_1,\bold e_1]=\bold 0,
&&[\bold e_1,\bold e_2]=\bold e_3,
&&[\bold e_1,\bold e_3]=-\bold e_2,
\qquad\quad\\
&[\bold e_2,\bold e_1]=-\bold e_3,
&&[\bold e_2,\bold e_2]=\bold 0,
&&[\bold e_2,\bold e_3]=\bold e_1,
\qquad\quad
\mytag{42.1}\\
&[\bold e_3,\bold e_1]=\bold e_2,
&&[\bold e_3,\bold e_2]=-\bold e_1,
&&[\bold e_3,\bold e_3]=\bold 0.
\qquad\quad
\endxalignat
$$\par
     Let's choose the first of the relationships \mythetag{42.1} and write its 
right hand side in the form of an expansion in the basis $\bold e_1,\,\bold e_2,
\,\bold e_3$:
$$
\hskip -2em
[\bold e_1,\bold e_1]=0\,\bold e_1+0\,\bold e_2+0\,\bold e_3.
\mytag{42.2}
$$
Let's compare the expansion \mythetag{42.2} with the expansion \mythetag{40.1}
written for the case $i=1$ and $j=1$:
$$
\hskip -2em
[\bold e_1,\bold e_1]=C^1_{11}\,\bold e_1+
C^2_{11}\,\bold e_2+C^3_{11}\,\bold e_3.
\mytag{42.3}
$$
Due to the uniqueness of the expansion of a vector in a basis (see
Theorem~\mythetheorem{19.1}) from \mythetag{42.2} and \mythetag{42.3}
we derive
$$
\xalignat 3
&C^1_{11}=0,
&&C^2_{11}=0,
&&C^3_{11}=0.
\qquad\quad
\mytag{42.4}
\endxalignat
$$\par
     Now let's choose the second relationship \mythetag{42.1} and write its 
right hand side in the form of an expansion in the basis $\bold e_1,
\,\bold e_2,\,\bold e_3$:
$$
\hskip -2em
[\bold e_1,\bold e_2]=0\,\bold e_1+0\,\bold e_2+1\,\bold e_3.
\mytag{42.5}
$$
Comparing \mythetag{42.5} with the expansion \mythetag{40.1} written for
the case $i=1$ and $j=2$, we get the values of the following constants:
$$
\xalignat 3
&C^1_{12}=0,
&&C^2_{12}=0,
&&C^3_{12}=1.
\qquad\quad
\mytag{42.6}
\endxalignat
$$
Repeating this procedure for all relationships \mythetag{42.1}, we can get the
complete set of relationships similar to \mythetag{42.4} and \mythetag{42.6}.
Then we can organize them into a single list: 
$$
\xalignat 3
&C^1_{11}=0,
&&C^2_{11}=0,
&&C^3_{11}=0,
\qquad\quad\\
&C^1_{12}=0,
&&C^2_{12}=0,
&&C^3_{12}=1,
\qquad\quad\\
&C^1_{13}=0,
&&C^2_{13}=-1,
&&C^3_{13}=0,
\qquad\quad\\
&C^1_{21}=0,
&&C^2_{21}=0,
&&C^3_{21}=-1,
\qquad\quad\\
&C^1_{22}=0,
&&C^2_{22}=0,
&&C^3_{22}=0,
\qquad\quad
\mytag{42.7}\\
&C^1_{23}=1,
&&C^2_{23}=0,
&&C^3_{23}=0,
\qquad\quad\\
&C^1_{31}=0,
&&C^2_{31}=1,
&&C^3_{31}=0,
\qquad\quad\\
&C^1_{32}=-1,
&&C^2_{32}=0,
&&C^3_{32}=1,
\qquad\quad\\
&C^1_{33}=0,
&&C^2_{33}=0,
&&C^3_{33}=0.
\qquad\quad
\endxalignat
$$
The formulas \mythetag{42.7} determine all of the 27 structural constants 
of the vector product in a right orthonormal basis. Let's write this result
as a theorem. 
\mytheorem{42.1} For any right orthonormal basis $\bold e_1,\,\bold e_2,
\,\bold e_3$ in the space \,$\Bbb E$ the structural constants of the vector
product are determined by the formulas \mythetag{42.7}. 
\endproclaim
\mytheorem{42.2} For any left orthonormal basis $\bold e_1,\,\bold e_2,
\,\bold e_3$ in the space \,$\Bbb E$ the structural constants of the vector
product are derived from \mythetag{42.7} by changing signs 
{\tencyr\char '074}$+${\tencyr\char '076} for 
{\tencyr\char '074}$-${\tencyr\char '076} and vise versa. 
\endproclaim
\myexercise{42.1} Draw a left orthonormal basis and, applying the 
definition~\mythedefinition{35.1} to the pairwise vector products of  the
basis vectors, derive the relationships analogous to \mythetag{42.1}. 
Then prove the theorem~\mythetheorem{42.2}.
\endproclaim
\head
\SectionNum{43}{99} Levi-Civita symbol.
\endhead
\rightheadtext{\S\,43. Levi-Civita symbol.}
     Let's examine the formulas \mythetag{42.7} for the structural
constants of the vector product in a right orthonormal basis. One can
easily observe the following pattern in them:
$$
\hskip -2em
C^k_{ij}=0 \ \vtop{\hsize 4.5cm\baselineskip=9pt\noindent 
if there are coinciding values of the indices $i,j,k$.}
\mytag{43.1}
$$
The condition \mythetag{43.1} describes all of the cases where the
structural constants in \mythetag{42.7} do vanish. The cases where
$C^k_{ij}=1$ are described by the condition
$$
C^k_{ij}=1 \ \vtop{\hsize 5.7cm\baselineskip=9pt\noindent 
if the indices $i,j,k$ take the values $(1,2,3)$, $(2,3,1)$,
or $(3,1,2)$.}
\mytag{43.2}
$$
Finally, the cases where $C^k_{ij}=-1$ are described by the
condition
$$
C^k_{ij}=-1 \ \vtop{\hsize 5.7cm\baselineskip=9pt\noindent 
if the indices $i,j,k$ take the values $(1,3,2)$, $(3,2,1)$,
or $(2,1,3)$.}
\mytag{43.3}
$$\par
     The triples of numbers in \mythetag{43.2} and \mythetag{43.3} 
constitute the complete set of various permutations of the numbers 
$1$, $2$, $3$:
$$
\hskip -2em
\aligned
&(1,2,3),\quad (2,3,1),\quad (3,1,2),\\
&(1,3,2),\quad (3,2,1),\quad (2,1,3).
\endaligned
\mytag{43.4}
$$
The first three permutations in \mythetag{43.4} are called {\it even 
permutations}. They are produced from the right order of the numbers 
$1$, $2$, $3$ by applying an even number of pairwise transpositions to 
them. Indeed, we have
$$
\align
&(1,2,3);\\
&(1,2,3)\overset{1}\to{\longrightarrow}(2,1,3)
\overset{2}\to{\longrightarrow}(2,3,1);\\
&(1,2,3)\overset{1}\to{\longrightarrow}(1,3,2)
\overset{2}\to{\longrightarrow}(3,1,2).
\endalign
$$
The rest three permutations in \mythetag{43.4} are called {\it odd
permutations}. In the case of these three permutations we have
$$
\align
&(1,2,3)\overset{1}\to{\longrightarrow}(1,3,2);\\
\displaybreak
&(1,2,3)\overset{1}\to{\longrightarrow}(3,2,1);\\
&(1,2,3)\overset{1}\to{\longrightarrow}(2,1,3).
\endalign
$$
\mydefinition{43.1} The permutations \mythetag{43.4} constitute
a set, which is usually denoted through $S_3$. If $\sigma\in S_3$, 
then $(-1)^\sigma$ means the parity of the permutation $\sigma$:
$$
(-1)^\sigma=\cases \ \ 1\text{\ if the permutation $\sigma$ 
is even;}\\
-1\text{\ \ if the permutation $\sigma$ is odd.}
\endcases
$$
\enddefinition
    Zeros, unities, and minus unities from \mythetag{43.1},
\mythetag{43.2}, and \mythetag{43.3} are usually united into a single
numeric array:
$$
\varepsilon^{ijk}=\varepsilon_{ijk}=
\cases \ \ 0&\vtop{\hsize=4.5cm\baselineskip 9pt
\lineskip=2pt\noindent if there are coinciding values 
of the indices $i,j,k$;}\\
\vspace{1ex}
\ \ 1&\vtop{\hsize=4.8cm\baselineskip 9pt\lineskip=2pt
\noindent if the values of the indices $i,j,k$ form an even
permutation of the numbers $1,2,3$;}\\
\vspace{1ex}
-1&\vtop{\hsize=4.8cm\baselineskip 9pt\lineskip=2pt
\noindent if the values of the indices $i,j,k$ form an odd
permutation of the numbers $1,2,3$.}
\endcases
\quad
\mytag{43.5}
$$
\mydefinition{43.2} The numeric array $\varepsilon$ determined by the
formula \mythetag{43.5} is called the Levi-Civita symbol. 
\enddefinition
     When writing the components of the Levi-Civita symbol either
three upper indices or three lower indices are used. Thus we emphasize
the equity of all these three indices. Placing indices on different levels 
in the Levi-Civita symbol is not welcome. Summarizing what was said above,
the formulas \mythetag{43.1}, \mythetag{43.2}, and \mythetag{43.3} are
written as follows:
$$
\hskip -2em
C^k_{ij}=\varepsilon_{ijk}.
\mytag{43.6}
$$
\mytheorem{43.1} For any \pagebreak right orthonormal basis $\bold e_1,
\,\bold e_2,\,\bold e_3$ the structural constants of the vector product 
in such a basis are determined by the equality \mythetag{43.6}.
\endproclaim
    In the case of a left orthonormal basis we have the 
theorem~\mythetheorem{42.2}. It yields the equality
$$
\hskip -2em
C^k_{ij}=-\varepsilon_{ijk}.
\mytag{43.7}
$$
\mytheorem{43.2} For any left orthonormal basis $\bold e_1,
\,\bold e_2,\,\bold e_3$ the structural constants of the vector product 
in such a basis are determined by the equality \mythetag{43.7}.
\endproclaim
    Note that the equalities \mythetag{43.6} and \mythetag{43.7} violate
the index setting rule given in the definition~\mythedefinition{24.9}. 
The matter is that the structural constants of the vector product 
$C^k_{ij}$ are the components of a geometric object. The places of their
indices are determined by the index setting convention, which is known
as Einstein's tensorial notation (see Definition~\mythedefinition{20.1}). 
As for the Levi-Civita symbol, it is an array of purely algebraic origin. 
\par
     The most important property of the Levi-Civita symbol is its 
{\it complete skew symmetry\/} or {\it complete antisymmetry}. It is 
expressed by the following equalities:
$$
\xalignat 3
&\hskip -2em
\varepsilon_{ijk}=-\varepsilon_{j\kern 0.3pt ik},
&&\varepsilon_{ijk}=-\varepsilon_{ikj},
&&\varepsilon_{ijk}=-\varepsilon_{kj\kern 0.3pt i},
\qquad
\\
\vspace{-1.5ex}
\mytag{43.8}\\
\vspace{-1.5ex}
&\hskip -2em
\varepsilon^{ijk}=-\varepsilon^{j\kern 0.3pt ik},
&&\varepsilon^{ijk}=-\varepsilon^{ikj},
&&\varepsilon^{ijk}=-\varepsilon^{kj\kern 0.3pt i}.
\qquad
\endxalignat
$$
The equalities \mythetag{43.8} mean, that under the transposition of 
any two indices the quantity $\varepsilon_{ijk}=\varepsilon^{ijk}$ changes
its sign. These equalities are easily derived from \mythetag{43.5}.
\head
\SectionNum{44}{102} Calculation of the vector product through the 
coordinates of vectors in an orthonormal basis.
\endhead
\rightheadtext{\S\,44. \dots in an orthonormal basis.}
     Let's recall that the term {\it skew-angular basis\/} in this book is 
used as a synonym of an arbitrary basis (see Definition~\mythedefinition{29.1}). 
Let $\bold e_1,\,\bold e_2,\,\bold e_3$ be a right orthonormal basis. It can be
treated as a special case of a skew-angular basis. Substituting \mythetag{43.6} 
into the formula \mythetag{41.5}, we obtain the formula 
$$
\hskip -2em
[\bold a,\bold b]=\sum^3_{i=1}\sum^3_{j=1}\sum^3_{k=1}a^i
\,b^{\kern 0.6pt j}\,\varepsilon_{ijk}\,\bold e_k.
\mytag{44.1}
$$
Here $a^i$ and $b^{\kern 0.6pt j}$ are the coordinates of the vectors $\bold a$ 
and $\bold b$ in the basis $\bold e_1,\,\bold e_2,\,\bold e_3$.\par
     In order to simplify the formulas \mythetag{44.1} note that the majority 
components of the Levi-Civita symbol are equal to zero. Only six of its twenty 
seven components are nonzero. Applying the formula \mythetag{43.5}, we can bring
the formula \mythetag{44.1} to 
$$
\hskip -2em
\gathered
[\bold a,\bold b]=a^1\,b^2\,\bold e_3+a^2\,b^3\,\bold e_1
+a^3\,b^1\,\bold e_2\,-\\
-\,a^2\,b^1\,\bold e_3-a^3\,b^2\,\bold e_1-a^1\,b^3\,\bold e_2.
\endgathered
\mytag{44.2}
$$
Upon collecting similar terms the formula \mythetag{44.2} yields
$$
\hskip -2em
\gathered
[\bold a,\bold b]=\bold e_1\,(a^2\,b^3-a^3\,b^2)\,-\\
-\,\bold e_2\,(a^1\,b^3-a^3\,b^1)+\bold e_3\,(a^1\,b^2-a^2\,b^1).
\endgathered
\mytag{44.3}
$$
Now from the formula \mythetag{44.3} we derive
$$
[\bold a,\bold b]=\bold e_1\,\vmatrix a^2 & a^3\\ 
\vspace{1ex}b^2 & b^3\endvmatrix-\bold e_2\,\vmatrix a^1 & a^3\\
\vspace{1ex}b^1 & b^3\endvmatrix+\bold e_3\,\vmatrix a^1 & a^2\\
\vspace{1ex}b^1 & b^2\endvmatrix.\quad
\mytag{44.4}
$$
In deriving the formula \mythetag{44.4} we used the formula for the determinant
of a $2\times 2$ matrix (see \mycite{7}).\par
      Note that the right hand side of the formula \mythetag{44.4} is the 
expansion of the determinant of a $3\times 3$ matrix by its first row (see 
\mycite{7}). Therefore this formula can be written as:
$$
\pagebreak
\hskip -2em
[\bold a,\bold b]=\vmatrix\bold e_1 & \bold e_2 & \bold e_3\\
\vspace{1ex} a^1 & a^2 & a^3\\ 
\vspace{1ex}b^1 & b^2 & b^3\endvmatrix.
\mytag{44.5}
$$
Here $a^1,\,a^2,\,a^3$ and $b^1,\,b^2,\,b^3$ are the coordinates of the vectors 
$\bold a$ and $\bold b$ in the basis $\bold e_1,\,\bold e_2,\,\bold e_3$. They
fill the second and the third rows in the determinant \mythetag{44.5}.
\mydefinition{44.1} The formulas \mythetag{44.1} and \mythetag{44.5} are called
the {\it formulas for calculating the vector product through the coordinates of
vectors in a right orthonormal basis}.
\enddefinition
     Let's proceed to the case of a left orthonormal basis. In this case the
structural constants of the vector product are given by the formula
\mythetag{43.7}. Substituting \mythetag{43.7} into \mythetag{41.5}, we get
$$
\hskip -2em
[\bold a,\bold b]=-\sum^3_{i=1}\sum^3_{j=1}\sum^3_{k=1}a^i
\,b^{\kern 0.6pt j}\,\varepsilon_{ijk}\,\bold e_k.
\mytag{44.6}
$$
Then from \mythetag{44.6} we derive the formula
$$
\hskip -2em
[\bold a,\bold b]=-\vmatrix\bold e_1 & \bold e_2 & \bold e_3\\
\vspace{1ex} a^1 & a^2 & a^3\\ 
\vspace{1ex}b^1 & b^2 & b^3\endvmatrix.
\mytag{44.7}
$$
\mydefinition{44.2} The formulas \mythetag{44.6} and \mythetag{44.7} are
called the {\it formulas for calculating the vector product through the coordinates 
of vectors in a left orthonormal basis}.
\enddefinition
\head
\SectionNum{45}{104} Mixed product.
\endhead
\rightheadtext{\S\,45. Mixed product.}
\mydefinition{45.1} The {\it mixed product\/} of three free vectors $\bold a$, 
$\bold b$, and $\bold c$ is a number obtained as the scalar product of the vector
$\bold a$ by the vector product of $\bold b$ and $\bold c$:
$$
\hskip -2em
(\bold a,\bold b,\bold c)=(\bold a,[\bold b,\bold c]).
\mytag{45.1}
$$
\enddefinition
     As we see in the formula \mythetag{45.1}, the mixed product has three
multiplicands. They are separated from each other by commas. Commas are the 
multiplication signs in writing the mixed product, not by themselves, 
but together with the round brackets surrounding the whole expression.
\par 
     Commas and brackets in writing the mixed product are natural delimiters
for multiplicands: the first multiplicand is an expression between the opening
bracket and the first comma, the second multiplicand is an expression placed 
between two commas, and the third multiplicand is an expression between the
second comma and the closing bracket. Therefore in complicated expressions
no auxiliary delimiters are required. Foe example, in 
$$
(\bold a+\bold b,\bold c+\bold d,\bold e+\bold f)
$$
the sums $\bold a+\bold b$, $\bold c+\bold d$, and $\bold e+\bold f$ are
calculated first, then the mixed product itself is calculated.
\par
\head
\SectionNum{46}{105} Calculation of the mixed product through 
the coordinates of vectors in an orthonormal basis.
\endhead
\rightheadtext{\S\,46. \dots in an orthonormal basis.}
     The formula \mythetag{45.1} reduces the calculation of the mixed product
to successive calculations of the vectorial and scalar products. In the case
of the vectorial and scalar products we already have rather efficient formulas
for calculating them through the coordinates of vectors in an orthonormal 
basis.\par
     Let $\bold e_1,\,\bold e_2,\,\bold e_3$ be a right orthonormal basis 
and let $\bold a$, $\bold b$, and $\bold c$ be free vectors given by their
coordinates in this basis:
$$
\xalignat 3
&\bold a=\Vmatrix a^1\\\vspace{0.4ex}a^2\\\vspace{0.4ex}a^3
\endVmatrix,
&&\bold b=\Vmatrix b^1\\\vspace{0.4ex}b^2\\\vspace{0.4ex}b^3
\endVmatrix,
&&\bold c=\Vmatrix c^1\\\vspace{0.4ex}c^2\\\vspace{0.4ex}c^3
\endVmatrix.
\qquad\quad
\mytag{46.1}
\endxalignat
$$  
Let's denote $\bold d=[\bold b,\bold c]$. Then the formula \mythetag{45.1}
is written as 
$$
\hskip -2em
(\bold a,\bold b,\bold c)=(\bold a,\bold d).
\mytag{46.2}
$$
In order to calculate the vector $\bold d=[\bold b,\bold c]$ \pagebreak 
we apply the formula \mythetag{44.1} which now is written as follows:
$$
\hskip -2em
\bold d=\sum^3_{k=1}\biggl(\,\sum^3_{i=1}\sum^3_{j=1}
b^{\kern 0.6pt i}\,c^{\kern 0.8pt j}\,\varepsilon_{ijk}
\!\biggr)\,\bold e_k.
\mytag{46.3}
$$
The formula \mythetag{46.3} is an expansion of the vector $\bold d$ in 
the basis $\bold e_1,\,\bold e_2,\,\bold e_3$. Hence the coefficients 
in this expansion should coincide with the coordinates of the vector 
$\bold d$:
$$
\hskip -2em
d^{\kern 0.3pt k}=\sum^3_{i=1}\sum^3_{j=1}b^{\kern 0.6pt i}
\,c^{\kern 0.8pt j}\,\varepsilon_{ijk}.
\mytag{46.4}
$$\par
     The next step consists in using the coordinates of the vector
$\bold d$ from \mythetag{46.4} for calculating the scalar product in
the right hand side of the formula \mythetag{46.2}. The formula
\mythetag{33.3} now is written as 
$$
\hskip -2em
(\bold a,\bold d)=\sum^3_{k=1}a^k\,d^{\kern 0.3pt k}.
\mytag{46.5}
$$
Let's substitute \mythetag{46.4} into \mythetag{46.5} and take into
account \mythetag{46.2}. As a result we get the formula
$$
\hskip -2em
(\bold a,\bold b,\bold c)=\sum^3_{i=1}a^k\biggl(\,\sum^3_{i=1}
\sum^3_{j=1}b^{\kern 0.6pt i}\,c^{\kern 0.8pt j}
\,\varepsilon_{ijk}\!\biggr).
\mytag{46.6}
$$
Expanding the right hand side of the formula \mythetag{46.6} and changing 
the order of summations in it, we bring it to 
$$
\hskip -2em
(\bold a,\bold b,\bold c)=\sum^3_{i=1}\sum^3_{j=1}\sum^3_{k=1}
b^{\kern 0.6pt i}\,c^{\kern 0.8pt j}\,\varepsilon_{ijk}\,a^k.
\mytag{46.7}
$$\par
     Note that the right hand side of the formula \mythetag{46.7} differs
from that of the formula \mythetag{44.1} by changing $a^i$ for 
$b^{\kern 0.6pt i}$, changing $b^{\kern 0.6pt j}$ for $c^{\kern 0.8pt j}$,
and changing $\bold e_k$ for $a^k$. For this reason the formula \mythetag{46.7} 
can be brought to the following form analogous to \mythetag{44.5}:
$$
\hskip -2em
(\bold a,\bold b,\bold c)
=\vmatrix a^1 & a^2 & a^3\\ 
\vspace{1ex}b^1 & b^2 & b^3\\
\vspace{1ex}c^1 & c^2 & c^3
\endvmatrix.
\mytag{46.8}
$$\par
     Another way for transforming the formula \mythetag{46.7} is the use of
the complete antisymmetry of the Levi-Civita symbol \mythetag{43.8}. Applying
this property, we derive the identity $\varepsilon_{ijk}=\varepsilon_{k\kern 
0.4pt ij}$. Due to this identity, upon changing the order of multiplicands 
and redesignating the summation indices in the right hand side of the formula
\mythetag{46.7}, we can bring this formula to the following form:
$$
\hskip -2em
(\bold a,\bold b,\bold c)=\sum^3_{i=1}\sum^3_{j=1}\sum^3_{k=1}
a^i\,b^{\kern 0.6pt j}\,c^{\kern 0.4pt k}\,\varepsilon_{ijk}.
\mytag{46.9}
$$\par
\mydefinition{46.1} The formulas \mythetag{46.8} and \mythetag{46.9} are called 
the {\it formulas for calculating the mixed product through the coordinates of 
vectors in a right orthonormal basis}.
\enddefinition
\noindent The coordinates of the vectors $\bold a$, $\bold b$, and $\bold c$ 
used in the formulas \mythetag{46.8} and \mythetag{46.9} are taken from 
\mythetag{46.1}.\par
     Let's proceed to the case of a left orthonormal basis. Analogs of the formulas 
\mythetag{46.8} and \mythetag{46.9} for this case are obtained by changing the sign 
in the formulas \mythetag{46.8} and \mythetag{46.9}:
$$
\gather
\hskip -2em
(\bold a,\bold b,\bold c)
=-\vmatrix a^1 & a^2 & a^3\\ 
\vspace{1ex}b^1 & b^2 & b^3\\
\vspace{1ex}c^1 & c^2 & c^3
\endvmatrix.
\mytag{46.10}\\
\vspace{2ex}
\hskip -2em
(\bold a,\bold b,\bold c)=-\sum^3_{i=1}\sum^3_{j=1}\sum^3_{k=1}
a^i\,b^{\kern 0.6pt j}\,c^{\kern 0.4pt k}\,\varepsilon_{ijk}.
\mytag{46.11}
\endgather
$$
\mydefinition{46.2} The formulas \mythetag{46.10} and \mythetag{46.11} are called
the {\it formulas for calculating the mixed product through the coordinates of 
vectors in a left orthonormal basis}.
\enddefinition
     The formulas \mythetag{46.10} and \mythetag{46.11} can be derived using
the theorem~\mythetheorem{42.2} or comparing the formulas \mythetag{43.6} and
\mythetag{43.7}. 
\head
\SectionNum{47}{108} Properties of the mixed product.
\endhead
\rightheadtext{\S\,47. Properties of the mixed product.}
\mytheorem{47.1} The mixed product possesses the following four properties 
which are fulfilled for any four vectors $\bold a$, $\bold b$, $\bold c$, $\bold d$
and for any number $\alpha$:
\roster
\rosteritemwd=5pt
\item"1)" $(\bold a,\bold b,\bold c)=-(\bold a,\bold c,\bold b)$,
\newline $(\bold a,\bold b,\bold c)=-(\bold c,\bold b,\bold a)$,
\newline $(\bold a,\bold b,\bold c)=-(\bold b,\bold a,\bold c)$;
\item"2)" $(\bold a+\bold b,\bold c,\bold d)=(\bold a,\bold c,
\bold d)+(\bold b,\bold c,\bold d)$;
\item"3)" $(\alpha\,\bold a,\bold c,\bold d)=\alpha\,(\bold a,
\bold c,\bold d)$;
\item"4)" $(\bold a,\bold b,\bold c)=0$ if and only if the vectors $\bold a$, 
$\bold b$, and $\bold c$ are coplanar.
\endroster
\endproclaim
\mydefinition{47.1} The property 1) expressed by three equalities in the 
theorem~\mythetheorem{47.1} is called the property of {\it complete skew 
symmetry} or {\it complete antisymmetry}, the properties 2) and 3) are
called the properties of {\it linearity with respect to the first 
multiplicand}, the property 4) is called the {\it vanishing condition}.
\enddefinition
\demo{Proof of the theorem~\mythetheorem{47.1}} The first of the
three equalities composing the property of complete antisymmetry 1)
follows from the formula \mythetag{45.1} and the theorems~\mythetheorem{39.1}
and \mythetheorem{28.1}:
$$
(\bold a,\bold b,\bold c)=(\bold a,[\bold b,\bold c])
=(\bold a,-[\bold c,\bold b])=-(\bold a,[\bold c,\bold b])
=-(\bold a,\bold c,\bold b). 
$$
The other two equalities entering the property 1) cannot be derived in this
way. Therefore we need to use the formula \mythetag{46.8}. Transposition of
two vectors in the left hand side of this formula corresponds to the transposition
of two rows in the determinant in the right hand side of this formula. It is well
known that the transposition of any two rows in a determinant changes its sign.
This observation proves all of the three equalities composing the property of
complete antisymmetry for the mixed product.\par
     The properties of linearity 2) and 3) of the mixed product in the
theorem~\mythetheorem{47.1} are derived from the corresponding properties of the 
scalar and vectorial products due to the formula \mythetag{45.1}:
$$
\align
&\aligned
(\bold a+\bold b,\bold c,\bold d)&=(\bold a+\bold b,[\bold c,
\bold d])=(\bold a,[\bold c,\bold d])\,+\\
&\quad\,\ +\,(\bold b,[\bold c,\bold d])=(\bold a,\bold c,\bold d)
+(\bold b,\bold c,\bold d),
\endaligned\\
\vspace{1ex}
&(\alpha\,\bold a,\bold c,\bold d)=(\alpha\,\bold a,[\bold c,
\bold d])=\alpha\,(\bold a,[\bold c,\bold d])=
\alpha\,(\bold a,\bold c,\bold d).
\endalign
$$\par
      Let's proceed to proving the fourth property of the mixed product in the
theorem~\mythetheorem{47.1}. Assume that the vectors $\bold a$, $\bold b$, and
$\bold c$ are coplanar. In this case they are parallel to some plane $\alpha$ 
in the space $\Bbb E$ and one can choose their geometric realizations lying
on this plane $\alpha$. If $\bold b\nparallel\bold c$, then the vector product
$\bold d=[\bold b,\bold c]$ is nonzero and perpendicular to the plane
$\alpha$. As for the vector $\bold a$, it is parallel to this plane. Hence
 $\bold d\perp\bold a$, which yields the equalities $(\bold a,\bold b,\bold c)
 =(\bold a,[\bold b,\bold c])=(\bold a,\bold d)=0$.\par
     If $\bold b\parallel\bold c$, then the vector product $[\bold b,\bold c]$ 
is equal to zero and the equality $(\bold a,\bold b,\bold c)=0$ is derived
from $[\bold b,\bold c]=0$ with use of the initial formula \mythetag{45.1} 
for the scalar product.\par
     Now, conversely, assume that $(\bold a,\bold b,\bold c)=0$. If $\bold b
\parallel\bold c$, then the vectors $\bold a$, $\bold b$, and $\bold c$ 
determine not more than two directions in the space $\Bbb E$. For any two
lines in this space always there is a plane to which these lines are parallel.
In this case the vectors $\bold a$, $\bold b$, and $\bold c$ are coplanar
regardless to the equality $(\bold a,\bold b,\bold c)=0$.\par
     If $\bold b\nparallel\bold c$, then $\bold d=[\bold b,\bold c]
\neq\bold 0$. Choosing geometric realizations of the vectors $\bold b$
and $\bold c$ with some common initial point $O$, we easily build a plane
$\alpha$ comprising both of these geometric realizations. The vector 
$\bold d\neq\bold 0$ in perpendicular to this plane $\alpha$. Then from
$(\bold a,\bold b,\bold c)=(\bold a,[\bold b,\bold c])=(\bold a,\bold d)=0$ 
we derive $\bold a\perp\bold d$, which yields $\bold a\parallel\alpha$. The
vectors $\bold b$ and $\bold c$ are also parallel to the plane $\alpha$ since
their geometric realizations lie on this plane. Hence all of the three vectors
$\bold a$, $\bold b$, and $\bold c$ are parallel to the plane $\alpha$, which
means that they are coplanar. The theorem~\mythetheorem{47.1} is completely 
proved.
\qed\enddemo
\mytheorem{47.2} Apart from the properties 1)--4), the mixed product
possesses the following four properties which are fulfilled for any four vectors 
$\bold a$, $\bold b$, $\bold c$, $\bold d$ and for any number $\alpha$:
\roster
\rosteritemwd=5pt
\item"5)" $(\bold c,\bold a+\bold b,\bold d)=(\bold c,
\bold a,\bold d)+(\bold c,\bold b,\bold d)$;
\item"6)" $(\bold c,\alpha\,\bold a,\bold d)
=\alpha\,(\bold c,\bold a,\bold d)$;
\item"7)" $(\bold c,\bold d,\bold a+\bold b)=(\bold c,\bold d,
\bold a)+(\bold c,\bold d,\bold b)$;
\item"8)" $(\bold c,\bold d,\alpha\,\bold a)
=\alpha\,(\bold c,\bold d,\bold a)$.
\endroster
\endproclaim
\mydefinition{47.2} The properties 5) and 6) in the theorem \mythetheorem{47.2}
are called the properties of {\it linearity with respect to the second 
multiplicand}, the properties 7) and 8) are called the properties of 
{\it linearity with respect to the third multiplicand}.
\enddefinition
     The property 5) is derived from the property 2) in the 
theorem~\mythetheorem{47.1} in the following way:
$$
\gathered
(\bold c,\bold a+\bold b,\bold d)=-(\bold a+\bold b,\bold c,\bold d)
=-((\bold a,\bold c,\bold d)+(\bold b,\bold c,\bold d))=\\
=-(\bold a,\bold c,\bold d)-(\bold b,\bold c,\bold d)
=(\bold c,\bold a,\bold d)+(\bold c,\bold b,\bold d).
\endgathered
$$
The property 1) from this theorem is also used in the above calculations. 
As for the properties 6), 7), and 8) in the theorem~\mythetheorem{47.2}, 
they are also easily derived from the properties 2) and 3) with the use 
of the property 1). Indeed, we have
$$
\align
&(\bold c,\alpha\,\bold a,\bold d)=-(\alpha\,\bold a,\bold c,\bold d)
=-\alpha\,(\bold a,\bold c,\bold d)=\alpha\,(\bold c,\bold a,\bold d),\\
\vspace{1ex}
&\gathered
(\bold c,\bold d,\bold a+\bold b)=-(\bold a+\bold b,\bold d,\bold c)
=-((\bold a,\bold d,\bold c)+(\bold b,\bold d,\bold c))=\\
=-(\bold a,\bold d,\bold c)-(\bold b,\bold d,\bold c)
=(\bold c,\bold d,\bold a)+(\bold c,\bold d,\bold b),
\endgathered\\
\vspace{1ex}
&(\bold c,\bold d,\alpha\,\bold a)=-(\alpha\,\bold a,\bold d,\bold c)
=-\alpha\,(\bold a,\bold d,\bold c)=\alpha\,(\bold c,\bold d,\bold a).
\endalign
$$
The calculations performed prove the theorem~\mythetheorem{47.2}.
\head
\SectionNum{48}{111} The concept of the oriented volume.
\endhead
\rightheadtext{\S\,48. The concept of the oriented volume.}
\parshape 3 0cm 10cm 0cm 10cm 5cm 5cm 
     Let $\bold a,\,\bold b,\,\bold c$ be a right triple of non-coplanar
vectors in the space $\Bbb E$. Let's consider their mixed product $(\bold a,
\bold b,\bold c)$. Due to the item 4) from the theorem~\mythetheorem{47.1} 
the non-coplanarity of the vectors $\bold a,\,\bold b,\,\bold c$ means
$(\bold a,\bold b,\bold c)\neq 0$, which in turn due to \mythetag{45.1} 
implies $[\bold b,\bold c]\neq\bold 0$.\par
\parshape 9 5cm 5cm 5cm 5cm 5cm 5cm 5cm 5cm 5cm 5cm 5cm 5cm 
5cm 5cm 5cm 5cm 0cm 10cm 
     Due to the item 4) from the theorem \mythetheorem{39.1} the non-vanishing
condition $[\bold b,\bold c]\neq\bold 0$ means $\bold b\nparallel\bold c$. 
\vadjust{\vskip 5pt\hbox to 0pt{\kern 0pt\includegraphics{angemeng28.eps}\hss}\vskip -5pt}Let's build the geometric realizations of 
the non-collinear vectors $\bold b$ and $\bold c$ at some common initial 
point $O$ and denote them $\bold b=\overrightarrow{OB\,\,}\!$ and 
$\bold c=\overrightarrow{OC\,\,}\!$. Then we build the geometric realization
of the vector $\bold a$ at the same initial point $O$ and denote it through
$\bold a=\overrightarrow{OA\,\,}\!$. Let's complement the vectors 
$\overrightarrow{OA\,\,}\!$, $\overrightarrow{OB\,\,}\!$, and  
$\overrightarrow{OC\,\,}\!$ up to a skew-angular parallelepiped
as shown in Fig\.~48.1.\par
     Let's denote $\bold d=[\bold b,\bold c]$. The vector $\bold d$ is
perpendicular to the base plane of the parallelepiped, its length is
calculated by the formula $|\bold d|=|\bold b|\,|\bold c|\,\sin\alpha$.
It is easy to see that the length of $\bold d$ coincides with the base area 
of our parallelepiped, i\.\,e\. with the area of the parallelogram
built on the vectors $\bold b$ and $\bold c$:
$$
\hskip -2em
S=|\bold d|=|\bold b|\,|\bold c|\,\sin\alpha.
\mytag{48.1}
$$ 
\mytheorem{48.1} For any two vectors the length of their vector product
coincides with the area of the parallelogram built on these two vectors. 
\endproclaim
     The fact formulated in the theorem~\mythetheorem{48.1} is known as
the {\it geometric interpretation of the vector product}.\par
     Now let's return to Fig\.~48.1. Applying the formula \mythetag{45.1},
for the mixed product $(\bold a,\bold b,\bold c)$ we derive
$$
\hskip -2em
(\bold a,\bold b,\bold c)=(\bold a,[\bold b,\bold c])
=(\bold a,\bold d)=|\bold a|\,|\bold d|\,\cos\varphi.
\mytag{48.2}
$$
Note that $|\bold a|\,\cos\varphi$ is the length of the segment $[OF]$,
which coincides with the length of $[AH]$. The segment $[AH]$ is parallel to
the segment $[OF]$ and to the vector $\bold d$, which is perpendicular to
the base plane of the skew-angular parallelepiped shown in Fig\.~48.1. Hence
the segment $[AH]$ represents the height of this parallelepiped ad we have 
the formula
$$
\hskip -2em
h=|AH|=|\bold a|\,\cos\varphi.
\mytag{48.3}
$$
Now from \mythetag{48.1}, \mythetag{48.3}, and \mythetag{48.2} we derive
$$
\hskip -2em
(\bold a,\bold b,\bold c)=S\,h=V,
\mytag{48.4}
$$
i\.\,e\. the mixed product $(\bold a,\bold b,\bold c)$ in our case coincides
with the volume of the skew-angular parallelepiped built on the vectors
$\bold a$, $\bold b$, and $\bold c$.\par
     In the general case the value of the mixed product of three non-coplanar vectors 
$\bold a$, $\bold b$, $\bold c$ can be either positive or negative, while the volume
of a parallelepiped is always positive. Therefore in the general case the formula 
\mythetag{48.4} should be written as
$$
\hskip -2em
(\bold a,\bold b,\bold c)=\cases \ \ V,&\vtop{\hsize=3.7cm\baselineskip 9pt
\lineskip=2pt\noindent if $\bold a,\,\bold b,\,\bold c$ is a right triple 
of vectors;}\\
\vspace{1ex}
-V,&\vtop{\hsize=3.5cm\baselineskip 9pt
\lineskip=2pt\noindent if $\bold a,\,\bold b,\,\bold c$ is a left triple 
of vectors.}\\
\endcases
\mytag{48.5}
$$
\mydefinition{48.1} The oriented volume of an ordered triple of non-coplanar
vectors is a quantity which is equal to the volume of the parallelepiped built
on these vectors in the case where these vectors form a right triple \pagebreak
and which is equal to the volume of this parallelepiped taken with the minus 
sign in the case where these vectors form a left triple.
\enddefinition
The formula \mythetag{48.5} can be written as a theorem.
\mytheorem{48.2} The mixed product of any triple of non-copla\-nar vectors 
coincides with their oriented volume.
\endproclaim
\mydefinition{48.2} If $\bold e_1,\,\bold e_2,\,\bold e_3$ is a basis in the
space $\Bbb E$, then the oriented volume of the triple of vectors $\bold e_1,
\,\bold e_2,\,\bold e_3$ is called the oriented volume of this basis. 
\enddefinition
\head
\SectionNum{49}{113} Structural constants of the mixed product.
\endhead
\rightheadtext{\S\,49. Structural constants \dots}
     Let $\bold e_1,\,\bold e_2,\,\bold e_3$ be some basis in the space $\Bbb E$. 
Let's consider various mixed products composed by the vectors of this basis:
$$
\hskip -2em
c_{ijk}=(\bold e_i,\bold e_j,\bold e_k).
\mytag{49.1}
$$
\mydefinition{49.1} For any basis $\bold e_1,\,\bold e_2,\,\bold e_3$ in the space
$\Bbb E$ the quantities $c_{ijk}$ given by the formula \mythetag{49.1} are called
the {\it structural constants of the mixed product\/} in this basis. 
\enddefinition
     The formula \mythetag{49.1} is similar to the formula \mythetag{29.6} for the
components of the Gram matrix. However, the structural constants of the mixed product 
$c_{ijk}$ in \mythetag{49.1} constitute a three index array which cannot be laid 
into a matrix.\par
     An important property of the structural constants $c_{ijk}$ is their {\it complete
skew symmetry\/} or {\it complete antisymmetry}. This property is expressed by the 
following equalities:
$$
\xalignat 3
&\hskip -2em
c_{ijk}=-c_{j\kern 0.3pt ik},
&&c_{ijk}=-c_{ikj},
&&c_{ijk}=-c_{kj\kern 0.3pt i},
\qquad
\mytag{49.2}
\endxalignat
$$
The relationships \mythetag{49.2} mean that under the transposition of any two indices
the quantity $c_{ijk}$ changes its sign. These relationships are easily derived from
\mythetag{43.5} by applying the item 1) from the theorem~\mythetheorem{47.1} to the
right hand side of \mythetag{49.1}.\par
     The following relationships are an immediate consequence of the property of 
complete antisymmetry of the structural constants of the mixed product $c_{ijk}$:
$$
\xalignat 3
&\hskip -2em
c_{i\kern 0.3pt ik}=-c_{i\kern 0.3pt ik},
&&c_{ijj}=-c_{ijj},
&&c_{ij\kern 0.3pt i}=-c_{ij\kern 0.3pt i},
\qquad
\mytag{49.3}
\endxalignat
$$
They are derived by substituting $j=i$, $k=j$, and $k=i$ into \mythetag{49.2}. 
From the relationships \mythetag{49.3} the vanishing condition for the structural
constants $c_{ijk}$ is derived:
$$
\hskip -2em
c_{ijk}=0, \ \vtop{\hsize 4.5cm\baselineskip=9pt\noindent 
if there are coinciding values of the indices $i,j,k$.}
\mytag{49.4}
$$
Now assume that the values of the indices $i,\,j,\,k$ do not coincide. In this case,
applying the relationships \mythetag{49.2}, we derive
$$
\align
&\hskip -5em
c_{ijk}=c_{123} \ \vtop{\hsize 5.7cm\baselineskip=9pt\noindent 
if the indices $i,j,k$ take the values $(1,2,3)$, $(2,3,1)$,
or $(3,1,2)$;}\hskip -1em
\mytag{49.5}\\
\vspace{1ex}
&\hskip -5em
c_{ijk}=-c_{123} \ \vtop{\hsize 5.7cm\baselineskip=9pt\noindent 
if the indices $i,j,k$ take the values $(1,3,2)$, $(3,2,1)$,
or $(2,1,3)$.}\hskip -1em
\mytag{49.6}
\endalign
$$\par
     The next step consists in comparing the relationships \mythetag{49.4}, 
\mythetag{49.5}, and \mythetag{49.6} with the formula \mythetag{43.5} that
determines the Levi-Civita symbol $\varepsilon_{ijk}$. Such a comparison
yields
$$
\hskip -2em
c_{ijk}=c_{123}\ \varepsilon_{ijk}.
\mytag{49.7}
$$
Note that $c_{123}=(\bold e_1,\bold e_2,\bold e_3)$. This formula follows
from \mythetag{49.1}. Therefore the formula \mythetag{49.7} can be written
as 
$$
\hskip -2em
c_{ijk}=(\bold e_1,\bold e_2,\bold e_3)\ \varepsilon_{ijk}.
\mytag{49.8}
$$
\mytheorem{49.1} In an arbitrary basis $\bold e_1,\,\bold e_2,\,\bold e_3$ 
the structural constants of the mixed product are expressed by the formula
\mythetag{49.8} through the only one constant --- the oriented volume of 
the basis. 
\endproclaim
\head
\SectionNum{50}{115} Calculation of the mixed product through 
the coordinates of vectors in a skew-angular basis.
\endhead
\rightheadtext{\S\,50. \dots in a skew-angular basis.}
     Let $\bold e_1,\,\bold e_2,\,\bold e_3$ be a skew-angular basis in the 
space $\Bbb E$. Let's recall that the term {\it skew-angular basis} in this 
book is used as a synonym of an arbitrary basis (see 
Definition~\mythedefinition{29.1}). Let $\bold a$, $\bold b$, and $\bold c$ 
be free vectors given by their coordinates in the basis \nolinebreak $\bold e_1,
\,\bold e_2,\,\bold e_3$:
$$
\xalignat 3
&\bold a=\Vmatrix a^1\\\vspace{0.4ex}a^2\\\vspace{0.4ex}a^3
\endVmatrix,
&&\bold b=\Vmatrix b^1\\\vspace{0.4ex}b^2\\\vspace{0.4ex}b^3
\endVmatrix,
&&\bold c=\Vmatrix c^1\\\vspace{0.4ex}c^2\\\vspace{0.4ex}c^3
\endVmatrix.
\qquad\quad
\mytag{50.1}
\endxalignat
$$  
The formulas \mythetag{50.1} mean that we have the expansions
$$
\xalignat 3
&\bold a=\sum^3_{i=1}a^i\,\bold e_i,
&&\bold b=\sum^3_{j=1}b^{\kern 0.6pt j}\,\bold e_j,
&&\bold c=\sum^3_{k=1}c^{\kern 0.4pt k}\,\bold e_k.
\qquad\quad
\mytag{50.2}
\endxalignat
$$
Let's substitute \mythetag{50.2} into the mixed product $(\bold a,
\bold b,\bold c)$:
$$
\hskip -2em
(\bold a,\bold b,\bold c)=\biggl(\,\sum^3_{i=1}a^i\,\bold e_i,
\,\sum^3_{j=1}b^{\kern 0.6pt j}\,\bold e_j,\,\sum^3_{k=1}
c^{\kern 0.4pt k}\,\bold e_k\biggr).
\mytag{50.3}
$$
In order to transform the formula \mythetag{50.3} we apply the properties
of the mixed product 2), 5), and 7) from the theorems~\mythetheorem{47.1} 
and \mythetheorem{47.2}. Due to these properties we can bring the summation 
signs over $i$, $j$, and $k$ outside the brackets of the mixed product:
$$
\hskip -2em
(\bold a,\bold b,\bold c)=\sum^3_{i=1}\sum^3_{j=1}
\sum^3_{k=1}\,(a^i\,\bold e_i,b^{\kern 0.6pt j}\,\bold e_j,
c^{\kern 0.4pt k}\,\bold e_k).
\mytag{50.4}
$$
Now we apply the properties 3), 6), 8) from the theorems~\mythetheorem{47.1}~and 
\mythetheorem{47.2}. Due to these properties we can bring the numeric
factors $a^i$, $b^{\kern 0.6pt j}$, $c^{\kern 0.4pt k}$ outside the brackets of 
the mixed product \mythetag{50.4}:
$$
\hskip -2em
(\bold a,\bold b,\bold c)=\sum^3_{i=1}\sum^3_{j=1}
\sum^3_{k=1}a^i\,b^{\kern 0.6pt j}\,c^{\kern 0.4pt k}
\,(\bold e_i,\bold e_j,\bold e_k).
\mytag{50.5}
$$\par
     The quantities $(\bold e_i,\bold e_j,\bold e_k)$ are structural constants
of the mixed product in the basis $\bold e_1,\,\bold e_2,\,\bold e_3$ (see 
\mythetag{49.1}). Therefore the formula \mythetag{50.5} can be written as
follows:
$$
\hskip -2em
(\bold a,\bold b,\bold c)=\sum^3_{i=1}\sum^3_{j=1}
\sum^3_{k=1}a^i\,b^{\kern 0.6pt j}\,c^{\kern 0.4pt k}
\,c_{ijk}.
\mytag{50.6}
$$
Let's substitute \mythetag{49.8} into \mythetag{50.6} and take into account
that the oriented volume $(\bold e_1,\bold e_2,\bold e_3)$ does not depend
on the summation indices $i$, $j$, and $k$. Therefore the oriented volume 
$(\bold e_1,\bold e_2,\bold e_3)$ can be brought outside the sums as a common 
factor:
$$
\hskip -2em
(\bold a,\bold b,\bold c)=(\bold e_1,\bold e_2,
\bold e_3)\sum^3_{i=1}\sum^3_{j=1}\sum^3_{k=1}
a^i\,b^{\kern 0.6pt j}\,c^{\kern 0.4pt k}
\,\varepsilon_{ijk}.
\mytag{50.7}
$$\par
      Note that the formula \mythetag{50.7} differs from the formula
\mythetag{46.9} only by the extra factor $(\bold e_1,\bold e_2,\bold e_3)$ 
in its right hand side. As for the formula \mythetag{46.9}, it is brought
to the form \mythetag{46.8} by applying the properties of the Levi-Civita
symbol $\varepsilon_{ijk}$ only. For this reason the formula \mythetag{50.7} 
can be brought to the following form:
$$
\hskip -2em
(\bold a,\bold b,\bold c)=(\bold e_1,\bold e_2,
\bold e_3)
\,\vmatrix a^1 & a^2 & a^3\\ 
\vspace{1ex}b^1 & b^2 & b^3\\
\vspace{1ex}c^1 & c^2 & c^3
\endvmatrix.
\mytag{50.8}
$$
\mydefinition{50.1} The formulas \mythetag{50.6}, \mythetag{50.7}, and
\mythetag{50.8} are called the {\it formulas for calculating the mixed
product through the coordinates of vectors in a skew-angular basis}.
\enddefinition
\head
\SectionNum{51}{116} The relation of structural constants of the 
vectorial and mixed products.
\endhead
\rightheadtext{\S\,51. The relation of structural constants \dots}
      The structural constants of the mixed product are determined by the
formula \mythetag{49.1}. Let's apply the formula \mythetag{45.1} \pagebreak
in order to transform \mythetag{49.1}. As a result we get the formula
$$
\hskip -2em
c_{ijk}=(\bold e_i,[\bold e_j,\bold e_k]).
\mytag{51.1}
$$
Now we can apply the formula \mythetag{40.2}. Let's write it as follows:
$$
\hskip -2em
[\bold e_j,\bold e_k]=\sum^3_{q=1}C^q_{jk}\,\bold e_q.
\mytag{51.2}
$$
Substituting \mythetag{51.2} into \mythetag{51.1} and taking into account
the properties 5) and 6) from the theorem~\mythetheorem{28.2}, we derive
$$
\hskip -2em
c_{ijk}=\sum^3_{q=1}C^q_{jk}\,(\bold e_i,\bold e_q).
\mytag{51.3}
$$
Let's apply the formulas \mythetag{29.6} and \mythetag{30.1} to \mythetag{51.3}
and bring the formula \mythetag{51.3} to the form
$$
\hskip -2em
c_{ijk}=\sum^3_{q=1}C^q_{jk}\,g_{q\kern 0.5pt i}.
\mytag{51.4}
$$
The following formula is somewhat more beautiful:
$$
\hskip -2em
c_{ijk}=\sum^3_{q=1}C^q_{ij}\,g_{qk}.
\mytag{51.5}
$$
In order to derive the formula \mythetag{51.5} we apply the identity $c_{ijk}
=c_{jk\kern 0.4pt i}$ to the left hand side of the formula \mythetag{51.4}. 
This identity is derived from \mythetag{49.2}. Then we perform the cyclic 
redesignation of indices $i\to k\to j\to i$.\par
     The formula \mythetag{51.5} is the first formula relating the structural
constants of the vectorial and mixed products. It is important from the
theoretical point of view, but this formula is of little use practically. 
Indeed, it expresses the structural constants of the mixed product through
the structural constants of the vector product. But for the structural constants 
of the mixed product we already have the formula \mythetag{49.8} which is rather 
efficient. As for the structural constants of the vector product, we have no 
formula yet, except the initial definition \mythetag{40.2}. For this reason
we need to invert the formula \mythetag{51.5} and express $C^q_{ij}$ through
$c_{ijk}$. In order to reach this goal we need some auxiliary information on the
Gram matrix. 
\mytheorem{51.1} The Gram matrix $G$ of any basis in the space $\Bbb E$ is 
non-degenerate, i\.\,e\. its determinant is nonzero: $\det G\neq 0$. 
\endproclaim
\mytheorem{51.2} For any basis $\bold e_1,\,\bold e_2,\,\bold e_3$ in the space
$\Bbb E$ the determinant of the Gram matrix $G$ is equal to the square of the
oriented volume of this basis: 
$$
\hskip -2em
\det G=(\bold e_1,\bold e_2,\bold e_3)^2.
\mytag{51.6}
$$
\endproclaim
     The theorem~\mythetheorem{51.1} follows from the theorem~\mythetheorem{51.2}. 
Indeed, a basis is a triple of non-coplanar vectors. From the non-coplanarity
of the vectors $\bold e_1,\,\bold e_2,\,\bold e_3$ due to item 4) of the 
theorem~\mythetheorem{47.1} we get $(\bold e_1,\bold e_2,\bold e_3)\neq 0$. 
Then the formula \mythetag{51.6} yields
$$
\hskip -2em
\det G>0,
\mytag{51.7}
$$
while the theorem~\mythetheorem{51.1} follows from the inequality 
\mythetag{51.7}.\par
     I will not prove the theorem~\mythetheorem{51.2} right now at this place. 
This theorem is proved below in \S\,56.\par
     Let's proceed to deriving consequences from the theorem~\mythetheorem{51.1}.
It is known that each non-degenerate matrix has an inverse matrix (see \mycite{7}). 
Let's denote through $G^{-1}$ the matrix inverse to the Gram matrix $G$. In writing
the components of the matrix $G^{-1}$ the following convention is used.
\mydefinition{51.1} For denoting the components of the matrix $G^{-1}$ inverse to
the Gram matrix $G$ \pagebreak the same symbol $g$ as for the components of the 
matrix $G$ itself is used, but the components of the {\it inverse Gram matrix\/} 
are enumerated with two upper indices:
$$
\hskip -2em
G^{-1}=\Vmatrix 
g^{11} & g^{12} & g^{13}\\
\vspace{1ex}
g^{21} & g^{22} & g^{23}\\
\vspace{1ex}
g^{31} & g^{32} & g^{33}
\endVmatrix
\mytag{51.8}
$$ 
\enddefinition
     The matrices $G$ and $G^{-1}$ are inverse to each other. Their product in any
order is equal to the unit matrix:
$$
\xalignat 2
&\hskip -2em
G\cdot G^{-1}=1,
&&G^{-1}\cdot G=1.
\mytag{51.9}
\endxalignat
$$
From the regular course of algebra we know that each of the equalities \mythetag{51.9}
fixes the matrix $G^{-1}$ uniquely once the matrix $G$ is given (see. \mycite{7}). 
Now we apply the matrix transposition operation to both sides of the matrix equalities
\mythetag{51.9}:
$$
\hskip -2em
(G\cdot G^{-1})^{\sssize\top}=1^{\sssize\top}=1.
\mytag{51.10}
$$
Then we use the identity $(A\kern 1pt \cdot B)^{\sssize
\top}=B^{\sssize\top}\cdot A^{\sssize\top}$ from the exercise~\mytheexercise{29.2} 
in order to transform the formula \mythetag{51.10} and take into account the symmetry
of the matrix $G$ (see Theorem~\mythetheorem{30.1}):
$$
\hskip -2em
(G^{-1})^{\sssize\top}\cdot G^{\sssize\top}=
(G^{-1})^{\sssize\top}\cdot G=1.
\mytag{51.11}
$$
The rest is to compare the equality \mythetag{51.11} with the second matrix 
equality \mythetag{51.9}. This yields
$$
\hskip -2em
(G^{-1})^{\sssize\top}=G^{-1}
\mytag{51.12}
$$
The formula \mythetag{51.12} can be written as a theorem.
\mytheorem{51.3} For any basis $\bold e_1,\,\bold e_2,\,\bold e_3$ in the
space $\Bbb E$ the matrix $G^{-1}$ inverse to the Gram matrix of this
basis is symmetric. 
\endproclaim
     In terms of the matrix components of the matrix \mythetag{51.7} the equality
\mythetag{51.12} is written as an equality similar to \mythetag{30.1}:
$$
\hskip -2em
g^{\kern 0.3pt ij}=g^{\kern 0.6pt j\kern 0.5pt i}.
\mytag{51.13}
$$
The second equality \mythetag{51.9} is written as the following relationships 
for the components of the matrices \mythetag{51.8} and \mythetag{29.7}:
$$
\hskip -2em
\sum^3_{k=1}g^{sk}\,g_{kq}=\delta^s_q.
\mytag{51.14}
$$
Here $\delta^s_q$ is the Kronecker symbol defined by the formula \mythetag{23.3}. 
Using the symmetry identities \mythetag{51.13} and \mythetag{30.1}, we write 
the relationship \mythetag{51.14} in the following form:
$$
\hskip -2em
\sum^3_{k=1}g_{qk}\,g^{ks}=\delta^s_q.
\mytag{51.15}
$$
Now let's multiply both ides of the equality \mythetag{51.5} by $g^{ks}$ and
then perform summation over $k$ in both sides of this equality:
$$
\hskip -2em
\sum^3_{k=1}c_{ijk}\,g^{ks}=\sum^3_{k=1}\sum^3_{q=1}C^q_{ij}
\,g_{qk}\,g^{ks}.
\mytag{51.16}
$$
If we take into account the identity \mythetag{51.15}, then we can 
bring the formula \mythetag{51.16} to the following form:
$$
\hskip -2em
\sum^3_{k=1}c_{ijk}\,g^{ks}=\sum^3_{q=1}C^q_{ij}
\,\delta^s_q.
\mytag{51.17}
$$
When summing over $q$ in the right hand side of the equality \mythetag{51.17} 
the index $q$ runs over three values $1,\,2,\,3$. Then the Kronecker symbol
$\delta^s_q$ takes the values $0$ and $1$, the value $1$ is taken only once
when $q=s$. This means that only one of three summands in the right hand side 
of \mythetag{51.17} is nonzero. This nonzero summand is equal to $C^s_{ij}$. 
Hence the formula \mythetag{51.17} can be written in the following form:
$$
\hskip -2em
\sum^3_{k=1}c_{ijk}\,g^{ks}=C^s_{ij}.
\mytag{51.18}
$$
Let's change the symbol $k$ for $q$ and the symbol $s$ for $k$ in 
\mythetag{51.18}. Then we transpose left and right hand sides of this 
formula:
$$
\hskip -2em
C^k_{ij}=\sum^3_{q=1}c_{ijq}\,g^{qk}.
\mytag{51.19}
$$\par
     The formula \mythetag{51.19} is the second formula relating the structural
constants of the vectorial and mixed products. On the base of the relationships
\mythetag{51.5} and \mythetag{51.19} we formulate a theorem. 
\mytheorem{51.4} For any basis $\bold e_1,\,\bold e_2,\,\bold e_3$ in the space
$\Bbb E$ the structural constants of the vectorial and mixed products in this 
basis are related to each other in a one-to-one manner by means of the formulas 
\mythetag{51.5} and \mythetag{51.19}.
\endproclaim
\head
\SectionNum{52}{121} Effectivization of the formulas for calculating
vectorial and mixed products.
\endhead
\rightheadtext{\S\,52. Effectivization of the formulas \dots}
\baselineskip=12pt plus 0.02pt minus 0.02pt
      Let's consider the formula \mythetag{29.8} for calculating the scalar
product in a skew-angular basis. Apart from the coordinates 
of vectors, this formula uses the components of the Gram matrix 
\mythetag{29.7}. In order to get the components of this matrix one should 
calculate the mutual scalar products of the basis vectors (see formula 
\mythetag{29.6}), for this purpose one should measure their lengths and
the angles between them (see Definition~\mythedefinition{26.1}). No other
geometric constructions are required. For this reason the formula 
\mythetag{29.8} is recognized to be effective.\par
      Now let's consider the formula \mythetag{41.5} for calculating the 
vector product in a skew-angular basis. This formula uses the structural
constants of the vector product which are defined by means of the formula
\mythetag{40.2}. According to this formula, in order to calculate the
structural constants one should calculate the vector products 
$[\bold e_i,\bold e_j]$ in its left hand side. For this purpose one should
construct the normal vectors (perpendiculars) to the planes given by various
pairs of basis vectors $\bold e_i$, $\bold e_j$. (see 
Definition~\mythedefinition{35.1}). Upon calculating the vector products
$[\bold e_i,\bold e_j]$ one should expand them in the basis $\bold e_1,
\,\bold e_2,\,\bold e_3$, which require some auxiliary geometric 
constructions (see formula \mythetag{18.4} and Fig\.~18.1). For this reason
the efficiency of the formula \mythetag{41.5} is much less than the efficiency
of the formula \mythetag{29.8}.\par
       And finally we consider the formula \mythetag{50.7} for calculating 
the mixed product in a skew-angular basis. In order to apply this formula
one should know the value of the mixed product of the three basis vectors
$(\bold e_1,\bold e_2,\bold e_3)$. It is called the oriented volume of a basis 
(see Definition~\mythedefinition{48.2}). Due to the theorem~\mythetheorem{48.2} 
and the definition~\mythedefinition{48.1} for this purpose one should 
calculate the volume of the skew-angular parallelepiped built on the basis
vectors $\bold e_1,\,\bold e_2,\,\bold e_3$. In order to calculate the volume
of this parallelepiped one should know its base area and its height. The area
of its base is effectively calculated by the lengths of two basis vectors and
the angle between them (see formula \mythetag{48.1}). As for the height of the
parallelepiped, in order to find it one should drop a perpendicular from one 
of its vertices to its base plane. Since we need such an auxiliary geometric
construction, the formula \mythetag{50.7} is less effective as compared to
the formula \mythetag{29.8} in the case of the scalar product.\par
\baselineskip=12pt
     In order to make the formulas \mythetag{41.5} and \mythetag{50.7}    
effective we use the formula \mythetag{51.6}. It leads to the following 
relationship:
$$
\hskip -2em
(\bold e_1,\bold e_2,\bold e_3)=\pm\sqrt{\det G}.
\mytag{52.1}
$$
The sign in \mythetag{52.1} is determined by the orientation
of a \pagebreak basis:
$$
(\bold e_1,\bold e_2,\bold e_3)
=\cases \ \ \sqrt{\det G}&\vtop{\hsize=3.5cm\baselineskip 9pt
\lineskip=2pt\noindent if the basis $\bold e_1,\,\bold e_2,\,\bold e_3$
is right;}\\
\vspace{1ex}
-\sqrt{\det G}&\vtop{\hsize=3.5cm\baselineskip 9pt
\lineskip=2pt\noindent if the basis $\bold e_1,\,\bold e_2,\,\bold e_3$
is left.}\\
\endcases
\quad
\mytag{52.2}
$$\par
     Let's substitute the expression \mythetag{52.1} into the formula
for the mixed product \mythetag{50.7}. As a result we get
$$
\hskip -2em
(\bold a,\bold b,\bold c)=\pm\sqrt{\det G}\,\sum^3_{i=1}\sum^3_{j=1}
\sum^3_{k=1}a^i\,b^{\kern 0.6pt j}\,c^{\kern 0.4pt k}
\,\varepsilon_{ijk}.
\mytag{52.3}
$$
Similarly, substituting \mythetag{52.1} into the formula \mythetag{50.8}, 
we get
$$
\hskip -2em
(\bold a,\bold b,\bold c)=\pm\sqrt{\det G}
\,\vmatrix a^1 & a^2 & a^3\\ 
\vspace{1ex}b^1 & b^2 & b^3\\
\vspace{1ex}c^1 & c^2 & c^3
\endvmatrix.
\mytag{52.4}
$$
\mydefinition{52.1} The formulas \mythetag{52.3} and \mythetag{52.4} are 
called the effectivized formulas for calculating the mixed product throu\-gh 
the coordinates of vectors in a skew-angular basis.
\enddefinition
     In the case of the formula \mythetag{41.5}, in order to make it effective
we need the formulas \mythetag{49.8} and \mythetag{51.19}. Substituting the
expression \mythetag{52.1} into these formulas, we obtain
$$
\align
&\hskip -2em
c_{ijk}=\pm\sqrt{\det G}\ \varepsilon_{ijk},
\mytag{52.5}\\
\vspace{1ex}
&\hskip -2em
C^k_{ij}=\pm\sqrt{\det G}\,\sum^3_{q=1}\varepsilon_{ijq}\,g^{qk}.
\mytag{52.6}
\endalign
$$
\mydefinition{52.2} The formulas \mythetag{52.5} and \mythetag{52.6} 
are called the effectivized formulas for calculating the structural 
constants of the mixed and vectorial products.
\enddefinition
     Now let's substitute the formula \mythetag{52.6} into \mythetag{41.5}.
This leads to the following formula for the vector \pagebreak product:
$$
[\bold a,\bold b]=\pm\sqrt{\det G}\,\sum^3_{i=1}\sum^3_{j=1}
\sum^3_{k=1}\sum^3_{q=1}a^i\,b^{\kern 0.6pt j}\,\varepsilon_{ijq}
\,g^{qk}\,\bold e_k.
\mytag{52.7}
$$
\mydefinition{52.3} The formula \mythetag{52.7} is called the effectivized 
formula for calculating the vector product through the coordinates of vectors 
in a skew-angular basis.
\enddefinition
\head
\SectionNum{53}{124} Orientation of the space.
\endhead
\rightheadtext{\S\,53. The orientation of the space.}
     Let's consider the effectivized formulas \mythetag{52.1}, \mythetag{52.3}, 
\mythetag{52.4}, \mythetag{52.5}, \mythetag{52.6}, and \mythetag{52.7} from
\S\,52. Almost all information on a basis in these formulas is given by the
Gram matrix $G$. The Gram matrix arising along with the choice of a basis
reflects an important property of the space $\Bbb E$ --- its metric. 
\mydefinition{53.1} The {\it metric\/} of the space $\Bbb E$ is its structure
(its feature) that consists in possibility to measure the lengths of segments 
and the numeric values of angles in it.
\enddefinition
     The only non-efficiency remaining in the effectivized formulas 
\mythetag{52.1}, \mythetag{52.3}, \mythetag{52.4}, \mythetag{52.5}, \mythetag{52.6}, 
and \mythetag{52.7} is the choice of sign in them. As was said above, the sign 
in these formulas is determined by the orientation of a basis (see formula
\mythetag{52.2}). There is absolutely no possibility to determine the orientation
of a basis through the information comprised in its Gram matrix. The matter is that
the mathematical space $\Bbb E$ described by Euclid's axioms (see \mycite{6}), 
comprises the possibility to distinguish a pair of bases with different orientations 
from a pair of bases with coinciding orientations. However, it does not contain 
any reasons for to prefer bases with one of two possible orientations.\par 
     The concept of right triple of vectors (see Definition~\mythedefinition{34.2}) 
and the possibility to distinguish right triples of vectors from left ones is due
to the presence of people, it is due to their ability to observe vectors and compare
their rotation with the rotation of clock hands. Is there a {\it fundamental 
asymmetry between left and right\/} not depending on the presence of people and on 
other non-fundamental circumstances? Is the space {\it fundamentally oriented}? This 
is a question on the nature of the physical space $\Bbb E$. Some research in the 
field of elementary particle physics says that such an asymmetry does exist. As for
me, as the author of this book I cannot definitely assert that this problem is
finally resolved. 
\head
\SectionNum{54}{125} Contraction formulas.
\endhead
\rightheadtext{\S\,54. Contraction formulas.}
     Contraction formulas is a collection of four purely algebraic identities 
relating the Levi-Civita symbol and the Kronecker symbol with each other.
They are proved successively one after another. Once the first identity is proved,
each next identity is derived from the previous one.
\mytheorem{54.1} The Levi-Civita symbol and the Kronecker symbol are related by
the first contraction formula 
$$
\hskip -2em
\varepsilon^{mnp}\,\varepsilon_{ijk}
=\vmatrix
\delta^{\kern 0.4pt m}_i & \delta^{\kern 0.4pt m}_j 
& \delta^{\kern 0.4pt m}_k \\
\vspace{1ex}
\delta^{\kern 0.4pt n}_i & \delta^{\kern 0.4pt n}_j 
& \delta^{\kern 0.4pt n}_k \\
\vspace{1ex}
\delta^{\kern 0.4pt p}_i & \delta^{\kern 0.4pt p}_j 
& \delta^{\kern 0.4pt p}_k 
\endvmatrix.
\mytag{54.1}
$$
\endproclaim
\demo{Proof} Let's denote through $f^{mnp}_{ijk}$ the right hand side of the
first contraction formula \mythetag{54.1}. The transposition of any two lower 
indices in $f^{mnp}_{ijk}$ is equivalent to the corresponding transposition 
of two columns in the matrix \mythetag{54.1}. It is known that the transposition
of two columns of a matrix results in changing the sign of its determinant. 
This yield the following relationships for the quantities $f^{mnp}_{ijk}$:
$$
f^{mnp}_{ijk}=-f^{mnp}_{j\kern 0.3pt ik},\quad
f^{mnp}_{ijk}=-f^{mnp}_{ikj},\quad
f^{mnp}_{ijk}=-f^{mnp}_{kj\kern 0.3pt i}.
\quad
\mytag{54.2}
$$
The relationships \mythetag{54.2} are analogous to the relationships \mythetag{49.2}.
Repeating the considerations used in \S\,49 when deriving the formulas 
\mythetag{49.4}, \mythetag{49.5}, \mythetag{49.6}, from the formula \mythetag{54.2}
we derive
$$
f^{mnp}_{ijk}=
\cases \ \ 0&\vtop{\hsize=4.5cm\baselineskip 9pt
\lineskip=2pt\noindent if there are coinciding values of the indices $i,j,k$;}\\
\vspace{1ex}
\ \ f^{mnp}_{123}&\vtop{\hsize=4.8cm\baselineskip 9pt\lineskip=2pt
\noindent if the values of the indices $i,j,k$ form an even
permutation of the numbers $1,2,3$;}\\
\vspace{1ex}
-f^{mnp}_{123},&\vtop{\hsize=4.8cm\baselineskip 9pt\lineskip=2pt
\noindent if the values of the indices $i,j,k$ form an odd
permutation of the numbers $1,2,3$.}
\endcases
\quad
\mytag{54.3}
$$
Let's compare the formula \mythetag{54.3} with the formula \mythetag{43.5} defining
the Levi-Civita symbol. Such a comparison yields
$$
\hskip -2em
f^{mnp}_{ijk}=f^{mnp}_{123}\,\varepsilon_{ijk}.
\mytag{54.4}
$$\par
     Like the initial quantities $f^{mnp}_{ijk}$, the factor $f^{mnp}_{123}$
in \mythetag{54.4} is defined as the determinant of a matrix:
$$
\hskip -2em
f^{mnp}_{123}
=\vmatrix
\delta^{\kern 0.4pt m}_1 & \delta^{\kern 0.4pt m}_2 
& \delta^{\kern 0.4pt m}_3 \\
\vspace{1ex}
\delta^{\kern 0.4pt n}_1 & \delta^{\kern 0.4pt n}_2 
& \delta^{\kern 0.4pt n}_3 \\
\vspace{1ex}
\delta^{\kern 0.4pt p}_1 & \delta^{\kern 0.4pt p}_2 
& \delta^{\kern 0.4pt p}_3 
\endvmatrix.
\mytag{54.5}
$$
Transposition of any pair of the upper indices in $f^{mnp}_{123}$ is equivalent
to the corresponding transposition of rows in the matrix \mythetag{54.5}. Again 
we know that the transposition of any two rows of a matrix changes the sign of
its determinant. Hence we derive the following relationships for the
quantities \mythetag{54.5}:
$$
f^{mnp}_{123}=-f^{nmp}_{123},\quad
f^{mnp}_{123}=-f^{mp\kern 0.2pt n}_{123},\quad
f^{mnp}_{123}=-f^{p\kern 0.2pt nm}_{123}
\quad
\mytag{54.6}
$$
The relationships \mythetag{54.6} are analogous to the relationships \mythetag{54.2},
which are analogous to \mythetag{49.2}. Repeating the considerations used in \S\,49,
from the relationships \mythetag{54.6} we immediately derive the formula
analogous to the formula \mythetag{54.3}:
$$
f^{mnp}_{123}=
\cases \ \ 0&\vtop{\hsize=4.8cm\baselineskip 9pt
\lineskip=2pt\noindent if there are coinciding values of the indices $m,n,p$;}\\
\vspace{1ex}
\ \ f^{123}_{123}&\vtop{\hsize=5.0cm\baselineskip 9pt\lineskip=2pt
\noindent if the values of the indices $m,n,p$ form an even
permutation of the numbers $1,2,3$;}\\
\vspace{1ex}
-f^{123}_{123}&\vtop{\hsize=5.1cm\baselineskip 9pt\lineskip=2pt
\noindent if the values of the indices $m,n,p$ form an odd
permutation of the numbers $1,2,3$.}
\endcases
\quad
\mytag{54.7}
$$
Let's compare the formula \mythetag{54.7} with the formula \mythetag{43.5}
defining the Levi-Civita symbol. This comparison yields
$$
\hskip -2em
f^{mnp}_{123}=f^{123}_{123}\,\varepsilon^{mnp}.
\mytag{54.8}
$$\par
     Let's combine the formulas \mythetag{54.4} and \mythetag{54.8}, i\.\,e\.
we substitute \mythetag{54.8} into \mythetag{54.4}. This substitution leads 
to the formula
$$
\hskip -2em
f^{mnp}_{ijk}=f^{123}_{123}
\,\varepsilon^{mnp}\,\varepsilon_{ijk}.
\mytag{54.9}
$$
Now we need to calculate the coefficient $f^{123}_{123}$ in the formula
\mythetag{54.9}. Let's recall that the quantity $f^{mnp}_{ijk}$ was
defined as the right hand side of the formula \mythetag{54.1}. Therefore
we can write
$$
\hskip -2em
f^{mnp}_{ijk}
=\vmatrix
\delta^{\kern 0.4pt m}_i & \delta^{\kern 0.4pt m}_j 
& \delta^{\kern 0.4pt m}_k \\
\vspace{1ex}
\delta^{\kern 0.4pt n}_i & \delta^{\kern 0.4pt n}_j 
& \delta^{\kern 0.4pt n}_k \\
\vspace{1ex}
\delta^{\kern 0.4pt p}_i & \delta^{\kern 0.4pt p}_j 
& \delta^{\kern 0.4pt p}_k 
\endvmatrix
\mytag{54.10}
$$
and substitute $i=m=1$, $j=n=2$, $k=p=3$ into \mythetag{54.10}:
$$
\pagebreak
\hskip -2em
f^{123}_{123}
=\vmatrix
\delta^1_1 & \delta^1_2 & \delta^1_3 \\
\vspace{1ex}
\delta^2_1 & \delta^2_2 & \delta^2_3 \\
\vspace{1ex}
\delta^3_1 & \delta^3_2 & \delta^3_3 
\endvmatrix=\vmatrix
1 & 0 & 0 \\
\vspace{0.5ex}
0 & 1 & 0 \\
\vspace{0.5ex}
0 & 0 & 1
\endvmatrix=1.
\mytag{54.11}
$$\par
     Taking into account the value of the coefficient $f^{123}_{123}$ given
by the formula \mythetag{54.11}, we can transform the formula \mythetag{54.9} 
to 
$$
\hskip -2em
f^{mnp}_{ijk}=\varepsilon^{mnp}\,\varepsilon_{ijk}.
\mytag{54.12}
$$
Now the required contraction formula \mythetag{54.1} is obtained as a consequence
of \mythetag{54.10} and \mythetag{54.12}. The theorem~\mythetheorem{54.1} is proved.
\qed\enddemo
\mytheorem{54.2} The Levi-Civita symbol and the Kronecker symbol are related by
the second contraction formula
$$
\hskip -2em
\sum^3_{k=1}\varepsilon^{mnk}\,\varepsilon_{ijk}
=\vmatrix
\delta^{\kern 0.4pt m}_i & \delta^{\kern 0.4pt m}_j\\
\vspace{1ex}
\delta^{\kern 0.4pt n}_i & \delta^{\kern 0.4pt n}_j
\endvmatrix.
\mytag{54.13}
$$
\endproclaim
\demo{Proof} The second contraction formula \mythetag{54.13} is derived from the
first contraction formula \mythetag{54.1}. For this purpose we substitute
$p=k$ into the formula \mythetag{54.1} and take into account that
$\delta^{\kern 0.4pt k}_k=1$. This yields the formula 
$$
\hskip -2em
\varepsilon^{mnk}\,\varepsilon_{ijk}
=\vmatrix
\delta^{\kern 0.4pt m}_i & \delta^{\kern 0.4pt m}_j 
& \delta^{\kern 0.4pt m}_k \\
\vspace{1ex}
\delta^{\kern 0.4pt n}_i & \delta^{\kern 0.4pt n}_j 
& \delta^{\kern 0.4pt n}_k \\
\vspace{1ex}
\delta^{\kern 0.4pt k}_i & \delta^{\kern 0.4pt k}_j & 1
\endvmatrix.
\mytag{54.14}
$$
Let's insert the summation over $k$ to both sides of the formula 
\mythetag{54.14}. As a result we get the formula
$$
\hskip -2em
\sum^3_{k=1}\varepsilon^{mnk}\,\varepsilon_{ijk}
=\sum^3_{k=1}\,\vmatrix
\delta^{\kern 0.4pt m}_i & \delta^{\kern 0.4pt m}_j 
& \delta^{\kern 0.4pt m}_k \\
\vspace{1ex}
\delta^{\kern 0.4pt n}_i & \delta^{\kern 0.4pt n}_j 
& \delta^{\kern 0.4pt n}_k \\
\vspace{1ex}
\delta^{\kern 0.4pt k}_i & \delta^{\kern 0.4pt k}_j & 1
\endvmatrix.
\mytag{54.15}
$$\par
    The left hand side of the obtained formula \mythetag{54.15} coincides 
with the left hand side of the second contraction formula \mythetag{54.13}
to be proved. For this reason below we transform the right hand side of the 
formula \mythetag{54.15} only. \pagebreak Let's expand the determinant in 
the formula \mythetag{54.15} by its last row: 
$$
\aligned
\sum^3_{k=1}\varepsilon^{mnk}\,&\varepsilon_{ijk}
=\sum^3_{k=1}\left(\delta^{\kern 0.4pt k}_i
\,\vmatrix
\delta^{\kern 0.4pt m}_j 
& \delta^{\kern 0.4pt m}_k \\
\vspace{1ex}
\delta^{\kern 0.4pt n}_j 
& \delta^{\kern 0.4pt n}_k
\endvmatrix\right.-\\
&-\left.\delta^{\kern 0.4pt k}_j\,\vmatrix
\delta^{\kern 0.4pt m}_i & \delta^{\kern 0.4pt m}_k\\
\vspace{1ex}
\delta^{\kern 0.4pt n}_i & \delta^{\kern 0.4pt n}_k
\endvmatrix
+\vmatrix
\delta^{\kern 0.4pt m}_i & \delta^{\kern 0.4pt m}_j\\
\vspace{1ex}
\delta^{\kern 0.4pt n}_i & \delta^{\kern 0.4pt n}_j
\endvmatrix\right).
\endaligned
\mytag{54.16}
$$
The summation over $k$ in the right hand side of the formula 
\mythetag{54.16} applies to all of the three terms enclosed in 
the round brackets. Expanding these brackets, we get 
$$
\aligned
\sum^3_{k=1}&\varepsilon^{mnk}\,\varepsilon_{ijk}
=\sum^3_{k=1}\delta^{\kern 0.4pt k}_i
\,\vmatrix
\delta^{\kern 0.4pt m}_j 
& \delta^{\kern 0.4pt m}_k \\
\vspace{1ex}
\delta^{\kern 0.4pt n}_j 
& \delta^{\kern 0.4pt n}_k
\endvmatrix-\\
&-\sum^3_{k=1}\delta^{\kern 0.4pt k}_j\,\vmatrix
\delta^{\kern 0.4pt m}_i & \delta^{\kern 0.4pt m}_k\\
\vspace{1ex}
\delta^{\kern 0.4pt n}_i & \delta^{\kern 0.4pt n}_k
\endvmatrix
+\sum^3_{k=1}\vmatrix
\delta^{\kern 0.4pt m}_i & \delta^{\kern 0.4pt m}_j\\
\vspace{1ex}
\delta^{\kern 0.4pt n}_i & \delta^{\kern 0.4pt n}_j
\endvmatrix.
\endaligned
\mytag{54.17}
$$\par
     The first sum in the right hand side of \mythetag{54.17} contains
the factor $\delta^{\kern 0.4pt k}_i$. In performing the summation cycle
over $k$ this factor appears to be nonzero only once when $k=i$. For this
reason only one summand of the first sum does actually survive. In this 
term $k=i$. Similarly, in the second sum also only one its term does 
actually survive, in this term $k=j$. As for the last sum in the right
hand side of \mythetag{54.17}, the expression being summed does not depend
on $k$. Therefore it triples upon calculating this sum. Taking into 
account that $\delta^{\kern 0.4pt i}_i=\delta^{\kern 0.4pt j}_j=1$, we get
$$
\sum^3_{k=1}\varepsilon^{mnk}\,\varepsilon_{ijk}
=\vmatrix
\delta^{\kern 0.4pt m}_j 
& \delta^{\kern 0.4pt m}_i\\
\vspace{1ex}
\delta^{\kern 0.4pt n}_j 
& \delta^{\kern 0.4pt n}_i
\endvmatrix
-\vmatrix
\delta^{\kern 0.4pt m}_i & \delta^{\kern 0.4pt m}_j\\
\vspace{1ex}
\delta^{\kern 0.4pt n}_i & \delta^{\kern 0.4pt n}_j
\endvmatrix
+3\,\vmatrix
\delta^{\kern 0.4pt m}_i & \delta^{\kern 0.4pt m}_j\\
\vspace{1ex}
\delta^{\kern 0.4pt n}_i & \delta^{\kern 0.4pt n}_j
\endvmatrix.
$$
The first determinant in the right hand side of the above formula
differs from two others by the transposition of its columns. If
we perform this transposition once more, it changes its sign and
we get a formula with three coinciding determinants. Then, collecting 
the similar terms, we derive
$$
\hskip -2em
\sum^3_{k=1}\varepsilon^{mnk}\,\varepsilon_{ijk}
=(-1-1+3)\,\vmatrix
\delta^{\kern 0.4pt m}_i & \delta^{\kern 0.4pt m}_j\\
\vspace{1ex}
\delta^{\kern 0.4pt n}_i & \delta^{\kern 0.4pt n}_j
\endvmatrix.
\mytag{54.18}
$$
Now it is easy to see that the formula \mythetag{54.18} leads to the
required formula \mythetag{54.13}. The theorem~\mythetheorem{54.2} is
proved.\qed\enddemo
\mytheorem{54.3} The Levi-Civita symbol and the Kronecker symbol are 
related by the third contraction formula
$$
\hskip -2em
\sum^3_{j=1}\sum^3_{k=1}\varepsilon^{mjk}\,\varepsilon_{ijk}
=2\,\delta^{\kern 0.4pt m}_i.
\mytag{54.19}
$$
\endproclaim
\demo{Proof} The third contraction formula \mythetag{54.19} is derived from
the second contraction formula \mythetag{54.13}. For this purpose let's 
substitute $n=j$ into the formula \mythetag{54.13}, take into account that
$\delta^{\kern 0.4pt j}_j=1$, and insert the summation over $j$ into both
sides of the obtained equality. This yields the formula
$$
\sum^3_{j=1}\sum^3_{k=1}\varepsilon^{mjk}\,\varepsilon_{ijk}
=\sum^3_{j=1}\vmatrix
\delta^{\kern 0.4pt m}_i & \delta^{\kern 0.4pt m}_j\\
\vspace{1ex}
\delta^{\kern 0.4pt j}_i & 1
\endvmatrix=\sum^3_{j=1}(\delta^{\kern 0.4pt m}_i
-\delta^{\kern 0.4pt m}_j\,\delta^{\kern 0.4pt j}_i).
$$
Upon expanding the brackets the sum over $j$ in the right hand side 
of the above formula can be calculated explicitly:
$$
\sum^3_{j=1}\sum^3_{k=1}\varepsilon^{mjk}\,\varepsilon_{ijk}
=\sum^3_{j=1}\delta^{\kern 0.4pt m}_i-\sum^3_{j=1}
\delta^{\kern 0.4pt m}_j\,\delta^{\kern 0.4pt j}_i
=3\,\delta^{\kern 0.4pt m}_i-\delta^{\kern 0.4pt m}_i
=2\,\delta^{\kern 0.4pt m}_i.
$$
It is easy to see that the calculations performed prove the 
formula~\mythetag{54.19} and the theorem~\mythetheorem{54.3} in whole. 
\qed\enddemo
\mytheorem{54.4} The Levi-Civita symbol and the Kronecker symbol are 
related by the fourth contraction formula
$$
\hskip -2em
\sum^3_{i=1}\sum^3_{j=1}\sum^3_{k=1}\varepsilon^{ijk}
\,\varepsilon_{ijk}=6.
\mytag{54.20}
$$
\endproclaim
\demo{Proof} The fourth contraction formula \mythetag{54.20} is derived 
from the third contraction formula \mythetag{54.19}. For this purpose
we substitute $m=i$ into \mythetag{54.19} and insert the summation over 
$i$ into both sides of the obtained equality. This yields
$$
\sum^3_{i=1}\sum^3_{j=1}\sum^3_{k=1}\varepsilon^{ijk}
\,\varepsilon_{ijk}=2\,\sum^3_{i=1}\delta^{\kern 0.4pt i}_i
=2\,\sum^3_{i=1}1=2\cdot 3=6.
$$
The above calculations prove the formula \mythetag{54.20}, which completes
the proof of the theorem~\mythetheorem{54.4}.
\qed\enddemo
\head
\SectionNum{55}{131} The triple product expansion formula and the 
Jacobi identity.
\endhead
\rightheadtext{\S\,55. \dots and the Jacobi identity.}
\mytheorem{55.1} For any triple of free vectors $\bold a$, $\bold b$, and
$\bold c$ in the space $\Bbb E$ the following identity is fulfilled:
$$
\hskip -2em
[\bold a,[\bold b,\bold c]]=\bold b\,(\bold a,\bold c)
-\bold c\,(\bold a,\bold b),
\mytag{55.1}
$$
which is known as the triple product expansion formula\footnotemark.
\endproclaim
\footnotetext{\ In Russian literature the triple product expansion 
formula is known as the double vectorial product formula or the 
{\tencyr\char '074}BAC minus CAB{\tencyr\char '076} formula.}
\adjustfootnotemark{-1}
\demo{Proof} In order to prove the identity \mythetag{55.1} we choose some
right orthonormal basis $\bold e_1,\,\bold e_2,\,\bold e_3$ in the space 
$\Bbb E$ and let's expand the vectors $\bold a$, $\bold b$, and $\bold c$ 
in this basis:
$$
\xalignat 3
&\bold a=\shave{\sum^3_{i=1}}a^i\,\bold e_i,
&&\bold b=\shave{\sum^3_{j=1}}b^{\kern 0.6pt j}\,\bold e_j,
&&\bold c=\shave{\sum^3_{k=1}}c^{\kern 0.4pt k}\,\bold e_k.
\qquad\quad
\mytag{55.2}
\endxalignat
$$
Let's denote $\bold d=[\bold b,\bold c]$ and use the formula \mythetag{44.1}
for to calculate the vector $\bold d$. We write this formula as 
$$
\hskip -2em
\bold d=\sum^3_{k=1}\biggl(\,\sum^3_{i=1}\sum^3_{j=1}
b^{\kern 0.6pt i}\,c^{\kern 0.8pt j}\,\varepsilon_{ijk}
\!\biggr)\,\bold e_k.
\mytag{55.3}
$$
The formula \mythetag{55.3} is the expansion of the vector $\bold d$ in
the basis $\bold e_1,\,\bold e_2,\,\bold e_3$. Therefore we can get its
coordinates:
$$
\hskip -2em
d^{\kern 0.3pt k}=\sum^3_{i=1}\sum^3_{j=1}b^{\kern 0.6pt i}
\,c^{\kern 0.8pt j}\,\varepsilon_{ijk}.
\mytag{55.4}
$$\par
      Now we again apply the formula \mythetag{44.1} in order to calculate
the vector $[\bold a,[\bold b,\bold c]]=[\bold a,\bold d]$. In this case
we write it as follows:
$$
\hskip -2em
[\bold a,\bold d]=\sum^3_{n=1}\sum^3_{m=1}\sum^3_{k=1}
a^m\,d^{\kern 0.3pt k}\,\varepsilon_{mkn}\,\bold e_n.
\mytag{55.5}
$$
Substituting \mythetag{55.4} into the formula \mythetag{55.5}, we get
$$
\gathered
[\bold a,\bold d]=\sum^3_{n=1}\sum^3_{m=1}\sum^3_{k=1}
\sum^3_{i=1}\sum^3_{j=1}a^m\,b^{\kern 0.6pt i}\,c^{\kern 0.8pt j}
\,\varepsilon_{ijk}\,\varepsilon_{mkn}\,\bold e_n=\\
=\sum^3_{n=1}\sum^3_{m=1}\sum^3_{i=1}\sum^3_{j=1}
a^m\,b^{\kern 0.6pt i}\,c^{\kern 0.8pt j}
\,\biggl(\,\sum^3_{k=1}\varepsilon_{ijk}\,\varepsilon_{mkn}
\!\biggr)\,\bold e_n.
\endgathered
\mytag{55.6}
$$
The upper or lower position of indices in the Levi-Civita symbol does not 
matter (see \mythetag{43.5}). Therefore, taking into account \mythetag{43.8},
we can write $\varepsilon_{mkn}=\varepsilon^{mkn}=-\varepsilon^{mnk}$ and 
bring the formula \mythetag{55.6} to the following form:
$$
[\bold a,\bold d]=-\sum^3_{n=1}\sum^3_{m=1}\sum^3_{i=1}
\sum^3_{j=1}a^m\,b^{\kern 0.6pt i}\,c^{\kern 0.8pt j}
\,\biggl(\,\sum^3_{k=1}\varepsilon^{mnk}\,\varepsilon_{ijk}
\!\biggr)\,\bold e_n.\quad
\mytag{55.7}
$$
The sum enclosed in the round brackets in \mythetag{55.7} coincides with the
left hand side of the second contraction formula \mythetag{54.13}. Applying 
\mythetag{54.13}, we continue transforming the formula \mythetag{55.7}:
$$
\gather
[\bold a,\bold d]=-\sum^3_{n=1}\sum^3_{m=1}\sum^3_{i=1}
\sum^3_{j=1}a^m\,b^{\kern 0.6pt i}\,c^{\kern 0.8pt j}
\,\vmatrix
\delta^{\kern 0.4pt m}_i & \delta^{\kern 0.4pt m}_j\\
\vspace{1ex}
\delta^{\kern 0.4pt n}_i & \delta^{\kern 0.4pt n}_j
\endvmatrix
\,\bold e_n=\\
=\sum^3_{n=1}\sum^3_{m=1}\sum^3_{i=1}
\sum^3_{j=1}a^m\,b^{\kern 0.6pt i}\,c^{\kern 0.8pt j}
\bigl(\delta^{\kern 0.4pt n}_i\,\delta^{\kern 0.4pt m}_j
-\delta^{\kern 0.4pt m}_i\,\delta^{\kern 0.4pt n}_j\bigr)
\,\bold e_n=\\
\aligned
=\sum^3_{n=1}\sum^3_{m=1}&\sum^3_{i=1}\sum^3_{j=1}
a^m\,b^{\kern 0.6pt i}\,c^{\kern 0.8pt j}
\,\delta^{\kern 0.4pt n}_i\,\delta^{\kern 0.4pt m}_j
\,\bold e_n\,-\\
&-\,\sum^3_{n=1}\sum^3_{m=1}\sum^3_{i=1}\sum^3_{j=1}
a^m\,b^{\kern 0.6pt i}\,c^{\kern 0.8pt j}
\,\delta^{\kern 0.4pt m}_i\,\delta^{\kern 0.4pt n}_j
\,\bold e_n=
\endaligned\\
=\sum^3_{i=1}\sum^3_{j=1}a^j\,b^{\kern 0.6pt i}
\,c^{\kern 0.8pt j}\,\bold e_i-
\sum^3_{i=1}\sum^3_{j=1}a^i\,b^{\kern 0.6pt i}
\,c^{\kern 0.8pt j}\,\bold	 e_j.
\endgather
$$
The result obtained can be written as follows:
$$
[\bold a,[\bold b,\bold c]]=\biggl(\,\sum^3_{j=1}a^j
\,c^{\kern 0.8pt j}\biggr)\biggl(\,\sum^3_{i=1}
b^{\kern 0.6pt i}\,\bold e_i\biggr)
-\biggl(\,\sum^3_{i=1}a^i\,b^{\kern 0.6pt i}\biggr)
\biggl(\,\sum^3_{j=1}
c^{\kern 0.8pt j}\,\bold e_j\biggr).
$$\par
     Now, let's recall that the scalar product in an orthonormal basis
is calculated by the formula \mythetag{33.3}. The formula \mythetag{33.3} 
is easily recognized within the above relationship. Taking into account 
this formula, we get the equality
$$
\hskip -2em
[\bold a,[\bold b,\bold c]]=(\bold a,\bold c)\,\sum^3_{i=1}
b^{\kern 0.6pt i}\,\bold e_i-(\bold a,\bold b)\,\sum^3_{j=1}
c^{\kern 0.8pt j}\,\bold e_j.
\mytag{55.8}
$$
In order to bring \mythetag{55.8} to the ultimate form \mythetag{55.1} 
it is sufficient to find the expansions of the form \mythetag{55.2} for
\pagebreak $\bold b$ and $\bold c$ within the above formula \mythetag{55.8}. 
As a result the formula \mythetag{55.8} takes the form $[\bold a,[\bold b,
\bold c]]=(\bold a,\bold c)\,\bold b-(\bold a,\bold b)\,\bold c$, which 
in essential coincides with \mythetag{55.1}. The theorem~\mythetheorem{55.1} 
is proved.
\qed\enddemo
\mytheorem{55.2} For any triple of free vectors $\bold a$, $\bold b$, and
$\bold c$ in the space $\Bbb E$ the {\it Jacobi identity\/} is fulfilled:
$$
\hskip -2em
[\bold a,[\bold b,\bold c]]+[\bold b,[\bold c,\bold a]]
+[\bold c,[\bold a,\bold b]]=\bold 0.
\mytag{55.9}
$$
\endproclaim
\demo{Proof} The Jacoby identity \mythetag{55.9} is easily derived with the
use of the triple product expansion formula \mythetag{55.1}:
$$
\hskip -2em
[\bold a,[\bold b,\bold c]]=\bold b\,(\bold a,\bold c)
-\bold c\,(\bold a,\bold b).
\mytag{55.10}
$$
Let's twice perform the cyclic redesignation of vectors $\bold a\to\bold b\to
\bold c\to\bold a$ in the above formula \mythetag{55.10}. This yields
$$
\align
&\hskip -2em
[\bold b,[\bold c,\bold a]]=\bold c\,(\bold b,\bold a)
-\bold a\,(\bold b,\bold c),
\mytag{55.11}\\
&\hskip -2em
[\bold c,[\bold a,\bold b]]=\bold a\,(\bold c,\bold b)
-\bold b\,(\bold c,\bold a).
\mytag{55.12}
\endalign
$$
By adding the above three equalities \mythetag{55.10}, \mythetag{55.11}, and
\mythetag{55.12} in the left hand side we get the required expression 
$[\bold a,[\bold b,\bold c]]+[\bold b,[\bold c,\bold a]]+[\bold c,[\bold a,
\bold b]]$, while th right hand side of the resulting equality vanishes. 
The theorem~\mythetheorem{55.2} is proved.\qed\enddemo
\head
\SectionNum{56}{134} The product of two mixed products.
\endhead
\rightheadtext{\S\,56. \dots of two mixed products.}
\mytheorem{56.1} For any six free vectors $\bold a$, $\bold b$, $\bold c$, 
$\bold x$, $\bold y$, and $\bold z$ in the space $\Bbb E$ the following 
formula for the product of two mixed products $(\bold a,\bold b,\bold c)$
and $(\bold x,\bold y,\bold z)$ is fulfilled:
$$
\hskip -2em
(\bold a,\bold b,\bold c)\,(\bold x,\bold y,\bold z)
=\vmatrix 
(\bold a,\bold x) & (\bold b,\bold x) & (\bold c,\bold x)\\
\vspace{1ex}
(\bold a,\bold y) & (\bold b,\bold y) & (\bold c,\bold y)\\
\vspace{1ex}
(\bold a,\bold z) & (\bold b,\bold z) & (\bold c,\bold z)
\endvmatrix.\quad
\mytag{56.1}
$$
\endproclaim
     In order to prove the formula \mythetag{56.1} we need two properties 
of the matrix determinants which are known as the {\it linearity with respect
to a row} and the {\it linearity with respect to a column} (see \mycite{7}).
The first of them can be expressed by the formula
$$
\sum^r_{i=1}\alpha_i\,
\vmatrix
x^1_1 & \hdots & x^1_n\\
\vdots & \vdots & \vdots\\
x^k_1(i) & \hdots & x^k_n(i)\\
\vdots & \vdots & \vdots\\
x^n_1 & \hdots & x^n_n\\
\endvmatrix
=\vmatrix
x^1_1 & \hdots & x^1_n\\
\vdots & \vdots & \vdots\\
\dsize\sum^r_{i=1}\alpha_i\,x^k_1(i) & \hdots & \dsize\sum^r_{i=1}
\alpha_i\,x^k_n(i)\\
\vdots & \vdots & \vdots\\
x^n_1 & \hdots & x^n_n\\
\endvmatrix.
$$
\mylemma{56.1} If the $k$-th row of a square matrix is a linear combination of
some $r$ rows (non necessarily coinciding with the other rows of this matrix), 
then the determinant of such a matrix is equal to a linear combination of\/ $r$ 
separate determinants. 
\endproclaim
     The linearity with respect to a column is formulated similarly. 
\mylemma{56.2} If the $k$-th column of a square matrix is a linear combination 
of some $r$ columns (non necessarily coinciding with the other columns of this 
matrix), then the determinant of such a matrix is equal to a linear combination 
of\/ $r$ separate determinants. 
\endproclaim
     I do not write the formula illustrating the lemma~\mythelemma{56.2} since
it does not fit the width of a page in this book. This formula can be obtained 
from the above formula by transposing the matrices in both its sides. 
\demo{Proof of the theorem~\mythetheorem{56.1}} Let $\bold e_1,\,
\bold e_2,\,\bold e_3$ be some right orthonormal basis in the space
$\Bbb E$. Assume that the vectors $\bold a$, $\bold b$, $\bold c$, 
$\bold x$, $\bold y$, and $\bold z$ are given by their coordinates
in this basis. Let's denote through $L$ the left hand side
of the formula \mythetag{56.1}:
$$
\hskip -2em
L=(\bold a,\bold b,\bold c)\,(\bold x,\bold y,\bold z).
\mytag{56.2}
$$
In order to calculate $L$ we apply the formula \mythetag{46.9}. In the
case of the vectors $\bold a$, $\bold b$, $\bold c$ the formula
\mythetag{46.9} yields
$$
\hskip -2em
(\bold a,\bold b,\bold c)=\sum^3_{i=1}\sum^3_{j=1}\sum^3_{k=1}
a^i\,b^{\kern 0.6pt j}\,c^{\kern 0.4pt k}\,\varepsilon_{ijk}.
\mytag{56.3}
$$
In the case of the vectors $\bold x$, $\bold y$, $\bold z$ the formula
\mythetag{46.9} is written as 
$$
\hskip -2em
(\bold x,\bold y,\bold z)=\sum^3_{m=1}\sum^3_{n=1}\sum^3_{p=1}
x^m\,y^n\,z^p\,\varepsilon^{mnp}.
\mytag{56.4}
$$
Note that raising indices of the Levi-Civita symbol in \mythetag{56.4} does 
not change its values (see  \mythetag{43.5}). Now, multiplying the formulas 
\mythetag{56.3} and \mythetag{56.4}, we obtain the formula for $L$:
$$
L=\sum^3_{i=1}\sum^3_{j=1}\sum^3_{k=1}\sum^3_{m=1}\sum^3_{n=1}
\sum^3_{p=1}a^i\,b^{\kern 0.6pt j}\,c^{\kern 0.4pt k}\,
x^m\,y^n\,z^p\,\varepsilon^{mnp}\,\varepsilon_{ijk}.
\quad
\mytag{56.5}
$$
The product $\varepsilon^{mnp}\,\varepsilon_{ijk}$ in \mythetag{56.5} 
can be replaced by the matrix determinant taken from the first contraction
formula \mythetag{54.1}:
$$
L=\msum{3}\Sb i=1\ j=1\ k=1\\m=1\ n=1\ p=1\endSb
\kern -9pt a^i\,b^{\kern 0.6pt j}\,c^{\kern 0.4pt k}\,
x^m\,y^n\,z^p\,
\vmatrix
\delta^{\kern 0.4pt m}_i & \delta^{\kern 0.4pt m}_j 
& \delta^{\kern 0.4pt m}_k \\
\vspace{1ex}
\delta^{\kern 0.4pt n}_i & \delta^{\kern 0.4pt n}_j 
& \delta^{\kern 0.4pt n}_k \\
\vspace{1ex}
\delta^{\kern 0.4pt p}_i & \delta^{\kern 0.4pt p}_j 
& \delta^{\kern 0.4pt p}_k 
\endvmatrix.
\quad
\mytag{56.6}
$$\par
     The next step consists in applying the lemmas~\mythelemma{56.1} and
\mythelemma{56.2} in order to transform the formula \mythetag{56.6}. 
Applying the lemma~\mythelemma{56.2}, we bring the sum over $i$ and the 
associated factor $a^i$ into the first column of the determinant. 
Similarly, we bring the sum over $j$ and the factor $b^{\kern 0.6pt j}$ 
into the second column, and finally, we bring the sum over $k$ and 
its associated factor $c^{\kern 0.4pt k}$ into the third column of
the determinant. Then we apply the lemma~\mythelemma{56.1} in order to
distribute the sums over $m$, $n$, and $p$ to the rows of the determinant.
Simultaneously, we distribute the associated factors $x^m$, $y^n$, and
$z^p$ to the rows of the determinant. As a result of our efforts we get
the following formula:
$$
L=\vmatrix
\dsize\sum^3_{i=1}\sum^3_{m=1}a^i\,x^m\,\delta^{\kern 0.4pt m}_i & 
\dsize\sum^3_{j=1}\sum^3_{m=1}b^{\kern 0.6pt j}\,x^m
\,\delta^{\kern 0.4pt m}_j & 
\dsize\sum^3_{k=1}\sum^3_{m=1}c^{\kern 0.4pt k}\,x^m
\,\delta^{\kern 0.4pt m}_k \\
\dsize\vspace{1ex}
\dsize\sum^3_{i=1}\sum^3_{n=1}a^i\,y^n\,\delta^{\kern 0.4pt n}_i & 
\dsize\sum^3_{j=1}\sum^3_{n=1}b^{\kern 0.6pt j}\,y^n
\,\delta^{\kern 0.4pt n}_j & 
\dsize\sum^3_{k=1}\sum^3_{n=1}c^{\kern 0.4pt k}\,y^n
\,\delta^{\kern 0.4pt n}_k \\
\dsize\vspace{1ex}
\dsize\sum^3_{i=1}\sum^3_{p=1}a^i\,z^p\,\delta^{\kern 0.4pt p}_i & 
\dsize\sum^3_{j=1}\sum^3_{p=1}b^{\kern 0.6pt j}\,z^p
\,\delta^{\kern 0.4pt p}_j & 
\dsize\sum^3_{k=1}\sum^3_{p=1}c^{\kern 0.4pt k}\,z^p
\,\delta^{\kern 0.4pt p}_k 
\endvmatrix.
$$
Due to the Kronecker symbols the double sums in the above formula are 
reduced to single sums:
$$
\hskip -2em
L=\vmatrix
\dsize\sum^3_{i=1}a^i\,x^i & 
\dsize\sum^3_{j=1}b^{\kern 0.6pt j}\,x^{\kern 0.3pt j} & 
\dsize\sum^3_{k=1}c^{\kern 0.4pt k}\,x^k\\
\dsize\vspace{1ex}
\dsize\sum^3_{i=1}a^i\,y^i & 
\dsize\sum^3_{j=1}b^{\kern 0.6pt j}\,y^{\kern 0.3pt j} & 
\dsize\sum^3_{k=1}c^{\kern 0.4pt k}\,y^k\\
\dsize\vspace{1ex}
\dsize\sum^3_{i=1}a^i\,z^i & 
\dsize\sum^3_{j=1}b^{\kern 0.6pt j}\,z^j & 
\dsize\sum^3_{k=1}c^{\kern 0.4pt k}\,z^k
\endvmatrix.
\mytag{56.7}
$$\par
     Let's recall that our basis $\bold e_1,\,\bold e_2,\,\bold e_3$ is
orthonormal. The scalar product in such a basis is calculated according to
the formula \mythetag{33.3}. Comparing \mythetag{33.3} with \mythetag{56.7}, 
we see that all of the sums within the determinant \mythetag{56.7} 
are scalar products of vectors:
$$
\pagebreak
\hskip -2em
L=\vmatrix 
(\bold a,\bold x) & (\bold b,\bold x) & (\bold c,\bold x)\\
\vspace{1ex}
(\bold a,\bold y) & (\bold b,\bold y) & (\bold c,\bold y)\\
\vspace{1ex}
(\bold a,\bold z) & (\bold b,\bold z) & (\bold c,\bold z)
\endvmatrix.
\mytag{56.8}
$$
Now the formula \mythetag{56.1} follows from the formulas \mythetag{56.2} 
and \mythetag{56.8}. The theorem~\mythetheorem{56.1} is proved. 
\qed\enddemo
     The formula \mythetag{56.1} is valuable for us not by itself, but due to 
its consequences. Assume that $\bold e_1,\,\bold e_2,\,\bold e_3$ is some arbitrary
basis. Let's substitute $\bold a=\bold e_1$, $\bold b=\bold e_2$, $\bold c=\bold e_3$, 
$\bold x=\bold e_1$, $\bold y=\bold e_2$, $\bold z=\bold e_3$ into the formula
\mythetag{56.1}. As a result we obtain
$$
(\bold e_1,\bold e_2,\bold e_3)^2
=\vmatrix 
(\bold e_1,\bold e_1) & (\bold e_2,\bold e_1) 
& (\bold e_3,\bold e_1)\\
\vspace{1.5ex}
(\bold e_1,\bold e_2) & (\bold e_2,\bold e_2) 
& (\bold e_3,\bold e_2)\\
\vspace{1.5ex}
(\bold e_1,\bold e_3) & (\bold e_2,\bold e_3) 
& (\bold e_3,\bold e_3)
\endvmatrix.\quad
\mytag{56.9}
$$\par 
     Comparing \mythetag{56.9} with \mythetag{29.7} and \mythetag{29.6}, we
     see that the matrix \mythetag{56.9} differs from the Gram matrix $G$ by
transposing. If we take into account the symmetry of the Gram matrix (see
theorem~\mythetheorem{30.1}), than we find that they do coincide. Therefore
the formula \mythetag{56.9} is written as follows:
$$
\hskip -2em
(\bold e_1,\bold e_2,\bold e_3)^2=\det G. 
\mytag{56.10}
$$
The formula \mythetag{56.10} coincides with the formula \mythetag{51.6}. This fact
proves the theorem~\mythetheorem{51.2}, which was unproved in \S\,51.\par
\newpage
%---------------------------------------------------------------
\setfirstpage
\topmatter
\title\chapter{2}
Geometry of lines and surfaces.
\endtitle
\leftheadtext{CHAPTER~\uppercase\expandafter{\romannumeral 2}.
GEOMETRY OF LINES AND SURFACES.}
\endtopmatter
\document
\chapternum=2
      In this Chapter the tools of the vector algebra developed in
Chapter~\uppercase\expandafter{\romannumeral 1} are applied for describing
separate points of the space $\Bbb E$ and for describing some geometric forms
which are composed by these points. 
\head
\SectionNum{1}{139} Cartesian coordinate systems.
\endhead
\rightheadtext{\S\,1. Cartesian coordinate systems.}
\mydefinition{1.1} \parshape 4 0cm 10cm 0cm 10cm 0cm 10cm 5cm 5cm 
A {\it Cartesian coordinate system\/} in the space $\Bbb E$ is a
basis $\bold e_1,\,\bold e_2,\,\bold e_3$ \vadjust{\vskip 5pt\hbox to 
0pt{\kern 0pt\includegraphics{angemeng29.eps}\hss}
\vskip -5pt}complemented by some fixed point $O$ of this space. 
The point $O$ being a part of the Cartesian coordinate system 
$O,\,\bold e_1,\,\bold e_2,\,\bold e_3$ is called the {\it origin}
of this coordinate system.
\enddefinition
\mydefinition{1.2} \parshape 9 5cm 5cm 5cm 5cm 5cm 5cm 5cm 5cm 
5cm 5cm 5cm 5cm 5cm 5cm 5cm 5cm 0cm 10cm 
The vector $\overrightarrow{OX\,\,}\!$ binding the origin $O$ of a 
Cartesian coordinate system $O$, $\bold e_1,\,\bold e_2,\,\bold e_3$
with a point $X$ of the space $\Bbb E$ is called the 
{\it radius vec\-tor\/} of the point $X$ in this coordinate system. 
\enddefinition
\parshape 3 5cm 5cm 5cm 5cm 0cm 10cm 
     The free vectors $\bold e_1,\,\bold e_2,\,\bold e_3$ being constituent
parts of a Cartesian coordinate system $O,\,\bold e_1,\,\bold e_2,\,\bold e_3$ 
remain free. However they are often represented by geometric realizations
attached to the origin $O$ (see Fig\.~1.1). These geometric realizations are 
extended up to the whole lines which are called the {\it coordinate axes} of
the Cartesian coordinate system $O,\,\bold e_1,\,\bold e_2,\,\bold e_3$. 
\mydefinition{1.3} The {\it coordinates\/} of a point $X$ in a Cartesian
coordinate system $O,\,\bold e_1,\,\bold e_2,\,\bold e_3$ are the coordinates
of its radius vector $\bold r_X=\overrightarrow{OX\,\,}\!$ in the basis
$\bold e_1,\,\bold e_2,\,\bold e_3$.
\enddefinition
      Like other vectors, radius vectors of points are covered by the index 
setting convention (see Definition~\mythedefinitionchapter{20.1}{1} in 
Chapter~\uppercase\expandafter{\romannumeral 1}). The coordinates of 
radius vectors are enumerated by upper indices, they are usually arranged into 
columns:
$$
\hskip -2em
\bold r_X=\Vmatrix x^1\\\vspace{0.4ex}x^2\\\vspace{0.4ex}x^3\endVmatrix.
\mytag{1.1}
$$
However, in those cases where a point $X$ is represented by its coordinates 
these coordinates are arranged into a comma-separated row and placed within 
round brackets just after the symbol denoting the point $X$ itself:
$$
\hskip -2em
X=X(x^1,x^2,x^3).
\mytag{1.2}
$$
The upper position of indices in the formula \mythetag{1.2} is inherited 
from the formula \mythetag{1.1}.\par
      Cartesian coordinate systems can be either {\it rectangular} or 
{\it skew-angular}  depending on the type of basis used for defining them. 
In this book I follow a convention similar to that of the 
definition~\mythedefinitionchapter{29.1}{1} in 
Chapter~\uppercase\expandafter{\romannumeral 1}. 
\mydefinition{1.4} In this book a {\it skew-angular coordinate system\/} is 
understood as an {\it arbitrary coordinate system} where no restrictions for
the basis are imposed.
\enddefinition
\mydefinition{1.5} A rectangular coordinate system $O,\,\bold e_1,\,\bold e_2$,
$\bold e_3$ whose basis is orthonormal is called a {\it rectangular coordinate
system with unit scales along the axes}. 
\enddefinition
{\bf A remark}. The radius vector of a point $X$ written as 
$\overrightarrow{OX\,\,}\!$ is a geometric vector whose position in the space is
fixed. The radius vector of a point $X$ written as $\bold r_X$ is a free vector
We can perform various operations of vector algebra with this vector: we can
add it with other free vectors, multiply it by numbers, compose scalar products 
etc. But the vector $\bold r_X$ has a special mission --- to be a pointer to the
point $X$. It can do this mission only if its initial point is placed to the
origin. 
\myexercise{1.1} Relying on the definition~\mythedefinition{1.1}, formulate
analogous definitions for Cartesian coordinate systems on a line and on a plane.
\endproclaim
\head
\SectionNum{2}{141} Equations of lines and surfaces.
\endhead
\rightheadtext{\S\,2. Equations of lines and surfaces.}
     Coordinates of a single fixed point $X$ in the space $\Bbb E$ are three
fixed numbers (three constants). If these coordinates are changing, we get a 
moving point that runs over some set within the space $\Bbb E$. In this book
we consider the cases where this set is some line or some surface. The case
of a line differs from the case of a surface by its {\it dimension\/} or, in
other words, by the {\it number of degrees of freedom}. Each surface is
{\it two-dimensional} --- a point on a surface has two degrees of freedom. 
Each line is {\it one-dimensional} --- a point on a line has one degree of 
freedom.\par
     Lines and surfaces contain infinitely many points. Therefore they cannot 
be described by enumeration where each point is described separately. Lines
and surfaces are described by means of equations. The equations of lines and
surfaces are subdivided into several types. If the radius vector of a point
enters an equation as a whole without subdividing it into separate coordinates,
such an equation is called a {\it vectorial equation}. If the radius vector 
of a point enters an equation through its coordinates, such an equation 
is called a {\it coordinate equation}.\par
      Another attribute of the equations is the method of implementing the 
degrees of freedom. One or two degrees of freedom can be implemented in an
equation explicitly when the radius vector of a point is given as a function 
of one or two variables, which are called {\it parameters}. In this case the
equation is called {\it parametric}. {\it Non-parametric equations} behave
as obstructions decreasing the number of degrees of freedom from the initial
three to one or two.
\head
\SectionNum{3}{142} A straight line on a plane.
\endhead
\rightheadtext{\S\,3. A straight line on a plane.}
\parshape 3 0cm 10cm 0cm 10cm 5cm 5cm 
       Assume that some plane $\alpha$ in the space $\Bbb E$ is chosen and 
fixed. Then the number of degrees of freedom of a point immediately decreases
from three to two. \vadjust{\vskip 5pt\hbox to 
0pt{\kern 0pt\includegraphics{angemeng30.eps}\hss}
\vskip -5pt}In order to study various forms of equations defining a straight
line on the plane $\alpha$ we choose some coordinate system $O,\,\bold e_1,
\,\bold e_2$ on this plane. Then we can define the points of the plane $\alpha$
and the points of a line on it by means of their radius-vectors.\par
\parshape 5 5cm 5cm 5cm 5cm  5cm 5cm 5cm 5cm 0cm 10cm 
     {\bf \hskip 0pt plus 5pt minus 5pt 1. Vectorial parametric equation} of a
line on a plane. Let's consider a line on a plane with some coordinate system
$O$, $\bold e_1,\,\bold e_2$. Let $X$ be some arbitrary point on this line
(see Fig\.~3.1) and let $A$ be some fixed point of this line. The position 
of the point $X$ relative to the point $A$ is marked by the vector
$\overrightarrow{AX\,\,}\!$, while the position of the point $A$ itself is
determined by its radius vector $\bold r_0=\overrightarrow{OA\,\,}\!$. 
Therefore we have
$$
\hskip -2em
\bold r=\overrightarrow{OX\,\,}\!=\bold r_0+\overrightarrow{AX\,\,}\!.
\mytag{3.1}
$$
Let's choose and fix some nonzero vector $\bold a\neq\bold 0$ directed along 
the line in question. The vector $\overrightarrow{AX\,\,}\!$ is expressed 
through $\bold a$ by means of the following formula:
$$
\hskip -2em
\overrightarrow{AX\,\,}\!=\bold a\cdot t.
\mytag{3.2}
$$
From the formulas \mythetag{3.1} and \mythetag{3.2} we immediately derive:
$$
\hskip -2em
\bold r=\bold r_0+\bold a\cdot t.
\mytag{3.3}
$$\vskip -1ex
\mydefinition{3.1} The equality \mythetag{3.3} is called the {\it 
vectorial parametric\/} equation of a line on a plane. The constant vector
$\bold a\neq\bold 0$ in it is a {\it directional vector} of the line, 
while the variable $t$ is a {\it parameter}. The constant vector 
$\bold r_0$ in \mythetag{3.3} is the {\it radius vector of an initial 
point}. 
\enddefinition
     Each particular value of the parameter $t$ corresponds to some definite point 
on the line. The initial point $A$ with the radius vector $\bold r_0$ is associated 
with the value $t=0$.\par
     {\bf 2. Coordinate parametric equations} of a line on a plane. Let's 
determine the vectors $\bold r$, $\bold r_0$, and $\bold a$ in the vectorial 
parametric equation \mythetag{3.3} through their coordinates:
$$
\xalignat 3
&\hskip -2em
\bold r=\Vmatrix x\\ \vspace{0.5ex} y\endVmatrix, 
&&\bold r_0=\Vmatrix x_0\\ \vspace{0.5ex} y_0\endVmatrix, 
&&\bold a=\Vmatrix a_x\\ \vspace{0.5ex} a_y\endVmatrix.
\quad
\mytag{3.4}
\endxalignat
$$
Due to \mythetag{3.4} the equation \mythetag{3.3} is written as two equations:
$$
\hskip -2em
\cases x=x_0+a_x\,t,\\ y=y_0+a_y\,t.\endcases
\mytag{3.5}
$$\vskip -1ex
\mydefinition{3.2} The equalities \mythetag{3.5} are called the 
{\it coordinate parametric \/} equations of a straight line on a plane. The 
constants $a_x$ and $a_y$ in them cannot vanish simultaneously.
\enddefinition
     {\bf 3. Normal vectorial equation} of a line on a plane. Let $\bold n\neq
\bold 0$ be a vector lying on the plane $\alpha$ in question and being 
perpendicular to the line in question (see Fig\.~3.1). Let's apply the scalar
multiplication by the vector $\bold n$ to both sides of the equation 
\mythetag{3.3}. As a result we get
$$
\hskip -2em
(\bold r,\bold n)=(\bold r_0,\bold n)+(\bold a,\bold n)\,t.
\mytag{3.6}
$$
But $\bold a\perp\bold n$. Fort this reason the second term in the right hand side
of \mythetag{3.6} vanishes and the equation actually does not have the parameter $t$. 
The resulting equation is usually written as 
$$
\hskip -2em
(\bold r-\bold r_0,\bold n)=0.
\mytag{3.7}
$$
The scalar product of two constant vectors $\bold r_0$ and $\bold n$ is a numeric 
constant. If we denote $D=(\bold r_0,\bold n)$, then the equation \mythetag{3.7} 
can be written in the following form:
$$
\hskip -2em
(\bold r,\bold n)=D.
\mytag{3.8}
$$\vskip -1ex
\mydefinition{3.3} Any one of the two equations \mythetag{3.7} and \mythetag{3.8} 
is called the {\it normal vectorial\/} equation of a line on a plane. The constant
vector $\bold n\neq\bold 0$ in these equations is called a {\it normal vector\/} 
of this line. 
\enddefinition
     {\bf 4. Canonical equation\/} of a line on a plane. Let's consider the case 
where $a_x\neq 0$ and $a_y\neq 0$ in the equations \mythetag{3.5}. In this case 
the parameter $t$ can be expressed through $x$ and $y$:
$$
\xalignat 2
&\hskip -2em
t=\frac{x-x_0}{a_x}, 
&&t=\frac{y-y_0}{a_y}. 
\mytag{3.9}
\endxalignat
$$
From the equations \mythetag{3.9} we derive the following equality:
$$
\hskip -2em
\frac{x-x_0}{a_x}=\frac{y-y_0}{a_y}. 
\mytag{3.10}
$$
If $a_x=0$, then the the first of the equations \mythetag{3.5} turns to
$$
\hskip -2em
x=x_0.
\mytag{3.11}
$$
If $a_y=0$, then the second equation \mythetag{3.5} turns to
$$
\pagebreak
\hskip -2em
y=y_0.
\mytag{3.12}
$$
Like the equation \mythetag{3.10}, the equations \mythetag{3.11} and
\mythetag{3.12} do not comprise the parameter $t$. 
\mydefinition{3.4} Any one of the three equalities \mythetag{3.10}, 
\mythetag{3.11}, and \mythetag{3.12} is called the  {\it canonical\/} 
equation of a line on a plane. The constants $a_x$ and $a_y$ in the
equation \mythetag{3.10} should be nonzero.
\enddefinition
     {\bf 5. The equation of a line passing through two given points\/} 
on a plane. Assume that two distinct points $A\neq B$ on a plane are 
given. We write their coordinates
$$
\xalignat 2
&\hskip -2em
A=A(x_0,y_0), 
&&B=B(x_1,y_1). 
\mytag{3.13}
\endxalignat
$$
The vector $\bold a=\overrightarrow{AB\,\,}\!$ can be used for the
directional vector of the line passing through the points $A$ and $B$
in \mythetag{3.13}. Then from \mythetag{3.13} we derive the coordinates
of $\bold a$:
$$
\hskip -2em
\bold a=\Vmatrix a_x\\ \vspace{0.5ex} a_y\endVmatrix
=\Vmatrix x_1-x_0\\ \vspace{0.5ex} y_1-y_0\endVmatrix.
\mytag{3.14}
$$
Due to \mythetag{3.14} the equation \mythetag{3.10} can be written as 
$$
\hskip -2em
\frac{x-x_0}{x_1-x_0}=\frac{y-y_0}{y_1-y_0}. 
\mytag{3.15}
$$
The equation \mythetag{3.15} corresponds to the case where $x_1\neq x_0$ and
$y_1\neq y_0$. If $x_1=x_0$, then we write the equation \mythetag{3.11}:
$$
\hskip -2em
x=x_0=x_1.
\mytag{3.16}
$$
If $y_1=y_0$, we write the equation \mythetag{3.12}:
$$
\hskip -2em
y=y_0=y_1.
\mytag{3.17}
$$
The conditions $x_1=x_0$ and $y_1=y_0$ cannot be fulfilled simultaneously since
the points $A$ and $B$ are distinct, i\.\,e\. $A\neq B$.
\mydefinition{3.5} Any one of the three equalities \mythetag{3.15}, \mythetag{3.16}, 
and \mythetag{3.17} is called the {\it equation of a line passing through two 
given points\/} \mythetag{3.13} on a plane.
\enddefinition
     {\bf 6. General equation\/} of a line on a plane. Let's apply the formula
\mythetagchapter{29.8}{1} from Chapter~\uppercase\expandafter{\romannumeral 1} in 
order to calculate the scalar product in \mythetag{3.8}. In this particular case 
it is written as 
$$
\hskip -2em
(\bold r,\bold n)=\sum^2_{i=1}\sum^2_{j=1}r^i\,n^{\kern 0.6pt j}
\,g_{ij},
\mytag{3.18}
$$
where $g_{ij}$ are the components of the Gram matrix for the basis $\bold e_1,
\,\bold e_2$ on our plane (see Fig\.~3.1). In Fig\.~3.1 the basis $\bold e_1,
\,\bold e_2$ is drawn to be rectangular. However, in general case it could be
skew-angular as well. Let's introduce the notations
$$
\hskip -2em
n_i=\sum^2_{j=1}n^{\kern 0.6pt j}\,g_{ij}.
\mytag{3.19}
$$
The quantities $n^1$ and $n^2$ in \mythetag{3.18} and in \mythetag{3.19} are the
coordinates of the normal vector $\bold n$ (see Fig\.~3.1). 
\mydefinition{3.6} The quantities $n_1$ and $n_2$ produced from the coordinates 
of the normal vector $\bold n$ by mens of the formula \mythetag{3.19} are
called the {\it covariant coordinates\/} of the vector $\bold n$. 
\enddefinition
     Taking into account the notations \mythetag{3.19}, the formula \mythetag{3.18} 
is written in the following simplified form:
$$
\hskip -2em
(\bold r,\bold n)=\sum^2_{i=1}r^i\,n_i.
\mytag{3.20}
$$
Let's recall the previous notations \mythetag{3.4} and introduce new ones:
$$
\pagebreak
\xalignat 2
&\hskip -2em
A=n_1, 
&&B=n_2. 
\mytag{3.21}
\endxalignat
$$
Due to \mythetag{3.4}, \mythetag{3.20}, and \mythetag{3.21} the equation 
\mythetag{3.8} is written as
$$
\hskip -2em
A\,x+B\,y-D=0.
\mytag{3.22}
$$
\mydefinition{3.7} The equation \mythetag{3.22} is called the {\it general
equation\/} of a line on a plane. 
\enddefinition
\parshape 14 0cm 10cm 0cm 10cm 5cm 5cm 5cm 5cm 5cm 5cm 5cm 5cm 
5cm 5cm 5cm 5cm 5cm 5cm 5cm 5cm 5cm 5cm 5cm 5cm 5cm 5cm 0cm 10cm 
     {\bf 7. Double intersect equation\/} of a line on a plane. Let's consider
a line on a plane that does not pass through the origin and intersects with both
of the coordinate axes. \vadjust{\vskip 5pt\hbox 
to 0pt{\kern 0pt\includegraphics{angemeng31.eps}\hss}
\vskip -5pt}These conditions mean that $D\neq 0$, $A\neq 0$, and $B\neq 0$ in 
the equation \mythetag{3.22} of this line. Through $X$ and $Y$ in Fig\.~3.2 
two intercept points are denoted:
$$
\hskip -2em
\aligned
&X=X(a,0),\\
&Y=Y(0,b).
\endaligned
\mytag{3.23}
$$
The quantities $a$ and $b$ in \mythetag{3.23} are expressed through
the constant parameters $A$, $B$, and $D$ of the equation \mythetag{3.22}
by means of the following formulas:
$$
\xalignat 2
&\hskip -2em
a=D/A, &&b=D/B.
\mytag{3.24}
\endxalignat
$$
The equation \mythetag{3.22} itself in our present case can be written as 
$$
\hskip -2em
\frac{x}{D/A}+\frac{y}{D/B}=1.
\mytag{3.25}
$$
If we take into account \mythetag{3.24}, the equality \mythetag{3.25} turns to
$$
\pagebreak
\hskip -2em
\frac{x}{a}+\frac{y}{b}=1.
\mytag{3.26}
$$
\mydefinition{3.8} The equality \mythetag{3.26} is called the {\it double intercept
equation\/} of a line on a plane.
\enddefinition
     The name of the equation \mythetag{3.26} is due to the parameters $a$ and $b$ 
being the lengths of the segments $[OX]$ and $[OY]$ which the line intercepts
on the coordinate axes. 
\head
\SectionNum{4}{148} A plane in the space.
\endhead
\rightheadtext{\S\,4. A plane in the space.}
\parshape 3 0cm 10cm 0cm 10cm 5cm 5cm 
       Assume that some plane $\alpha$ in the space $\Bbb E$ is chosen and fixed. 
\vadjust{\vskip 5pt\hbox to 
0pt{\kern 0pt\includegraphics{angemeng32.eps}\hss}
\vskip -5pt}In order to study various equations determining this plane we choose 
some coordinate system $O,\,\bold e_1,\,\bold e_2,\,\bold e_3$ in the
space $\Bbb E$. Then we can describe the points of this plane $\alpha$ by their
radius vectors.\par
\parshape 11 5cm 5cm 5cm 5cm 5cm 5cm 5cm 5cm 5cm 5cm 5cm 5cm 5cm 5cm 
5cm 5cm	 5cm 5cm 5cm 5cm 0cm 10cm 
     {\bf \hskip 0pt plus 5pt minus 5pt 1. Vectorial parametric equation\/} of a
plane in the spa\-ce. Let's denote through $A$ some fixed point of the plane $\alpha$ 
(see Fig\.~4.1) and denote through $X$ some arbitrary point of this pla\-ne. 
The position of the point $X$ relative to the point $A$ is marked by the vector
$\overrightarrow{AX\,\,}\!$, while the position of the point $A$ is determined by 
its radius vector $\bold r_0=\overrightarrow{OA\,\,}\!$. For this reason
$$
\hskip -2em
\bold r=\overrightarrow{OX\,\,}\!=\bold r_0+\overrightarrow{AX\,\,}\!.
\mytag{4.1}
$$\par 
     Let's choose and fix some pair of non-collinear vectors $\bold a\nparallel
\bold b$ lying on the plane $\alpha$. Such vectors constitute a basis on this 
plane (see Definition~\mythedefinitionchapter{17.1}{1} from
Chapter~\uppercase\expandafter{\romannumeral 1}). The vector 
$\overrightarrow{AX\,\,}\!$ is expanded in the basis of the vectors $\bold a$ 
and $\bold b$:
$$
\hskip -2em
\overrightarrow{AX\,\,}\!=\bold a\cdot t+\bold b\cdot\tau.
\mytag{4.2}
$$
Since $X$ is an arbitrary point of the plane $\alpha$, the numbers 
$t$ and $\tau$ in \mythetag{4.2} are two variable parameters. Upon substituting
\mythetag{4.2} into the formula \mythetag{4.1}, we get the following equality:
$$
\hskip -2em
\bold r=\bold r_0
+\bold a\cdot t+\bold b\cdot\tau.
\mytag{4.3}
$$
\mydefinition{4.1} The equality \mythetag{4.3} is called the {\it vectorial
parametric equation\/} of a plane in the space. The non-collinear vectors 
$\bold a$ and $\bold b$ in it are called {\it directional vectors} of a plane, 
while $t$ and $\tau$ are called {\it parameters}. The fixed vector $\bold r_0$ 
is the {\it radius vector of an initial point}. 
\enddefinition
     {\bf 2. Coordinate parametric equation} of a plane in the space. Let's 
determine the vectors $\bold r$, $\bold r_0$, $\bold a$, $\bold b$, in the 
vectorial parametric equation \mythetag{4.3} through their coordinates:     
$$
\xalignat 4
&\bold r=\Vmatrix x\\ \vspace{0.5ex} y\\ \vspace{0.5ex} z\endVmatrix, 
&&\bold r_0=\Vmatrix x_0\\ \vspace{0.5ex} y_0\\ \vspace{0.5ex} z_0
\endVmatrix, 
&&\bold a=\Vmatrix a_x\\ \vspace{0.5ex} a_y\\ \vspace{0.5ex} a_z
\endVmatrix,
&&\bold b=\Vmatrix b_x\\ \vspace{0.5ex} b_y\\ \vspace{0.5ex} b_z
\endVmatrix.
\qquad\quad
\mytag{4.4}
\endxalignat
$$
Due to \mythetag{4.4} the equation \mythetag{4.3} is written as three equations
$$
\hskip -2em
\cases x=x_0+a_x\,t+b_x\,\tau,
\\ y=y_0+a_y\,t+b_y\,\tau,
\\ z=z_0+a_z\,t+b_z\,\tau.
\endcases
\mytag{4.5}
$$
\mydefinition{4.2} The equalities \mythetag{4.5} are called the {\it coordinate 
parametric equations\/} of a plane in the space. The triples of constants $a_x,
\,a_y,\,a_z$ and $b_x,\,b_y,\,b_z$ in these equations cannot be proportional to 
each other.
\enddefinition
     {\bf 3. Normal vectorial equation} of a plane in the space. Let $\bold n\neq
\bold 0$ be a vector perpendicular to the plane $\alpha$ (see Fig\.~4.1). Let's 
apply the scalar multiplication by the vector $\bold n$ to both sides of the 
equation \mythetag{4.3}. As a result we get
$$
\hskip -2em
(\bold r,\bold n)=(\bold r_0,\bold n)+(\bold a,\bold n)\,t
+(\bold b,\bold n)\,\tau.
\mytag{4.6}
$$
But $\bold a\perp\bold n$ and $\bold b\perp\bold n$. For this reason the second 
and the third terms in the right hand side of \mythetag{4.6} vanish and the
equation actually does not have the parameters $t$ and $\tau$. The resulting 
equation is usually written as 
$$
\hskip -2em
(\bold r-\bold r_0,\bold n)=0.
\mytag{4.7}
$$
The scalar product of two constant vectors $\bold r_0$ and $\bold n$ is a numeric
constant. If we denote $D=(\bold r_0,\bold n)$, then the equation \mythetag{4.7} 
can be written in the following form:
$$
\hskip -2em
(\bold r,\bold n)=D.
\mytag{4.8}
$$
\mydefinition{4.3} Any one of the two equalities \mythetag{4.7} and \mythetag{4.8} 
is called the {\it normal vectorial equation\/} of a plane in the space. The constant
vector $\bold n\neq\bold 0$ in these equations is called a {\it normal vector\/} 
of this plane.
\enddefinition
     {\bf 4. Canonical equation\/} of a plane in the space. The vector product of
two non-coplanar vectors $\bold a\nparallel\bold b$ lying on the plane $\alpha$ can 
be chosen for the normal vector $\bold n$ in \mythetag{4.7}. Substituting 
$\bold n=[\bold a,\bold b]$ into the equation \mythetag{4.7}, we get the relationship
$$
\hskip -2em
(\bold r-\bold r_0,[\bold a,\bold b])=0.
\mytag{4.9}
$$
Applying the definition of the mixed product (see formula \mythetagchapter{45.1}{1} in
Chapter~\uppercase\expandafter{\romannumeral 1}), the scalar product in the formula
\mythetag{4.9} can be transformed to a mixed product:
$$
\hskip -2em
(\bold r-\bold r_0,\bold a,\bold b)=0.
\mytag{4.10}
$$\par
     Let's transform the equation \mythetag{4.10} into a coordinate form. For this
purpose we use the coordinate presentations of the vectors $\bold r$, $\bold r_0$, 
$\bold a$, and $\bold b$ taken from \mythetag{4.4}. Applying the formula
\mythetagchapter{50.8}{1} from Chapter~\uppercase\expandafter{\romannumeral 1}) 
\pagebreak and taking into account the fact that the oriented volume of a basis 
$(\bold e_1,\bold e_2,\bold e_3)$ is nonzero, we derive
$$
\hskip -2em
\vmatrix x-x_0 & y-y_0 & z-z_0\\ 
\vspace{0.5ex}a_x & a_y & a_z\\
\vspace{1ex}b_x & b_y & b_z
\endvmatrix=0.
\mytag{4.11}
$$\par
     If we use not the equation \mythetag{4.7}, but the equation \mythetag{4.8}, then
instead of the equation \mythetag{4.10} we get the following relationship:
$$
\hskip -2em
(\bold r,\bold a,\bold b)=D.
\mytag{4.12}
$$
In the coordinate form the relationship \mythetag{4.12} is written as
$$
\hskip -2em
\vmatrix x & y & z\\ 
\vspace{0.5ex}a_x & a_y & a_z\\
\vspace{1ex}b_x & b_y & b_z
\endvmatrix=D.
\mytag{4.13}
$$
\mydefinition{4.4} Any one of the two equalities \mythetag{4.11} and \mythetag{4.13} 
is called the {\it canonical equation\/} of a plane in the space. The triples of
constants $a_x,\,a_y,\,a_z$ and $b_x,\,b_y,\,b_z$ in these equations should not be
proportional to each other.
\enddefinition
\mydefinition{4.5} The equation \mythetag{4.10} and the equation \mythetag{4.12}, 
where $\bold a\nparallel\bold b$, are called the {\it vectorial forms of the 
canonical equation\/} of a plane in the space.
\enddefinition
     {\bf 5. The equation of a plane passing through three given points\/}. Assume 
that three points $A$, $B$, and $C$ in the space $\Bbb E$ nit lying on a straight
line are given. We write their coordinates
$$
\align
&\hskip -2em
A=A(x_0,y_0,z_0),\\
&\hskip -2em
B=B(x_1,y_1,z_1),
\mytag{4.14}\\
&\hskip -2em
C=C(x_2,y_2,z_2). 
\endalign
$$
The vectors $\bold a=\overrightarrow{AB\,\,}\!$ and 
$\bold b=\overrightarrow{AC\,\,}\!$ can be chosen for the directional vectors
of a plane \pagebreak passing through the points $A$, $B$, and $C$. Then from 
the formulas \mythetag{4.14} we derive
$$
\bold a=\Vmatrix a_x\\ \vspace{0.5ex} a_y\\ \vspace{0.5ex} a_z
\endVmatrix
=\Vmatrix x_1-x_0\\ \vspace{0.5ex} y_1-y_0\\ \vspace{0.5ex} 
z_1-z_0\endVmatrix,\quad
\bold b=\Vmatrix b_x\\ \vspace{0.5ex} b_y\\ \vspace{0.5ex} b_z
\endVmatrix
=\Vmatrix x_2-x_0\\ \vspace{0.5ex} y_2-y_0\\ \vspace{0.5ex} 
z_2-z_0\endVmatrix.\quad
\mytag{4.15}
$$
Due to \mythetag{4.15} the equation \mythetag{4.11} can be written as 
$$
\hskip -2em
\vmatrix x-x_0 & y-y_0 & z-z_0\\ 
\vspace{0.5ex}x_1-x_0 & y_1-y_0 & z_1-z_0\\
\vspace{1ex}x_2-x_0 & y_2-y_0 & z_2-z_0
\endvmatrix=0.
\mytag{4.16}
$$\par
     If we denote through $\bold r_1$ and $\bold r_2$ the radius vectors 
of the points $B$ and $C$ from \mythetag{4.14}, then \mythetag{4.15} can be
written as
$$
\xalignat 2
&\hskip -2em
\bold a=\bold r_1-\bold r_0,
&&\bold b=\bold r_2-\bold r_0.
\mytag{4.17}
\endxalignat
$$
Substituting \mythetag{4.17} into the equation \mythetag{4.10}, we get
$$
\hskip -2em
(\bold r-\bold r_0,\bold r_1-\bold r_0,\bold r_2-\bold r_0)=0.
\mytag{4.18}
$$
\mydefinition{4.6} The equality \mythetag{4.16} is called the {\it equation
of a plane passing through three given points}.
\enddefinition
\mydefinition{4.7} The equality \mythetag{4.18} is called the {\it vectorial 
form of the equation of a plane passing through three given points}.
\enddefinition
     {\bf 6. General equation\/} of a plane in the space. Let's apply the
formula \mythetagchapter{29.8}{1} from 
Chapter~\uppercase\expandafter{\romannumeral 1} in order to calculate
the scalar product $(\bold r,\bold n)$ in the equation \mythetag{4.8}. 
This yields
$$
\hskip -2em
(\bold r,\bold n)=\sum^3_{i=1}\sum^3_{j=1}r^i\,n^{\kern 0.6pt j}
\,g_{ij},
\mytag{4.19}
$$
where $g_{ij}$ are the components of the Gram matrix for the basis $\bold e_1,
\,\bold e_2,\,\bold e_3$ (see Fig\.~4.1). In Fig\.~4.1 the basis $\bold e_1,
\,\bold e_2,\,\bold e_3$ is drawn to be rectangular. However, in general case
it could be skew-angular as well. Let's introduce the notations
$$
\hskip -2em
n_i=\sum^3_{j=1}n^{\kern 0.6pt j}\,g_{ij}.
\mytag{4.20}
$$
The quantities $n^1$, $n^2$, and $n^3$ in \mythetag{4.19} and \mythetag{4.20} are
the coordinates of the normal vector $\bold n$ (see Fig\.~4.1). 
\mydefinition{4.8} The quantities $n_1$, $n_2$, and $n_3$ produced from the 
coordinates of the normal vector $\bold n$ by mens of the formula \mythetag{4.20} 
are called the {\it covariant coordinates\/} of the vector $\bold n$. 
\enddefinition
     Taking into account the notations \mythetag{4.20}, the formula \mythetag{4.19} 
is written in the following simplified form:
$$
\hskip -2em
(\bold r,\bold n)=\sum^2_{i=1}r^i\,n_i.
\mytag{4.21}
$$
Let's recall the previous notations \mythetag{4.4} and introduce new ones:
$$
\xalignat 3
&\hskip -2em
A=n_1, 
&&B=n_2, 
&&C=n_3. 
\quad
\mytag{4.22}
\endxalignat
$$
Due to \mythetag{4.4}, \mythetag{4.21}, and \mythetag{4.22} the equation
\mythetag{4.8} is written as
$$
\hskip -2em
A\,x+B\,y+C\,y-D=0.
\mytag{4.23}
$$
\mydefinition{4.9} The equation \mythetag{4.23} is called the {\it general
equation\/} of a plane in the space.
\enddefinition
\parshape 19 0cm 10cm 0cm 10cm 0cm 10cm 0cm 10cm 0cm 10cm 0cm 10cm 
0cm 10cm 5cm 5cm 5cm 5cm 5cm 5cm 5cm 5cm 5cm 5cm 5cm 5cm 5cm 5cm 
5cm 5cm 5cm 5cm 5cm 5cm 5cm 5cm 
0cm 10cm 
{\bf 7. Triple intercept equation\/} of a plane in the space. Let's 
consider a plane in the space that does not pass through the origin and 
intersects with each of the three coordinate axes. These conditions mean
that $D\neq 0$, $A\neq 0$, $B\neq 0$, and $C\neq 0$ in \mythetag{4.23}. 
\pagebreak Through $X$, $Y$, and $Z$ in Fig\.~4.2 below we denote three 
intercept points of the plane. Here are the coordinates of these three 
intercept points $X$, $Y$, and $Z$ produced by our plane:
$$
\aligned
&X=X(a,0,0),\\
&Y=Y(0,b,0),\\
&Z=Y(0,0,c).
\endaligned
\mytag{4.24}
$$
The quantities $a$, $b$, and $c$ \vadjust{\vskip 5pt\hbox 
to 0pt{\kern 0pt\includegraphics{angemeng33.eps}\hss}
\vskip -5pt}in \mythetag{4.24} are expressed through the constant 
parameters $A$, $B$, $C$, and $D$ of the equation \mythetag{4.23} 
by means of the formulas
$$
\aligned
&a=D/A,\\
&b=D/B,\\
&c=D/C.
\endaligned
\mytag{4.25}
$$
The equation of the plane \mythetag{4.23} itself can be written as
$$
\hskip -2em
\frac{x}{D/A}+\frac{y}{D/B}+\frac{z}{D/C}=1.
\mytag{4.26}
$$
If we take into account \mythetag{4.25}, the equality \mythetag{4.26} turns
to
$$
\hskip -2em
\frac{x}{a}+\frac{y}{b}+\frac{z}{c}=1.
\mytag{4.27}
$$
\mydefinition{4.10} The equality \mythetag{4.27} is called the {\it triple
intercept  equation} of a plane in the space.
\enddefinition
\head
\SectionNum{5}{154} A straight line in the space.
\endhead
\rightheadtext{\S\,5. A straight line in the space.}
       Assume that some straight line $a$ in the space $\Bbb E$ is chosen 
and fixed. In order to study various equations determining this line we 
choose some coordinate system $O,\,\bold e_1,\,\bold e_2,\,\bold e_3$ 
in the space $\Bbb E$. Then we can describe the points of the line
\pagebreak by their radius vectors relative to the origin $O$.\par
\parshape 17 0cm 10cm 0cm 10cm 5cm 5cm 5cm 5cm 5cm 5cm 5cm 5cm 
5cm 5cm 5cm 5cm 5cm 5cm 5cm 5cm 5cm 5cm 5cm 5cm 5cm 5cm 5cm 5cm 
5cm 5cm 5cm 5cm 
0cm 10cm 
     {\bf 1. Vectorial parametric equation\/} of a line in the space. 
\vadjust{\vskip 5pt\hbox to 
0pt{\kern 0pt\includegraphics{angemeng34.eps}\hss}
\vskip -5pt}Let's denote through $A$ some fixed point on the line (see 
Fig\.~5.1) and denote through $X$ an arbitrary point of this line. The 
position of the point $X$ relative to the point $A$ is marked by the vector
$\overrightarrow{AX\,\,}\!$, while the position of the point $A$ itself is
determined by its radius vector $\bold r_0=\overrightarrow{OA\,\,}\!$. 
Therefore we have 
$$
\hskip -2em
\bold r=\bold r_0+\overrightarrow{AX\,\,}\!.
\mytag{5.1}
$$
Let's choose and fix some non-zero vector $\bold a\neq\bold 0$ directed
along the line in question. The vector $\overrightarrow{AX\,\,}\!$ is
expressed through the vector $\bold a$ by means of the formula
$$
\hskip -2em
\overrightarrow{AX\,\,}\!=\bold a\cdot t.
\mytag{5.2}
$$
From the formulas \mythetag{5.1} and \mythetag{5.2} we immediately derive
$$
\hskip -2em
\bold r=\bold r_0+\bold a\cdot t.
\mytag{5.3}
$$
\mydefinition{5.1} The equality \mythetag{5.3} is called the {\it vectorial
parametric equation\/} of the line in the space. The constant vector $\bold a
\neq\bold 0$ in this equation is called a {\it directional vector} of the line, 
while $t$ is a {\it parameter}. The constant vector $\bold r_0$ is the 
{\it radius vector of an initial point}. 
\enddefinition
     Each particular value of the parameter $t$ corresponds to some definite 
point on the line. The initial point $A$ with the radius vector $\bold r_0$ is 
associated with the value $t=0$.\par
     {\bf 2. Coordinate parametric equations} of a line in the space. Let's 
determine the vectors $\bold r$, $\bold r_0$, and $\bold a$ in the vectorial 
parametric equation \mythetag{5.3} through their coordinates:
$$
\xalignat 3
&\hskip -2em
\bold r=\Vmatrix x\\ \vspace{0.5ex} y\\ \vspace{0.5ex} z
\endVmatrix, 
&&\bold r_0=\Vmatrix x_0\\ \vspace{0.5ex} y_0\\ \vspace{0.5ex} 
z_0\endVmatrix, 
&&\bold a=\Vmatrix a_x\\ \vspace{0.5ex} a_y\\ \vspace{0.5ex} 
a_z\endVmatrix.
\quad
\mytag{5.4}
\endxalignat
$$
Due to \mythetag{5.4} the equation \mythetag{5.3} is written as three 
equations:
$$
\hskip -2em
\cases x=x_0+a_x\,t,\\ y=y_0+a_y\,t,\\ z=z_0+a_z\,t.\endcases
\mytag{5.5}
$$
\mydefinition{5.2} The equalities \mythetag{5.5} are called the
{\it coordinate parametric equations\/} of a line in the space. 
The constants $a_x$, $a_y$, $a_z$ in them cannot vanish simultaneously
\enddefinition
     {\bf 3. Vectorial equation} of a line in the space. Let's apply the
vectorial multiplication by the vector $\bold a$ to both sides of the 
equation \mythetag{5.3}. As a result we get
$$
\hskip -2em
[\bold r,\bold a]=[\bold r_0,\bold a]+[\bold a,\bold a]\,t.
\mytag{5.6}
$$
Due to the item 4 of the theorem~\mythetheoremchapter{39.1}{1} in 
Chapter~\uppercase\expandafter{\romannumeral 1} the vector product
$[\bold a,\bold a]$ in \mythetag{5.6} is equal to zero. For this reason 
the equation \mythetag{5.6} actually does not contain the parameter $t$. 
This equation is usually written as follows:
$$
\hskip -2em
[\bold r-\bold r_0,\bold a]=0.
\mytag{5.7}
$$\par 
     The vector product of two constant vectors $\bold r_0$ and $\bold a$ 
is a constant vector. If we denote this vector $\bold b=[\bold r_0,\bold a]$, 
then the equation of the line \mythetag{5.7} can be written as
$$
\hskip -2em
[\bold r,\bold a]=\bold b\text{, \ where \ }\bold b\perp\bold a.
\mytag{5.8}
$$
\mydefinition{5.3} Any one of the two equalities \mythetag{5.7} and 
\mythetag{5.8} is called the {\it vectorial\footnotemark\ equation\/} of a 
line in the space. The constant vector $\bold b$ in the equation \mythetag{5.8} 
should be perpendicular to the directional vector $\bold a$. 
\enddefinition
\footnotetext{ \ The terms {\eightcyr\char '074}vectorial equation{\eightcyr
\char '076} and {\eightcyr\char '074}vectorial parametric equation{\eightcyr
\char '076} are often confused and the term {\eightcyr\char '074}vector 
equation{\eightcyr\char '076} is misapplied to \mythetag{5.3}.}
\adjustfootnotemark{-1}
     {\bf 4. Canonical equation \/} of a line in the space. Let's consider 
the case where $a_x\neq 0$, $a_y\neq 0$, and $a_z\neq 0$ in the equations 
\mythetag{5.5}. Then the parameter $t$ can be expressed through $x$, $y$, and
$z$:
$$
\xalignat 3
&\hskip -2em
t=\frac{x-x_0}{a_x}, 
&&t=\frac{y-y_0}{a_y},
&&t=\frac{z-z_0}{a_z}.
\quad 
\mytag{5.9}
\endxalignat
$$
From the equations \mythetag{5.9} one can derive the equalities
$$
\hskip -2em
\frac{x-x_0}{a_x}=\frac{y-y_0}{a_y}=\frac{z-z_0}{a_z}. 
\mytag{5.10}
$$\par
     If $a_x=0$, then instead of the first equation \mythetag{5.9} from
\mythetag{5.5} we derive $x=x_0$. Therefore instead of \mythetag{5.10} we
write
$$
\xalignat 2
&\hskip -2em
x=x_0,
&&\frac{y-y_0}{a_y}=\frac{z-z_0}{a_z}.
\mytag{5.11}
\endxalignat
$$
If $a_y=0$, then instead of the second equation \mythetag{5.9} from
\mythetag{5.5} we derive $y=y_0$. Therefore instead of \mythetag{5.10} we
write
$$
\xalignat 2
&\hskip -2em
y=y_0,
&&\frac{x-x_0}{a_x}=\frac{z-z_0}{a_z}.
\mytag{5.12}
\endxalignat
$$
If $a_z=0$, then instead of the third equation  \mythetag{5.9} from
\mythetag{5.5} we derive$z=z_0$. Therefore instead of \mythetag{5.10}
we write
$$
\xalignat 2
&\hskip -2em
z=z_0,
&&\frac{x-x_0}{a_x}=\frac{y-y_0}{a_y}.
\mytag{5.13}
\endxalignat
$$\par
     If $a_x=0$ and $a_y=0$, then instead of \mythetag{5.10} we write
$$
\xalignat 2
&\hskip -2em
x=x_0,
&&y=y_0.
\mytag{5.14}
\endxalignat
$$
If $a_x=0$ and $a_z=0$, then instead of \mythetag{5.10} we write
$$
\xalignat 2
&\hskip -2em
x=x_0,
&&z=z_0.
\mytag{5.15}
\endxalignat
$$
If $a_y=0$ and $a_z=0$, then instead of \mythetag{5.10} we write
$$
\xalignat 2
&\hskip -2em
y=y_0,
&&z=z_0.
\mytag{5.16}
\endxalignat
$$
\mydefinition{5.4} Any one of the seven pairs of equations \mythetag{5.10}, 
\mythetag{5.11}, \mythetag{5.12}, \mythetag{5.13}, \mythetag{5.14},
\mythetag{5.15}, and \mythetag{5.16} is called the {\it canonical equation\/} 
of a line in the space.
\enddefinition
     {\bf 5. The equation of a line passing through two given points\/} in the
space. Assume that two distinct points $A\neq B$ in the space are given.
We write their coordinates
$$
\xalignat 2
&\hskip -2em
A=A(x_0,y_0,z_0), 
&&B=B(x_1,y_1,z_1). 
\quad
\mytag{5.17}
\endxalignat
$$
The vector $\bold a=\overrightarrow{AB\,\,}\!$ can be used for the
directional vector of the line passing through the points $A$ and $B$
in \mythetag{5.17}. Then from \mythetag{5.17} we derive the coordinates
of $\bold a$:
$$
\hskip -2em
\bold a=\Vmatrix a_x\\ \vspace{0.5ex} a_y\\ \vspace{0.5ex} a_z
\endVmatrix
=\Vmatrix x_1-x_0\\ \vspace{0.5ex} y_1-y_0
\\ \vspace{0.5ex} z_1-z_0\endVmatrix.
\mytag{5.18}
$$
Due to \mythetag{5.18} the equations \mythetag{5.10} can be written as
$$
\hskip -2em
\frac{x-x_0}{x_1-x_0}=\frac{y-y_0}{y_1-y_0}=\frac{z-z_0}{z_1-z_0}. 
\mytag{5.19}
$$
The equations \mythetag{5.19} correspond to the case where the inequalities
$x_1\neq x_0$, $y_1\neq y_0$, and $z_1\neq z_0$ are fulfilled.\par
     If $x_1=x_0$, then instead of \mythetag{5.19} we write the equations
$$
\xalignat 2
&\hskip -2em
x=x_0=x_1,
&&\frac{y-y_0}{y_1-y_0}=\frac{z-z_0}{z_1-z_0}.
\qquad
\mytag{5.20}
\endxalignat 
$$
If $y_1=y_0$, then instead of \mythetag{5.19} we write the equations
$$
\xalignat 2
&\hskip -2em
y=y_0=y_1,
&&\frac{x-x_0}{x_1-x_0}=\frac{z-z_0}{z_1-z_0}.
\qquad
\mytag{5.21}
\endxalignat 
$$
If $z_1=z_0$, then instead of \mythetag{5.19} we write the equations
$$
\xalignat 2
&\hskip -2em
z=z_0=z_1,
&&\frac{x-x_0}{x_1-x_0}=\frac{y-y_0}{y_1-y_0}.
\qquad
\mytag{5.22}
\endxalignat 
$$\par
     If $x_1=x_0$ and $y_1=y_0$, then instead of \mythetag{5.19} we write 
$$
\xalignat 2
&\hskip -2em
x=x_0=x_1,
&&y=y_0=y_1.
\qquad
\mytag{5.23}
\endxalignat 
$$
If $x_1=x_0$ and $z_1=z_0$, then instead of \mythetag{5.19} we write 
$$
\xalignat 2
&\hskip -2em
x=x_0=x_1,
&&z=z_0=z_1.
\qquad
\mytag{5.24}
\endxalignat 
$$
If $y_1=y_0$ and $z_1=z_0$, then instead of \mythetag{5.19} we write 
$$
\xalignat 2
&\hskip -2em
y=y_0=y_1,
&&z=z_0=z_1.
\qquad
\mytag{5.25}
\endxalignat 
$$
The conditions $x_1=x_0$, $y_1=y_0$, and $z_1=z_0$ cannot be fulfilled
simultaneously since $A\neq B$.\par
\mydefinition{5.5} Any one of the seven pairs of equalities \mythetag{5.19}, 
\mythetag{5.20}, \mythetag{5.21}, \mythetag{5.22}, \mythetag{5.23}, 
\mythetag{5.24},and \mythetag{5.25} is called the equation of a line
{\it passing through two given points} $A$ and $B$ with the coordinates
\mythetag{5.17}.
\enddefinition
     {\bf 6. The equation\/} of a line in the space {\bf as the intersection 
of two planes}. In the vectorial form the equations of two intersecting planes
can be written as \mythetag{4.8}:
$$
\xalignat 2
&\hskip -2em
(\bold r,\bold n_1)=D_1,
&&(\bold r,\bold n_2)=D_2.
\mytag{5.26}
\endxalignat
$$
For the planes given by the equations \mythetag{5.26} do actually intersect
their normal vectors should be non-parallel: $\bold n_1\nparallel\bold n_2$.\par
     In the coordinate form the equations of two intersecting planes
can be written as \mythetag{4.23}:
$$
\hskip -2em
\aligned
&A_1\,x+B_1\,y+C_1\,y-D_1=0,\\
&A_2\,x+B_2\,y+C_2\,y-D_2=0.
\endaligned
\mytag{5.27}
$$
\mydefinition{5.6} Any one of the two pairs of equalities \mythetag{5.26} and
\mythetag{5.27} is called the equation of a line in the space obtained {\it as 
the intersection of two planes}.
\enddefinition
\head
\SectionNum{6}{160} Ellipse. Canonical equation of an ellipse.
\endhead
\rightheadtext{\S\,6. Ellipse.}
\mydefinition{6.1} An {\it ellipse\/} is a set of points on some plane the sum
of distances from each of which to some fixed points $F_1$ and $F_2$ of this
plane is a constant which is the same for all points of the set. The points
$F_1$ and $F_2$ are called the {\it foci\/} of the ellipse.
\enddefinition
\parshape 11 5cm 5cm 5cm 5cm 5cm 5cm 5cm 5cm 5cm 5cm 5cm 5cm 5cm 5cm 5cm 5cm 
5cm 5cm 5cm 5cm 0cm 10cm 
     Assume that an ellipse with the foci $F_1$ and $F_2$ is given. Let's draw
the line connecting the points $F_1$ and $F_2$ and choose this line for the 
$x$-axis of a coordinate system. \vadjust{\vskip 5pt\hbox to 
0pt{\kern 0pt\includegraphics{angemeng35.eps}\hss}
\vskip -5pt}Let's denote through $O$ the midpoint of the segment $[F_1F_2]$ and
choose it for the origin. We choose the second coordinate axis (the $y$-axis) 
to be perpendicular to the $x$-axis on the ellipse plane (see Fig\.~6.1). We
choose the unity scales along the axes. This means that the basis of the 
coordinate system constructed is orthonormal.\par
     Let $M=M(x,y)$ be some arbitrary point of the ellipse. According to the 
definition \mythedefinition{6.1} the sum $|MF_1|+|MF_2|$ is equal to a constant 
which does not depend on the position of $M$ on the ellipse. Let's denote this 
constant through $2\,a$ and write
$$
\hskip -2em
|MF_1|+|MF_2|=2\,a.
\mytag{6.1}
$$\par
The length of the segment $[F_1F_2]$ is also a constant. Let's denote this constant
through $2\,c$. This yields the relationships
$$
\xalignat 2
&\hskip -2em
|F_1O|=|OF_2|=c,
&&|F_1F_2|=2\,c.
\mytag{6.2}
\endxalignat
$$
From the triangle inequality $|F_1F_2|\leqslant |MF_1|+|MF_2|$ we derive the
following inequality for the constants $c$ and $a$:
$$
\hskip -2em
c\leqslant a.
\mytag{6.3}
$$\par
     The case $c=a$ in \mythetag{6.3} corresponds to a degenerate ellipse. 
In this case the triangle inequality $|F_1F_2|\leqslant |MF_1|+|MF_2|$ turns 
to the equality $|F_1F_2|=|MF_1|+|MF_2|$, while the $MF_1F_2$ itself collapses 
to the segment $[F_1F_2]$. Since $M\in [F_1F_2]$, a degenerate ellipse with the 
foci $F_1$ and $F_2$ is the segment $[F_1F_2]$. The case of a degenerate ellipse 
is usually excluded by setting
$$
\hskip -2em
0\leqslant c<a.
\mytag{6.4}
$$\par
     The formulas \mythetag{6.2} determine the foci $F_1$ and $F_2$ of the 
ellipse in the coordinate system we have chosen:
$$
\xalignat 2
&\hskip -2em
F_1=F_1(-c,0),
&&F_2=F_2(c,0).
\mytag{6.5}
\endxalignat
$$
Having determined the coordinates of the points $F_1$ and $F_2$ and knowing the
coordinates of the point $M=M(x,y)$, \pagebreak we derive the following 
relationships for the segments $[MF_1]$ and $[MF_2]$:
$$
\hskip -2em
\aligned
&|MF_1|=\sqrt{y^2+(x+c)^2},\\
&|MF_2|=\sqrt{y^2+(x-c)^2}.
\endaligned
\mytag{6.6}
$$\par
     {\bf Derivation of the canonical equation of an ellipse}. Let's substitute
\mythetag{6.6} into the equality \mythetag{6.1} defining the ellipse:
$$
	\hskip -2em
\sqrt{y^2+(x+c)^2}+\sqrt{y^2+(x-c)^2}=2\,a.
\mytag{6.7}
$$
Then we move one of the square roots to the right hand side of the formula
\mythetag{6.7}. As a result we derive
$$
\hskip -2em
\sqrt{y^2+(x+c)^2}=2\,a-\sqrt{y^2+(x-c)^2}.
\mytag{6.8}
$$
Let's square both sides of the equality \mythetag{6.8}:
$$
\hskip -2em
\gathered
y^2+(x+c)^2=4\,a^2\,-\\
-\,4\,a\,\sqrt{y^2+(x-c)^2}+y^2+(x-c)^2.
\endgathered
\mytag{6.9}
$$
Upon expanding brackets and collecting similar terms the equality
\mythetag{6.9} can be written in the following form:
$$
\hskip -2em
4\,a\,\sqrt{y^2+(x-c)^2}=4\,a^2-4\,x\,c.
\mytag{6.10}
$$
Let's cancel the factor four in \mythetag{6.10} and square both sides 
of this equality. As a result we get the formula
$$
\hskip -2em
a^2\,(y^2+(x-c)^2)=a^4-2\,a^2\,x\,c+x^2\,c^2.
\mytag{6.11}
$$
Upon expanding brackets and recollecting similar terms the equality
\mythetag{6.11} can be written in the following form:
$$
\hskip -2em
x^2\,(a^2-c^2)+y^2\,a^2=a^2\,(a^2-c^2).
\mytag{6.12}
$$\par
     Both sides of the equality \mythetag{6.12} contain the quantity 
$a^2-c^2$. Due to the inequality \mythetag{6.4} this quantity is positive. 
For this reason it can be written as the square of some number $b>0$:
$$
\hskip -2em
a^2-c^2=b^2.
\mytag{6.13}
$$
Due to \mythetag{6.13} the equality \mythetag{6.12} can be written as 
$$
\hskip -2em
x^2\,b^2+y^2\,a^2=a^2\,b^2.
\mytag{6.14}
$$
Since $b>0$ and $a>0$ (see the inequalities \mythetag{6.4}), the equality
\mythetag{6.14} transforms to the following one:
$$
\hskip -2em
\frac{x^2}{a^2}+\frac{y^2}{b^2}=1.
\mytag{6.15}
$$
\mydefinition{6.2} The equality \mythetag{6.15} is called the {\it canonical\/} 
equation of an ellipse. 
\enddefinition
\mytheorem{6.1} For each point $M=M(x,y)$ of the ellipse determined by the
initial equation \mythetag{6.7} its coordinates obey the canonical 
equation \mythetag{6.15}.  
\endproclaim
     The proof of the theorem~\mythetheorem{6.1} consists in the above 
calculations that lead to the canonical equation of an ellipse \mythetag{6.15}. 
This canonical equation yields the following important inequalities:
$$
\xalignat 2
&\hskip -2em
|x|\leqslant a,
&&|y|\leqslant b.
\mytag{6.16}
\endxalignat
$$
\mytheorem{6.2} The canonical equation of an ellipse \mythetag{6.15} is 
equivalent to the initial equation \mythetag{6.7}.
\endproclaim
\demo{Proof} In order to prove the theorem~\mythetheorem{6.2} we transform
the expressions \mythetag{6.6} relying on \mythetag{6.15}. From \mythetag{6.15} 
we derive
$$
\hskip -2em
y^2=b^2-\frac{b^2}{a^2}\,x^2.
\mytag{6.17}
$$
Substituting \mythetag{6.17} into the first formula \mythetag{6.6}, we get
$$
\gathered
|MF_1|=\sqrt{b^2-\frac{b^2}{a^2}\,x^2+x^2+2\,x\,c+c^2}=\\
=\sqrt{\frac{a^2-b^2}{a^2}\,x^2+2\,x\,c+(b^2+c^2)}.
\endgathered
\mytag{6.18}
$$
Now we take into account the relationship \mythetag{6.13} and write the 
equality \mythetag{6.18} in the following form:
$$
|MF_1|=\sqrt{a^2+2\,x\,c+\frac{c^2}{a^2}\,x^2}
=\sqrt{\Bigl(\frac{a^2+c\,x}{a}\Bigr)^2}.
\quad
\mytag{6.19}
$$
Upon calculating the square root the formula \mythetag{6.19} yields
$$
\hskip -2em
|MF_1|=\frac{|a^2+c\,x|}{a}.
\mytag{6.20}
$$
From the inequalities \mythetag{6.4} and \mythetag{6.16} we derive
$a^2+c\,x>0$. Therefore the modulus signs in \mythetag{6.20} can be omitted:
$$
\hskip -2em
|MF_1|=\frac{a^2+c\,x}{a}=a+\frac{c\,x}{a}.
\mytag{6.21}
$$
In the case of the second formula \mythetag{6.6} the considerations similar
to the above ones yield the following result:
$$
\hskip -2em
|MF_2|=\frac{a^2-c\,x}{a}=a-\frac{c\,x}{a}.
\mytag{6.22}
$$
Now it is easy to see that the equation \mythetag{6.7} written as 
$|MF_1|+|MF_2|=2\,a$ due to \mythetag{6.6} is immediate from \mythetag{6.21} 
and \mythetag{6.22}. The theorem~\mythetheorem{6.2} is proved. 
\qed\enddemo
     Let's consider again the inequalities \mythetag{6.16}. The 
coordinates of any point $M$ of the ellipse obey these inequalities. 
The first of the inequalities \mythetag{6.16} turns to an equality 
if $M$ coincides with $A$ or if $M$ coincides with $C$ (see
Fig\.~6.1). The second inequality \mythetag{6.16} turns to an equality
if $M$ coincides with $B$ or if $M$ coincides $D$.\par
\mydefinition{6.3} The points $A$, $B$, $C$, and $D$ in Fig\.~6.1 are
called the {\it vertices\/} of an ellipse. The segments $[AC]$ and $[BD]$ 
are called the {\it axes\/} of an ellipse, while the segments $[OA]$, $[OB]$,
$[OC]$, and $[OD]$ are called its {\it semiaxes}.
\enddefinition
     The constants $a$, $b$, and $c$ obey the relationship \mythetag{6.13}. 
From this relationship and from the inequalities \mythetag{6.4} we derive 
$$
\hskip -2em
0<b\leqslant a.
\mytag{6.23}
$$
\mydefinition{6.4} Due to the inequalities \mythetag{6.23} the semiaxis 
$[OA]$ in Fig\.~6.1 is called the {\it major semiaxis\/} of the ellipse,
while the semiaxis $[OB]$ is called its {\it minor semiaxis}. 
\enddefinition
\mydefinition{6.5} A coordinate system $O,\,\bold e_1,\,\bold e_2$
with an orthonormal basis $\bold e_1,\,\bold e_2$ where an ellipse
is given by its canonical equation \mythetag{6.15} and where the 
inequalities \mythetag{6.23} are fulfilled is called a 
{\it canonical coordinate system\/} of this ellipse. 
\enddefinition
\head
\SectionNum{7}{165} The eccentricity and directrices of an ellipse.
The property of directrices.
\endhead
\rightheadtext{\S\,7. The eccentricity and directrices\dots}
     The shape and sizes of an ellipse are determined by two constants $a$
and $b$ in its canonical equation \mythetag{6.15}. Due to the relationship
\mythetag{6.13} the constant $b$ can be expressed through the constant
$c$. Multiplying both constants $a$ and $c$ by the same number, we change
the sizes of an ellipse, but do not change its shape. The ratio of these two
constants 
$$
\hskip -2em
\varepsilon=\frac{c}{a}.
\mytag{7.1}
$$
is responsible for the shape of an ellipse. 
\mydefinition{7.1} The quantity $\varepsilon$ defined by the relationship
\mythetag{7.1}, where $a$ is the major semiaxis and $c$ is the half of the 
interfocal distance, is called the {\it eccentricity\/} of an ellipse. 
\enddefinition
     The eccentricity \mythetag{7.1} is used in order to define 
one more numeric parameter of an ellipse. It is usually denoted through $d$:
$$
\hskip -2em
d=\frac{a}{\varepsilon}=\frac{a^2}{c}.
\mytag{7.2}
$$
\mydefinition{7.2} \parshape 4 0cm 10cm 0cm 10cm 0cm 10cm 5cm 5cm
On the plane of an ellipse there are two lines parallel
to its minor axis and placed at the distance $d$ given by the formula 
\mythetag{7.2} from its center. These lines are called {\it directrices\/} 
of an ellipse. 
\enddefinition
\parshape 11 5cm 5cm 5cm 5cm 5cm 5cm 5cm 5cm 5cm 5cm 
5cm 5cm 5cm 5cm 5cm 5cm 5cm 5cm 5cm 5cm 0cm 10cm
     Each ellipse has two foci and two directrices. Each directrix has the
corresponding focus of it.
\vadjust{\vskip 5pt\hbox to 
0pt{\kern 0pt\includegraphics{angemeng36.eps}\hss}
\vskip -5pt}This is that of two foci which is more close to the directrix in
question. Let $M=M(x,y)$ be some arbitrary point of an ellipse. Let's connect
this point with the left focus of the ellipse $F_1$ and drop the perpendicular 
from it to the left directrix of the ellipse. Let's denote through $H_1$ the
base of such a perpendicular and calculate its length:
$$
\hskip -2em
|MH_1|=|x-(-d)|=|d+x|.
\mytag{7.3}
$$
Taking into account \mythetag{7.2}, the formula \mythetag{7.3} can be brought
to
$$
\hskip -2em
|MH_1|=\Bigl|\frac{a^2}{c}+x\Bigr|=\frac{|a^2+c\,x|}{c}.
\mytag{7.4}
$$
The length of the segment $MF_1$ was already calculated above. Initially it 
was given by one of the formulas \mythetag{6.6}, but later the more simple 
expression \mythetag{6.20} was derived for it: 
$$
\hskip -2em
|MF_1|=\frac{|a^2+c\,x|}{a}.
\mytag{7.5}
$$
If we divide \mythetag{7.5} by \mythetag{7.4}, we obtain the following 
relationship: 
$$
\hskip -2em
\frac{|MF_1|}{|MH_1|}=\frac{c}{a}=\varepsilon.
\mytag{7.6}
$$
The point $M$ can change its position on the ellipse. Then the numerator and the 
denominator of the fraction \mythetag{7.6} are changed, but its value remains
unchanged. This fact is known as the property of directrices. 
\mytheorem{7.1} The ratio of the distance from some arbitrary point $M$ of an 
ellipse to its focus and the distance from this point to the corresponding
directrix is a constant which is equal to the eccentricity of the ellipse.
\endproclaim
\head
\SectionNum{8}{167} The equation of a tangent line to an 
ellipse.
\endhead
\rightheadtext{\S\,8. The equation of a tangent line \dots}
\parshape 4 0cm 10cm 0cm 10cm 0cm 10cm 5cm 5cm 
     Let's consider an ellipse given by its canonical equation 
\mythetag{6.15} in its canonical coordinate system (see
Definition~\mythedefinition{6.5}). \vadjust{\vskip 5pt\hbox to 
0pt{\kern 0pt\includegraphics{angemeng37.eps}\hss}
\vskip -5pt}Let's draw a tangent line to this ellipse and denote through
$M=M(x_0,y_0)$ its tangency point (see Fig\.~8.1). Our goal is to write
the equation of a tangent line to an ellipse. 
\par
\parshape 10 5cm 5cm 5cm 5cm 5cm 5cm 5cm 5cm 5cm 5cm 5cm 5cm 5cm 5cm 5cm 5cm 
5cm 5cm 0cm 10cm 
     An ellipse is a curve composed by two halves --- the upper half and the
lower half. Any one of these two halves can be treated as a graph of a 
function
$$
\hskip -2em
y=f(x)
\mytag{8.1}
$$
with the domain $(-a,a)$. The equation of a tangent line to the graph of a
function is given by the following well-known formula (see \mycite{9}):
$$
\hskip -2em
y=y_0+f'(x_0)\,(x-x_0). 
\mytag{8.2}
$$\par
     In order to apply the formula \mythetag{8.2} to an ellipse we need to 
calculate the derivative of the function \mythetag{8.1}. Let's substitute 
the expression \mythetag{8.1} for $y$ into the formula \mythetag{6.15}:
$$
\hskip -2em
\frac{x^2}{a^2}+\frac{(f(x))^2}{b^2}=1.
\mytag{8.3}
$$
The equality \mythetag{8.3} is fulfilled identically in $x$. Let's 
differentiate the equality \mythetag{8.3} with respect to $x$. This yields
$$
\hskip -2em
\frac{2\,x}{a^2}+\frac{2\,f(x)\,f'(x)}{b^2}=0.
\mytag{8.4}
$$
Let's apply the formula \mythetag{8.4} for to calculate the derivative $f'(x)$:
$$
\hskip -2em
f'(x)=-\frac{b^2\,x}{a^2\,f(x)}.
\mytag{8.5}
$$
In order to substitute \mythetag{8.5} into the equation \mythetag{8.2} we 
change $x$ for $x_0$ and $f(x)$ for $f(x_0)=y_0$. As a result we get
$$
\hskip -2em
f'(x_0)=-\frac{b^2\,x_0}{a^2\,y_0}.
\mytag{8.6}
$$\par
     Let's substitute \mythetag{8.6} into the equation of the tangent line 
\mythetag{8.2}. This yields the following relationship
$$
\hskip -2em
y=y_0-\frac{b^2\,x_0}{a^2\,y_0}\,(x-x_0). 
\mytag{8.7}
$$
Eliminating the denominator, we write the equality \mythetag{8.7} as 
$$
\hskip -2em
a^2\,y\,y_0+b^2\,x\,x_0=a^2\,y_0^2+b^2\,x_0^2.
\mytag{8.8}
$$
Now let's divide both sides of the equality \mythetag{8.8} by $a^2\,b^2$:
$$
\hskip -2em
\frac{x\,x_0}{a^2}+\frac{y\,y_0}{b^2}=\frac{x_0^2}{a^2}
+\frac{y_0^2}{b^2}.
\mytag{8.9}
$$\par
     Note that the point $M=M(x_0,y_9)$ is on the ellipse. Therefore its 
coordinates satisfy the equation \mythetag{6.15}:
$$
\hskip -2em
\frac{x_0^2}{a^2}+\frac{y_0^2}{b^2}=1.
\mytag{8.10}
$$
Taking into account \mythetag{8.10}, we can transform \mythetag{8.9} to
$$
\hskip -2em
\frac{x\,x_0}{a^2}+\frac{y\,y_0}{b^2}=1.
\mytag{8.11}
$$
This is the required equation of a tangent line to an ellipse. 
\mytheorem{8.1} For an ellipse determined by its canonical equation 
\mythetag{6.15} its tangent line that touches this ellipse at the 
point $M=M(x_0,y_0)$ is given by the equation \mythetag{8.11}.
\endproclaim
     The equation \mythetag{8.11} is a particular case of the equation 
\mythetag{3.22} where the constants $A$, $B$, and $D$ are given by the
formulas
$$
\xalignat 3
&\hskip -2em
A=\frac{x_0}{a^2},
&&B=\frac{y_0}{b^2},
&&D=1.
\quad
\mytag{8.12}
\endxalignat
$$
According to the definition~\mythedefinition{3.6} and the formulas 
\mythetag{3.21} the constants $A$ and $B$ in \mythetag{8.12} are
the covariant components of the normal vector for the tangent line 
to an ellipse. The tangent line equation \mythetag{8.11} is written
in a canonical coordinate system of an ellipse. The basis of such a
coordinate system is orthonormal. Therefore the formula \mythetag{3.19} 
and the formula \mythetagchapter{32.4}{1} from 
Chapter~\uppercase\expandafter{\romannumeral 1} yield the following
relationships:
$$
\xalignat 2
&\hskip -2em
A=n_1=n^1, 
&&B=n_2=n^2. 
\mytag{8.13}
\endxalignat
$$
\mytheorem{8.2} The quantities $A$ and $B$ in \mythetag{8.12} are the
coordinates of the normal vector $\bold n$ for the tangent line
to an ellipse which is given by the equation \mythetag{8.11}.
\endproclaim
\noindent The relationships \mythetag{8.13} prove the
theorem~\mythetheorem{8.2}.
\head
\SectionNum{9}{170} Focal property of an ellipse.
\endhead
\rightheadtext{\S\,9. Focal property of an ellipse.}
     The term {\it focus\/} is well known in optics. It means a point
where light rays converge upon refracting in lenses or upon reflecting 
in curved mirrors. In the case of an ellipse let's assume 
that it is manufactured of a thin strip of some flexible material and 
assume that its inner surface is covered with a light reflecting layer. 
For such a materialized ellipse one can formulate the \nolinebreak
following fo\-cal \nolinebreak property. 
\mytheorem{9.1} A light ray emitted from one of the foci of an ellipse 
upon reflecting on its inner surface passes through the other focus of
this ellipse. 
\endproclaim
\mytheorem{9.2} The perpendicular to a tangent line of an ellipse
drawn at the tangency point is a bisector in the triangle composed
by the tangency point and two foci of the ellipse.
\endproclaim
\parshape 1 5cm 5cm 
     The theorem~\mythetheorem{9.2} is a geometric version 
of the theorem~\mythetheorem{9.1}. 
\vadjust{\vskip 5pt\hbox to 
0pt{\kern 0pt\includegraphics{angemeng38.eps}\hss}
\vskip -5pt}These theorems are equivalent due to the reflection law saying 
that the angle of reflection is equal to the angle of incidence.
\demo{Proof of the theorem~\mythetheorem{9.2}}\parshape 5 4.9cm 5.1cm
5cm 5cm 5cm 5cm 5cm 5cm 0cm 10cm \linebreak 
Let's choose an arbitrary point $M=M(x_0,y_0)$ on the ellipse and 
draw the tangent line to the ellipse through this point as shown in 
Fig\.~9.1. Then we draw the perpendicular $[MN]$ to the tangent line 
through the point $M$. The segment $[MN]$ is directed along the normal
vector of the tangent line. In order to prove that this segment is 
a bisector of the triangle $F_1MF_2$ it is sufficient to prove the
equality
$$
\hskip -2em
\cos(\widehat{F_1MN})=\cos(\widehat{F_2MN}). 
\mytag{9.1}
$$
The cosine equality \mythetag{9.1} \pagebreak is equivalent to the 
following  equality for the scalar products of vectors:
$$
\hskip -2em
\frac{(\overrightarrow{MF_1\,}\!,\bold n)}{|MF_1|}
=\frac{(\overrightarrow{MF_2\,}\!,\bold n)}{|MF_2|}.
\mytag{9.2}
$$\par
     The coordinates of the points $F_1$ and $F_2$ in a canonical 
coordinate system of the ellipse are known (see formulas \mythetag{6.5}).
The coordinates of the point $M=M(x_0,y_0)$ are also known. Therefore we 
can find the coordinates of the vectors $\overrightarrow{MF_1\,}\!$ and
$\overrightarrow{MF_2\,}\!$ used in the above formula \mythetag{9.2}:
$$
\xalignat 2
&\hskip -2em
\overrightarrow{MF_1\,}\!
=\Vmatrix -c-x_0\\ -y_0\endVmatrix,
&&\overrightarrow{MF_2\,}\!
=\Vmatrix c-x_0\\ -y_0\endVmatrix.
\mytag{9.3}
\endxalignat
$$
The coordinates of the normal vector $\bold n$ of the tangent line in 
Fig\.~9.1 are given by the formulas \mythetag{8.12} and \mythetag{8.13}:
$$
\hskip -2em
\bold n=\Vmatrix
\dfrac{x_0}{a^2}\\
\vspace{2ex}
\dfrac{y_0}{b^2}
\endVmatrix.
\mytag{9.4}
$$
Using \mythetag{9.3} and \mythetag{9.4}, we apply the formula 
\mythetagchapter{33.3}{1} from Chapter~\uppercase\expandafter{\romannumeral 1} 
for calculating the scalar products in \mythetag{9.2}:
$$
\hskip -2em
\aligned
&(\overrightarrow{MF_1\,}\!,\bold n)
=\frac{-c\,x_0-x_0^2}{a^2}-\frac{y_0^2}{b^2}
=\frac{-c\,x_0}{a^2}-\frac{x_0^2}{a^2}-\frac{y_0^2}{b^2},\\
\vspace{1ex}
&(\overrightarrow{MF_2\,}\!,\bold n)
=\frac{c\,x_0-x_0^2}{a^2}-\frac{y_0^2}{b^2}
=\frac{c\,x_0}{a^2}-\frac{x_0^2}{a^2}-\frac{y_0^2}{b^2}.
\endaligned
\mytag{9.5}
$$\par
     The point $M=M(x_0,y_0)$ lies on the ellipse. Therefore its coordinates
should satisfy the equation of the ellipse \mythetag{6.15}:
$$
\hskip -2em
\frac{x_0^2}{a^2}+\frac{y_0^2}{b^2}=1.
\mytag{9.6}
$$
Due to \mythetag{9.6} the formulas \mythetag{9.5} simplify and take the form
$$
\hskip -2em
\aligned
&(\overrightarrow{MF_1\,}\!,\bold n)
=\frac{-c\,x_0}{a^2}-1=-\frac{a^2+c\,x_0}{a^2},\\
\vspace{1ex}
&(\overrightarrow{MF_2\,}\!,\bold n)
=\frac{c\,x_0}{a^2}-1
=-\frac{a^2-c\,x_0}{a^2}.
\endaligned
\mytag{9.7}
$$\par
     In order to calculate the denominators in the formula \mythetag{9.2} we use
the formulas \mythetag{6.21} and \mythetag{6.22}. In this case, when applied to
the point $M=M(x_0,y_0)$, they yield
$$
\xalignat 2
&\hskip -2em
|MF_1|=\frac{a^2+c\,x_0}{a},
&&|MF_2|=\frac{a^2-c\,x_0}{a}.
\quad
\mytag{9.8}
\endxalignat
$$
From the formulas \mythetag{9.7} and \mythetag{9.8} we easily derive 
the equalities
$$
\xalignat 2
&\hskip -2em
\frac{(\overrightarrow{MF_1\,}\!,\bold n)}{|MF_1|}
=-\frac{1}{a},
&&\frac{(\overrightarrow{MF_2\,}\!,\bold n)}{|MF_2|}
=-\frac{1}{a}
\endxalignat
$$
which proves the equality \mythetag{9.2}. As a result the 
theorem~\mythetheorem{9.2} and the theorem~\mythetheorem{9.1} equivalent 
to it both are proved.\qed\enddemo
\head
\SectionNum{10}{172} Hyperbola. Canonical equation of 
a hyperbola.
\endhead
\rightheadtext{\S\,10. Hyperbola.}
\mydefinition{10.1} A {\it hyperbola\/} is a set of points on some plane the 
modulus of the difference of distances from each of which to some fixed points 
$F_1$ and $F_2$ of this plane is a constant which is the same for all points 
of the set. The points $F_1$ and $F_2$ are called the {\it foci\/} of the 
hyperbola.
\enddefinition
\parshape 7 0cm 10cm 0cm 10cm 0cm 10cm 0cm 10cm 0cm 10cm 0cm 10cm 
5cm 5cm 
     Assume that a hyperbola with the foci $F_1$ and $F_2$ is given. Let's 
draw the line connecting the points $F_1$ and $F_2$ and choose this line for 
the $x$-axis of a coordinate system. Let's denote through $O$ the midpoint of 
the segment $[F_1F_2]$ and choose it for the origin. We choose the second 
coordinate axis (the $y$-axis) to be perpendicular to the $x$-axis on the 
hyperbola plane (see Fig\.~10.1). We choose the unity scales along the axes. 
This means that \vadjust{\vskip 5pt\hbox to 
0pt{\kern 0pt\includegraphics{angemeng39.eps}\hss}
\vskip -5pt}the basis of the coordinate system constructed is orthonormal.
\par
\parshape 8 5cm 5cm 5cm 5cm 5cm 5cm 5cm 5cm 5cm 5cm 5cm 5cm 5cm 5cm 0cm 10cm 
     Let $M=M(x,y)$ be some arbitrary point of the hyperbola in Fig\.~10.1. 
According to the definition of a hyperbola \mythedefinition{10.1} the modulus
of the difference $\bigl||MF_1|-|MF_2|\bigr|$ is a constant which does not
depend on a location of the point $M$ on the hyperbola. We denote this constant
through $2\,a$ and write the formula 
$$
\hskip -2em
\bigl||MF_1|-|MF_2|\bigr|=2\,a.
\mytag{10.1}
$$
According to \mythetag{10.1} the points of the hyperbola are subdivided into two 
subsets. For the points of one of these two subsets the condition \mythetag{10.1} 
is written as 
$$
\hskip -2em
|MF_1|-|MF_2|=2\,a.
\mytag{10.2}
$$
These points constitute the right branch of the hyperbola. For the 
points of the other subset the condition \mythetag{10.1} is written 
as
$$
\hskip -2em
|MF_2|-|MF_1|=2\,a.
\mytag{10.3}
$$
These points constitute the left branch of the hyperbola.\par
     The length of the segment $[F_1F_2]$ is also a constant. Let's 
denote this constant through $2\,c$. This yields 
$$
\xalignat 2
&\hskip -2em
|F_1O|=|OF_2|=c,
&&|F_1F_2|=2\,c.
\mytag{10.4}
\endxalignat
$$
From the triangle inequality $|F_1F_2|\geqslant \bigl||MF_1|-|MF_2|\bigr|$, 
\pagebreak from \mythetag{10.2}, from \mythetag{10.3}, and from \mythetag{10.4} 
we derive
$$
\hskip -2em
c\geqslant a.
\mytag{10.5}
$$\par
     The case $c=a$ in \mythetag{10.5} corresponds to a degenerate hyperbola. 
In this case the inequality $|F_1F_2|\geqslant \bigl||MF_1|-|MF_2|\bigr|$ turns
to the equality $|F_1F_2|=\bigl||MF_1|-|MF_2|\bigr|$ which distinguishes two 
subcases. In the first subcase the equality $|F_1F_2|=\bigl||MF_1|-|MF_2|\bigr|$ 
is written as
$$
|F_1F_2|=|MF_1|-|MF_2|.
\mytag{10.6}
$$
The equality \mythetag{10.6} means that the triangle $F_1MF_2$ collapses into
the segment $[F_1M]$, i\.\,e\. the point $M$ lies on the ray going along
the $x$-axis to the right from the point $F_2$ (see Fig\.~10.1).\par
     In the second subcase of the case $c=a$ the equality $|F_1F_2|=\break
=\bigl||MF_1|-|MF_2|\bigr|$ is written as follows
$$
|F_1F_2|=|MF_2|-|MF_1|.
\mytag{10.7}
$$
The equality \mythetag{10.7} means that the triangle $F_1MF_2$ collapses 
into the segment $[F_1M]$, i\.\,e\. the point $M$ lies on the ray going
along the $x$-axis to the left from the point $F_1$ (see Fig\.~10.1).
\par
     As we see considering the above two subcases, if $c=a$ the degenerate
hyperbola is the union of two non-intersecting rays lying on the $x$-axis 
and being opposite to each other.\par
     Another form of a degenerate hyperbola arises if $a=0$. In this case 
the branches of the hyperbola straighten and, gluing with each other, lie
on the $y$-axis.\par
     Both cases of degenerate hyperbolas are usually excluded from the 
consideration by means of the following inequalities:
$$
\hskip -2em
c>a>0.
\mytag{10.8}
$$\par
     The formulas \mythetag{10.4} determine the coordinates \pagebreak
of the foci $F_1$ and $F_2$ of our hyperbola in the chosen coordinate 
system:
$$
\xalignat 2
&\hskip -2em
F_1=F_1(-c,0),
&&F_2=F_2(c,0).
\mytag{10.9}
\endxalignat
$$
Knowing the coordinates of the points $F_1$ and $F_2$ and knowing the 
coordinates of the point $M=M(x,y)$, we write the formulas 
$$
\hskip -2em
\aligned
&|MF_1|=\sqrt{y^2+(x+c)^2},\\
&|MF_2|=\sqrt{y^2+(x-c)^2}.
\endaligned
\mytag{10.10}
$$\par
     {\bf Derivation of the canonical equation of a hyperbola}. As we 
have seen above, the equality \mythetag{10.1} defining a hyperbola breaks
into two equalities \mythetag{10.2} and \mythetag{10.3} corresponding to
the right and left branches of the hyperbola. Let's unite them back into 
a single equality of the form
$$
\hskip -2em
|MF_1|-|MF_2|=\pm\,2\,a.
\mytag{10.11}
$$     
Let's substitute the formulas \mythetag{10.10} into the equality 
\mythetag{10.11}:
$$
\sqrt{y^2+(x+c)^2}-\sqrt{y^2+(x-c)^2}=\pm\,2\,a.
\mytag{10.12}
$$
Then we move one of the square roots to the right hand side of the formula
\mythetag{10.12}. As a result we derive
$$
\sqrt{y^2+(x+c)^2}=\pm\,2\,a+\sqrt{y^2+(x-c)^2}.
\mytag{10.13}
$$
Squaring both sides of the equality \mythetag{10.13}, we get
$$
\hskip -2em
\gathered
y^2+(x+c)^2=4\,a^2\,\pm\\
\pm\,4\,a\,\sqrt{y^2+(x-c)^2}+y^2+(x-c)^2.
\endgathered
\mytag{10.14}
$$
Upon expanding brackets and collecting similar terms the equality
\mythetag{10.14} can be written in the following form:
$$
\hskip -2em
\mp 4\,a\,\sqrt{y^2+(x-c)^2}=4\,a^2-4\,x\,c.
\mytag{10.15}
$$
Let's cancel the factor four in \mythetag{10.15} and square both sides 
of this equality. As a result we get the formula
$$
\hskip -2em
a^2\,(y^2+(x-c)^2)=a^4-2\,a^2\,x\,c+x^2\,c^2.
\mytag{10.16}
$$
Upon expanding brackets and recollecting similar terms the equality
\mythetag{10.16} can be written in the following form:
$$
\hskip -2em
x^2\,(a^2-c^2)+y^2\,a^2=a^2\,(a^2-c^2).
\mytag{10.17}
$$\par
     An attentive reader can note that the above calculations are almost
literally the same as the corresponding calculations for the case of an 
ellipse. As for the resulting formula \mythetag{10.17}, it coincides exactly
with the formula \mythetag{6.12}. But, nevertheless, there is a difference. 
It consists in the inequalities \mythetag{10.8}, which are different from the
inequalities \mythetag{6.4} for an ellipse.\par
     Both sides of the equality \mythetag{10.17} comprises the quantity
$a^2-c^2$. Due to the inequalities \mythetag{10.8} this quantity is negative. 
For this reason the quantity $a^2-c^2$ can be written as the square of some 
positive quantity $b>0$ taken with the minus sign:
$$
\hskip -2em
a^2-c^2=-b^2.
\mytag{10.18}
$$
Due to \mythetag{10.18} the equality \mythetag{10.17} can be written as 
$$
\hskip -2em
-x^2\,b^2+y^2\,a^2=-a^2\,b^2.
\mytag{10.19}
$$
Since $b>0$ and $a>0$ (see inequalities \mythetag{10.8}), the above
equality \mythetag{10.19} transforms to the following one:
$$
\hskip -2em
\frac{x^2}{a^2}-\frac{y^2}{b^2}=1.
\mytag{10.20}
$$
\mydefinition{10.2} The equality \mythetag{10.20} is called the 
{\it canonical equation\/} of a hyperbola.
\enddefinition
\mytheorem{10.1} For each point $M(x,y)$ of the hyperbola determined
by the initial equation \mythetag{10.12} its coordinates satisfy the
canonical equation \mythetag{10.20}.  
\endproclaim
     The above derivation of the canonical equation \mythetag{10.20}
of a hyperbola proves the theorem~\mythetheorem{10.1}. The canonical	
equation leads \mythetag{10.20} to the following important inequality:
$$
\hskip -2em
|x|\geqslant a.
\mytag{10.21}
$$
\mytheorem{10.2} The canonical equation of a hyperbola \mythetag{10.20}
is equivalent to the initial equation \mythetag{10.12}.
\endproclaim
\demo{Proof} The proof of the theorem~\mythetheorem{10.2} is analogous
to the proof of the theorem~\mythetheorem{6.2}. In order to prove 
the theorem~\mythetheorem{10.2} we calculate the expressions \mythetag{10.10} 
relying on the equation \mythetag{10.20}. The equation \mythetag{10.20} 
itself can be written as follows:
$$
\hskip -2em
y^2=\frac{b^2}{a^2}\,x^2-b^2.
\mytag{10.22}
$$
Substituting \mythetag{10.22} into the first formula \mythetag{10.10}, we
get
$$
\gathered
|MF_1|=\sqrt{\frac{b^2}{a^2}\,x^2-b^2+x^2+2\,x\,c+c^2}=\\
=\sqrt{\frac{a^2+b^2}{a^2}\,x^2+2\,x\,c+(c^2-b^2)}.
\endgathered
\mytag{10.23}
$$
Now we take into account the relationship \mythetag{10.18} and write 
the equality \mythetag{10.23} in the following form:
$$
|MF_1|=\sqrt{a^2+2\,x\,c+\frac{c^2}{a^2}\,x^2}
=\sqrt{\Bigl(\frac{a^2+c\,x}{a}\Bigr)^2}.
\quad
\mytag{10.24}
$$
Upon calculating the square root the formula \mythetag{10.24} yields
$$
\hskip -2em
|MF_1|=\frac{|a^2+c\,x|}{a}.
\mytag{10.25}
$$
From the inequalities \mythetag{10.8} and \mythetag{10.21} we derive
$|c\,x|>a^2$. Therefore the formula \mythetag{10.25} can be written as
follows:
$$
|MF_1|=\frac{c\,|x|+\sign(x)\,a^2}{a}=\frac{c\,|x|}{a}
+\sign(x)\,a.
\mytag{10.26}
$$
In the case of the second formula \mythetag{10.10} the considerations, 
analogous to the above ones, yield the following result:
$$
|MF_2|=\frac{c\,|x|-\sign(x)\,a^2}{a}=\frac{c\,|x|}{a}
-\sign(x)\,a.
\mytag{10.27}
$$
Let's subtract the equality \mythetag{10.27} from the equality 
\mythetag{10.26}. Then we get the following relationships:
$$
\hskip -2em
|MF_1|-|MF_2|=2\,\sign(x)\,a=\pm\,2\,a.
\mytag{10.28}
$$
The plus sign in \mythetag{10.28} corresponds to the case $x>0$, which 
corresponds to the right branch of the hyperbola in Fig\.~10.1. The minus
sign corresponds  to the left branch of the hyperbola. Due to what was 
said the equality \mythetag{10.28} is equivalent to the equality 
\mythetag{10.11}, which in turn is equivalent to the equality 
\mythetag{10.12}. The theorem~\mythetheorem{10.2} is proved. 
\qed\enddemo
     Let's consider again the inequality \mythetag{10.21} which should 
be fulfilled for the $x$-coordinate of any point $M$ on the hyperbola. 
The inequality \mythetag{10.21} turns to an equality if $M$ coincides
with $A$ or if $M$ coincides with $C$ (see Fig\.~10.1).\par
\mydefinition{10.3} The points $A$ and $C$ in Fig\.~10.1 are called the
{\it vertices\/} of the hyperbola. The segment $[AC]$ is called the 
{\it transverse axis\/} or the {\it real axis\/} of the hyperbola, while 
the segments $[OA]$ and $[OC]$ are called its {\it real semiaxes}.
\enddefinition
     The constant $a$ in the equation of the hyperbola \mythetag{10.20} 
is the length of the segment $[OA]$ in Fig\.~10.1 (the length of the
real semiaxis of a hyperbola). \pagebreak As for the constant $b$, there 
is no segment of the length $b$ in Fig\.~10.1. For this reason the constant 
$b$ is called the length of the {\it imaginary semiaxis\/} of a hyperbola, 
i\.\,e\. the semiaxis which does not actually exist. 
\mydefinition{10.4} A coordinate system $O,\,\bold e_1,\,\bold e_2$
with an orthonormal basis $\bold e_1,\,\bold e_2$ where a hyperbola
is given by its canonical equation \mythetag{10.20} is called a 
{\it canonical coordinate system\/} of\linebreak  this hyperbola. 
\enddefinition
\head
\SectionNum{11}{179} The eccentricity and directrices of a hyperbola.
The property of directrices.
\endhead
\rightheadtext{\S\,11. The eccentricity and directrices \dots}
     The shape and sizes of a hyperbola are determined by two constants 
$a$ and $b$ in its canonical equation \mythetag{10.20}.. Due to the 
relationship \mythetag{10.18} the constant $b$ can be expressed through 
the constant $c$. Multiplying both constants $a$ and $c$ by the same number, 
we change the sizes of a hyperbola, but do not change its shape. The ratio 
of these two constants 
$$
\hskip -2em
\varepsilon=\frac{c}{a}.
\mytag{11.1}
$$
is responsible for the shape of a hyperbola. 
\mydefinition{11.1} The quantity $\varepsilon$ defined by the relationship
\mythetag{11.1}, where $a$ is the real semiaxis and $c$ is the half of the 
interfocal distance, is called the {\it eccentricity\/} of a hyperbola. 
\enddefinition
     The eccentricity \mythetag{11.1} is used in order to define 
one more parameter of a hyperbola. It is usually denoted through $d$:
$$
\hskip -2em
d=\frac{a}{\varepsilon}=\frac{a^2}{c}.
\mytag{11.2}
$$
\mydefinition{11.2}On the plane of a hyperbola there are two lines 
perpendicular to its real axis and placed at the distance $d$ given by 
the formula \mythetag{7.2} from its center. These lines are called 
{\it directrices\/} of a hyperbola. 
\enddefinition
\parshape 14 5cm 5cm 5cm 5cm 5cm 5cm 5cm 5cm 5cm 5cm 5cm 5cm
5cm 5cm 5cm 5cm 5cm 5cm 5cm 5cm 5cm 5cm 5cm 5cm 5cm 5cm 0cm 10cm
     Each hyperbola has two foci and two directrices. Each directrix has 
the corresponding focus of it. \vadjust{\vskip 5pt\hbox to 
0pt{\kern 0pt\includegraphics{angemeng40.eps}\hss}
\vskip -5pt}This is that of the two foci which is more close to the directrix 
in question. Let $M(x,y)$ be some arbitrary point of a hyperbola. Let's 
connect this point with the left focus of the hyperbola $F_1$ and drop the 
perpendicular from it to the left directrix of the hyperbola. Let's denote 
through $H_1$ the base of such a perpendicular and calculate its length
$|MH_1|$:
$$
\hskip -2em
|MH_1|=|x-(-d)|=|d+x|.
\mytag{11.3}
$$
Taking into account \mythetag{11.2}, the formula \mythetag{11.3} can be 
brought to
$$
\hskip -2em
|MH_1|=\Bigl|\frac{a^2}{c}+x\Bigr|=\frac{|a^2+c\,x|}{c}.
\mytag{11.4}
$$
The length of the segment $MF_1$ was already calculated above. Initially 
it was given by one of the formulas \mythetag{10.10}, but later the more 
simple expression \mythetag{10.25} was derived for it: 
$$
\hskip -2em
|MF_1|=\frac{|a^2+c\,x|}{a}.
\mytag{11.5}
$$
If we divide \mythetag{11.5} by \mythetag{11.4}, we obtain the following 
relationship: 
$$
\hskip -2em
\frac{|MF_1|}{|MH_1|}=\frac{c}{a}=\varepsilon.
\mytag{11.6}
$$
The point $M$ can change its position on the hyperbola. Then the numerator 
\pagebreak and the denominator of the fraction \mythetag{11.6} are changed, 
but its value remains unchanged. This fact is known as the property of 
directrices. 
\mytheorem{11.1} The ratio of the distance from some arbitrary point $M$ of 
a hyperbola to its focus and the distance from this point to the corresponding
directrix is a constant which is equal to the eccentricity of the hyperbola.
\endproclaim
\head
\SectionNum{12}{181} The equation of a tangent line to a hyperbola.
\endhead
\rightheadtext{\S\,12. \dots of a tangent line to a hyperbola.}
\parshape 3 0cm 10cm 0cm 10cm 5cm 5cm 
     Let's consider a hyperbola given by its canonical equation 
\mythetag{10.20} in its canonical coordinate system (see
Definition~\mythedefinition{10.4}). \vadjust{\vskip 5pt\hbox to 
0pt{\kern 0pt\includegraphics{angemeng41.eps}\hss}
\vskip -5pt}Let's draw a tangent line to this hyperbola and denote through
$M=M(x_0,y_0)$ its tangency\linebreak point (see Fig\.~12.1). Our goal is 
to write the equation of a tangent line to a hyperbola. 
\par
\parshape 7 5cm 5cm 5cm 5cm 5cm 5cm 5cm 5cm 5cm 5cm 5cm 5cm 0cm 10cm 
     A hyperbola consists of two branches, each branch being a curve 
composed by two halves --- the upper half and the lower half. The upper 
halves of the hyperbola branches can be treated as a graph 
of some function of the form 
$$
\hskip -2em
y=f(x)
\mytag{12.1}
$$
defined in the union of two intervals $(-\infty,-a)\cup(a,+\infty)$. 
Lower halves of a hyperbola can also be treated as a graph of some function
of the form \mythetag{12.1} with the same domain. The equation of a tangent 
line to the graph of the function \mythetag{12.1} is given by the following 
well-known formula (see \mycite{9}):
$$
\hskip -2em
y=y_0+f'(x_0)\,(x-x_0). 
\mytag{12.2}
$$\par
     In order to apply the formula \mythetag{12.2} to a hyperbola \pagebreak
we need to calculate the derivative of the function \mythetag{12.1}. Let's 
substitute the function \mythetag{12.1} into the equation \mythetag{10.20}:
$$
\hskip -2em
\frac{x^2}{a^2}-\frac{(f(x))^2}{b^2}=1.
\mytag{12.3}
$$
The equality \mythetag{12.3} is fulfilled identically in $x$. Let's 
differentiate the equality \mythetag{12.3} with respect to $x$. This yields
$$
\hskip -2em
\frac{2\,x}{a^2}-\frac{2\,f(x)\,f'(x)}{b^2}=0.
\mytag{12.4}
$$
Let's apply \mythetag{12.4} for to calculate the derivative $f'(x)$:
$$
\hskip -2em
f'(x)=\frac{b^2\,x}{a^2\,f(x)}.
\mytag{12.5}
$$
In order to substitute \mythetag{12.5} into the equation \mythetag{12.2} we 
change $x$ for $x_0$ and $f(x)$ for $f(x_0)=y_0$. As a result we get
$$
\hskip -2em
f'(x_0)=\frac{b^2\,x_0}{a^2\,y_0}.
\mytag{12.6}
$$\par
     Let's substitute \mythetag{12.6} into the equation of the tangent line 
\mythetag{12.2}. This yields the following relationship
$$
\hskip -2em
y=y_0+\frac{b^2\,x_0}{a^2\,y_0}\,(x-x_0). 
\mytag{12.7}
$$
Eliminating the denominator, we write the equality \mythetag{12.7} as 
$$
\hskip -2em
a^2\,y\,y_0-b^2\,x\,x_0=a^2\,y_0^2-b^2\,x_0^2.
\mytag{12.8}
$$
Now let's divide both sides of the equality \mythetag{12.8} by $a^2\,b^2$:
$$
\pagebreak
\hskip -2em
\frac{x\,x_0}{a^2}-\frac{y\,y_0}{b^2}=\frac{x_0^2}{a^2}
-\frac{y_0^2}{b^2}.
\mytag{12.9}
$$\par
     Note that the point $M=M(x_0,y_9)$ is on the hyperbola. Therefore its 
coordinates satisfy the equation \mythetag{10.20}:
$$
\hskip -2em
\frac{x_0^2}{a^2}-\frac{y_0^2}{b^2}=1.
\mytag{12.10}
$$
Taking into account \mythetag{12.10}, we can transform \mythetag{12.9} to
$$
\hskip -2em
\frac{x\,x_0}{a^2}-\frac{y\,y_0}{b^2}=1.
\mytag{12.11}
$$
\mytheorem{12.1} For a hyperbola determined by its canonical equation 
\mythetag{10.20} its tangent line that touches this hyperbola at the 
point $M=M(x_0,y_0)$ is given by the equation \mythetag{12.11}.
\endproclaim
     The equation \mythetag{12.11} is a particular case of the equation 
\mythetag{3.22} where the constants $A$, $B$, and $D$ are given by the
formulas
$$
\xalignat 3
&\hskip -2em
A=\frac{x_0}{a^2},
&&B=-\frac{y_0}{b^2},
&&D=1.
\quad
\mytag{12.12}
\endxalignat
$$
According to the definition~\mythedefinition{3.6} and the formulas 
\mythetag{3.21} the constants $A$ and $B$ in \mythetag{12.12} are
the covariant components of the normal vector for the tangent line 
to a hyperbola. The tangent line equation \mythetag{12.11} is written
in a canonical coordinate system of a hyperbola. The basis of such a
coordinate system is orthonormal. Therefore the formula \mythetag{3.19} 
and the formula \mythetagchapter{32.4}{1} from 
Chapter~\uppercase\expandafter{\romannumeral 1} yield the following
relationships:
$$
\xalignat 2
&\hskip -2em
A=n_1=n^1, 
&&B=n_2=n^2. 
\mytag{12.13}
\endxalignat
$$
\mytheorem{12.2} The quantities $A$ and $B$ in \mythetag{12.12} are the
coordinates of the normal vector $\bold n$ for the tangent line
to a hyperbola which is given by the equation \mythetag{12.11}.
\endproclaim
\noindent The relationships \mythetag{12.13} prove the
theorem~\mythetheorem{12.2}.
\head
\SectionNum{13}{184} Focal property of a hyperbola.
\endhead
\rightheadtext{\S\,13. Focal property of a hyperbola.}
     Like in the case of an ellipse, assume that a hyperbola is manufactured
of a thin strip of some flexible material and assume that its surface is 
covered with a light reflecting layer. For such a hyperbola we can formulate 
the following focal property. 
\mytheorem{13.1} A light ray emitted in one focus of a hyperbola 
upon reflecting on its surface goes to infinity so that its backward 
extension passes through the other focus of this hyperbola. 
\endproclaim
\mytheorem{13.2} The tangent line of a hyperbola is a bisector in the triangle 
composed by the tangency point and two foci of the hyperbola.
\endproclaim
\parshape 3 0cm 10cm 0cm 10cm 5cm 5cm 
     The theorem~\mythetheorem{13.2} is a purely geometric version 
of the theorem~\mythetheorem{13.1}. 
\vadjust{\vskip 5pt\hbox to 
0pt{\kern 0pt\includegraphics{angemeng42.eps}\hss}
\vskip -5pt}These theorems are equivalent due to the reflection law saying 
that the angle of reflection is equal to the angle of incidence.
\demo{Proof}\parshape 10 5cm 5cm 5cm 5cm 5cm 5cm 5cm 5cm 
5cm 5cm 5cm 5cm 5cm 5cm 5cm 5cm 5cm 5cm 0cm 10cm 
Let's consider some point $M=M(x_0,y_0)$ on the hyperbola and draw the 
tangent line to the hyperbola through this point as shown in Fig\.~13.1. 
Let $N$ be the $x$-intercept of this tangent line. Due to $M$ and $N$, 
we have the segment $[MN]$. This segment is perpendicular to the normal vector
$\bold n$ of the tangent line. In order to prove that $[MN]$ is a bisector 
in the triangle $F_1MF_2$ it is sufficient to prove that $\bold n$ is 
directed along the bisector for the external angle of this triangle at its 
vertex $M$. This condition can be written as
$$
\hskip -2em
\frac{(\overrightarrow{F_1M\,\,}\!\!,\bold n)}{|F_1M|}
=\frac{(\overrightarrow{MF_2\,}\!,\bold n)}{|MF_2|}.
\mytag{13.1}
$$\par
     The coordinates of the points $F_1$ and $F_2$ in a canonical 
coordinate system of the hyperbola are known (see formulas \mythetag{10.9}).
The coordinates of the point $M=M(x_0,y_0)$ are also known. Therefore we 
can find the coordinates of the vectors $\overrightarrow{F_1M\,\,}\!$ and
$\overrightarrow{MF_2\,}\!$ used in the above formula \mythetag{13.1}:
$$
\xalignat 2
&\hskip -2em
\overrightarrow{F_1M\,\,}\!
=\Vmatrix x_0+c\\ y_0\endVmatrix,
&&\overrightarrow{MF_2\,}\!
=\Vmatrix c-x_0\\ -y_0\endVmatrix.
\mytag{13.2}
\endxalignat
$$
The tangent line that touches the hyperbola at the point $M=\break =M(x_0,y_0)$ 
is given by the equation \mythetag{12.11}. The coordinates of the normal vector 
$\bold n$ of this tangent line in Fig\.~13.1 are given by the formulas 
\mythetag{12.12} and \mythetag{12.13}:
$$
\hskip -2em
\bold n=\Vmatrix
\dfrac{x_0}{a^2}\\
\vspace{2ex}
-\dfrac{y_0}{b^2}
\endVmatrix
\mytag{13.3}
$$
Relying upon \mythetag{13.2} and \mythetag{13.3}, we apply the formula
\mythetagchapter{33.3}{1} from Chapter~\uppercase\expandafter{\romannumeral 1} 
in order to calculate the scalar products in \mythetag{13.1}:
$$
\aligned
&(\overrightarrow{F_1M\,\,}\!\!,\bold n)
=\frac{c\,x_0+x_0^2}{a^2}-\frac{y_0^2}{b^2}
=\frac{c\,x_0}{a^2}+\frac{x_0^2}{a^2}-\frac{y_0^2}{b^2},\\
\vspace{1ex}
&(\overrightarrow{MF_2\,}\!,\bold n)
=\frac{c\,x_0-x_0^2}{a^2}+\frac{y_0^2}{b^2}
=\frac{c\,x_0}{a^2}-\frac{x_0^2}{a^2}+\frac{y_0^2}{b^2}.
\endaligned
\quad
\mytag{13.4}
$$\par
     The coordinates of the point $M$ satisfy the equation \mythetag{10.20}:
$$
\hskip -2em
\frac{x_0^2}{a^2}-\frac{y_0^2}{b^2}=1.
\mytag{13.5}
$$
Due to \mythetag{13.5} the formulas \mythetag{13.4} simplify to 
$$
\hskip -2em
\aligned
&(\overrightarrow{F_1M\,\,}\!\!,\bold n)
=\frac{c\,x_0}{a^2}+1=\frac{c\,x_0+a^2}{a^2},\\
\vspace{1ex}
&(\overrightarrow{MF_2\,}\!,\bold n)
=\frac{c\,x_0}{a^2}-1
=\frac{c\,x_0-a^2}{a^2}.
\endaligned
\mytag{13.6}
$$\par
     In order to calculate the denominators in the formula 
\mythetag{13.1} we use the formulas \mythetag{10.26} and 
\mythetag{10.27}. In this case, when applied to
the point $M=M(x_0,y_0)$, they yield
$$
\hskip -2em
\aligned
&|MF_1|=\frac{c\,|x_0|+\sign(x_0)\,a^2}{a},\\
&|MF_2|=\frac{c\,|x_0|-\sign(x_0)\,a^2}{a}.
\endaligned
\mytag{13.7}
$$
Due to the purely numeric identity $|x_0|=\sign(x_0)\,x_0$ we can
write \mythetag{13.7} in the following form:
$$
\hskip -2em
\aligned
&|MF_1|=\frac{c\,x_0+a^2}{a}\,\sign(x_0),\\
&|MF_2|=\frac{c\,x_0-a^2}{a}\,\sign(x_0).
\endaligned
\mytag{13.8}
$$
From the formulas \mythetag{13.6} and \mythetag{13.8} we easily derive
the equalities 
$$
\xalignat 2
&\frac{(\overrightarrow{F_1M\,\,}\!\!,\bold n)}{|MF_1|}
=\frac{\sign(x_0)}{a},
&&\frac{(\overrightarrow{MF_2\,}\!,\bold n)}{|MF_2|}
=\frac{\sign(x_0)}{a}
\endxalignat
$$
that prove the equality \mythetag{13.1}. The theorem~\mythetheorem{13.2} 
is proved.
\qed\enddemo
As we noted above the theorem~\mythetheorem{13.1} is equivalent to the
theorem~\mythetheorem{13.2} due to the light reflection law. 
Therefore the theorem~\mythetheorem{13.1} is also proved. 
\head
\SectionNum{14}{186} Asymptotes of a hyperbola.
\endhead
\rightheadtext{\S\,14. Asymptotes of a hyperbola.}
\parshape 6 0cm 10cm 0cm 10cm 0cm 10cm 0cm 10cm 0cm 10cm 5cm 5cm 
      {\it Asymptotes \/} are usually understood as some straight
lines to which some points of a given curve come unlimitedly close 
along some unlimitedly long fragments of this curve.
Each hyperbola has two asymptotes (see Fig\.~14.1). In a canonical 
coordinate system the asymptotes of the hyperbola given by the 
equation \mythetag{10.20} are determined by the following equations:
$$
y=\pm\frac{b}{a}\,x.
\mytag{14.1}
$$
One of the asymptotes is associated with the plus sign in the
formula~\mythetag{14.1}, 
\vadjust{\vskip 5pt\hbox to 
0pt{\kern 0pt\includegraphics{angemeng43.eps}\hss}
\vskip -5pt}the other asymptote is associated with the opposite
minus sign.
\par
\parshape 3 5cm 5cm 5cm 5cm 0cm 10cm
     The theory of asymptotes is closely related to the theory of
limits. This theory is usually studied within the course of mathematical
analysis (see \mycite{9}). For this reason I do not derive the equations
\mythetag{14.1} in this book. 
\head
\SectionNum{15}{187} Parabola. Canonical equation of a parabola.
\endhead
\rightheadtext{\S\,15. Parabola.}
\mydefinition{15.1}\parshape 1 5cm 5cm 
{\it A parabola\/} is a set of points on some plane each of which is
equally distant from some fixed point $F$ of this plane and from some
straight line $d$ lying on this plane. \vadjust{\vskip 5pt\hbox to 
0pt{\kern -5pt\includegraphics{angemeng44.eps}\hss}
\vskip -5pt}The point $F$ is called the {\it focus\/} of this parabola, 
while the line $d$ is called its {\it directrix}.
\enddefinition
\parshape 8 5cm 5cm 5cm 5cm 5cm 5cm 5cm 5cm 5cm 5cm 5cm 5cm 5cm 5cm 
0cm 10cm 
     Assume that a parabola with the focus $F$ and with the 
directrix $d$ is given. Let's drop the perpendicular from the point 
$F$ onto the line $d$ and denote through $D$ the base of such a 
perpendicular. Let's choose the line $DF$ for the $x$-axis of a coordinate
system. Then we denote through $O$ the midpoint of the segment $[DF]$ and
choose this point $O$ for the origin. And finally, we draw the $y$-axis 
of a coordinate system through the point $O$ perpendicular to the line
$DF$ (see Fig\.~15.1). Choosing the unity scales along the axes, we 
ultimately fix a coordinate system with an orthonormal basis on the parabola 
plane. \par
     Let's denote through $p$ the distance from the focus $F$ of the parabola
to its directrix $d$, i\.\,e\. we set 
$$
\hskip -2em
|DF|=p.
\mytag{15.1}
$$
The point $O$ is the midpoint of the segment $[DF]$. Therefore the equality
\mythetag{15.1} leads to the following equalities:
$$
|DO|=|OF|=\frac{p}{2}.
\mytag{15.2}
$$
The relationships \mythetag{15.2} determine the coordinates of $D$ and $F$:
$$
\xalignat 2
&\hskip -2em
D=D(-p/2,0),
&&F=F(p/2,0).
\mytag{15.3}
\endxalignat
$$\par
    Let $M=M(x,y)$ be some arbitrary point of the parabola. According to
the definition~\mythedefinition{15.1}, the following equality is fulfilled:
$$
\hskip -2em
|MF|=|MH|
\mytag{15.4}
$$
(see Fig\.~15.1). Due to \mythetag{15.3} the length of the segment $[MF]$
in the chosen coordinate system is given by the formula
$$
\hskip -2em
|MF|=\sqrt{y^2+(x-p/2)^2)}.
\mytag{15.5}
$$
The formula for the length of $[MH]$ is even more simple:
$$
\hskip -2em
|MH|=x+p/2.
\mytag{15.6}
$$
Substituting \mythetag{15.5} and \mythetag{15.6} into \mythetag{15.4}, we
get the equation
$$
\hskip -2em
\sqrt{y^2+(x-p/2)^2)}=x+p/2.
\mytag{15.7}
$$\par
     Let's square both sides of the equation \mythetag{15.7}:
$$
\hskip -2em
y^2+(x-p/2)^2)=(x+p/2)^2.
\mytag{15.8}
$$
Upon expanding brackets and collecting similar terms in \mythetag{15.8},
we bring this equation to the following form:
$$
y^2=2\,p\,x.
\mytag{15.9}
$$
\mydefinition{15.2} The equality \mythetag{15.9} is called the 
{\it canonical equation\/} of a parabola.
\enddefinition
\mytheorem{15.1} For each point $M(x,y)$ of the parabola determined 
by the initial equation \mythetag{15.7} its coordinates satisfy the 
canonical equation \mythetag{15.9}.  
\endproclaim
     Due to \mythetag{15.1} the constant $p$ in the equation \mythetag{15.9} 
is a non-negative quantity. The case $p=0$ corresponds to the degenerate
parabola. From the definition~\mythedefinition{15.1} it is easy to derive
that in this case the parabola turns to the straight line coinciding with the 
$x$-axis in Fig\.~15.1. The case of the degenerate parabola is excluded by
means of the inequality 
$$
\hskip -2em
p>0.
\mytag{15.10}
$$
Due to the inequality \mythetag{15.10} from the equation \mythetag{15.9}
we derive
$$
\hskip -2em
x\geqslant 0.
\mytag{15.11}
$$
\mytheorem{15.2} The canonical equation of the parabola \mythetag{15.9}
is equivalent to the initial equation \mythetag{15.7}.
\endproclaim
\demo{Proof} In order to prove the theorem~\mythetheorem{15.2} \pagebreak
it is sufficient to invert the calculations performed in deriving the 
equality \mythetag{15.9} from \mythetag{15.7}. Note that the passage from 
\mythetag{15.8} to \mythetag{15.9} is invertible. The passage from 
\mythetag{15.7} to \mythetag{15.8} is also invertible due to the inequality 
\mythetag{15.11}, which follows from the equation \mythetag{15.9}. This 
observation completes the proof of the theorem~\mythetheorem{15.2}.
\qed\enddemo
\mydefinition{15.3} The point $O$ in Fig\.~15.1 is called the 
{\it vertex\/} of the parabola, the line $DF$ coinciding with the $x$-axis
is called the {\it axis\/} of the parabola.
\enddefinition
\mydefinition{15.4} A coordinate system $O,\,\bold e_1,\,\bold e_2$
with an orthonormal basis $\bold e_1,\,\bold e_2$ where a parabola
is given by its canonical equation \mythetag{15.9} and where the 
inequality \mythetag{15.10} is fulfilled is called a 
{\it canonical coordinate system\/} of this parabola. 
\enddefinition
\head
\SectionNum{16}{190} The eccentricity of a parabola.
\endhead
\rightheadtext{\S\,16. The eccentricity of a parabola.}
     The definition of a parabola \mythedefinition{15.1} is substantially 
different from the definition of an ellipse \mythedefinition{6.1} and from
the definition of a hyperbola \mythedefinition{10.1}. But it is similar to
the property of directrices of an ellipse in the theorem~\mythetheorem{7.1}
and to the property of directrices of a hyperbola in the 
theorem~\mythetheorem{11.1}. Comparing the definition of a parabola 
\mythedefinition{15.1} with these theorems, we can formulate the following
definition. 
\mydefinition{16.1} The eccentricity of a parabola is postulated to be
equal to the unity: $\varepsilon=1$. 
\enddefinition
\head
\SectionNum{17}{190} The equation of a tangent line to a parabola.
\endhead
\rightheadtext{\S\,17. \dots tangent line to a parabola.}
     Let's consider a parabola given by its canonical equation \mythetag{15.9} 
in its canonical coordinate system (see Definition~\mythedefinition{15.4}). 
Let's draw a tangent line to this parabola and denote through $M=M(x_0,y_0)$
the tangency point (see Fig\.~17.1). Our goal is to write the equation of the
tangent line to the parabola through the point $M=M(x_0,y_0)$.\par
\parshape 16 0cm 10cm 0cm 10cm 5cm 5cm 5cm 5cm 5cm 5cm 5cm 5cm 
5cm 5cm 5cm 5cm 5cm 5cm 5cm 5cm 5cm 5cm 5cm 5cm 5cm 5cm 5cm 5cm 
5cm 5cm 0cm 10cm
     An parabola is a curve composed by two halves --- the upper half and the
lower half. Any one of these two halves of a parabola can be treated as a graph 
of some function
$$
\hskip -2em
y=f(x)
\mytag{17.1}
$$
with the domain $(0,+\infty)$. \vadjust{\vskip 5pt\hbox to 
0pt{\kern 0pt\includegraphics{angemeng45.eps}\hss}
\vskip -5pt}The equation of a tangent line to the graph of
a function \mythetag{17.1} is given by the well-known formula
$$
\gathered
y-y_0=\\
=f'(x_0)\,(x-x_0). 
\endgathered
\quad
\mytag{17.2}
$$
(see \mycite{9}). In order to apply the formula \mythetag{17.2} 
to a parabola one should calculate the derivative of the function 
\mythetag{17.1}. Let's substitute \mythetag{17.1} into the equation 
\mythetag{15.9}:
$$
\hskip -2em
(f(x))^2=2\,p\,x.
\mytag{17.3}
$$
The equality \mythetag{17.3} is fulfilled identically in $x$. Let's 
differentiate the equality \mythetag{17.3} with respect to $x$. This yields
$$
\hskip -2em
2\,f(x)\,f'(x)=2\,p.
\mytag{17.4}
$$
Let's apply the formula \mythetag{17.4} for to calculate the derivative 
$$
\hskip -2em
f'(x)=\frac{p}{f(x)}.
\mytag{17.5}
$$
In order to substitute \mythetag{17.5} into the equation \mythetag{17.2} we 
change $x$ for $x_0$ and $f(x)$ for $f(x_0)=y_0$. As a result we get
$$
\hskip -2em
f'(x_0)=\frac{p}{y_0}.
\mytag{17.6}
$$\par
     Let's substitute \mythetag{17.6} \pagebreak into the equation of the 
tangent line \mythetag{17.2}. This yields the following relationship:
$$
\hskip -2em
y-y_0=\frac{p}{y_0}\,(x-x_0). 
\mytag{17.7}
$$
Eliminating the denominator, we write the equality \mythetag{17.7} as 
$$
\hskip -2em
y\,y_0-y_0^2=p\,x-p\,x_0.
\mytag{17.8}
$$
Note that the point $M=M(x_0,y_9)$ is on the parabola. Therefore its 
coordinates satisfy the equality \mythetag{15.9}:
$$
\hskip -2em
y_0^2=2\,p\,x_0.
\mytag{17.9}
$$
Taking into account \mythetag{17.9}, we can transform \mythetag{17.8} to
$$
\hskip -2em
y\,y_0=p\,x+p\,x_0.
\mytag{17.10}
$$
This is the required equation of a tangent line to a parabola. 
\mytheorem{17.1} For a parabola determined by its canonical equation 
\mythetag{15.9} the tangent line that touches this parabola at the 
point $M=M(x_0,y_0)$ is given by the equation \mythetag{17.10}.
\endproclaim
     Let's write the equation of a tangent line to a parabola in the 
following slightly transformed form:
$$
\hskip -2em
p\,x-y\,y_0+p\,x_0=0.
\mytag{17.11}
$$
The equation \mythetag{17.11} is a special instance of the equation 
\mythetag{3.22} where the constants $A$, $B$, and $D$ are given by 
the formulas
$$
\xalignat 3
&\hskip -2em
A=p,
&&B=-y_0,
&&D=p\,x_0.
\quad
\mytag{17.12}
\endxalignat
$$
According to the definition~\mythedefinition{3.6} and the formulas
\mythetag{3.21}, the constants $A$ and $B$ in \mythetag{17.12} are
the covariant components of the normal vector of a tangent line
to a parabola. \pagebreak The equation \mythetag{17.11} is written 
in a canonical coordinate system of the parabola. The basis of a 
canonical system is orthonormal (see Definition~\mythedefinition{15.4}). 
In the case of an orthonormal basis the formula \mythetag{3.19} and 
the formula \mythetagchapter{32.4}{1} from  
Chapter~\uppercase\expandafter{\romannumeral 1} yield
$$
\xalignat 2
&\hskip -2em
A=n_1=n^1, 
&&B=n_2=n^2. 
\mytag{17.13}
\endxalignat
$$
\mytheorem{17.2} The quantities $A$ and $B$ in \mythetag{17.12} 
are the coordinates of the normal vector $\bold n$ of a tangent 
line to a parabola in the case where this tangent line is given
by the equation \mythetag{17.10}.
\endproclaim
\noindent The relationships \mythetag{17.13} prove the
theorem~\mythetheorem{17.2}.\par
\head
\SectionNum{18}{193} Focal property of a parabola.
\endhead
\rightheadtext{\S\,18. Focal property of a parabola.}
     Assume that we have a parabola manufactured of a thin strip
of some flexible material covered with a light reflecting layer.
For such a parabola the following focal property is formulated. 
\mytheorem{18.1} A light ray emitted from the focus of a para\-bola 
upon reflecting on its surface goes to infinity parallel to the
axis of this parabola. 
\endproclaim
\mytheorem{18.2}\parshape 1 5cm 5cm 
For any tangent line of a parabola \vadjust{\vskip 5pt\hbox to 
0pt{\kern 0pt\includegraphics{angemeng46.eps}\hss}
\vskip -5pt}the triangle formed by the tangency point $M$,
its focus $F$, and  by the point $N$ at which this tangent 
line intersects the axis of the parabola is an isosceles
triangle, i\.\,e\. the equality $|MF|=|NF|$ holds.
\endproclaim
\parshape 6 5cm 5cm 5cm 5cm 5cm 5cm 5cm 5cm 5cm 5cm 0cm 10cm 
     As we see in Fig\.~18.1, the theorem \mythetheorem{18.1}
follows from the theorem \mythetheorem{18.2} due to the reflection
law saying that the angle of reflection is equal to the angle 
of incidence and due to the equality of inner crosswise lying angles
in the intersection of two parallel lines by a third line
(see \mycite{6}). 
\demo{Proof of the theorem~\mythetheorem{18.2}} For to prove
this theorem we choose some tangency point $M(x_0,y_0)$ and write
the equation of the tangent line to the parabola in the form of
\mythetag{17.10}:
$$
\hskip -2em
y\,y_0=p\,x+p\,x_0.
\mytag{18.1}
$$
Let's find the intersection point of the tangent line \mythetag{18.1} 
with the axis of the parabola. In a canonical coordinate system the
axis of the parabola coincides with the $x$-axis (see 
Definition~\mythedefinition{15.3}). Substituting $y=0$ into the
equation \mythetag{18.1}, we get $x=-x_0$, which determines the coordinates
of the point $N$:
$$
\hskip -2em
N=N(-x_0,0).
\mytag{18.2}
$$
From \mythetag{18.2} and \mythetag{15.3} we derive the length of the segment
$[NF]$:
$$
\hskip -2em
|NF|=p/2-(-x_0)=p/2+x_0.
\mytag{18.3}
$$
In the case of a parabola the length of the segment $[MF]$ coincides with the 
length of the segment $[MH]$ (see Definition~\mythedefinition{15.1}). Therefore
from \mythetag{15.3} we derive
$$
\hskip -2em
|MF|=|MH|=x_0-(-p/2)=x_0+p/2.
\mytag{18.4}
$$
Comparing \mythetag{18.3} and \mythetag{18.4} we get the required equality 
$|MF|=|NF|$. As a result the theorem~\mythetheorem{18.2} is proved. As for the
theorem~\mythetheorem{18.1}, it equivalent to the theorem~\mythetheorem{18.2}.
\qed\enddemo
\head
\SectionNum{19}{194} The scale of eccentricities.
\endhead
\rightheadtext{\S\,19. The scale of eccentricities.}
     The eccentricity of an ellipse is determine by the formula \mythetag{7.1},
where the parameters $c$ and $a$ are related by the inequalities \mythetag{6.4}. 
Hence the eccentricity of an ellipse obeys the inequalities
$$
\hskip -2em
0\leqslant\varepsilon<1.
\mytag{19.1}
$$\par
     The eccentricity of a parabola is equal to unity by definition. Indeed, the
definition~\mythedefinition{16.1} yields 
$$
\hskip -2em
\varepsilon=1.
\mytag{19.2}
$$\par
     The eccentricity of a hyperbola is defined by the formula \mythetag{11.1}, 
where the parameters $c$ and $a$ obey the inequalities \mythetag{10.8}. Hence
the eccentricity of a hyperbola obeys the inequalities
$$
\hskip -2em
1<\varepsilon\leqslant+\infty.
\mytag{19.3}
$$\par
     The formulas \mythetag{19.1}, \mythetag{19.2}, and \mythetag{19.3} show
that the eccentricities of ellipses, parabolas, and hyperbolas fill the interval
from $0$ to $+\infty$ without omissions, i\.\,e\. we have the {\it continuous
scale of eccentricities}.
\head
\SectionNum{20}{195} Changing a coordinate system.
\endhead
\rightheadtext{\S\,20. Changing a coordinate system.}
     Let $O,\,\bold e_1,\,\bold e_2,\,\bold e_3$ and $\tilde O,
\,\tilde\bold e_1,\,\tilde\bold e_2,\,\tilde\bold e_3$ be two Cartesian 
coordinate systems in the space $\Bbb E$. They consist of the bases 
$\bold e_1,\,\bold e_2,\,\bold e_3$ and $\tilde\bold e_1,\,\tilde\bold e_2,
\,\tilde\bold e_3$ complemented with two points $O$ and $\tilde O$, which
are called origins (see Definition~\mythedefinition{1.1}). The transition
from the basis $\bold e_1,\,\bold e_2,\,\bold e_3$ to the basis $\tilde\bold e_1,
\,\tilde\bold e_2,\,\tilde\bold e_3$ and vice versa is described by two transition
matrices $S$ and $T$ whose components are in the following transition formulas:
$$
\xalignat 2
&\hskip -2em
\tilde\bold e_j=\sum^3_{i=1}S^{\kern 0.5pt i}_j\,\bold e_i,
&&\bold e_j=\sum^3_{i=1}T^{\kern 0.5pt i}_j\,\tilde\bold e_i
\mytag{20.1}
\endxalignat
$$
(see formulas \mythetagchapter{22.4}{1} and \mythetagchapter{22.9}{1}
in Chapter~\uppercase\expandafter{\romannumeral 1}).\par
     In order to describe the transition from the coordinate system $O,\,\bold e_1,
\,\bold e_2,\,\bold e_3$ to the coordinate system $\tilde O,\,\tilde\bold e_1,
\,\tilde\bold e_2,\,\tilde\bold e_3$ and vice versa two auxiliary parameters 
are employed. These are the vectors $\bold a=\overrightarrow{O\tilde O\,\,}\!$ 
and $\tilde\bold a=\overrightarrow{\tilde OO\,\,}\!$. 
\mydefinition{20.1} The vectors $\bold a=\overrightarrow{O\tilde O\,\,}\!$ 
and $\tilde\bold a=\overrightarrow{\tilde OO\,\,}\!$ are called the
{\it origin displacement vectors}. 
\enddefinition
     The origin displacement vector $\bold a=\overrightarrow{O\tilde O\,\,}\!$
is usually expanded in the basis $\bold e_1,\,\bold e_2,\,\bold e_3$, while the
other origin displacement vector $\tilde\bold a=\overrightarrow{\tilde OO\,\,}\!$ 
is expanded in the basis $\tilde\bold e_1,\,\tilde\bold e_2,\,\tilde\bold e_3$:
$$
\xalignat 2
&\hskip -2em
\bold a=\sum^3_{i=1}a^i\,\bold e_i,
&&\tilde\bold a=\sum^3_{i=1}\tilde a^i\,\tilde\bold e_i.
\mytag{20.2}
\endxalignat
$$
The vectors $\bold a=\overrightarrow{O\tilde O\,\,}\!$ are
$\tilde\bold a=\overrightarrow{\tilde OO\,\,}\!$ opposite to each other, i\.\,e\. 
the following relationships are fulfilled:
$$
\xalignat 2
&\hskip -2em
\bold a=-\tilde\bold a,
&&\tilde\bold a=-\bold a.
\mytag{20.3}
\endxalignat
$$
Their coordinates in the expansions \mythetag{20.2} are related with each other
by means of the formulas
$$
\xalignat 2
&\hskip -2em
\tilde a^i=-\sum^3_{j=1}
T^{\kern 0.5pt i}_j\,a^{\kern 0.2pt j},
&&a^i=-\sum^3_{j=1}
S^{\kern 0.5pt i}_j\,\tilde a^{\kern 0.2pt j}.
\mytag{20.4}
\endxalignat
$$
The formulas \mythetag{20.4} are derived from the formulas \mythetag{20.3} 
with the use of the formulas \mythetagchapter{25.4}{1} and
\mythetagchapter{25.5}{1} from Chapter~\uppercase\expandafter{\romannumeral 1}.
\head
\SectionNum{21}{196} Transformation of the coordinates of a point
under a change of a coordinate system.
\endhead
\rightheadtext{\S\,21. Transformation of the coordinates of a point \dots}
     Let $O,\,\bold e_1,\,\bold e_2,\,\bold e_3$ and $\tilde O,
\,\tilde\bold e_1,\,\tilde\bold e_2,\,\tilde\bold e_3$ be two Cartesian 
coordinate systems in the space $\Bbb E$ and let $X$ be some arbitrary point
in the space $\Bbb E$. Let's denote through 
$$
\xalignat 2
&\hskip -2em
X=X(x^1,x^2,x^3),
&&X=X(\tilde x^1,\tilde x^2,\tilde x^3)
\quad
\mytag{21.1}
\endxalignat
$$
the presentations of the point $X$ in these two coordinate systems. The radius
vectors of the point $X$ in these systems are related with each other by means
of the relationships
$$
\xalignat 2
&\hskip -2em
\bold r_X=\bold a+\tilde\bold r_X,
&&\tilde\bold r_X=\tilde\bold a+\bold r_X, 
\mytag{21.2}
\endxalignat
$$
The coordinates of $X$ in \mythetag{21.1} are the coordinates of its 
radius vectors in the bases of the corresponding coordinate systems:
$$
\xalignat 2
&\hskip -2em
\bold r_X=\sum^3_{j=1}x^{\kern 0.4pt j}\,\bold e_j,
&&\tilde\bold r_X=\sum^3_{j=1}\tilde x^{\kern 0.4pt j}
\,\tilde\bold e_j.
\mytag{21.3}
\endxalignat
$$
From \mythetag{21.2}, \mythetag{21.3}, and \mythetag{20.2}, applying the 
formulas \mythetagchapter{25.4}{1} and \mythetagchapter{25.5}{1} from
Chapter~\uppercase\expandafter{\romannumeral 1}, we easily derive the
following relationships:
$$
\xalignat 2
&\hskip -2em
\tilde x^i=\sum^3_{j=1}
T^{\kern 0.5pt i}_j\,x^{\kern 0.4pt j}+\tilde a^i,
&&x^i=\sum^3_{j=1}
S^{\kern 0.5pt i}_j\,\tilde x^{\kern 0.4pt j}+a^i.
\quad
\mytag{21.4}
\endxalignat
$$
\mytheorem{21.1} Under a change of coordinate systems in the space $\Bbb E$
determined by the formulas \mythetag{20.1} and \mythetag{20.2} the coordinates
of points are transformed according to the formulas \mythetag{21.4}.
\endproclaim
     The formulas \mythetag{21.4} are called the {\it direct\/} and 
{\it inverse transformation formulas} for the coordinates of a point {\it under
a change of a Cartesian coordinate system}.\par
\head
\SectionNum{22}{197} Rotation of a rectangular coordinate system
on a plane. The rotation matrix.
\endhead
\rightheadtext{\S\,22. Rotation of a coordinate system\dots}
    Let $O,\,\bold e_1,\,\bold e_2$ and $\tilde O,\,\tilde\bold e_1,
\,\tilde\bold e_2$ be two Cartesian coordinate system on a plane. The formulas
\mythetag{20.1}, \mythetag{20.2}, and \mythetag{20.4} in this case are written
as follows:
$$
\allowdisplaybreaks
\xalignat 2
\hskip -2em
\tilde\bold e_j&=\sum^2_{i=1}S^{\kern 0.5pt i}_j\,\bold e_i,
&\bold e_j&=\sum^2_{i=1}T^{\kern 0.5pt i}_j\,\tilde\bold e_i,
\mytag{22.1}\\
\hskip -2em
\bold a&=\sum^2_{i=1}a^i\,\bold e_i,
&\tilde\bold a&=\sum^2_{i=1}\tilde a^i\,\tilde\bold e_i,
\mytag{22.2}\\
\hskip -2em
\tilde a^i&=-\sum^2_{j=1}
T^{\kern 0.5pt i}_j\,a^{\kern 0.2pt j},
&a^i&=-\sum^2_{j=1}
S^{\kern 0.5pt i}_j\,\tilde a^{\kern 0.2pt j}.
\quad
\mytag{22.3}
\endxalignat
$$
Let $X$ be some arbitrary point of the plane. Its coordinates in the
coordinate systems $O,\,\bold e_1,\,\bold e_2$ and $\tilde O,\,\tilde\bold e_1,
\,\tilde\bold e_2$ are transformed according to the formulas similar to 
\mythetag{21.4}:
$$
\xalignat 2
&\hskip -2em
\tilde x^i=\sum^2_{j=1}
T^{\kern 0.5pt i}_j\,x^{\kern 0.4pt j}+\tilde a^i,
&&x^i=\sum^2_{j=1}
S^{\kern 0.5pt i}_j\,\tilde x^{\kern 0.4pt j}+a^i.
\quad
\mytag{22.4}
\endxalignat
$$\par
\parshape 13 0cm 10cm 0cm 10cm 5cm 5cm 5cm 5cm 5cm 5cm 5cm 5cm 
5cm 5cm 5cm 5cm 5cm 5cm 5cm 5cm 5cm 5cm 5cm 5cm 0cm 10cm
     Assume that the origins $O$ and $\tilde O$ do coincide. In this case
the parameters $a^1,\,a^2$ and $\tilde a^1,\,\tilde a^2$ in the formulas
\mythetag{22.2}, \mythetag{22.3}, \mythetag{22.3}, and \mythetag{22.4}
\vadjust{\vskip 5pt\hbox to 
0pt{\kern -5pt\includegraphics{angemeng47.eps}\hss}
\vskip -5pt}do vanish. Under the assumption $O=\tilde O$ we consider 
the special case where the bases $\bold e_1,\,\bold e_2$ and $\tilde\bold e_1,
\,\tilde\bold e_2$ both are orthonormal and where one of them is produced from
the other by means of the rotation by some angle $\varphi$ (see Fig\.~22.1).
For two bases on a plane the transition matrices $S$ and $T$ are square matrices
$2\times 2$. The components of the direct transition matrix $S$ are taken
from the following formulas:
$$
\hskip -2em
\aligned
&\tilde\bold e_1=\cos\varphi\cdot\bold e_1+\sin\varphi\cdot\bold e_2,\\
&\tilde\bold e_2=-\sin\varphi\cdot\bold e_1+\cos\varphi\cdot\bold e_2.
\endaligned
\mytag{22.5}
$$
The formulas \mythetag{22.5} are derived on the base of Fig\.~22.1.\par
     Comparing the formulas \mythetag{22.5} with the first relationship 
in \mythetag{22.1}, we get $S^1_1=\cos\varphi$, $S^2_1=\sin\varphi$, 
$S^2_1=-\sin\varphi$, $S^2_2=\cos\varphi$. Hence we have the following formula
for the matrix $S$:
$$
\hskip -2em
S=\Vmatrix \cos\varphi & -\sin\varphi\\
\vspace{1.2ex}
\sin\varphi & \cos\varphi
\endVmatrix.
\mytag{22.6}
$$
\mydefinition{22.1} The square matrix of the form \mythetag{22.6} is called 
the {\it rotation matrix by the angle\/} $\varphi$.
\enddefinition
     The inverse transition matrix $T$ is the inverse matrix for $S$ (see
theorem~\mythetheoremchapter{23.1}{1} in 
Chapter~\uppercase\expandafter{\romannumeral 1}). Due to this fact and due
to the formula \mythetag{22.6} we can calculate the matrix $T$:
$$
\hskip -2em
T=\Vmatrix \cos(-\varphi) & -\sin(-\varphi)\\
\vspace{1.2ex}
\sin(-\varphi) & \cos(-\varphi)
\endVmatrix.
\mytag{22.7}
$$
The matrix \mythetag{22.7} is also a rotation matrix by the angle $\varphi$. 
But the angle $\varphi$ in it is taken with the minus sign, which means that
the rotation is performed in the opposite direction.\par 
     Let's write the relationships \mythetag{22.4} taking into account that 
$\bold a=0$ and $\tilde\bold a=0$, which follows from $O=\tilde O$, and taking 
into account the formulas \mythetag{22.6} and \mythetag{22.7}: 
$$
\align
&\hskip -2em
\aligned
&\tilde x^1=\cos(\varphi)\,x^1+\sin(\varphi)\,x^2,\\
&\tilde x^2=-\sin(\varphi)\,x^1+\cos(\varphi)\,x^2,
\endaligned
\mytag{22.8}\\
\vspace{2ex}
&\hskip -2em
\aligned
&x^1=\cos(\varphi)\,\tilde x^1-\sin(\varphi)\,\tilde x^2,\\
&x^2=\sin(\varphi)\,\tilde x^1+\cos(\varphi)\,\tilde x^2.
\endaligned
\mytag{22.9}
\endalign
$$
The formulas \mythetag{22.8} and \mythetag{22.9} are the transformation formulas
for the coordinates of a point under the rotation of the rectangular coordinate
system shown in Fig\.~22.1.\par
\head
\SectionNum{23}{199} Curves of the second order.
\endhead
\rightheadtext{\S\,23. Curves of the second order.}
\mydefinition{23.1} A {\it curve of the second order\/} or a \pagebreak
{\it quadric\/} on a plane is a curve which is given by a polynomial equation 
of the second order in some Cartesian coordinate system
$$
A\,x^2+2\,B\,x\,y+C\,y^2+2\,D\,x+2\,E\,y+F=0.
\quad
\mytag{23.1}
$$
\enddefinition
    Here $x=x^1$ and $y=x^2$ are the coordinates of a point on a plane. 
Note that the transformation of these coordinates under a change of one 
coordinate system for another is given by the functions of the first order 
in $x^1$ and $x^2$ (see formulas \mythetag{22.4}). For this reason the general
form of the equation of a quadric \mythetag{23.1} does not change under 
a change of a coordinate system though the values of the parameters 
$A$, $B$, $C$, $D$, $E$, and $F$ in it can change. From what was said 
we derive the following theorem.
\mytheorem{23.1} For any curve of the second order on a plane, i\.\,e\. for 
any quadric, there is some rectangular coordinate system with an orthonormal 
basis such that this curve is given by an equation of the form \mythetag{23.1}
in this coordinate system.
\endproclaim
\head
\SectionNum{24}{200} Classification of curves of the second order.
\endhead
\rightheadtext{\S\,24. Classification of curves \dots}
     Let $\Gamma$ be a curve of the second order on a plane given by an equation 
of the form \mythetag{23.1} in some rectangular coordinate system with the 
orthonormal basis. Passing from one of such coordinate systems to another, one can
change the constant parameters $A$, $B$, $C$, $D$, $E$, and $F$ in \mythetag{23.1},
and one can always choose a coordinate system in which the equation \mythetag{23.1}
takes its most simple form. 
\mydefinition{24.1} The problem of finding a rectangular coordinate system with the 
orthonormal basis in which the equation of a curve of the second order $\Gamma$ 
takes its most simple form is called the problem of {\it bringing\/} the equation
of a curve $\Gamma$ to its {\it canonical form}.
\enddefinition
     An ellipse, a hyperbola, and parabola are examples of curves of the second 
order on a plane. \pagebreak The canonical forms of the equation \mythetag{23.1} 
for these curves are already known to us (see formulas \mythetag{6.15}, 
\mythetag{10.20}, and \mythetag{15.9}). 
\mydefinition{24.2} The problem of grouping curves by the form of their canonical
equations is called the problem of {\it classification\/} of curves of the second
order.
\enddefinition
\mytheorem{24.1} For any curve of the second order $\Gamma$ there is a rectangular
coordinate system with an orthonormal basis where the constant parameter $B$ of 
the equation \mythetag{23.1} for $\Gamma$ is equal to zero: $B=0$. 
\endproclaim
\demo{Proof} Let $O,\,\bold e_1,\,\bold e_2$ be rectangular coordinate system with 
an orthonormal basis $\bold e_1,\,\bold e_2$ where the equation of the curve $\Gamma$ 
has the form \mythetag{23.1} (see Theorem~\mythetheorem{23.1}). If $B=0$ in
\mythetag{23.1}, then $O,\,\bold e_1,\,\bold e_2$ is a required coordinate system.
\par
     If $B\neq 0$, then we perform the rotation of the coordinate system 
$O,\,\bold e_1,\,\bold e_2$ about the point $O$ by some angle $\varphi$. Such a 
rotation is equivalent to the change of variables
$$
\hskip -2em
\aligned
&x=\cos(\varphi)\,\tilde x-\sin(\varphi)\,\tilde y,\\
&y=\sin(\varphi)\,\tilde x+\cos(\varphi)\,\tilde y
\endaligned
\mytag{24.1}
$$
in the equation \mythetag{23.1} (see formulas \mythetag{22.9}). Upon substituting
\mythetag{24.1} into the equation \mythetag{23.1} we get the analogous equation
$$
\tilde A\,\tilde x^2+2\,\tilde B\,\tilde x\,\tilde y+\tilde C
\,\tilde y^2+2\,\tilde D\,\tilde x+2\,\tilde E\,\tilde y
+\tilde F=0
\quad
\mytag{24.2}
$$
whose parameters $\tilde A$, $\tilde B$, $\tilde C$, $\tilde D$, $\tilde E$, and
$\tilde F$ are expressed through the parameters $A$, $B$, $C$, $D$, $E$, and $F$ 
of the initial equation \mythetag{23.1}. For the parameter $\tilde B$ in
\mythetag{24.2} we derive the formula
$$
\hskip -2em
\gathered
\tilde B=(C-A)\,\cos(\varphi)\,\sin(\varphi)
+B\,\cos^2(\varphi)\,-\\
\vspace{1ex}
-\,B\,\sin^2(\varphi)=\frac{C-A}{2}\,\sin(2\,\varphi)
+B\,\cos(2\,\varphi). 
\endgathered
\mytag{24.3}
$$
Since $B\neq 0$, we determine $\varphi$ as a solution of the
equation
$$
\hskip -2em
\cot(2\,\varphi)=\frac{A-C}{2\,B}. 
\mytag{24.4}
$$
The equation \mythetag{24.4} is always solvable in the form of 
$$
\hskip -2em
\varphi=\frac{\pi}{4}-\frac{1}{2}\,\arctan\Bigl(\frac{A-C}{2\,B}
\,\Bigr). 
\mytag{24.5}
$$
Comparing \mythetag{24.4} and \mythetag{24.3}, we see that having 
determined the angle $\varphi$ by means of the formula \mythetag{24.5}, 
we provide the vanishing of the parameter $\tilde B=0$ in the rotated
coordinate system. The theorem~\mythetheorem{24.1} is proved.
\qed\enddemo
     Let's apply the theorem~\mythetheorem{24.1} and write the equation of
the curve $\Gamma$ in the form of the equation \mythetag{23.1} with $B=0$:
$$
\hskip -2em
A\,x^2+C\,y^2+2\,D\,x+2\,E\,y+F=0.
\mytag{24.6}
$$
The equation \mythetag{24.6} provides the subdivision of all curves of the
second order on a plane into three types:
\roster
\rosteritemwd=0pt
\item"--" the {\bf elliptic type} where $A\neq 0$, $C\neq 0$ and the parameters
          $A$ and $C$ are of the same sign, i\.\,e\. $\sign(A)=\sign(C)$;
\item"--" the {\bf hyperbolic type} where $A\neq 0$, $C\neq 0$ and the parameters
          $A$ and $C$ are of different signs, i\.\,e\. $\sign(A)\neq\sign(C)$;
\item"--" the {\bf parabolic type} where $A=0$ or $C=0$.
\endroster
The parameters $A$ and $C$ in \mythetag{24.6} cannot vanish simultaneously since in
this case the degree of the polynomial in \mythetag{24.6} would be lower than two,
which would contradict the definition~\mythedefinition{23.1}.
\par
     {\bf Curves of the elliptic type}. If the conditions $A\neq 0$, $C\neq 0$, and
$\sign(A)=\sign(C)$ are fulfilled, without loss of generality we can assume that
$A>0$ and $C>0$. In this case the equation \mythetag{24.6} can be written as
follows:
$$
A\,\Bigl(x+\frac{D}{A}\,\Bigr)^2
+C\,\Bigl(y+\frac{E}{C}\,\Bigr)^2
+\Bigl(F-\frac{D^2}{A}-\frac{E^2}{C}\,\Bigr)=0.
\quad
\mytag{24.7}
$$
Let's perform the following change of variables in \mythetag{24.7}:
$$
\xalignat 2
&\hskip -2em
x=\tilde x-\frac{D}{A},
&&y=\tilde y-\frac{E}{C}.
\mytag{24.8}
\endxalignat
$$
The change of variables \mythetag{24.8} corresponds to the displacement of the origin 
without rotation (the case of the unit matrices $S=1$ and $T=1$ in \mythetag{22.4}). 
In addition to \mythetag{24.8} we denote
$$
\hskip -2em
\tilde F=F-\frac{D^2}{A}-\frac{E^2}{C}.
\mytag{24.9}
$$
Taking into account \mythetag{24.8} and \mythetag{24.9}, we write \mythetag{24.7} 
as 
$$
\hskip -2em
A\,\tilde x^2+C\,\tilde y^2+\tilde F=0,
\mytag{24.10}
$$
where the coefficients $A$ and $C$ are positive: $A>0$ and $C>0$.\par
    The equation \mythetag{24.10} provides the subdivision of curves of the elliptic 
type into three subtypes:
\roster
\rosteritemwd=10pt
\item"--" the case of an {\bf ellipse} where $\tilde F<0$;
\item"--" the case of an {\bf imaginary ellipse} where $\tilde F>0$;
\item"--" the case of a {\bf point} where $\tilde F=0$.
\endroster\par
     In the case of an ellipse the equation \mythetag{24.10} is brought to the form
\mythetag{6.15} in the variables $\tilde x$ and $\tilde y$. As we know, this equation 
describes an ellipse.\par
     In the case of an imaginary ellipse the equation \mythetag{24.10} reduces to
the equation which is similar to the equation of an ellipse:
$$
\hskip -2em
\frac{\tilde x^2}{a^2}+\frac{\tilde y^2}{b^2}=-1.
\mytag{24.11}
$$
The equation \mythetag{24.11} has no solutions. Such an equation describes the empty
set of points.\par 
    The case of a point is sometimes called the case of a {\it pair of imaginary
intersecting lines}, which is somewhat not exact. In this case the equation \mythetag{24.10} 
describes a single point with the coordinates $\tilde x=0$ and $\tilde y=0$.\par
     {\bf Curves of the hyperbolic type}. If the conditions $A\neq 0$, $C\neq 0$, and
$\sign(A)\neq\sign(C)$ are fulfilled, without loss of generality we can assume that
$A>0$ and $C<0$. In this case the equation \mythetag{24.6} can be written as
$$
A\,\Bigl(x+\frac{D}{A}\,\Bigr)^2
-\tilde C\,\Bigl(y-\frac{E}{\tilde C}\,\Bigr)^2
+\Bigl(F-\frac{D^2}{A}+\frac{E^2}{\tilde C}\,\Bigr)=0,
\quad
\mytag{24.12}
$$
where $-C=\tilde C>0$. Let's denote
$$
\hskip -2em
\tilde F=F-\frac{D^2}{A}+\frac{E^2}{\tilde C}
\mytag{24.13}
$$
and then perform the following change of variables, which corresponds to a displacement
of the origin:
$$
\xalignat 2
&\hskip -2em
x=\tilde x-\frac{D}{A},
&&y=\tilde y+\frac{E}{\tilde C}.
\mytag{24.14}
\endxalignat
$$
Due to \mythetag{24.13} and \mythetag{24.14} the equation \mythetag{24.12} is written
as
$$
A\,\tilde x^2-\tilde C\,\tilde y^2+\tilde F=0,
\quad
\mytag{24.15}
$$
where the coefficients $A$ and $\tilde C$ are positive: $A>0$ and $\tilde C>0$.\par
    The equation \mythetag{24.15} provides the subdivision of curves of the hyperbolic
type into two subtypes:
\roster
\rosteritemwd=0pt
\item"--" the case of a {\bf hyperbola} where $\tilde F\neq 0$;
\item"--" the case of a {\bf pair of intersecting lines} where $\tilde F=0$.
\endroster\par
     In the case of a hyperbola the equation \mythetag{24.15} is brought to the form 
\mythetag{10.20} in the variables $\tilde x$ and $\tilde y$. As we know, it describes 
a hyperbola.\par
     In the case of a pair of intersecting lines the left hand side of the equation
\mythetag{24.15} \pagebreak is written as a product of two multiplicands and the 
equation \mythetag{24.15} is brought to
$$
\hskip -2em
(\sqrt{A}\,\tilde x+\sqrt{\tilde C}\,\tilde y)
(\sqrt{A}\,\tilde x-\sqrt{\tilde C}\,\tilde y)=0.
\mytag{24.16}
$$
The equation \mythetag{24.16} describes two line on a plane that intersect at the point
with the coordinates $\tilde x=0$ and$\tilde y=0$.\par
     {\bf Curves of the parabolic type}. For curves of this type 
there are two options in the equation \mythetag{24.6}: $A=0$, $C\neq 0$
or $C=0$, $A\neq 0$. But the second option reduces to the first one upon
changing variables $x=-\tilde y$, $y=\tilde x$, which corresponds to the
rotation by the angle $\varphi=\pi/2$. Therefore without loss of generality
we can assume that $A=0$ and $C\neq 0$. Then the equation \mythetag{24.6} is brought
to
$$
\hskip -2em
y^2+2\,D\,x+2\,E\,y+F=0.
\mytag{24.17}
$$
In order to transform the equation \mythetag{24.17} we apply the change of variables
$y=\tilde y-E$, which corresponds to the displacement of the origin along the $y$-axis.
Then we denote $\tilde F=F+E^2$. As a result the equation \mythetag{24.17} is written
as
$$
\hskip -2em
\tilde y^2+2\,D\,x+\tilde F=0.
\mytag{24.18}
$$\par
    The equation \mythetag{24.18} provides the subdivision of curves of the parabolic
type into four subtypes:    
\roster
\rosteritemwd=10pt
\item"--" the case of a {\bf parabola} where $D\neq 0$;
\item"--" the case of a {\bf pair of parallel lines}\newline 
where $D=0$ and $\tilde F<0$;
\item"--" the case of a {\bf pair of coinciding lines}\newline 
where $D=0$ and $\tilde F=0$;
\item"--" the case of a {\bf pair of imaginary parallel lines}\newline 
where $D=0$ and $\tilde F>0$.
\endroster\par
    In the case of a parabola the equation \mythetag{24.18} reduces to the equation
\mythetag{15.9} and describes a parabola.\par
    In the case of a pair of parallel lines we introduce the notation 
$\tilde F=-y_0^2$ into the equation \mythetag{24.18}. As a result the equation 
\mythetag{24.18} is written in the following form:
$$
\hskip -2em
(\tilde y+y_0)(\tilde y-y_0)=0.
\mytag{24.19}
$$
The equation \mythetag{24.19} describes a pair of lines parallel to the $y$-axis
and being at the distance $2\,y_0$ from each other.\par
     In the case of a pair of coinciding lines the equation \mythetag{24.18} reduces
to the form $\tilde y^2=0$, which describes a single line coinciding with the 
$y$-axis.\par 
     In the case of a pair of imaginary parallel lines the equation \mythetag{24.18} 
has no solutions. It describes the empty set.\par
\head
\SectionNum{25}{206} Surfaces of the second order.
\endhead
\rightheadtext{\S\,25. Surfaces of the second order.}
\mydefinition{25.1} A {\it surface of the second order\/} or a {\it quadric\/} 
in the space $\Bbb E$ is a surface which in some Cartesian coordinate system 
is given by a polynomial equation of the second order:
$$
\gathered
A\,x^2+2\,B\,x\,y+C\,y^2+2\,D\,x\,z+2\,E\,y\,z\,+\\
+\,F\,z^2+2\,G\,x+2\,H\,y+2\,I\,z+J=0.
\endgathered
\quad
\mytag{25.1}
$$
\enddefinition
    Here $x=x^1$, $y=x^2$, $z=x^3$ are the coordinates of points of the space
$\Bbb E$. Note that the transformation of the coordinates of points under a change
of a coordinate system is given by functions of the first order in $x^1$, $x^2$, 
$x^3$ (see formulas \mythetag{21.4}). For this reason the general form of a quadric 
equation \mythetag{25.1} remains unchanged under a change of a coordinate system, 
though the values of the parameters $A$, $B$, $C$, $D$, $E$, $F$, $G$, $H$, $I$, and
$J$ can change. The following theorem is immediate from what was said. 
\mytheorem{25.1} For any surface of the second order in the space $\Bbb E$, i\.\,e\.
for any quadric, there is some rectangular coordinate system with an orthonormal 
basis \pagebreak such that this surface is given by an equation of the form 
\mythetag{25.1} in this coordinate system.
\endproclaim
\head
\SectionNum{26}{207} Classification of surfaces of the second order.
\endhead
\rightheadtext{\S\,26. Classification of surfaces \dots}
     The problem of classification of surfaces of the second order in $\Bbb E$ is 
solved below following the scheme explained in \S\,2 of 
Chapter~\uppercase\expandafter{\romannumeral 6} in the book \mycite{1}. Let $S$ be
surface of the second order given by the equation \mythetag{25.1} on some rectangular
coordinate system with an orthonormal basis (see Theorem~\mythetheorem{25.1}). Let's
arrange the parameters $A$, $B$, $C$, $D$, $E$, $F$, $G$, $H$, $I$ of the equation
\mythetag{25.1} into two matrices
$$
\xalignat 2
\hskip -2em
\Cal F=\Vmatrix A & B & D\\
\vspace{1ex}
B & C & E\\
\vspace{1ex}
D & E & F
\endVmatrix,
&&\Cal D=\Vmatrix G\\ 
\vspace{1ex} H\\ 
\vspace{1ex} I
\endVmatrix.
\quad
\mytag{26.1}
\endxalignat
$$
The matrices \mythetag{26.1} are used in the following theorem.
\mytheorem{26.1} For any surface of the second order $S$ there is a rectangular
coordinate system with an orthonormal basis such that the matrix $\Cal F$ in
\mythetag{26.1} is diagonal, while the matrix $\Cal D$ is related to the matrix
$\Cal F$ by means of the formula $\Cal F\cdot\Cal D=0$.
\endproclaim
     The proof of the theorem~\mythetheorem{26.1} can be found in \mycite{1}. 
Applying the theorem~\mythetheorem{26.1}, we can write the equation
\mythetag{25.1} as
$$
A\,x^2+C\,y^2+F\,z^2+2\,G\,x+2\,H\,y+2\,I\,z+J=0.
\quad
\mytag{26.2}
$$
The equation \mythetag{26.2} and the theorem~\mythetheorem{26.1} provide the
subdivision of all surfaces of the second order in $\Bbb E$ into four types:
\roster
\rosteritemwd=0pt
\item"--" the {\bf elliptic type} where $A\neq 0$, $C\neq 0$, $F\neq 0$ and
          the quantities $A$, $C$, and $F$ are of the same sign; 
\item"--" the {\bf hyperbolic type} where $A\neq 0$, $C\neq 0$, $F\neq 0$ and
          the quantities $A$, $C$ and $F$ are of different signs; 
\item"--" the {\bf parabolic type} where exactly one of the quantities 
          $A$, $C$. and $F$ is equal to zero \pagebreak and exactly one of the 
          quantities $G$, $H$, and $I$ is nonzero.
\item"--" the {\bf cylindrical type} in all other cases. 
\endroster\par
     {\bf Surfaces of the elliptic type}. From the conditions
$A\neq 0$, $C\neq 0$, $F\neq 0$ in \mythetag{26.2} and from the condition
$\Cal F\cdot\Cal D=0$ in the theorem~\mythetheorem{26.1} we derive
$G=0$, $H=0$, and $I=0$. Since the quantities $A$, $C$, and $F$ are of the same
sign, without loss of generality we can assume that all of them are positive. 
Hence for all surfaces of the elliptic type we can write \mythetag{26.2} as 
$$
\hskip -2em
A\,x^2+C\,y^2+F\,z^2+J=0,
\mytag{26.3}
$$
where $A>0$, $C>0$, and $F>0$. The equation \mythetag{26.3} provides the subdivision 
of surfaces of the elliptic type into three subtypes:
\roster
\rosteritemwd=10pt
\item"--" the case of an {\bf ellipsoid} where $J<0$;
\item"--" the case of an {\bf imaginary ellipsoid} where $J>0$;
\item"--" the case of a {\bf point} where $J=0$.
\endroster\par
     The case of an ellipsoid is the most non-trivial of the three. In this case
the equation \mythetag{26.3} is brought to
$$
\hskip -2em
\frac{x^2}{a^2}+\frac{y^2}{b^2}+\frac{z^2}{c^2}=1.
\mytag{26.4}
$$
The equation \mythetag{26.4} describes the surface which is called 
\vadjust{\vskip 5pt\hbox to 
0pt{\kern -5pt\includegraphics{angemeng48.eps}\hss}
\vskip 134pt}\pagebreak an {\it ellipsoid}. This surface is shown in Fig\.~26.1.
\par
     In the case of an imaginary ellipsoid the equation \mythetag{26.3} has no
solutions. It describes the empty set.\par
     In the case of a point the equation \mythetag{26.3} can be written in the
form very similar to the equation of an ellipsoid \mythetag{26.4}:
$$
\hskip -2em
\frac{x^2}{a^2}+\frac{y^2}{b^2}+\frac{z^2}{c^2}=0.
\mytag{26.5}
$$
The equation \mythetag{26.5} describes a single point in the space $\Bbb E$ 
with the coordinates $x=0$, $y=0$, $z=0$.\par
     {\bf Surfaces of the hyperbolic type}. From the three conditions
$A\neq 0$, $C\neq 0$, $F\neq 0$ in \mythetag{26.2} and from the condition
$\Cal F\cdot\Cal D=0$ in the theorem~\mythetheorem{26.1} we derive $G=0$, $H=0$,
and $I=0$. The quantities $A$, $C$, and $F$ are of different signs. Without 
loss of generality we can assume that two of them are positive and one of them
is negative. By exchanging axes, which preserves the orthogonality of coordinate
systems and orthonormality of their bases, we can transform the equation 
\mythetag{26.2} so that we would have $A>0$, $C>0$, and $F<0$. As a result we 
conclude that for all surfaces of the hyperbolic type the initial equation 
\mythetag{26.2} can be brought to the form
$$
\hskip -2em
A\,x^2+C\,y^2+F\,z^2+J=0,
\mytag{26.6}
$$
where $A>0$, $C>0$, $F<0$. The equation \mythetag{26.6} provides the subdivision 
of surfaces of the hyperbolic type into three subtypes:
\roster
\rosteritemwd=5pt
\item"--" the case of a {\bf hyperboloid of one sheet} where $J<0$;
\item"--" the case of a {\bf hyperboloid of two sheets} where $J>0$;
\item"--" the case of a {\bf cone} where $J=0$.
\endroster\par
     In the case of a hyperboloid of one sheet the equation \mythetag{26.6} can be
written in the following form:
$$
\hskip -2em
\frac{x^2}{a^2}+\frac{y^2}{b^2}-\frac{z^2}{c^2}=1.
\mytag{26.7}
$$
The equation \mythetag{26.7} describes a surface which is called the {\it hyperboloid 
of one sheet}. It is shown in Fig\.~26.2.\par
     In the case of a hyperboloid of two sheets the equation \mythetag{26.6} can be
written in the following form:
$$
\hskip -2em
\frac{x^2}{a^2}+\frac{y^2}{b^2}-\frac{z^2}{c^2}=-1.
\mytag{26.8}
$$
The equation \mythetag{26.8} describes a surface which is called the 
\vadjust{\vskip 15pt\hbox to 
0pt{\kern -5pt\includegraphics{angemeng49.eps}\hss}
\vskip 135pt}{\it hyperboloid of two sheets}. This surface is shown in 
Fig\.~26.3.\par
     In the case of a cone the equation \mythetag{26.6} is transformed to the
equation which is very similar to the equations \mythetag{26.7} and \mythetag{26.8}, 
but with zero in the right hand side:  
$$
\hskip -2em
\frac{x^2}{a^2}+\frac{y^2}{b^2}-\frac{z^2}{c^2}=0.
\mytag{26.9}
$$
The equation \mythetag{26.9} describes a surface which is called the {\it cone}. 
This surface is shown in Fig\.~26.4.\par
     {\bf Surfaces of the parabolic type}. For this type of surfaces exactly one
of the three quantities $A$, $C$, and $F$ is equal to zero. By exchanging axes, 
which preserves the orthogonality of coordinate systems and orthonormality of 
their bases, we can transform the equation \mythetag{26.2} so that we would have 
$A\neq 0$, $C\neq 0$, and $F=0$. From $A\neq 0$, $C\neq 0$, and from the condition
$\Cal F\cdot\Cal D=0$ in the theorem~\mythetheorem{26.1} we derive $G=0$ and $H=0$. 
The value of $I$ is not determined by the condition $\Cal F\cdot\Cal D=0$. However,
according to the definition of surfaces of the parabolic type exactly one of the
three quantities $G$, $H$, $I$ should be nonzero. Due to $G=0$ and $H=0$ we
conclude that $I\neq 0$. As a result the equation \mythetag{26.2} is written as 
$$
\hskip -2em
A\,x^2+C\,y^2+2\,I\,z+J=0,
\mytag{26.10}
$$
where $A\neq 0$, $C\neq 0$, and $I\neq 0$. The condition $I\neq 0$ means that
we can perform the displacement of the origin along the $z$-axis equivalent
to the change of variables
$$
\hskip -2em
z\to z-\frac{J}{2\,I}. 
\mytag{26.11}
$$
Upon applying the change of variables \mythetag{26.11} to the equation 
\mythetag{26.10} this equation is written as 
$$
\hskip -2em
A\,x^2+C\,y^2+2\,I\,z=0,
\mytag{26.12}
$$
where $A\neq 0$, $C\neq 0$, $I\neq 0$. The equation \mythetag{26.12}
provides the subdivision of surfaces of the parabolic type into two 
subtypes:
\roster
\rosteritemwd=5pt
\item"--" the case of an {\bf elliptic paraboloid} where the quantities
$A\neq 0$ and $C\neq 0$ are of the same sign;
\item"--" the case of a {\bf hyperbolic paraboloid}, where the quantities 
$A\neq 0$ and $C\neq 0$ are of different signs.
\endroster\par
     In the case of an elliptic paraboloid the equation \mythetag{26.12}
can be written in the following form:
$$
\hskip -2em
\frac{x^2}{a^2}+\frac{y^2}{b^2}=2\,z.
\mytag{26.13}
$$
The equation \mythetag{26.13} describes a surface which is called 
\vadjust{\vskip 5pt\hbox to 
0pt{\kern 5pt\includegraphics{angemeng50.eps}\hss}
\vskip 145pt}the {\it elliptic paraboloid}. This surface is shown in
Fig\.~26.5.\par
     In the case of a hyperbolic paraboloid the equation \mythetag{26.12}
can be transformed to the following form:
$$
\hskip -2em
\frac{x^2}{a^2}-\frac{y^2}{b^2}=2\,z.
\mytag{26.14}
$$
The equation \mythetag{26.14} describes a saddle surface which is called 
the {\it hyperbolic paraboloid}. This surface is shown in Fig\.~26.6.\par
     {\bf Surfaces of the cylindrical type}. According to the results 
from \S\,2 in Chapter~\uppercase\expandafter{\romannumeral 6} of the book
\mycite{1}, in the cylindrical case the dimension reduction occurs. This 
means that there is a rectangular coordinate system with an orthonormal
basis where the variable $z$ drops from the equation \mythetag{26.2}:
$$
\hskip -2em
A\,x^2+C\,y^2+2\,G\,x+2\,H\,y+J=0.
\mytag{26.15}
$$
The classification of surfaces of the second order described by the equation 
\mythetag{26.15} is equivalent to the classification of curves of the second 
order on a plane described by the equation \mythetag{24.6}. The complete 
type list of such surfaces contains nine cases:
\roster
\rosteritemwd=10pt
\item"--" the case of an {\bf elliptic cylinder};
\item"--" the case of an {\bf imaginary elliptic cylinder};
\item"--" the case of a {\bf straight line}\vadjust{\vskip 5pt\hbox to 
0pt{\kern 5pt\includegraphics{angemeng51.eps}\hss}
\vskip 145pt};
\item"--" the case of a {\bf hyperbolic cylinder};
\item"--" the case of a {\bf pair of intersecting planes};
\item"--" the case of a {\bf parabolic cylinder};
\item"--" the case of a {\bf pair of parallel planes};
\item"--" the case of a {\bf pair of coinciding planes};
\item"--" the case of a {\bf pair of imaginary parallel planes}.
\endroster\par
\parshape 5 0cm 10cm 0cm 10cm 5cm 5cm 5cm 5cm 5cm 5cm 
     In the case of an elliptic cylinder the equation \mythetag{26.15}
can be transformed to the following form:
$$
\hskip 2.5em
\frac{x^2}{a^2}+\frac{y^2}{b^2}=1.
\mytag{26.16}
$$
The equation \mythetag{26.16} coincides with the equation of an ellipse
on a plane \mythetag{6.15}. \vadjust{\vskip 5pt\hbox to 
0pt{\kern 5pt\includegraphics{angemeng52.eps}\hss}
\vskip -5pt}In the space $\Bbb E$ it describes a surface which is called
the {\it elliptic cylinder}. This surface is shown in Fig\.~26.7.\par
\parshape 3 5cm 5cm 5cm 5cm 0cm 10cm 
     In the case of an imaginary elliptic \ cylinder \ \pagebreak the \ 
equation \mythetag{26.15} can be transformed to the following form:
$$
\frac{x^2}{a^2}+\frac{y^2}{b^2}=-1.
\mytag{26.17}
$$
The equation \mythetag{26.17} has no solutions. Such an equation describes
the empty set.\par 
     In the case of a straight line the equation \mythetag{26.15} can be 
brought to the form similar to \mythetag{26.16} and \mythetag{26.17}:
$$
\hskip -2em
\frac{x^2}{a^2}+\frac{y^2}{b^2}=0.
\mytag{26.18}
$$
The equation \mythetag{26.18} describes a straight line in the space coinciding 
with the $z$-axis. In the canonical form this line is given by the equations 
$x=0$ and $y=0$ (see \mythetag{5.14}).\par
     In the case of a hyperbolic cylinder the equation \mythetag{26.15} can be
transformed to the following form:
$$
\hskip -2em
\frac{x^2}{a^2}-\frac{y^2}{b^2}=1.
\mytag{26.19}
$$
The equation \mythetag{26.19} coincides with the equation of a hyperbola on a 
plane \mythetag{6.15}. In the spatial case it describes a surface which is 
called the {\it hyperbolic cylinder}. It is shown in Fig\.~26.8.\par
     The next case in the list is the case of a pair of intersecting planes. 
In this case the equation \mythetag{26.15} can be brought to 
$$
\hskip -2em
\frac{x^2}{a^2}-\frac{y^2}{b^2}=0.
\mytag{26.20}
$$
The equation \mythetag{26.20} describes the union of two intersecting planes 
in the space given by the equations
$$
\xalignat 2
&\hskip -2em
\frac{x}{a}-\frac{y}{b}=0,
&&\frac{x}{a}+\frac{y}{b}=0.
\mytag{26.21}
\endxalignat
$$
The planes \mythetag{26.21} intersect along a line which coincides with the
$z$-axis. This line is given by the equations $x=0$ and $y=0$.\par
     In the case of a parabolic cylinder the equation \mythetag{26.15} reduces
to the equation coinciding with the equation of a parabola 
$$
\hskip -2em
y^2=2\,p\,x.
\mytag{26.22}
$$
In the space the equation \mythetag{26.22} describes a surface which is called 
the {\it parabolic cylinder}. This surface is shown in Fig\.~26.9.\par
     In the case of a pair of parallel planes the equation \mythetag{26.15} is
brought to $y^2-y_0^2=0$, where $y_0\neq 0$. It describes two parallel planes 
given by the equations
$$
\xalignat 2
&\hskip -2em
y=y_0,
&&y=-y_0.
\mytag{26.23}
\endxalignat
$$\par
     In the case of a pair of coinciding planes the equation \mythetag{26.15}
is also brought to $y^2-y_0^2=0$, but the parameter $y_0$ in it is equal to zero. 
Due to $y_0=0$ two planes \mythetag{26.23} are glued into a single plane which is
perpendicular to the $y$-axis and is given by the equation $y=0$.\par
     In the case of a pair of imaginary parallel planes the equation \mythetag{26.15} 
is brought to $y^2+y_0^2=0$, where $y_0\neq 0$. Such an equation has no solutions. 
For this reason it describes the empty set.\par
\newpage
%---------------------------------------------------------------
\setfirstpage
\topmatter
\title
REFERENCES.
\endtitle
\endtopmatter
\document
\rightheadtext{References.}
\leftheadtext{References.}

\par
\newpage
%---------------------------------------------------------------
\setfirstpage
\topmatter
\title
Contacts
\endtitle
\endtopmatter
\document
\line{\vtop{\hsize 5cm
{\bf Address: }
\medskip\noindent
Ruslan A. Sharipov,\newline
Dep\. of Mathematics\newline 
and Information Techn\.,\newline
Bashkir State University,\newline
32 Zaki Validi street,\newline
Ufa 450074, Russia
\medskip
{\bf Phone:}\medskip
\noindent
+7-(347)-273-67-18 (Office)\newline
+7-(917)-476-93-48 (Cell)
}\hss
\vtop{\hsize 4.3cm
{\bf Home address:}\medskip\noindent
Ruslan A. Sharipov,\newline
5 Rabochaya street,\newline
Ufa 450003, Russia
\vskip 1cm
{\bf E-mails:}\medskip
\noindent
r-sharipov\@mail.ru\newline
R\hskip 0.5pt\_\hskip 1.5pt Sharipov\@ic.bashedu.ru}}
\bigskip
{\bf URL's:}\medskip
\noindent
\myhref{http://ruslan-sharipov.ucoz.com}
{http:/\negskp/ruslan-sharipov.ucoz.com}\newline
\myhref{http://freetextbooks.narod.ru/}
{http:/\negskp/freetextbooks.narod.ru}\newline
\myhref{http://sovlit2.narod.ru/}
{http:/\negskp/sovlit2.narod.ru}\newline
\par
\newpage
%---------------------------------------------------------------
\setfirstpage
\topmatter
\title
Appendix
\endtitle
\leftheadtext{List of publications.}
\endtopmatter
\document
\rightheadtext{List of publications.}

\par\newpage
%--------------------------------------------
\setfirstpage
\topmatter
\endtopmatter
\document
\vphantom{k}
\vfill
\centerline{\tencyrit Uchebnoe izdanie}
\vskip 0.7cm
\centerline{\tencyrbf SHARIPOV Ruslan Abdulovich}
\vskip 0.7cm
\centerline{\tencyrbf KURS ANALITICHESKO\wnIshort\ GEOMETRII}
\vskip 1.0cm
\centerline{\tencyrbf Uchebnoe posobie}
\vskip 1.0cm
\centerline{{\tencyrit Redaktor\/} \tencyr G\. G\. Sina{\ae}skaya}
\centerline{{\tencyrit Korrektor\/} \tencyr A\. I\. Nikolaeva}
\vskip 1.0cm
\centerline{\tencyrit Litsenziya na izdatel\wnsoft skuyu deyatelnost\wnsoft}
\centerline{\tencyrit LR {\tencyr\No} \tencyr 021319 ot 05.01.1999\,g.}
\vskip 0.3cm
\centerline{\tencyr Podpisano v pechat{\wnsoft} 02.11.2010\,\ g.}
\centerline{\tencyr Format 60$\times$84/16. Usl\. pech\. l\. 13,11. 
Uch\.-izd\kern 0.3pt\. l\. 11,62.}
\centerline{\tencyr Tirazh 100. Izd\kern 0.3pt\. \No\ 243. 
Zakaz 69\,a.}
\vskip 0.3cm
\centerline{\tencyrit Redaktsionno-izdatel{\wnsoft}ski{\ae} centr}
\centerline{\tencyrit Bashkirskogo gosudarstvennogo universiteta}
\centerline{\tencyrit 450074, RB, g\. Ufa, ul\. Zaki Validi, 32.}
\vskip 0.3cm
\centerline{\tencyrit Otpechatano na mnozhitel{\wnsoft}nom uchastke}
\centerline{\tencyrit Bashkirskogo gosudarstvennogo universiteta}
\centerline{\tencyrit 450074, RB, g\. Ufa, ul\. Zaki Validi, 32.}
\par\newpage
%---------------------------------------------------------------
\enddocument
\end